\documentclass[a4paper,10pt]{article}

\usepackage{amsmath,amssymb,amsthm} 
\usepackage[english]{babel}

\newtheorem{theo}{Theorem}[section]
\newtheorem{defi}[theo]{Definition}
\newtheorem{lemm}[theo]{Lemma}
\newtheorem{prop}[theo]{Proposition}
\newtheorem{coro}[theo]{Corollary}

\newtheorem{exam}[theo]{Example}

\newcommand{\Na}{\mathbb N}

\newcommand{\Ra}{\mathbb R}

\newcommand{\Ca}{\mathbb C}

\newcommand{\scal}[1]{\langle #1 \rangle}

\newcommand{\finpreuve}{\hfill $\Box$}

\newcommand{\name}{$\underline{\qquad \qquad}$}

\newcommand{\curv}{c }

\newcommand{\bott}{\gamma_n}

\newcommand{\diffeo}{ \Psi }

\newcommand{\ginfty}{ g }

\newcommand{\cutoff}{ \kappa }

\begin{document}

\author{  Jean-Marc Bouclet}
\title{{\bf \sc Strichartz estimates on asymptotically hyperbolic manifolds}}

\maketitle

\begin{abstract} We prove local in time Strichartz estimates without loss for the
 restriction of the solution of the Schr\"odinger equation, outside a large compact set,
 on a class of asymptotically hyperbolic manifolds.
\end{abstract}

\section{The results}
\setcounter{equation}{0}
 
Let $ ({\mathcal M} , G) $ be a riemannian manifold of dimension
$ n \geq 2 $. Denote by $ \Delta_G $ the associated
Laplace-Beltrami operator and by $ d G $ the riemannian volume
density. The Strichartz estimates for the Schr\"odinger equation
\begin{eqnarray}
 i \partial_t u + \Delta_G u = 0, \qquad u_{| t= 0} = u_0 ,
 \label{Schrodingerlibre}
\end{eqnarray}
are basically estimates of
$$ ||u||_{L^p ([0,1],L^q ({\mathcal M},dG))} :=
\left( \int_0^1 || u (t,.) ||_{L^q ({\mathcal M},dG)}^p d t \right)^{1/p} ,
 $$
in terms of certain $ L^2 $ quantities of $ u_0 $, when the pair
of exponents $ (p,q) $ satisfies the so called admissibility
conditions
\begin{eqnarray}
\frac{2}{p} + \frac{n}{q} = \frac{n}{2}, \qquad p \geq 2 , \ \
(p,q) \ne (2,\infty) . \label{paireadmissible}
\end{eqnarray}
We recall that Strichartz inequalities play an important role in the proof of
local existence results for non linear Schr\"odinger equations. We won't consider such applications in this paper and will
only focus on Strichartz estimates.

Let us review some classical results. If $
{\mathcal M} = \Ra^n $ with the flat metric $ G_{\rm Eucl} $, it is
well known \cite{Strichartz,GinibreVelo,KeelTao} that
\begin{eqnarray}
|| u ||_{ L^p ([0,1],L^q (\Ra^n))} \lesssim ||u_0||_{L^2(\Ra^n)} .
\label{globalinspace}
\end{eqnarray}
In this model case, the time interval $ [0,1] $ can  be
replaced by $ \Ra $ and the Strichartz estimates are then said to
be global in time\footnote{In this paper, however, we will only
consider local in time Strichartz estimates, ie with $ t \in [0,1]
$.}. Furthermore, the conditions (\ref{paireadmissible}) are seen
to be natural by considering the action of the scaling $ u (t,x)
\mapsto u (t/\lambda^2, x/\lambda) $ on both Schr\"odinger
equation and Strichartz estimates.

In more general situations, estimates of the form
(\ref{globalinspace}) have sometimes to be replaced by
\begin{eqnarray}
|| u ||_{ L^p ([0,1],L^q ({\mathcal M},dG))} \lesssim
||u_0||_{H^s({\mathcal M},dG)} , \qquad s \geq 0,
\label{globalinspaceloss}
\end{eqnarray}
where 
$$  ||u_0||_{H^s({\mathcal M},dG)} := ||(1-\Delta_G)^{s/2}u_0||_{L^2({\mathcal M},dG)}  , $$ 
is the natural $L^2 $ Sobolev norm. 
 If $ s > 0 $,
(\ref{globalinspaceloss}) are called Strichartz estimates with
loss (of $s$ derivatives). Notice that, under fairly general assumptions on $ ({\mathcal M},G) $, we have the Sobolev embeddings, namely $ H^s({\mathcal M},dG) \subset L^q ({\mathcal M},dG) $
for $ s > \frac{n}{2} - \frac{n}{q} $. They show that (\ref{globalinspaceloss}) holds automatically if $s$ is large enough and the point of Strichartz estimates with loss (and a fortiori without loss) is to consider smaller $s$ than those given by Sobolev embeddings.

Such inequalities have been proved by
Bourgain \cite{Bourgain} for the flat torus $ {\mathbb T}^n $ , $ n = 1,2 $, with any $s > 0$, ie
with 'almost no loss' (for certain $ (p,q) $) and by
Burq-G\'erard-Tzvetkov \cite{BGT} for any compact manifold with $
s = 1/p$. The techniques of \cite{BGT} are actually very robust
and can be applied to prove the same results on many non
compact manifolds. The estimates of \cite{BGT} are
known to be sharp, by considering for instance $ {\mathcal M}
= {\mathbb S}^3 $ with $ p = 2 $ and certain subsequences of
eigenfunctions of the Laplacian. This counterexample
can then be used to construct quasi-modes and show that
(\ref{globalinspaceloss}) can not hold in general with $s=0$, even for non compact manifolds.

A natural question is therefore to find (sufficient)
conditions leading to estimates with no loss. 

 A classical one
is the non trapping condition. We recall that $ ({\mathcal M},G )
$ is non trapping if all geodesics escape to infinity (which forces $ {\mathcal M} $ to be non compact). It was for
instance shown in \cite{StTa,RZ,BoucletTzvetkov} that, for non
trapping perturbations of the flat metric on $ \Ra^n $, (\ref{globalinspaceloss}) holds with $ s=0$. By
perturbation, we mean that $ G - G_{\rm Eucl} $ is small
near infinity and we refer to these papers for more details. In
\cite{HTW}, the more general case of non trapping asymptotically
conic manifolds was considered. To emphasize the difference with
 asymptotically hyperbolic manifolds studied in this paper,
we simply recall that $ ({\mathcal M},G) $ is asymptotically conic if $ G $ is close to 
$ d r^2 + r^2 g $, in a neighborhood of infinity diffeomorphic to $ (R,+\infty) \times S $, for some fixed metric $ g $ on a compact manifold $ S$. The asymptotically Euclidean case corresponds to the case where $ S = {\mathbb S}^{n-1} $.

The non trapping condition has however several drawbacks. First it is a non generic property and second there is no simple criterion to check whether a manifold is trapping or not. Furthermore, it is not clearly a {\it necessary} condition to get Strichartz estimates without loss.

In \cite{BoucletTzvetkov}, we partially got rid off this condition by considering Strichartz estimates localized near spatial infinity. For long range perturbations of the Euclidean metric\footnote{ie, for some $ \tau > 0 $, $ \partial_x^{\alpha}(G(x)-G_{\rm Eucl}) = {\mathcal O}(\scal{x}^{-\tau-|\alpha|}) $} $ G $ on $ {\mathcal M} = \Ra^n $, trapping or not, we proved the existence of $ R > 0 $ large enough such that, if $ \chi \in C_0^{\infty}(\Ra^n) $ satisfies $ \chi \equiv 1 $ for $ |x| \leq R $, then
\begin{eqnarray}
 || (1-\chi) u ||_{L^p ([0,1];L^q(\Ra^n,dG))} \lesssim ||u_0 ||_{L^2 (\Ra^n,dG)}. \label{localalinfini}
 \end{eqnarray}
 This shows that the possible loss in Strichartz estimates can only come from a bounded region where the metric is essentially arbitrary (recall that being asymptotically Euclidean is only a condition at infinity). One can loosely interpret this result as follows: as long as the metric is close to a model one, for which one has Strichartz estimates without loss, the solution to the Schr\"odinger equation satisfies Strichartz estimates without loss too. 
 
The first goal of the present paper is to show that the same result holds in negative curvature, more precisely for asymptotically hyperbolic (AH) manifolds. Let us however point out that, even if our result (Theorem \ref{theoprincipal} below) is formally the same as in the asymptotically Euclidean case \cite[Theorem 1]{BoucletTzvetkov}, its proof involves new arguments using the negative curvature. Slightly more precisely, one of the messages of this paper is that, by taking advantage of certain curvature effects described at the end of this Section, we prove Strichartz estimates using long time (microlocal) parametrices of the Schr\"odinger group which are localized in very narrow regions of the phase space, much smaller than those considered in the asymptotically Euclidean situation \cite{BoucletTzvetkov}.

As far as the Schr\"odinger equation is concerned, Strichartz estimates on
negatively curved spaces have been studied by Banica \cite{Banica} and  Pierfelice \cite{Pier1,Pier2}
 (see \cite{Tatahype} for the wave equation). In
\cite{Pier1}, the author considers perturbations of the Schr\"odinger equation
on the hyperbolic space
 $ {\mathbb H}^n $ by singular time dependent radial potentials, with radial
 initial data (and also radial source terms) and derives some weighted
 Strichartz estimates without loss. The non radial case for the free
 Schr\"odinger equation on $ {\mathbb H}^n $ is studied in \cite{Banica} where
 weighted Strichartz estimates are obtained too. The more general case of
 certain Lie groups, namely Damek-Ricci spaces, was considered in \cite{Pier2}
 for global in time estimates (see also \cite{BanicaCarles} for the 2-dimensional case)
 and further
generalized in \cite{BanicaDuyckaerts}. In these last papers, only radial data are considered.

In this article, we give a proof of Strichartz estimates at infinity which is
purely (micro)local and therefore, to many extents, stable under
perturbation. In particular, we do not use any Lie group structure nor any
spherical symmetry. We won't assume either any non trapping condition.  We
refer to Definition \ref{defAH} below for precise statements and simply quote
here that our class of manifolds contains $ {\mathbb H}^n $, some of its
quotients and  perturbations thereof.  
In particular, we do not assume that the curvature is constant, even near infinity. 

Powerful microlocal techniques for AH manifolds have already been developed
by Melrose and its school (see \cite{MaMe1} and the references in \cite{Melr0}). These geometric methods, based on compactification and  blowup considerations, are perfectly designed for conformally compact manifolds (with boundary) but do not clearly apply to the more general manifolds we shall study here.

Let us finally  mention that Theorem \ref{theoprincipal} reduces the proof of potential improvements of Burq-G\'erard-Tzvetkov inequalities to local in space estimates of the form 
$$ || \chi u ||_{ L^p ([0,1],L^q ({\mathcal M},dG))} \lesssim
||u_0||_{H^s({\mathcal M},dG)} ,  $$
with $ 0 \leq s < 1/p $. It would be interesting to know if such inequalities holds for some trapping AH manifolds. 

\medskip

Before stating and discussing precisely our result, let us fix the framework.

\bigskip

\begin{defi}[AH manifold] \label{defAH}  $ ({\mathcal M}^n,G )  $ is asymptotically hyperbolic if there exist a compact set $ {\mathcal K}
 \Subset {\mathcal M} $, a real number $ R_{\mathcal K} > 0 $, a compact manifold without boundary $ S $ and a function
\begin{eqnarray}
r \in C^{\infty}({\mathcal M},\Ra) , \qquad r (m) \rightarrow + \infty , \qquad m \rightarrow \infty, \label{rglobal}
\end{eqnarray}
which is a coordinate near $ \overline{ {\mathcal M} \setminus {\mathcal K} } $  such that: we have an isometry
\begin{eqnarray}
\diffeo : ( {\mathcal M} \setminus {\mathcal K} , G ) \rightarrow  \left( (R_{\mathcal K},+\infty)_r \times S , d r^2 + e^{2r} g (r) \right) , \label{diffeomorphisminfini}
\end{eqnarray}
 where $ g(r) $ is a family of metrics on $ S $ depending smoothly on $ r  $ such that, for some $ \tau > 0 $
 and some fixed metric $ \ginfty $ on $ S $, we have
\begin{eqnarray}
|| \partial_r^k \left( g (r) - \ginfty \right) ||_{C^{\infty}(S
,T^* S \otimes T^* S)} \lesssim r^{- \tau - k}, \qquad r >
R_{\mathcal K},  \label{longueportee}
\end{eqnarray}
for all $ k \geq 0 $ and all semi-norms $ || \cdot ||_{C^{\infty}(S ,T^* S \otimes T^* S)} $ in the space of smooth sections of $ T^* S \otimes T^* S $.  \end{defi}

With no loss of generality, we can assume that the decay rate $ \tau $ in (\ref{longueportee}) satisfies
\begin{eqnarray}
0 < \tau < 1 . \label{tauxlongueportee}
\end{eqnarray}
Therefore, by analogy with the standard terminology in Euclidean
scattering, $ d r^2 + e^{2r} g(r) $ can be considered as a  {\it
long range} perturbation of the metric $ d r^2 + e^{2r} \ginfty $. Notice that the conformally compact case quoted above corresponds to the special situation where $ g(r) $ is of the form $ \tilde{g}(e^{-r}) $, for some family of metrics $ (\tilde{g}(x))_{0 \leq x \ll 1} $
depending smoothly on $ x \in [0, x_0 ) $ ($ x_0 $ small enough) up to $ x = 0 $. In that case, $ g (r) $ is an exponentially small perturbation of  $ g = \tilde{g}(0) $. The assumption (\ref{longueportee}) is therefore more general. 

\medskip

We next denote by $ \Delta_{G} $ the Laplace-Beltrami operator associated to this metric.
It is classical that this operator is essentially self-adjoint on $ L^2 ({\mathcal M},dG) $, from $ C_0^{\infty} ({\mathcal M}) $,  and therefore generates a unitary group $ e^{i t \Delta_G} $.

Our main result is the following.

\begin{theo} \label{theoprincipal} There exists $ \chi \in C_0^{\infty}({\mathcal M}) $, with $  \chi \equiv 1 $ on a
sufficiently large compact set, such that, for all  pair $ (p,q) $
satisfying (\ref{paireadmissible}),
\begin{eqnarray}
 || (1-\chi) e^{it\Delta_G} u_0 ||_{L^p ([0,1];L^q({\mathcal M},dG))} \lesssim ||u_0||_{L^2({\mathcal M},dG)}, \qquad u_0 \in C_0^{\infty}({\mathcal M}) . \label{vraiStrichartz}
\end{eqnarray}
\end{theo}


This theorem  is the AH analogue of Theorem 1 of \cite{BoucletTzvetkov} in the asymptotically Euclidean case. 

To be more complete, let us point out that the analysis contained in this paper and a classical argument due to \cite{StTa} (see also 
\cite[Section 5]{BoucletTzvetkov}), using the local smoothing effect \cite{Doi}, would give
 the following global in space estimates.

\begin{theo} If in addition $ ({\mathcal M},G) $ is non trapping, then we have global in space Strichartz estimates with no loss: 
for all  pair $ (p,q) $
satisfying (\ref{paireadmissible}),
\begin{eqnarray}
 ||  e^{it\Delta_G} u_0 ||_{L^p ([0,1];L^q({\mathcal M},dG))} \lesssim ||u_0||_{L^2({\mathcal M},dG)}, \qquad u_0 \in C_0^{\infty}({\mathcal M}) . 
 \nonumber
\end{eqnarray}
\end{theo}

We state this result as a theorem although we won't explicitly prove it. The techniques are fairly well known and don't involve any new argument in the present context. We simply note that resolvent estimates implying the local smoothing effect can be found in \cite{CardosoVodev}.

\bigskip

We now describe, quite informally, the key points of the analysis developed
in this paper. Assuming for simplicity that $ S = {\mathbb S}^1 $ we consider
the model case where the principal symbol of the Laplacian is
$$  p = \rho^2 + e^{-2r} \eta^2 .  $$
For convenience, we introduce
$$ P := - e^{(n-1)r/2} \Delta_G e^{-(n-1)r/2} $$ 
which is  self-adjoint with respect to $ dr d \theta $, instead of $
e^{(n-1)r}dr d\theta $ for the Laplacian itself. 

Recall first that, by the Keel-Tao $ TT^* $ Theorem \cite{KeelTao}, proving
Strichartz estimates (without loss) is 
mainly reduced to prove certain dispersion estimates.
Using the natural semi-classical time scaling $ t \mapsto h t $,  this basically requires to control the propagator $ e^{-ithP} $ for semi-classical times of order $ h^{-1} $. Such a control on the full propagator
is out of reach (basically because of trapped trajectories) but, fortunately, studying some of its cutoffs will be sufficient.

After fairly classical reductions, we will work with semi-classical
pseudo-differential operators localized where $ r \gg 1 $ and $ p \in I $, $ I
$ being a (relatively) compact interval of $ (0,+\infty) $. We can split the latter region
 into two areas defined by
$$ \Gamma^+ = \left\{  r \gg 1 , \ p \in I , \  \rho > -\frac{ p^{1/2}}{2} \right\} , \qquad \Gamma^- = 
\left\{  r \gg 1 , \ p \in 1 , \  \rho < \frac{ p^{1/2} }{2} \right\} , $$
respectively called outgoing and incoming areas. In the asymptotically Euclidean case, it turns out that one can give accurate approximations of $ e^{-ithP} \chi^{\pm} $ for times $t$ such that $ 0 \leq \pm t \lesssim h^{-1} $, if $ \chi^{\pm} $
are pseudo-differential cutoffs localized in $ \Gamma^{\pm} $. This is not the AH case: here we are only able to approximate 
$ e^{-ithP} \chi^{\pm}_{\rm s} $
for cutoffs $ \chi^{\pm}_{\rm s} $ localized in much smaller areas, namely
$$  \Gamma^+_{\rm s} (\epsilon) = \left\{  r \gg 1 , \ p \in I , \  \frac{\rho}{  p^{1/2} } > 1-\epsilon^2 \right\} , \qquad \Gamma^-_{\rm s} (\epsilon) = 
\left\{  r \gg 1 , \ p \in I , \   \frac{\rho}{ - p^{1/2} } < \epsilon^2 -1 \right\} , $$
which we call {\it strongly} outgoing/incoming areas. Here $ \epsilon $ will be a fixed small real number.
 We then obtain approximations of the form
\begin{eqnarray}
  e^{-ith P} \chi^{\pm}_{\rm s} = J_{S^{\pm}}(a^{\pm}) e^{-ithD_r^2} J_{S^{\pm}}(b^{\pm})^* + {\mathcal O}(h^N), \qquad 0 \leq \pm t \lesssim h^{-1} . 
  \label{IsozakiKitadainformel} 
\end{eqnarray} 
Here $ e^{-ithD_r^2} $ is the semi-classical group associated to the radial part $ D_r^2 $ of $ P $. $  J_{S^{\pm}}(a^{\pm})$ and $ J_{S^{\pm}}(b^{\pm}) $
are Fourier integral operators with phase essentially of the form
$$ S^{\pm} \approx r \rho + \theta \eta + \frac{ e^{-2r}\eta^2 }{4 \rho} , $$
that is the sum of the {\it free phase} $ r \rho + \theta \eta $ and of a term whose Hessian is non degenerate in $ \eta $, which will be crucial for the final stationary phase argument (the small factor $ e^{-2r} $  will be eliminated by a change of variable). The non degeneracy of the full phase of the parametrix (\ref{IsozakiKitadainformel}) in $\rho$ will come of course from $ e^{-ithD_r^2} $.

The approximation (\ref{IsozakiKitadainformel}) is the AH Isozaki-Kitada parametrix and it will be used very similarly to the usual Euclidean one as in \cite{BoucletTzvetkov}. Its main interest is to give microlocal approximations of the propagator for times of size $ h^{-1} $. Recall however the big difference with the 
asymptotically Euclidean case where one is able to consider cutoffs supported in $ \Gamma^{\pm} $ rather than $ \Gamma^{\pm}_{\rm s}(\epsilon) $ in the AH case. We therefore have to consider the left parts, namely
$$ \Gamma^{\pm}_{\rm inter} = \Gamma^{\pm} \setminus \Gamma^{\pm}_{\rm s} (\epsilon) , $$
which we call {\it intermediate} areas. These areas  will only contribute to the dispersion  estimates for small times using the following argument.
By choosing $ \delta $ small enough and by splitting the interval $ (-1/2,1-\epsilon^2) $ into small intervals of size $ \delta $, we can write 
$$  \Gamma^{\pm}_{\rm inter} = \cup_{l \lesssim \delta^{-1}}  \left\{  r \gg 1 , \ p \in I , \  \pm \frac{\rho}{ p^{1/2}} \in (\sigma_l , \sigma_l + \delta) \right\}
=  \cup_{l \lesssim \delta^{-1}} \Gamma^{\pm}_{\rm inter} (l,\epsilon,\delta) . $$
By looking carefully at the Hamiltonian flow $ \Phi_p^t $ of $ p $, it turns out that, for any fixed (small) time $ t_0 $, we can choose $ \delta $ (which depends also on $ \epsilon $) such that
\begin{eqnarray}
 \Phi_p^t \left( \Gamma^{\pm}_{\rm inter} (l,\epsilon,\delta) \right) \cap \Gamma^{\pm}_{\rm inter} (l,\epsilon,\delta)
 = \emptyset, \qquad \pm t \geq t_0 . \label{effetdecourbure}
\end{eqnarray}
By semi-classical propagation, this implies that, for pseudo-differential operators $ \chi_{\rm inter}^{\pm} $ localized in $ \Gamma^{\pm}_{\rm inter} (l,\epsilon,\delta) $,
$$ \chi_{\rm inter}^{\pm} e^{-ithP}\chi_{\rm inter}^{\pm * } = {\mathcal O}(h^{\infty}), \qquad \pm t \geq t_0 . $$
Such operators typically appear in the  $ T T^* $ argument and the estimate above reduces the proof of dispersion estimates to times $ |t| \leq t_0 $. The latter range of times can then be treated by fairly standard geometric optics approximation.

We interpret (\ref{effetdecourbure}) as a {\it negative curvature effect} on the geodesic flow which we can roughly describe as follows, say in the outgoing case. For initial conditions $ (r,\theta,\rho,\eta) $ in $ \Gamma^{+}_{\rm inter} (l,\epsilon,\delta) $, the  bounds $ 1/2 < \rho \leq (1 - \epsilon^2) p^{1/2} $ yields the following lower bound,  
$$ \dot{\rho}^t = 2 e^{-2r^t}(\eta^t)^2 \gtrsim \epsilon^2 , $$
over a sufficiently long  time, if we set $  (r^t , \theta^t , \rho^t , \eta^t):= \Phi_p^t $. This ensures that $ \rho^t / p^{1/2} $ increases fast enough to leave the interval $ (\sigma_l , \sigma_l + \delta) $ before $ t = t_0 $ and give
(\ref{effetdecourbure}). In the asymptotically flat case, ie with $ r^{-2} $ instead of $ e^{-2r} $, we have
$ \dot{\rho}^2 = 2 (r^t)^{-3}(\eta^t)^2  $ and its control from below is not as good, basically because of the `extra' third power of $ (r^t)^{-1}  $. 
\bigskip

This paper is organized as follows. 

\medskip

In Section \ref{ideepreuvetheoreme}, we introduce all the necessary definitions, and some additional results, needed to
prove Theorem \ref{theoprincipal}. The latter proof is given in Subsection \ref{soussection4} using microlocal approximations which will be proved in Sections \ref{IsozakiKitada}, \ref{WKBsection} and \ref{sectiondispersion}.

In Section \ref{sectionflot}, we study the properties of the geodesic flow in outgoing/incoming areas required to construct the phases involved in the Isozaki-Kitada parametrix. This parametrix is then constructed in Section \ref{IsozakiKitada}.

In Section \ref{WKBsection} we prove two results:  the small semi-classical time approximation of the Schr\"odinger group
by the WKB method and  the propagation of the microlocal support (Egorov theorem). These results are essentially well known. We need however to check that all the symbols and phases belong to the natural classes (for AH geometry) of Definition \ref{definhyp}  below.   Furthermore, we use our Egorov theorem to obtain a propagation property in a time scale of size $ h^{-1} $ which is not quite standard.

Finally, in Section \ref{sectiondispersion}, we prove dispersion estimates using basically stationary phase estimates
in the parametrices obtained in  Sections \ref{IsozakiKitada} and \ref{WKBsection}.

\medskip

Up to the semi-classical functional calculus, which is taken from \cite{BoucletLpCF,BoucletLP} and whose results are recalled in subsection \ref{soussection2}, this paper is essentially self contained. This is not only for the reader's convenience, but also due to the fact that the results of Section \ref{WKBsection} 
do require proofs in the AH setting, although they are in principle well known. The construction of Section \ref{IsozakiKitada} is new.

\section{The strategy of the proof of Theorem \ref{theoprincipal}} \label{ideepreuvetheoreme}
\setcounter{equation}{0}





\subsection{The setup}

Before discussing the proof of Theorem \ref{theoprincipal}, we give the form of the Laplacian, volume densities and related
objects on AH manifolds.

The isometry (\ref{rglobal}) defines polar coordinates: $r$ is the radial coordinate and $S$ will be called the
 angular manifold. Coordinates on $S$ will be denoted by $ \theta_1 , \ldots , \theta_{n-1} $.

A finite atlas on $ {\mathcal M} \setminus {\mathcal K} $ is
obtained as follows. By (\ref{diffeomorphisminfini}), we have a
natural `projection' $ \pi_S : ( {\mathcal M} \setminus {\mathcal
K} , G ) \rightarrow S $ defined as the second component of $ \Psi $, ie
\begin{eqnarray}
 \Psi (m) = (r(m), \pi_{S}(m)) \in (R_{\mathcal K},+\infty) \times S, \qquad m \in {\mathcal M} \setminus {\mathcal K} .
 \label{projectionS}
\end{eqnarray}
 Choosing a finite cover of the angular manifold by
coordinate patches $ U_{\iota} $, ie
\begin{eqnarray}
 S = \cup_{\iota \in {\mathcal I}} U_\iota , \label{cover}
\end{eqnarray}
with corresponding diffeomorphisms
\begin{eqnarray}
 \psi_{\iota} : U_{\iota} \rightarrow \psi_{\iota} (U_{\iota}) \subset \Ra^{n-1} , \label{diffeocover}
\end{eqnarray}
we consider the open sets
$$ {\mathcal U}_{\iota} := \diffeo^{-1} \left( (R_{\mathcal K},+\infty) \times U_{\iota} \right) \subset {\mathcal M}
\setminus {\mathcal K } , $$
and then define diffeomorphisms
\begin{eqnarray}
 \Psi_{\iota} : {\mathcal U}_{\iota} \rightarrow (R_{\mathcal K},+\infty) \times \psi_{\iota}(U_{\iota}) \subset \Ra^n , \label{Diffeoaveciota} 
\end{eqnarray}
by
$$ \Psi_{\iota}(m) = \left( r(m) , \psi_{\iota} (\pi_S(m)) \right) . $$
The collection $ ( {\mathcal U}_{\iota} , \Psi_{\iota} )_{\iota \in {\mathcal I}} $ is then an atlas on $ {\mathcal M} \setminus {\mathcal K} $. If $ \theta_1 , \ldots , \theta_{n-1} $ are the coordinates in $ U_{\iota}$, ie $ \psi_{\iota} = (\theta_1 , \ldots , \theta_{n-1}) $, the coordinates in $ {\mathcal U}_{\iota} $ are then $ (r,\theta_1 , \ldots , \theta_{n-1}) $.

We now give formulas for  the riemannian measure $ dG $ and the Laplacian $ \Delta_G $ on $ {\mathcal M} \setminus {\mathcal K} $.
In local coordinates $ \theta=(\theta_1 , \ldots , \theta_{n-1}) $ on $ S $, the riemannian density associated to $ g (r) $ reads
$$ d g(r) := \mbox{det} \left( g (r,\theta) \right)^{1/2} |d \theta_1 \wedge \cdots \wedge d \theta_{n-1} | , $$
where $ \mbox{det} ( g (r,\theta) ) = \mbox{det} (g_{jk} (r,\theta)) $ if $ g (r) = g_{jk} (r,\theta) d \theta_j d \theta_k $ (using the summation convention). Then, in local coordinates on $ {\mathcal M} \setminus {\mathcal K} $, the riemannian density is
\begin{eqnarray}
d G = e^{(n-1)r} \mbox{det} \left( g (r,\theta) \right)^{1/2} |d r \wedge d \theta_1 \wedge \cdots \wedge d \theta_{n-1} | . \label{formededG}
\end{eqnarray}

Let us now consider the Laplacian. Slightly abusing the notation, we set
\begin{eqnarray}
 \curv (r,s) = \frac{1}{2} \frac{\partial_r \mbox{det}( g(r,s) )}{\mbox{det}( g(r,s) )}, \qquad r > R_{\mathcal K}, \ s \in  S . \label{courbure}
\end{eqnarray}
This means, for fixed $r$, that the quotient of $ \partial_r \mbox{det} ( g_{jk}(r,\theta) ) $
by $ 2 \mbox{det} ( g_{jk}(r,\theta) ) $ defines a function on $S$ (ie is independent of the coordinates). We then have
\begin{eqnarray}
\Delta_{G} =  \partial_r^2  + e^{-2r} \Delta_{g(r)} +  \curv (r,s) \partial_r + (n-1) \partial_r  . \nonumber
\end{eqnarray}

\bigskip

It will turn out be convenient to work with the following density
\begin{eqnarray}
\widehat{dG} = e^{(1-n)r} d G ,
\end{eqnarray}
rather than $ d G $ itself. In particular, we will use the following elementary property: for all relatively compact subset
$ V_{\iota}^{\prime} \Subset \psi_{\iota}(U_{\iota}) $, all $ R > R_{\mathcal K} $ and all $ 1 \leq q \leq \infty $, we have the equivalence of norms
\begin{eqnarray}
||u||_{L^q ({\mathcal M}, \widehat{dG})} \approx || u \circ \Psi_{\iota}^{-1} ||_{L^q (\Ra^n)
 } , \qquad \mbox{supp}(u) \subset (R,+\infty) \times V_{\iota}^{\prime} , \label{equivalenceLp}
\end{eqnarray}
 $ L^q( \Ra^n ) $ being the usual Lebesgue space. This is  a simple consequence of (\ref{longueportee}) and (\ref{formededG}) (we consider $ R > R_{\mathcal K} $ since (\ref{longueportee}) gives an upper bound for $ \mbox{det} \left( g (r,\theta) \right) $
as $ r \rightarrow R_{\mathcal K} $ , not a lower bound).

 We then have a unitary isomorphism
\begin{eqnarray}
L^2 ({\mathcal M},\widehat{dG}) \ni u \rightarrow e^{- \frac{(n-1)r}{2}} u \in L^2 ({\mathcal M},dG) , \label{unitaire}
\end{eqnarray}
and $ \Delta_G $ is  unitarily equivalent to the operator
\begin{eqnarray}
 \widehat{\Delta}_G := e^{ \bott r} \Delta_G e^{- \bott r} , \qquad \bott = \frac{n-1}{2}, \label{defDeltachapeau}
\end{eqnarray}
on $ L^2 ({\mathcal M},\widehat{dG}) $. This operator reads
\begin{eqnarray}
\widehat{\Delta}_G = \partial_r^2 + e^{-2r} \Delta_{g(r)} + \curv(r,\theta) \partial_r - \bott \curv (r,\theta) - \bott^2 , \label{coefficientsconcrets}
\end{eqnarray}
and we will work the following one
\begin{eqnarray}
P = - \widehat{\Delta}_G - \bott^2 . \label{defP}
\end{eqnarray}
If $ q_{\iota} (r,.,.) $ is the principal symbol of $ - \Delta_{g(r)} $ in the chart $ U_{\iota} $, namely
\begin{eqnarray}
  q_{\iota} (r,\theta,\xi) = \sum_{1 \leq k,l \leq n-1} g^{kl}(r,\theta) \xi_k \xi_l , \label{notationLaplacieninfini}
\end{eqnarray}
the principal symbol of $ P $ in the chart $ {\mathcal U}_{\iota}  $ is then
\begin{eqnarray}
 p_{\iota} & = &   \rho^2 + e^{-2r} q_{\iota} (r,\theta,\eta)  , \label{symbolprincipaliota} \\
  & = & \rho^2 + q_{\iota}(r,\theta,e^{-r}\eta) , \nonumber
\end{eqnarray}
The full symbol of $ P $ is of the form $ p_{\iota} + p_{\iota,1} + p_{\iota,0} $ with
\begin{eqnarray}
p_{\iota,j} & = & \sum_{k + |\beta| = j} a_{\iota,k \beta}(r,\theta) \rho^k (e^{-r}\eta)^{\beta} , \qquad j = 0,1 . \label{termesunzero}
\end{eqnarray}
The terms of degree $1$ in $ \eta $ come from the first order terms of the symbol of $ - \Delta_{g(r)} $. In the expression of $ \Delta_G $ they carry a factor $ e^{-2r} $ and therefore, if $ j=1$, $ k = 0 $ and $ | \beta | = 1 $ above, we could write
$ a_{\iota,k \beta}(r,\theta) = e^{-r} b_{\iota,k \beta}(r,\theta)  $ for some  function $b_{\iota,k \beta}$ bounded as $r \rightarrow \infty $. This remark and (\ref{longueportee}) show more precisely that, for all $ V \Subset \psi_{\iota} (U_{\iota})  $, 
the coefficients in (\ref{termesunzero}) decay as 
\begin{eqnarray}
| \partial_r^j \partial_{\theta}^{\alpha} a_{\iota, k \beta}  (r,\theta) | \leq C_{Vj\alpha} \scal{r}^{-\tau-1-j}, \qquad \theta \in V , \ r \geq R_{\mathcal K} + 1 . \label{borneshortrange}
\end{eqnarray}
 The decay rate $ -\tau - 1 -j  $  will be important to solve  transport equations for the Isozaki-Kitada parametrix. This is the main reason of the long range assumption (\ref{longueportee}).

\subsection{Pseudo-differential operators and the spaces $ {\mathcal B}_{\rm hyp} (\Omega) $} \label{soussection1}

We will consider $h$-pseudo-differential operators ($h$-$\Psi$Dos)  in a neighborhood of infinity and the `calculus' will be rather elementary. For instance, we will only consider compositions of operators with symbols supported in the same coordinate patch and no invariance result under diffeomorphism will be necessary.

The first step is to construct a suitable partition of unity near infinity.
Using the cover (\ref{cover}) and the related diffeomorphisms (\ref{diffeocover}), we consider a partition of unity on $ S $ of the form
\begin{eqnarray}
 \sum_{\iota \in {\mathcal I}} \kappa_{\iota} \circ \psi_{\iota} = 1 , \qquad \mbox{with} \qquad
 \kappa_{\iota} \in C_0^{\infty}(\Ra^{n-1}) , \ \ \ \mbox{supp}( \kappa_{\iota} ) \Subset \psi_{\iota} (U_{\iota}) , \label{partitionangulaire}
\end{eqnarray}
and a function $ \cutoff \in C^{\infty}(\Ra) $ such that
\begin{eqnarray}
 \mbox{supp}(\cutoff) \subset [ R_{\mathcal K}+1,+\infty), \qquad \cutoff
 \equiv 1 \ \ \ \mbox{on} \ \ [ R_{\mathcal K}+2, +  \infty) . \label{cutoffradial}
 \end{eqnarray}
Then, the functions $ (\kappa \otimes \kappa_{\iota}) \circ \Psi_{\iota} \in C^{\infty}({\mathcal M}) $ satisfy
\begin{eqnarray}
\sum_{\iota \in {\mathcal I}} (\kappa \otimes \kappa_{\iota}) \circ \Psi_{\iota} (m) = \begin{cases} 1 & \mbox{if} \ \  r (m) \geq R_{\mathcal K} + 2 , \\
0 & \mbox{if} \ \  r (m) \leq R_{\mathcal K} + 1 ,
\end{cases}
 \label{presquepartition}
\end{eqnarray}
which means that they define a partition of unity near infinity. We could obtain a partition of unity on ${\mathcal M}$ by adding a finite number of compactly supported functions (in coordinate patches) be we won't need it since the whole analysis in this paper will be localized near infinity.

We also consider $ \widetilde{\cutoff} \in C^{\infty}(\Ra) $ and $ \widetilde{\kappa}_{\iota} \in C_0^{\infty}(\Ra^{n-1}) $, for all $ \iota \in {\mathcal I} $, with the following properties
\begin{eqnarray}
 \widetilde{\cutoff} \equiv 1 \ \ \ \mbox{on} \ \ (R_{\mathcal K}+1/2,+ \infty), \qquad  \widetilde{\kappa}_{\iota} \equiv 1 \ \ \ \mbox{near} \ \
 \mbox{supp}(\kappa_{\iota})
, \\
 \mbox{supp}(\widetilde{\cutoff}) \subset (R_{\mathcal K}+1/4,+\infty), \qquad \mbox{supp}(\widetilde{\kappa}_{\iota}) \Subset \psi_{\iota}(U_{\iota})   . \label{pourLp}
\end{eqnarray}
We next choose, for each $ \iota \in {\mathcal I} $, two relatively compact open subsets $ V_{\iota} $ and $ V_{\iota}^{\prime} $ such that
\begin{eqnarray}
\mbox{supp}(\kappa_{\iota}) \Subset V_{\iota} \Subset V_{\iota}^{\prime} \Subset \mbox{supp}(\widetilde{\kappa}_{\iota}) \qquad \mbox{and} \qquad \widetilde{\kappa}_{\iota} \equiv 1 \ \ \ \mbox{near} \ \
 V_{\iota}^{\prime} . \label{margeouverts}
\end{eqnarray}

We are now ready to define our $\Psi $DOs. In the following definition, we will say that $ a \in C^{\infty}(\Ra^{2n}) $
is a symbol if either $ a \in C^{\infty}_b (\Ra^{2n}) $, ie bounded with all derivatives bounded, or
$$ a (r,\theta,\rho,\eta) = \sum  a_{k \beta} (r,\theta) \rho^k \eta^{\beta} , $$
with $  a_{k \beta} \in C^{\infty}_b (\Ra^n) $, the sum being finite. We shall give examples below. Notice that throughout  this paper, $ \rho $ and $ \eta $ will denote respectively the dual variables to $ r $ and $ \theta $.
\begin{defi} \label{definitOpchapeau} For $ \iota \in {\mathcal I} $, all $ h \in (0,1] $ and all symbol $ a $ such that
\begin{eqnarray}
 \emph{supp}(a) \subset [ R_{\mathcal K}+1,+\infty) \times V^{\prime}_{\iota} \times \Ra^n , \label{symbolpseudo}
\end{eqnarray}
we define
$$ \widehat{O \! p }_{\iota}(a) : C_0^{\infty}({\mathcal M}) \rightarrow C^{\infty}({\mathcal M}), $$
by
\begin{eqnarray}
 \left( \widehat{O \! p }_{\iota}(a) u \right) \circ \Psi_{\iota}^{-1} (r,\theta) =  a (r,\theta,h D_r , h D_{\theta})\left(  \widetilde{\cutoff}(r)  \widetilde{\kappa}_{\iota} (\theta) ( u \circ \Psi_{\iota}^{-1})(r,\theta) \right) . \label{defpseudo}
\end{eqnarray}
\end{defi}
Note the cutoff $  \widetilde{\kappa} \otimes \widetilde{\kappa}_{\iota} $ in the right hand side of (\ref{defpseudo}). It makes the Schwartz kernel of $ \widehat{O \! p }_{\iota}(a)  $ supported in a closed subset of $ {\mathcal M}^2 $ strictly contained in the patch $ {\mathcal U}_{\iota}^2  $ so that $ \widehat{O \! p }_{\iota}(a) $ is fully defined by the prescription of
$  \Psi_{\iota *} \widehat{O \! p }_{\iota}(a) \Psi_{\iota}^* $. For future reference, we recall that the kernel of the latter operator is
\begin{eqnarray}
 (2 \pi h)^{-n} \int \! \! \int e^{\frac{i}{h} (r-r^{\prime})\rho
 + \frac{i}{h} (\theta - \theta^{\prime})\cdot \eta} a (r,\theta, \rho , \eta) d \rho d \eta \widetilde{\chi} (r^{\prime}) \widetilde{\chi}_{\iota}(\theta^{\prime})  . \label{noyauiota}
\end{eqnarray}

The notation $ \widehat{ O \! p }_{\iota} $ refers to the following relation with the measure $ \widehat{dG} $: if $ a \in C^{\infty}_b (\Ra^{2n}) $ satisfies (\ref{symbolpseudo}), then
\begin{eqnarray}
||  \widehat{O \! p }_{\iota}(a) ||_{L^2 ({\mathcal M},\widehat{dG}) \rightarrow L^2 ({\mathcal M} , \widehat{dG})} \lesssim 1, \qquad h \in (0,1] .
\label{borneL2gratuite}
\end{eqnarray}
This is a direct consequence of the Calder\'on-Vaillancourt theorem using (\ref{equivalenceLp}) with $q=2$, (\ref{pourLp}) and (\ref{margeouverts}). In the `gauge' defined by $ d G $, the latter gives
\begin{eqnarray}
 ||  e^{ - \bott r  } \widehat{O \! p }_{\iota}(a) e^{ \bott r } ||_{L^2 ({\mathcal M}, dG ) \rightarrow L^2 ({\mathcal M} , dG)} \lesssim 1, \qquad h \in (0,1] . \label{borneL2apoids}
\end{eqnarray}
Working with the measure $ \widehat{dG} $ is therefore more transparent and avoid to deal with exponential weights.

\medskip

Let us now describe the typical symbols  
we shall use in this paper. Using (\ref{partitionangulaire}), (\ref{cutoffradial}), (\ref{presquepartition}) and  (\ref{margeouverts}), we can write
\begin{eqnarray}
 h^2 P = \sum_{\iota \in {\mathcal I}} \widehat{O \! p }_{\iota} \left( ( \kappa \otimes \kappa_{\iota} ) \times (  p_{\iota} + h p_{\iota,1} + h^2 p_{\iota,0} ) \right) , \qquad r > R_{\mathcal K} + 2 , \label{partitionP}
\end{eqnarray}
using (\ref{notationLaplacieninfini}), (\ref{symbolprincipaliota}) and (\ref{termesunzero}). 
One observes that the symbols  involved in (\ref{partitionP}) are of the form
\begin{eqnarray}
a_{\iota}(r,\theta,\rho,\eta) = \tilde{a}_{\iota}(r,\theta,\rho,e^{-r}\eta), \label{formsymbol}
\end{eqnarray}
with $ \tilde{a}_{\iota} \in S^{2} (\Ra^n \times \Ra^n) $.
It will turn out that the functional calculus of $h^2 P$ (or $ h^2 \Delta_G $) will involve more generally symbols of this form with $ a_{\iota} \in S^{-\infty}(\Ra^n \times \Ra^n) $. For instance, if $ f \in C_0^{\infty}(\Ra) $, the semi-classical principal symbol of $ f (h^2 P) $ or $ f (-h^2 \Delta_{G}) $ will be
\begin{eqnarray}
 f (\rho^2 + q_{\iota} (r,\theta,e^{-r}\eta)) , \label{symboleprincipalcalculfonctionnel}
\end{eqnarray}
which, once multiplied by the cutoff $ \kappa \otimes \kappa_{\iota} $, is of the form (\ref{formsymbol}) with $ \tilde{a}_{\iota} \in S^{- \infty}(\Ra^n \times \Ra^n) $. This type of symbols is the model of functions described in  Definition \ref{definhyp} below.
To state this definition, we introduce the notation
$$ D_{\rm hyp}^{j \alpha k \beta} := e^{r|\beta|}
\partial_{\eta}^{\beta} \partial_r^j \partial_{\theta}^{\alpha} \partial_{\rho}^k, $$
for all $ j,k \in \Na $ and $ \alpha,\beta \in \Na^{n-1} . $
\begin{defi} \label{definhyp} Given an open set $ \Omega \subset  T^* \Ra^n_+ = (0,+\infty)_r \times \Ra^{n-1}_{\theta} \times \Ra_{\rho} \times \Ra^{n-1}_{\theta} $, we define
$$ {\mathcal B}_{\rm hyp} (\Omega) = \{ a \in C^{\infty}(\Omega) \ | \ D_{\rm hyp}^{j \alpha k \beta} a \in L^{\infty}(
\Omega), \ \ \mbox{for all} \ j,k \in \Na, \ \alpha,\beta \in
\Na^{n-1} \} ,
$$
and
$$ {\mathcal S}_{\rm hyp}(\Omega) = \{ a \in C^{\infty}(\Ra^{2n}) \ | \ \emph{supp}(a) \subset \Omega \ \mbox{and} \ a \in {\mathcal B}_{\rm hyp}(\Omega) \} . $$
 A family $ (a_{\nu})_{\nu \in  \Lambda} $ is  bounded in $
{\mathcal B}_{\rm hyp}(\Omega) $ if, for all $ j,k,\alpha,\beta $,
 $ ( D_{\rm hyp}^{j \alpha k \beta} a_{\nu} )_{\nu \in \Lambda} $
is bounded in $ L^{\infty}(\Omega) $.
\end{defi}

Note that considering $ \Omega \subset  T^* \Ra^n_+ $ is not necessary but, since we shall work only in the region where $ r \gg 1 $, this will be sufficient.

\begin{exam} \label{exemple} Consider the following  diffeomorphism from $ \Ra^{2n} $ onto itself
\begin{eqnarray}
 F_{\rm hyp} : (r,\theta,\rho,\eta) \mapsto (r,\theta,\rho,e^{-r} \eta) . \label{Fhyp}
\end{eqnarray}
If $ a_{\iota} \in S^{0}(\Ra^n \times \Ra^n) $ is supported in $ F_{\rm hyp} (\Omega) $, with $ \Omega \subset T^* \Ra^{n}_+ $,
 then (\ref{formsymbol}) belongs to $ {\mathcal S}_{\rm hyp}( \Omega ) $.
\end{exam}

\noindent {\it Proof.} We only need to check that (\ref{formsymbol}) belongs to $ {\mathcal B}_{\rm hyp}( \Omega ) $. We have
\begin{eqnarray*}
\partial_r \left( \tilde{a}_{\iota}(r,\theta,\rho,e^{-r}\eta) \right) & = &  (\partial_r \tilde{a}_{\iota})(r,\theta,\rho,e^{-r}\eta) +
e^{-r} \eta \cdot (\partial_{\xi } \tilde{a}_{\iota})(r,\theta,\rho,\xi)_{| \xi = e^{-r}\eta},
\end{eqnarray*}
which is bounded since $ \xi \cdot \partial_{\xi} a_{\iota} $ is bounded. Similarly
\begin{eqnarray*}
e^r \partial_{\eta  } \left( \tilde{a}_{\iota}(r,\theta,\rho,e^{-r}\eta) \right) & = &
  (\partial_{\xi } \tilde{a}_{\iota})(r,\theta,\rho,\xi)_{| \xi = e^{-r}\eta} ,
\end{eqnarray*}
is bounded too. Derivatives with respect to $ \rho, \theta $ are harmless and higher order derivatives in $r , \eta $ are treated similarly. \finpreuve

\bigskip

In the following lemma, we give a characterization of functions in $ {\mathcal B}_{\rm hyp}(\Omega) $.
\begin{lemm} \label{lemmesymbole} 
Let $ \Omega \subset T^* \Ra^n_+  $ be an open subset and assume that
\begin{eqnarray}
 F_{\rm hyp} (\Omega) \subset \Ra^n_+ \times B, \qquad \mbox{with} \ \ B \ \ \mbox{bounded} . \label{ensembleborne}
\end{eqnarray}
Then,  a function $ a \in C^{\infty} (\Omega)$ is of the form
\begin{eqnarray}
 a (r,\theta,\rho,\eta) = \tilde{a}(r,\theta,\rho,e^{-r}\eta) \qquad \mbox{with} \ \ \tilde{a} \in C^{\infty}_b \left(
 F_{\rm hyp}(\Omega) \right) , \label{point1}
\end{eqnarray}
if and only if, for all $  j,k,\alpha,\beta$,
\begin{eqnarray}
 D_{\rm hyp}^{j \alpha k \beta} a \in L^{\infty}(\Omega) . \label{point2}
\end{eqnarray}
Here $ C^{\infty}_b (\Omega) $ (resp. $ C^{\infty}_b(
 F_{\rm hyp}(\Omega) )$) is the space of smooth functions bounded with all derivatives bounded on $ \Omega $ (resp. $ F_{\rm hyp}(\Omega) $).
\end{lemm}

\noindent {\it Proof.} That (\ref{point1}) implies (\ref{point2}) is proved in the same way as Example \ref{exemple}: the boundedness of $ \xi \cdot \partial_{\xi} \tilde{a} $ follows from the boundedness of $ \xi = e^{-r} \eta $ in $ F_{\rm hyp} (\Omega) $ by (\ref{ensembleborne}) and the fact that $ \tilde{a} \in C^{\infty}_b \left( F_{\rm hyp}(\Omega) \right) $. Conversely, one checks by induction that
$$ \tilde{a}(r,\theta,\rho,\xi) := a (r,\theta,\rho,e^{r} \xi) , $$
belongs to $ C^{\infty}_b \left( F_{\rm hyp}(\Omega) \right) $, using again the boundedness of $\xi$ on $ F_{\rm hyp}(\Omega) $. \finpreuve

\bigskip
\begin{exam} For all
$ f \in C_0^{\infty}(\Ra^n) $, all $ R > R_{\mathcal K} $ and all $ V \Subset \Psi_{\iota}(U_{\iota}) $,  (\ref{symboleprincipalcalculfonctionnel}) satisfies the conditions of this lemma with $ \Omega = (R,+\infty) \times V \times \Ra^n$.
\end{exam}

\noindent {\it Proof.} By (\ref{longueportee}),  there exists $ C > 1 $ such that
\begin{eqnarray}
C^{-1}|\xi|^2 \lesssim q_{\iota} (r,\theta,\xi) \lesssim C |\xi|^2 , \qquad r > R, \ \ \theta \in V , \ \ \xi \in \Ra^{n-1}, \label{ellipticiteuniforme}
\end{eqnarray}
and, using the notation (\ref{notationLaplacieninfini}),
\begin{eqnarray}
| \partial_r^j \partial_{\theta}^{\alpha} g^{kl}(r,\theta) | \leq C_{jk}, \qquad r > R, \ \ \theta \in V . \label{estimeeuniforme}
\end{eqnarray}
Therefore,  (\ref{ellipticiteuniforme}) and the compact support of $f$ ensure that $ e^{-r} \eta $ and $ \rho $ are bounded, ie that (\ref{ensembleborne}) holds on the support of (\ref{symboleprincipalcalculfonctionnel}). Then, (\ref{estimeeuniforme}) implies that
$ f (\rho^2 + q_{\iota}(r,\theta,\xi)) $ belongs to $ C^{\infty}_b (F_{\rm hyp} (\Omega)) $ (notice that
here $F_{\rm hyp} (\Omega)=(R,+\infty) \times V \times \Ra^n $). \finpreuve

\bigskip

We conclude this subsection with the following useful remarks. If $a,b \in {\mathcal S}_{\rm hyp}(\Omega) $ for some $ \Omega $ 
(such $a,b$ satisfy (\ref{symbolpseudo})), 
we have the composition rule
\begin{eqnarray}
\widehat{O \! p }_{\iota}(a)\widehat{O \! p }_{\iota}(b) =\widehat{O \! p }_{\iota}( (a \# b)(h) ) , \label{compositionpseudodifferentiel}
\end{eqnarray}
if $ (a \# b) (h) $ denotes the full symbol of $ a (r,\theta,h D_r , h D_{\theta}) b (r,\theta,h D_r , h D_{\theta}) $. In particular all the terms
of the expansion of $ (a \# b) (h)  $ belong to $ {\mathcal S}_{\rm hyp}(\Omega) $ and are supported in $ \mbox{supp}(a) \cap \mbox{supp}(b)  $. 
Similarly, for all $ N \geq 0$, 
we have
\begin{eqnarray}
\widehat{O \! p }_{\iota}(a)^* = \widehat{O \! p }_{\iota}(a^*_0 + \cdots + h^N a_N^* ) + h^{N+1} R_N (a,h) \label{adjointpseudodifferentiel} 
\end{eqnarray}
with $ a^*_0 , \ldots , a_N^* \in  {\mathcal S}_{\rm hyp}(\Omega) $ supported in $ \mbox{supp}(a)  $ and $ || R_N (a,h) ||_{L^2({\mathcal M},\widehat{dG}) \rightarrow L^2 ({\mathcal M} , \widehat{dG})} \lesssim 1  $ for $ h \in (0,1] $.

\subsection{The functional calculus} \label{soussection2}
In Proposition \ref{propcalculufonctionnel} below, we give two pseudo-differential approximations of $ f (h^2 P) $ near infinity of $ {\mathcal M} $, when $ f \in C_0^{\infty}(\Ra) $. The first approximation, namely (\ref{developpementL2}), is given in terms of the `quantization' $ \widehat{O \! p}_{\iota} $ defined in the previous subsection. This is the one we shall mostly use in this paper. However, at some crucial points, we shall need another approximation, (\ref{developpementpropre}), which uses {\it properly supported} $ \Psi $DO.

To define such properly supported operators, we  need a function
$$ \zeta \in C_0^{\infty}(\Ra^n), \qquad \zeta \equiv 1 \ \ \mbox{near 0} , \qquad \mbox{supp}(\zeta) \ \ \mbox{small enough} , $$
which will basically be used as a cutoff near the diagonal. The smallness of the support will be fixed in the following definition.
\begin{defi}
For $ \iota \in {\mathcal I} $, all $ h \in (0,1] $ and all symbol $ a $ satisfying (\ref{symbolpseudo}),
we define
$$ O \! p _{\iota, {\rm pr}}(a) : C_0^{\infty}({\mathcal M}) \rightarrow C^{\infty}({\mathcal M}), $$
as the unique operator with kernel supported in $ {\mathcal U}_{\iota}^2 $ and such that the kernel of $  \Psi_{\iota}^* \widehat{O \! p }_{\iota}(a) \Psi_{\iota *} $ is
\begin{eqnarray}
   (2 \pi h)^{-n} \int \! \! \int e^{\frac{i}{h} (r-r^{\prime})\rho
 + \frac{i}{h} (\theta - \theta^{\prime})\cdot \eta} a (r,\theta, \rho , \eta) d \rho d \eta \zeta (r-r^{\prime},\theta-\theta^{\prime})  . \label{defpseudopr}
\end{eqnarray}
\end{defi}
The interest of choosing the support of $ \zeta $ small enough is that, using (\ref{symbolpseudo}), we can assume that,  on the support of (\ref{defpseudopr}),
$ r^{\prime}  $ belongs to a neighborhood of $ [R_{\mathcal K} + 1 , + \infty ) $ and $ \theta^{\prime} $ belongs to a neighborhood of
$ V^{\prime}_{\iota} $. For instance, we may assume that $ r^{\prime} \in \widetilde{\kappa}^{-1}(1) $ and $ \theta^{\prime} \in
\widetilde{\kappa}_{\iota}^{-1}(1) $ so that we can put  a factor $ \widetilde{\kappa}(r^{\prime}) \widetilde{\kappa}_{\iota}(\theta^{\prime}) $ for free to the right hand side of (\ref{defpseudopr}). The latter implies, using (\ref{equivalenceLp}), (\ref{noyauiota}), (\ref{defpseudopr}), the standard off diagonal fast decay of kernels
 of $ \Psi $DO and the Calder\'on-Vaillancourt theorem that, for all $ a \in C^{\infty}_b (\Ra^{2n}) $ satisfying (\ref{symbolpseudo}) and all $ N \in \Na $, we have
\begin{eqnarray}
||  \widehat{O \! p }_{\iota}(a) - O \! p _{\iota, {\rm pr}}(a)
||_{L^2 ({\mathcal M},\widehat{dG}) \rightarrow L^2 ({\mathcal M}
, \widehat{dG})} \lesssim h^N, \qquad h \in (0,1] .
\label{pasproprepropre}
\end{eqnarray}
This shows that, up to remainders of size $ h^{\infty} $, $ \widehat{O \! p }_{\iota} (a) $ and $ O \! p _{\iota, {\rm pr}}(a) $ coincide as bounded operators
on $ L^2 ({\mathcal M} , \widehat{dG}) $.
Under the same assumptions on $a$, we also have
\begin{eqnarray}
 ||   O \! p_{\iota, {\rm pr}}(a)  ||_{L^2 ({\mathcal M}, dG ) \rightarrow L^2 ({\mathcal M} , dG)} \lesssim 1, \qquad h \in (0,1] , \label{borneL2pr}
\end{eqnarray}
which is a first difference with $ \widehat{O \! p }_{\iota}(a) $ for which we have only (\ref{borneL2apoids}) in general. The estimate
(\ref{borneL2pr}) is equivalent to the uniform boundedness\footnote{for $ h \in (0,1] $}
of $ e^{\bott r }O \! p_{\iota,pr}(a) e^{- r} $ on $ L^2 ({\mathcal M},\widehat{dG}) $. The latter is obtained similarly to (\ref{borneL2gratuite}), using the Calder\'on-Vaillancourt theorem, for we only have to consider the kernel obtained by multiplying (\ref{defpseudopr})
 by  $e^{\bott(r-r^{\prime})}$ which is bounded (as well as its derivatives) on the support of $ \zeta (r-r^{\prime},\theta-\theta^{\prime})$.

 In other words, (\ref{borneL2pr}) can be interpreted as a boundedness result between (exponentially) weighted $L^2$ spaces. Similar properties holds for $ L^q $ spaces (under suitable assumptions on the symbol $a$) and they are the main reason for considering properly supported operators. In particular, they lead to following proposition where we collect the estimates we shall need in this paper.
 We refer to \cite{BoucletLpCF} for the proof.

 \begin{prop} \label{propcalculufonctionnel} Let $f \in C_0^{\infty}(\Ra)$ and $ I \Subset (0,\infty) $ be an open interval containing $ \emph{supp}(f) $. Let $ \chi_{\mathcal K} \in C_0^{\infty}({\mathcal M}) $ and $ R > R_{\mathcal K} + 1 $ be such that
 $$ \chi_{\mathcal K}(m) = 1 \qquad \mbox{if} \ \ r (m) \leq R + 1 . $$
 Then, for all $ N \geq 0 $ and all $ \iota \in {\mathcal I} $, we can find symbols
\begin{eqnarray}
  a_{\iota,0} (f) , \ldots , a_{\iota,N} (f) \in {\mathcal S}_{\rm hyp} \left( ( R , \infty ) \times V_{\iota} \times \Ra^n \cap p_{\iota}^{-1} (I)   \right) , \label{choixsymboles}
\end{eqnarray}
(where $ p_{\iota} $ is the principal symbol of $P$ in the chart $ {\mathcal U}_{\iota} $) such that, if we set
 $$ a_{\iota}^{(N)}(f,h) =  a_{\iota,0} (f) + h a_{\iota,1}(f) + \cdots + h^N a_{\iota,N}(f), $$
 we have
\begin{eqnarray}
(1-\chi_{\mathcal K}) f (h^2 P) & = &  \sum_{\iota \in {\mathcal I}} \widehat{O \! p}_{\iota}( a_{\iota}^{(N)}(f,h) ) + h^{N+1} \widehat{R}_{N}(f,h) ,  \label{developpementL2} \\
& = & \sum_{\iota \in {\mathcal I}} O \! p_{\iota, {\rm pr}}( a_{\iota}^{(N)}(f,h) ) + h^{N+1} R_{N,{\rm pr}}(f,h) ,  \label{developpementpropre}
\end{eqnarray}
where,  for each $ q \in [2,\infty] $,
\begin{eqnarray}
|| e^{- \bott r} R_{N,{\rm pr}} (f,h) ||_{L^2 ({\mathcal M},\widehat{dG}) \rightarrow L^q({\mathcal M},dG)} & \lesssim & h^{-n \left( \frac{1}{2} - \frac{1}{q} \right)}, \qquad h \in (0,1] , \label{restepr}
\end{eqnarray}
and
\begin{eqnarray}
||  \widehat{R}_{N} (f,h) ||_{L^2 ({\mathcal M},\widehat{dG}) \rightarrow L^2({\mathcal M},\widehat{dG})} & \lesssim & 1 , \qquad h \in (0,1] . \label{restechapeau}
\end{eqnarray}
In addition, for all $ \iota \in {\mathcal I} $ and all $ q \in [2,\infty] $, we have
\begin{eqnarray}
|| e^{- \bott r} O \! p_{\iota, {\rm pr}}( a_{\iota}^{(N)}(f,h) ) ||_{L^2 ({\mathcal M},\widehat{dG}) \rightarrow L^{q}({\mathcal M},dG)} & \lesssim & h^{- n \left( \frac{1}{2} - \frac{1}{q} \right) }, \qquad h \in (0,1] ,
\end{eqnarray}
and, for all $ q \in [1,\infty] $ and all $ \gamma \in \Ra $,
\begin{eqnarray}
|| e^{- \gamma r} O \! p_{\iota, {\rm pr}}( a_{\iota}^{(N)}(f,h) ) e^{\gamma r}
 ||_{L^q ({\mathcal M},\widehat{dG}) \rightarrow L^q({\mathcal M},\widehat{dG})}
  & \lesssim & 1 , \qquad h \in (0,1] . \label{expondec}
\end{eqnarray}
\end{prop}

To make (\ref{choixsymboles}) more explicit, let us  quote for instance that
$$ a_{\iota,0}(f) (r,\theta,\rho,\eta) = \kappa (r) \kappa_{\iota} (\theta) f (\rho^2 + q_{\iota}(r,\theta,e^{-r}\eta)) \times (1 - \chi_{\mathcal K})(\Psi_{\iota}^{-1}(r,\theta)) . $$
More generally, (\ref{choixsymboles}) and  Lemma \ref{lemmesymbole} show that $ a_{\iota,0} (f) , \ldots , a_{\iota,N} (f) $ are of the form (\ref{formsymbol}) with $ \tilde{a}_{\iota}(r,\theta,\rho,\xi) $
compactly supported with respect to $ (\rho , \xi ) $.

The estimate (\ref{expondec}) basically means that $ O \!
p_{\iota, {\rm pr}}( a_{\iota}^{(N)}(f,h) ) $ preserves all $ L^q
$ spaces with any exponential weights. In particular, since $ L^q
({\mathcal M},dG) = e^{- \bott r/ q} L^q ({\mathcal M},\widehat{dG})
$, replacing $ \widehat{dG} $ by $ dG $ in (\ref{expondec}) would
give a completely equivalent statement. This estimate is the main
reason for introducing properly supported operators. Of course,
 (\ref{expondec}) holds for other symbols than those  involved in
the functional calculus of $ P $. We have more generally (see
\cite{BoucletLpCF}) for all $ \gamma \in \Ra $,
\begin{eqnarray}
|| e^{- \gamma r} O \! p_{\iota, {\rm pr}}( a_{\iota} ) e^{\gamma
r}
 ||_{L^q ({\mathcal M},\widehat{dG}) \rightarrow L^q({\mathcal M},\widehat{dG})}
  & \lesssim & 1 , \qquad h \in (0,1], \label{proprepropre}
\end{eqnarray}
for any $ q \in [1,\infty] $ and any
$$ a_{\iota} \in {\mathcal S}_{\rm hyp} \left( (R_{\mathcal K} + 1 , + \infty) \times V_{\iota}^{\prime} \times \Ra^n \cap
p_{\iota}^{-1}(I^{\prime}) \right) , $$ provided $ I^{\prime} $ is
bounded.

By the unitary equivalence of $ P $ and $ - \Delta_G - \bott^2 $, we would get a very similar pseudo-differential expansion for $ f (-h^2 \Delta_G) $. Note that, here,  we have only described $ (1-\chi_{\mathcal K}) f (h^2 P) $ since this will be sufficient for our present purpose, but of course there is a completely analogous result for the compactly supported part $ \chi_{\mathcal K} f (h^2 P) $ (see \cite{BoucletLpCF}). Such an approximation of $ f (-h^2 \Delta_G)$ was used in \cite{BoucletLP} to prove the next two propositions.
\begin{prop} Consider a dyadic partition of unit
$$ 1 = f_0 (\lambda) + \sum_{k \geq 0} f (2^{-k} \lambda), $$
for $ \lambda $ in a neighborhood of $ [0,+\infty) $, with
\begin{eqnarray}
f_0 \in C_0^{\infty}(\Ra), \qquad f \in C_0^{\infty} \left( [1/4, 4] \right) .
\end{eqnarray}
Then, for all $ \chi \in C_0^{\infty} (\mathcal M) $, for all $ q \in [2,\infty ) $, we have
$$ || (1-\chi)  u ||_{L^q ({\mathcal M},dG)} 
\lesssim \left( \sum_{k \geq 0} ||(1-\chi) f(-h^2 \Delta_G) u ||_{L^q ({\mathcal M},dG)}^2 \right)^2 + ||u||_{L^2 ({\mathcal M},d G)} . $$
\end{prop}
This proposition leads to the following classical reduction.
\begin{prop} \label{localisationfrequence} Let $ \chi \in C_0^{\infty} ({\mathcal M}) $ and $ (p,q) $ be an admissible pair. Then (\ref{vraiStrichartz}) holds true  if and only if  there exists $ C $ such that
\begin{eqnarray}
 || (1-\chi) e^{it \Delta_G} f (-h^2 \Delta_G) u_0 ||_{L^p ([0,1]; L^q ({\mathcal M},dG))} \leq C || u_0 ||_{L^2 ({\mathcal M},dG)}, \label{theolocalise}
 \end{eqnarray}
 for all $ h \in (0,1 ] $ and $ u_0 \in C_0^{\infty}({\mathcal M}) $.
\end{prop}

This result is essentially well known and proved in \cite{BoucletLP} for a class of non compact manifolds. We simply recall here that the $ L^q \rightarrow L^q $ boundedness of the spectral cutoffs $ f (-h^2 \Delta_G) $ is not necessary to prove this result, although the latter slightly simplifies the proof when it is available.


\subsection{Outgoing and incoming areas} \label{soussection3}
Propositions \ref{propcalculufonctionnel} and \ref{localisationfrequence} lead to a microlocalization of Theorem
\ref{theoprincipal}: as we shall see more precisely in subsection \ref{soussection4}, they allow to reduce the proof of (\ref{vraiStrichartz}) to the same estimate in which $ (1-\chi) $ is replaced by $h$-$\Psi $DOs. This microlocalization, ie the support of the symbols in (\ref{choixsymboles}), is however still too rough to simplify the proof of Theorem \ref{theoprincipal} in a significant way. The purpose of this subsection is to describe  convenient regions which will refine this localization.

\begin{defi} \label{definitsortantentrantiota} Fix $ \iota \in {\mathcal I} $. Let $ R > R_{\mathcal K} +1 $, $ V \Subset V^{\prime}_{\iota} $ be an open subset (see \ref{margeouverts}), $ I \Subset (0,+\infty) $ be an open interval and
 $ \sigma \in (-1,1)
$. We define
\begin{eqnarray}
\Gamma^{\pm}_{\iota} (R,V,I,\sigma) = \{ (r,\theta,\rho,\eta) \in \Ra^{2n}
\ | \ r > R , \ \theta \in V , \ p_{\iota} \in I , \ \pm \rho > - \sigma
p_{\iota}^{1/2} \} , \nonumber
\end{eqnarray}
where $ p_{\iota} $ is the principal symbol of $P$ in the chart $ {\mathcal U}_{\iota} $ given by (\ref{symbolprincipaliota}).
 The open
set $ \Gamma^+_{\iota} (R,V,I,\sigma) $ (resp. $ \Gamma^-_{\iota} (R,V,I,\sigma)
$) is called an outgoing (resp. incoming) area.
\end{defi}

We note in passing that, except from the localization in $ \theta $, these areas are defined using only the variable $r$, its dual $ \rho $ and the principal symbol of $ P $. In particular, up to the choice of the coordinate $ r $, the conditions $r > R $, $ p_{\iota} \in I $ and $ \pm \rho > - \sigma  p_{\iota}^{1/2} $ define invariant  subsets of $ T^* {\mathcal M} $.  However the whole analysis in this paper will be localized in charts and we will not use this invariance property.

Let us record some useful properties of outgoing/incoming areas. First, they decrease with respect to $ V , I , \sigma $ and $ R^{-1} $:
\begin{eqnarray}
 R_1 \geq R_2, \ \ \ V_1 \subset V_2, \ \ \ I_1 \subset I_2, \ \ \ \sigma_1 \leq \sigma_2 \ \ \ \Rightarrow \ \ \
 \Gamma^{\pm}_{\iota} (R_1,V_1,I_1,\sigma_1) \subset \Gamma^{\pm}_{\iota} (R_2,V_2,I_2,\sigma_2) . \label{decroissancezones}
\end{eqnarray}
Second, we have
\begin{eqnarray}
 \Gamma^{+}_{\iota} (R,V,I,1/2) \cup \Gamma^{-}_{\iota} (R,V,I,1/2) = (R,+\infty) \times V \times \Ra^n \cap p_{\iota}^{-1}(I) . \label{recouvrementzones}
\end{eqnarray}
Here we have chosen $ \sigma =1/2 $ but any $ \sigma \in (0,1) $ would work as well.

We will use the following elementary property.
\begin{prop} \label{partitionsiota}  
Any symbol $ a \in {\mathcal S}_{\rm hyp} \left( (R,+\infty) \times V \times \Ra^n \cap p_{\iota}^{-1}(I) \right) $ can be written
$$ a = a^+ + a^-, \qquad \mbox{with} \ \ \ a^{\pm} \in {\mathcal S}_{\rm hyp} \left(  \Gamma^{\pm}_{\iota} (R,V,I,1/2) \right) . $$
\end{prop}

\noindent {\it Proof.} See part {\it ii)} of Proposition \ref{partitions}. \finpreuve

\bigskip

 This splitting into outgoing/incoming areas was sufficient  to use the Isozaki-Kitada parametrix in the
asymptotically Euclidean case; in the AH case, we will only be able to construct this parametrix
in much smaller areas, called {\it strongly} outgoing/incoming areas, which we now introduce.

  We first describe briefly the meaning of such areas, say in the outgoing case. Basically, being in an outgoing area means that $ \rho $ is {\it not too close to} $ - p^{1/2} $; the aim of strongly outgoing areas is to guarantee that $ \rho $ is {\it  very close to} $ p^{1/2} $, which is of course a much stronger restriction. This amounts essentially to chose $ \sigma  $ close to $ - 1 $ in the definition of outgoing areas. We will measure this closeness in term of a small parameter $ \epsilon $. It will actually be convenient to have the other parameters, namely $ R, V, I  $, depending also on $ \epsilon $, so we introduce

\begin{eqnarray}
 R (\epsilon)  =  1/ \epsilon,   \qquad
 V_{\iota,\epsilon}  =   \{ \theta \in \Ra^{n-1} \ | \ \mbox{dist}(\theta,V_{\iota}) < 
 \epsilon^2 \}, \qquad I (\epsilon) =  (1/4-\epsilon, 4+\epsilon) , \label{RVepsiloniota}
\end{eqnarray}
where we recall that $ V_{\iota} $ is defined in (\ref{margeouverts}).
\begin{defi} \label{zonesstrong} For all $ \epsilon > 0 $ small enough, we set
$$ \Gamma^{\pm}_{\iota,{\rm s}} (\epsilon) := \Gamma^{\pm} (R(\epsilon),V_{\iota,\epsilon},I(\epsilon),\epsilon^2-1)
. $$ The open set $ \Gamma^{+}_{\iota,{\rm s}}(\epsilon) $ (resp. $
\Gamma^{-}_{\iota,{\rm s}} (\epsilon) $) is called a strongly outgoing
(resp. incoming) area.
\end{defi}

The main interest of such areas is to ensure that $ e^{-r}|\eta| $
is small if $ \epsilon $ is small. Indeed, if $ q \in [0,+\infty) $
and $ - 1 < \sigma < 0 $, we have the equivalence
\begin{eqnarray}
   \pm \rho > - \sigma
(\rho^2 + q )^{1/2}  \qquad  \Leftrightarrow  \qquad  \pm
\rho
> 0 \qquad \mbox{and} \qquad q < \sigma^{-2} (1-\sigma^2)
\rho^2 . \label{equivalencepratique}
\end{eqnarray}
Therefore, there exists $ C $ such that, for all $  \epsilon   $ small enough and $  (r,\theta,\rho,\eta) \in \Gamma^{\pm}_{\iota,{\rm s}} (\epsilon)
$,
\begin{eqnarray}
  q_{\iota} (r,\theta, e^{-r}\eta) \leq C \epsilon^2 , \nonumber
\end{eqnarray}
which, by (\ref{ellipticiteuniforme}), is equivalent to
\begin{eqnarray}
 |e^{-r} \eta| \lesssim \epsilon . \label{crucialdiffeo} 
\end{eqnarray}
Note also that, by (\ref{decroissancezones}), strongly outgoing/incoming areas decrease with $ \epsilon $.

\bigskip

We now quote a result which motivates, at least partially, the introduction of strongly outgoing/incoming areas. 

Denote by $ \Phi_{\iota}^t $ the Hamiltonian flow of $ p_{\iota} $. This is of course the geodesic flow written in the chart $ \Psi_{\iota}({\mathcal U}_{\iota}) \times \Ra^n $ of $ T^* {\mathcal M} $.

\begin{prop} \label{verszonesortante} Let 
$ \sigma \in (-1,1) $. Then, for all $ \epsilon > 0 $ small enough, there exists $ T_{\epsilon} > 0 $ such that,
$$ \Phi^t_{\iota} \left( \Gamma^{\pm}_{\iota} (R(\epsilon),V_{\iota},I(\epsilon),\sigma) \right) \subset \Gamma^{\pm}_{ \iota,{\rm \epsilon}} (\epsilon),
\qquad \pm t \geq T_{\epsilon} .  $$
\end{prop}

\noindent {\it Proof.} Follows from Corollary \ref{dependencedudomaine} and Proposition \ref{sortantversfortementsortant}. \finpreuve

\bigskip

Note that, since $ p_{\iota} $ is only defined in the chart $ \Psi_{\iota}({\mathcal U}_{\iota}) \times \Ra^n $, its flow is not complete. We shall however see in Section \ref{sectionflot} that, for any initial data  $ (r,\theta,\rho,\eta) \in  \Gamma^{\pm}_{\iota} (R(\epsilon),V_{\iota},I,\sigma) $, $ \Phi_{\iota}^t (r,\theta,\rho,\eta) $ is well defined for all $ \pm t \geq 0 $, ie $ \Phi_{\iota}^t (r,\theta,\rho,\eta) \in \Psi_{\iota}({\mathcal U}_{\iota}) \times \Ra^n $ for all $ \pm t \geq 0 $.

Proposition \ref{verszonesortante} essentially states that the forward (resp. backward flow) sends outgoing (resp. incoming) areas into strongly outgoing (resp. incoming) areas in finite positive (resp. negative) time.

\bigskip

\bigskip

The last type of regions we need to consider are the {\it intermediate} areas. They should have two properties: firstly they should essentially cover the complement of strongly outgoing/incoming areas in outgoing/incoming areas
and, secondly, be  small enough.

To define them we need the following. For all $ \epsilon > 0 $ and all $ \delta > 0 $, we can find $ L + 1 $ real numbers, $ \sigma_0, \ldots , \sigma_L $,
\begin{eqnarray}
\left( \frac{\epsilon}{2} \right)^2 - 1 = \sigma_0 <  \sigma_2 <\ldots < \sigma_L = 1/2 , \label{conditionsurlessigmas}
\end{eqnarray}
   such that
\begin{eqnarray}
 ((\epsilon/2)^2 - 1 , 1/2) = \cup_{l = 1}^{L-1} (\sigma_{l-1},\sigma_{l + 1})  , \label{recouvrementinterval}
\end{eqnarray}
and
\begin{eqnarray}
 |\sigma_{l+1} - \sigma_{l-1}| \leq \delta . \label{tailleintervalle}
\end{eqnarray}
Note that the intervals overlap in (\ref{recouvrementinterval}), since $(\sigma_{l-1},\sigma_{l + 1})$ always contains $ \sigma_{l} $.

\begin{defi} \label{definitionintermediaireoriginale} The intermediate outgoing/incoming area associated to the cover (\ref{recouvrementinterval}) are
$$ \Gamma^{\pm}_{\iota,{\rm inter}} (\epsilon,\delta;l) := \left\{ (r,\theta,\rho,\eta) \in \Ra^{2n} \ | \ r > R(\epsilon), \ \theta \in V_{\iota}
, \ p_{\iota} \in I(\epsilon), \  \pm \frac{\rho}{p_{\iota}^{1/2}} \in (-\sigma_{l+1}, -\sigma_{l-1}) \right\}  , $$
for $ 1 \leq l \leq L-1 $.
\end{defi}

Notice that, by definition,
\begin{eqnarray}
\Gamma^{\pm}_{\iota,{\rm inter}} (\epsilon,\delta;l) \subset
\Gamma^{\pm}_{\iota} (R(\epsilon), V_{\iota} , I(\epsilon) , 1/2 ) .
\label{inclusionlogique}
\end{eqnarray}

In the notation,  we only specify  the parameters which are relevant for our analysis, namely $ \epsilon,\delta $,  but,
of course, intermediate areas depend on the choice of $ \sigma_1 , \ldots , \sigma_L $. Here $ \delta $ measures the
smallness and Proposition \ref{dynamiquedelta} below
will explain how to choose this parameter.

 We first give the following result.
\begin{prop} \label{deuxiemedecoupage} Fix $ \epsilon > 0 $ small enough, $ \delta > 0 $ and $ \sigma_0 , \ldots , \sigma_L $  satisfying (\ref{conditionsurlessigmas}), (\ref{recouvrementinterval}) and (\ref{tailleintervalle}).
Then, any symbol
$$ a^{\pm} \in {\mathcal S}_{\rm hyp} (  \Gamma^{\pm}_{\iota}(R(\epsilon),V_{\iota},I(\epsilon),1/2) )  $$
can be written
\begin{eqnarray}
 a^{\pm} = a^{\pm}_{\rm s} + a_{1,{\rm inter}}^{\pm} + \cdots + a_{L-1,{\rm inter}}^{\pm}, \nonumber
\end{eqnarray}
with
\begin{eqnarray}
 a^{\pm}_{\rm s} \in {\mathcal S}_{\rm hyp} (  \Gamma^{\pm}_{\iota,{\rm s}} ( \epsilon) )  ,
\qquad a_{l,{\rm inter}}^{\pm} \in {\mathcal S}_{\rm hyp} (  \Gamma^{\pm}_{\iota,{\rm inter}} (\epsilon,\delta;l) ) \nonumber .
\end{eqnarray}
\end{prop}

\noindent {\it Proof.} Follows from Proposition \ref{deuxiemedecoupagesansiota}. \finpreuve

\bigskip

We conclude this subsection with the following proposition which will be crucial for the proof of Theorem \ref{theoprincipal} and motivates the introduction of intermediate areas.

\begin{prop} \label{dynamiquedelta} Fix   $ \underline{t} > 0 $. Then, for all $ \epsilon > 0 $ small enough, we can find $ \delta > 0 $ small enough such that, for any choice of $ \sigma_0, \ldots , \sigma_L $ satisfying (\ref{conditionsurlessigmas}), (\ref{recouvrementinterval}) and (\ref{tailleintervalle}), we have, for all $ 1 \leq l \leq L-1 $,
\begin{eqnarray}
\Phi^t_{\iota} \left( \Gamma^{\pm}_{\iota,{\rm inter}} (\epsilon,\delta;l) \right) \cap \Gamma^{\pm}_{\iota,{\rm inter}} (\epsilon,\delta;l) = \emptyset, \nonumber
\end{eqnarray}
provided that
\begin{eqnarray}
\pm t \geq \underline{t}. \nonumber
\end{eqnarray}

\end{prop}

\noindent {\it Proof.} Follows from Corollary \ref{dependencedudomaine} and Proposition \ref{dynamiquedeltasansiota}.

\subsection{The main steps of the proof of Theorem \ref{theoprincipal}} \label{soussection4}
We already know from Proposition \ref{localisationfrequence} that we only have to find $
\chi \in C_0^{\infty}({\mathcal M}) $ such that
(\ref{theolocalise}) holds, which is equivalent to
\begin{eqnarray}
 || e^{- \bott r} (1-\chi) f (h^2 P) e^{-it P}  u_0 ||_{L^p ([0,1]; L^q ({\mathcal M},dG))} \leq C || u_0 ||_{L^2 ({\mathcal M},\widehat{dG})},
 \label{reformulationmesure}
 \end{eqnarray}
using the unitary map (\ref{unitaire}) and
(\ref{coefficientsconcrets}), (\ref{defP}).  

Before choosing $ \chi $, we introduce the following operators. Choose  a cutoff $
\tilde{f} \in C_0^{\infty}((0,+\infty)) $ such that $ \tilde{f} f = f $. 
\begin{lemm} For all $ \chi \in C_0^{\infty}({\mathcal M})  $, we can write
$$ (1-\chi) \tilde{f}(h^2 P) = (1-\chi) A_{\rm pr}(h) + R (h) $$
with $ R (h) $ satisfying, for all $ q \in [2,\infty] $,
\begin{eqnarray}
|| e^{-\bott r} R(h) ||_{L^2 ({\mathcal M},\widehat{dG}) \rightarrow L^q({\mathcal M},dG)} & \lesssim & 1 , \label{SobolevLq}
\end{eqnarray}
and $ A_{\rm pr}(h) $  such that, for all $ q \in [2,\infty] $,
\begin{eqnarray}
|| e^{-\bott r} A_{\rm pr}(h) ||_{L^2 ({\mathcal M},\widehat{dG}) \rightarrow  L^{q}({\mathcal M},dG)} & \lesssim &
 h^{-n\left( \frac{1}{2} - \frac{1}{q} \right)} , \label{aadjoindre} \\
|| e^{-\bott r} A_{\rm pr}(h) e^{\bott r} ||_{L^{\infty} ({\mathcal
M},dG) \rightarrow L^{\infty}({\mathcal M},dG)} & \lesssim & 1 , 
\label{propreLinfini} \\
||  A_{\rm pr}(h)^* e^{- \bott r} ||_{L^1 ({\mathcal M},\widehat{dG}) \rightarrow L^2({\mathcal M},\widehat{dG})} &
 \lesssim & h^{-n/2}, \label{aadjoindrefait} \\
|| e^{\bott r} A_{\rm pr}(h)^* e^{- \bott r} ||_{L^{1} ({\mathcal
M},\widehat{dG}) \rightarrow L^{1}({\mathcal M},\widehat{dG})} &
\lesssim &  1 . \label{propreL1}
\end{eqnarray}
\end{lemm}

\noindent {\it Proof.} It is an immediate consequence of Proposition \ref{propcalculufonctionnel}. Using (\ref{developpementpropre}), with $ N $
such that $ N + 1 \geq n/2 $, we  define $A_{\rm pr}(h)$ as the sum
of the properly supported pseudo-differential operators. We thus have (\ref{SobolevLq}), (\ref{aadjoindre}) and (\ref{propreLinfini}). The estimates (\ref{aadjoindrefait}) and (\ref{propreL1}) are obtained by
 taking the adjoints (with $q = \infty$ in (\ref{aadjoindre}))
with respect to $ \widehat{dG} $. \finpreuve
 
 \bigskip
 
Basically, the operators $  e^{- \bott r} A_{\rm pr}(h) $ and $  A_{\rm pr}(h)^* e^{- \bott r} $ will be used as 'ghost cutoffs' to deal with remainder terms of parametrices which will be $ {\mathcal O}(h^N) $ in $ {\mathcal L}(L^2({\mathcal M},\widehat{dG})) $, using the Sobolev embeddings (\ref{aadjoindre}) and (\ref{aadjoindrefait}). They will be 'transparent' for the principal terms of the parametrices by (\ref{propreLinfini}) and (\ref{propreL1}),  which use crucially that they are properly supported.

\medskip

For $ \epsilon $ to be fixed below, we choose $ \chi \in C_0^{\infty}({\mathcal M}) $ such that
$$ \chi \equiv 1 \ \ \ \mbox{for} \ \ \ r (m) \leq 3 \epsilon^{-1} . $$
\begin{prop} \label{propositioninversezero} To prove (\ref{reformulationmesure}), it is sufficient to show that, for some $ \epsilon $ small enough and all 
\begin{eqnarray}
 a_{\iota} \in {\mathcal S}_{\rm hyp}  \left( R (\epsilon)
\times V_{\iota} \times \Ra^n \cap p_{\iota}^{-1}(I(\epsilon)) \right) , \label{proprietesaiota}
\end{eqnarray}
 where we
recall that $ R (\epsilon) = \epsilon^{-1} $ and $ I(\epsilon)= (1/4 - \epsilon , 4 + \epsilon) $, we have
\begin{eqnarray}
 || e^{-\bott r} A_{\rm pr} (h) \widehat{O \! p}_{\iota} (a_{\iota}) e^{-itP} u_0 ||_{L^p ([0,1];  L^q ({\mathcal M},dG))}
  \leq C || u_0 ||_{L^2 ({\mathcal M},\widehat{dG})} . \label{preuvereduitepropre}
\end{eqnarray}
\end{prop}

\noindent {\it Proof.} Choose $ \chi_0 \in
C_0^{\infty}({\mathcal M}) $ such that
$$ \chi_0 \equiv 1 \ \ \ \mbox{for} \ \ \ r (m) \leq \epsilon^{-1} , \qquad \chi_0 \equiv 0 \ \ \ \mbox{for} \ \ \ r (m) \geq 2 \epsilon^{-1} . $$
We then have $ (1-\chi_0) \equiv 1 $ near $ \mbox{supp} (1 - \chi)
$ so, by the proper support of the kernel of $ A_{\rm pr} (h) $,
we also have
$$ (1-\chi) A_{\rm pr}(h) = (1-\chi) A_{\rm pr} (h) (1-\chi_0) , $$
at least for $ \epsilon $ small enough. The latter and
(\ref{SobolevLq}) reduces the proof of (\ref{reformulationmesure}) to the study of
\begin{eqnarray*}
  e^{-\bott r} A_{\rm pr} (h) (1-\chi_0) f (h^2 P) e^{-it P}  .
 \end{eqnarray*}
By splitting $ (1-\chi_0) f (h^2 P) $ using (\ref{developpementL2}) with $ N + 1 \geq n/2 $, we obtain the result using  (\ref{restechapeau}) and (\ref{aadjoindre}). \finpreuve

\bigskip

We now introduce a second small parameter $ \delta > 0 $. By Propositions \ref{partitionsiota}
 and \ref{deuxiemedecoupage},  for
all $ \delta
> 0 $, any $ a_{\iota} $ satisfying (\ref{proprietesaiota}) can be written
\begin{eqnarray}
 a_{\iota} =  a^{+}_{\rm s} + a^{-}_{\rm s} + \sum_{l=1}^{L-1} a_{l,{\rm inter}}^{+} + a_{l,{\rm inter}}^{+} ,  \label{decoupagesymbole}
\end{eqnarray}
with
\begin{eqnarray}
    a^{\pm}_{\rm s} \in {\mathcal S}_{\rm hyp} (
\Gamma^{\pm}_{\iota,{\rm s}} ( \epsilon) )  , \qquad
 a_{l,{\rm inter}}^{\pm} \in {\mathcal S}_{\rm hyp} (  \Gamma^{\pm}_{\iota,{\rm inter}} (\epsilon,\delta;l) )
 . \label{pourdelta}
\end{eqnarray}
 
\begin{prop} \label{propositioninverse} To prove (\ref{preuvereduitepropre}), it is sufficient to show that, for some $ \epsilon $ and $ \delta  $ small enough, we have:
\begin{eqnarray}
||e^{- \bott r} A_{\rm pr} (h) \widehat{O \! p}_{\iota} (a_{\rm
s}^{\pm}) e^{-ithP} \widehat{O \! p}_{\iota} (a_{\rm s}^{\pm})^*
A_{\rm pr} (h)^* e^{- \bott r} ||_{L^1 (\widehat{dG}) \rightarrow
L^{\infty}(dG)} \leq C_{\epsilon}|ht|^{-n/2} ,
\label{dispersionfortesemiclassique}
\end{eqnarray}
and
\begin{eqnarray}||e^{- \bott r} A_{\rm pr} (h) \widehat{O \! p}_{\iota} (a_{l,{\rm inter}}^{\pm})
e^{-ithP} \widehat{O \! p}_{\iota} (a_{l,{\rm inter}}^{\pm})^*
A_{\rm pr} (h)^* e^{- \bott r} ||_{L^1 (\widehat{dG}) \rightarrow
L^{\infty}(dG)} \leq C_{\epsilon,\delta}|ht|^{-n/2} ,
\label{dispersionintersemiclassique}
\end{eqnarray}
for
\begin{eqnarray}
h \in (0,1 ] \qquad 0 \leq \pm t \leq 2 h^{-1} .
\label{timeorientation}
\end{eqnarray}
\end{prop}

Recall that the important point in this lemma is
(\ref{timeorientation}), ie that we can take $ t \geq 0 $
for outgoing localizations, and $ t \leq 0 $ for incoming
ones. 

\finpreuve

\noindent {\it Proof.} Let us define
$$ T_{\rm s}^{\pm} (t,h,\epsilon) = e^{- \bott r} A_{\rm pr} (h) \widehat{O \! p}_{\iota} (a_{\rm s}^{\pm}) e^{-itP} , $$
and
$$ T_{l,{\rm inter}}^{\pm}(t,h,\epsilon,\delta) = e^{- \bott r} A_{\rm pr} (h) \widehat{O \! p}_{\iota} (a_{l,{\rm inter}}) e^{-itP} . $$
By (\ref{borneL2gratuite}) and (\ref{aadjoindre}) (with $ q=2 $), we have,
$$ ||T_{\rm s}^{\pm} (t,h,\epsilon) ||_{L^2 (\widehat{dG}) \rightarrow L^2(dG)} +
|| T_{l,{\rm inter}}^{\pm} (t,h,\epsilon,\delta)||_{L^2 (\widehat{dG}) \rightarrow L^2(dG)} \leq C_{\epsilon,\delta}, \qquad h \in (0,1] , \ t \in \Ra , $$
hence by the Keel-Tao Theorem \cite{KeelTao}, (\ref{preuvereduitepropre}) would follow from the estimates
\begin{eqnarray}
||T_{\rm s}^{\pm} (t,h,\epsilon) T_{\rm s}^{\pm} (s,h,\epsilon)^*
||_{L^1 (dG) \rightarrow L^{\infty}(dG)} \leq
C_{\epsilon}|t-s|^{-n/2} , \label{dispersionforte}
\end{eqnarray}
and
\begin{eqnarray}
||T_{l,{\rm inter}}^{\pm} (t,h,\epsilon) T_{l,{\rm inter}}^{\pm}
(s,h,\epsilon)^* ||_{L^1 (dG) \rightarrow L^{\infty}(dG)} \leq
C_{\epsilon,\delta}|t-s|^{-n/2} , \label{dispersioninter}
\end{eqnarray}
for $ h \in (0,1]$ and $ t ,s \in [0,1] $.
Using the time rescaling $ t \mapsto ht $, the fact that $ L^1
(dG) = e^{-2 \bott r} L^1 (\widehat{dG}) $ and that the adjoint of
(\ref{unitaire}) is given by $ e^{ \bott r} $,
(\ref{dispersionforte}) and (\ref{dispersioninter}) are respectively equivalent to (\ref{dispersionfortesemiclassique}) and 
(\ref{dispersionintersemiclassique}), for $ h \in (0,1] $ and $ |t| \leq 2 h^{-1} $. 
The reduction (\ref{timeorientation}) to $ \pm t \geq 0 $ is obtained similarly to \cite[Lemma 4.3]{BoucletTzvetkov}. We only recall here
that it is based on the simple observation that the operators $ T (t) T (s)^* $ considered above are of the form $ B e^{-i(t-s)P} B^* $
so that $ L^{\infty} $ bounds on their Schwartz kernel for $ \pm (t-s) \geq 0 $ give automatically bounds for $ \pm (t-s) \leq 0 $
by taking the adjoints. \finpreuve

\bigskip

As we shall see, there basically be two reasons for choosing $ \epsilon $ small enough.
The next result is the first condition.

\begin{prop} \label{sousIK} For all $ \epsilon > 0 $ small enough and all $ a^{\pm}_{\rm s} \in {\mathcal S}_{\rm hyp} \left(
\Gamma^{\pm}_{\iota,{\rm s}} ( \epsilon) \right) $, we can write
$$ e^{-ithP} \widehat{O \! p}_{\iota} (a_{\rm s}^{\pm})^* = E^{\pm}_{\rm IK}(t,h)  + h^n R_{\rm IK}^{\pm}(t,h) , $$
with
\begin{eqnarray}
||e^{- \bott r}E_{\rm IK}^{\pm}(t,h) e^{- \bott r}||_{L^1(\widehat{dG})
\rightarrow L^\infty(dG)} & \lesssim & |ht|^{-n/2} , \label{dispersionIsozakiKitada}  \\
||R_{\rm IK}^{\pm}(t,h)||_{L^2(\widehat{dG}) \rightarrow
L^2(\widehat{dG})} & \lesssim & 1 , \label{resteL2IsozakiKitada}
\end{eqnarray}
for
$$  h \in (0,1] , \qquad 0 \leq
\pm t \leq 2 h^{-1} . $$
\end{prop}

\noindent {\it Proof.} By $  (\ref{adjointpseudodifferentiel}) $, the result follows 
  from Theorem \ref{IsozakiKitadaansatz}  and Section \ref{sectiondispersion}.  \finpreuve

\bigskip

Proposition \ref{sousIK} is mainly an application of the Isozaki-Kitada parametrix. It has the following consequence.

\begin{prop} \label{commentutiliserIK} For all  $ \epsilon > 0 $ small enough, (\ref{dispersionfortesemiclassique})
holds for all $ h, t $ satisfying (\ref{timeorientation}).
\end{prop}

\noindent {\it Proof.} We first replace $ \widehat{O \! p}_{\iota}
(a_{\rm s}^{\pm}) $ by $ O \! p_{\iota, {\rm pr}} (a_{\rm
s}^{\pm}) $ to the left of $ e^{-ithP} $ in
(\ref{dispersionfortesemiclassique}). The  remainder term, which
is $ {\mathcal O}(h^{\infty}) $ in $ {\mathcal L}(L^2
(\widehat{dG})) $ by (\ref{pasproprepropre}), produces a term of
size $ {\mathcal O}(h^{\infty}) $ in $ {\mathcal
L}(L^1(\widehat{dG}),L^{\infty}(dG)) $ using (\ref{aadjoindre})
(with $q=\infty$) and (\ref{aadjoindrefait}). We then use
Proposition \ref{sousIK}: the remainder term satisfies
$$ ||e^{- \bott r} A_{\rm pr} (h) O \! p_{\iota,{\rm pr}} (a_{\rm
s}^{\pm}) e^{-ithP} h^n R_{\rm IK}^{\pm}(t,h) A_{\rm pr} (h)^*
e^{- \bott r} ||_{L^1 (\widehat{dG}) \rightarrow L^{\infty}(dG)}
\lesssim 1 \lesssim |ht|^{-d/2} , $$ and the main term $
E^{\pm}_{\rm IK}(t,h) $ gives the expected contribution using
(\ref{propreLinfini}), (\ref{propreL1}) and (\ref{proprepropre})
for $ O \! p_{\iota, {\rm pr}} (a_{\rm
s}^{\pm}) $. \finpreuve

\bigskip

The second condition on $ \epsilon $ will come from Proposition \ref{dynamiquedelta}. The latter proposition depends on some fixed small time
which will be given by the following result.
\begin{prop} \label{sousWKB} There exists $ t_{\rm WKB} > 0 $ such that, for all $ \epsilon > 0 $ small enough and   all symbol $ a^{\pm}
\in {\mathcal S}_{\rm hyp} \left( \Gamma^{\pm}_{\iota} (
R(\epsilon),V_{\iota},I,1/2) \right) $, we can write
$$ e^{-ithP} \widehat{O \! p}_{\iota} (a^{\pm})^* = E^{\pm}_{\rm WKB}(t,h)  + h^n R_{\rm WKB}^{\pm}(t,h) , $$
with
\begin{eqnarray}
||e^{- \bott r}E_{\rm WKB}^{\pm}(t,h) e^{- \bott r}||_{L^1(\widehat{dG})
\rightarrow L^\infty(dG)} & \lesssim & |ht|^{-n/2} , \label{weightWKB}  \\
||R_{\rm WKB}^{\pm}(t,h)||_{L^2(\widehat{dG}) \rightarrow
L^2(\widehat{dG})} & \lesssim & 1 , \nonumber
\end{eqnarray}
for
\begin{eqnarray}
  h \in (0,1] , \qquad 0 \leq \pm t \leq t_{\rm WKB} . \label{shortrangeoftime}
\end{eqnarray}
\end{prop}

\noindent {\it Proof.} See Sections \ref{WKBsection} and \ref{sectiondispersion}.

\bigskip

The first consequence of this proposition are the following small times dispersion estimates.

\begin{prop} \label{commentutiliserWKB}
For all $ \epsilon > 0 $, all $ \delta > 0 $ and all $ a^{\pm}_{l,{\rm inter}}$ satisfying (\ref{pourdelta}), the estimate
 (\ref{dispersionintersemiclassique})
holds for all $ h, t $ satisfying (\ref{shortrangeoftime}).
\end{prop}

\noindent {\it Proof.} It is completely similar to the proof of
Proposition \ref{commentutiliserIK}. \finpreuve

\bigskip
 
 We can now give the second condition on $ \epsilon $, also giving the choice of $ \delta $.

\begin{prop} \label{propositionfinale} If $ \epsilon $ is small enough, we can choose $ \delta > 0 $ small enough such that,
for all $ 1 \leq l \leq L-1 $, all 
$$ b_{l,{\rm inter}}^{\pm} \in {\mathcal S}_{\rm hyp} (  \Gamma^{\pm}_{\iota,{\rm inter}} (\epsilon,\delta;l) ) , $$
 and all $ N \geq 0 $, we have
\begin{eqnarray}
|| \widehat{O \! p}_{\iota} (b_{l,{\rm inter}}^{\pm}) e^{-ithP}
\widehat{O \! p}_{\iota} (b_{l,{\rm inter}}^{\pm})^* ||_{L^2
(\widehat{dG}) \rightarrow L^2 (\widehat{dG})} \leq C_{l,N} h^N \label{longueorthogonalite}
\end{eqnarray}
for
$$ h \in (0,1], \qquad   t_{\rm WKB} \leq \pm t \leq 2 h^{-1} . $$
\end{prop}

\noindent {\it Proof.} See Section \ref{WKBsection}.

\bigskip

This is, at least intuitively, a consequence of Proposition
\ref{dynamiquedelta} with $ \underline{t} = t_{\rm WKB} $ and of the Egorov Theorem which states that $ e^{-ithP}
\widehat{O \! p}_{\iota} (b_{l,{\rm inter}}^{\pm})^* $ lives semi-classically in the region 
$ \Phi^t_{\iota} \left( \mbox{supp}(b_{l,{\rm inter}}^{\pm}) \right) $.

\bigskip

We summarize the above reasoning as follows.

\medskip

\noindent {\bf Proof of Theorem \ref{theoprincipal}.}   Using Proposition \ref{commentutiliserIK}, we choose first $ \epsilon_0 > 0 $ small enough so that, for all $ \epsilon \in (0, \epsilon_0 ] $, (\ref{dispersionfortesemiclassique}) holds for $ 0 \leq \pm t \leq 2 h^{-1} $. By possibly decreasing $ \epsilon_0 $, we then choose $ t_{\rm WKB} $ according to Proposition \ref{sousWKB}, uniformly with respect to $ \epsilon \in (0,\epsilon_0 ] $.
Next, according to Proposition \ref{propositionfinale}, we fix $  \epsilon \in (0, \epsilon_0 ] $ and $ \delta > 0 $ small enough such that (\ref{longueorthogonalite}) holds for $ t_{\rm WKB} \leq \pm t \leq 2 h^{-1} $.  Using (\ref{aadjoindre}), (\ref{aadjoindrefait}) and
Proposition \ref{propositionfinale} with $ N = n $ and $ b^{\pm}_{l,{\rm inter}} = a^{\pm}_{l,{\rm inter}} $ defined by (\ref{decoupagesymbole}), 
we have
\begin{eqnarray}||e^{-\bott r} A_{\rm pr} (h) \widehat{O \! p}_{\iota} (a_{l,{\rm inter}}^{\pm})
e^{-ithP} \widehat{O \! p}_{\iota} (a_{l,{\rm inter}}^{\pm})^*
A_{\rm pr} (h)^* e^{- \bott r} ||_{L^1 (\widehat{dG}) \rightarrow
L^{\infty}(dG)} & \leq & \ C_{\epsilon,\delta} \nonumber \\ & \lesssim & |ht|^{-n/2} ,
\nonumber
\end{eqnarray}
for $   t_{\rm WKB} \leq \pm t \leq 2 h^{-1} $. On the other hand, (\ref{dispersionintersemiclassique}) holds for $ 0 \leq \pm t \leq t_{\rm WKB} $, using Proposition \ref{sousWKB}. Therefore 
(\ref{dispersionintersemiclassique}) holds for $ 0 \leq \pm t \leq 2 h^{-1} $. By Proposition \ref{propositioninverse}, this proves (\ref{preuvereduitepropre}) for all $ a_{\iota} $ satisfying  (\ref{proprietesaiota}). By Proposition \ref{propositioninversezero}, this implies (\ref{reformulationmesure}) which, by Proposition \ref{localisationfrequence}, implies Theorem \ref{theoprincipal}. 
\finpreuve

\section{Estimates on the geodesic flow near infinity} \label{sectionflot}
\setcounter{equation}{0}
In this section, we describe some properties of the Hamiltonian flow   of functions of the form
\begin{eqnarray}
 p (r,\theta,\rho,\eta) = \rho^2 + w(r)  q (r,\theta,\eta)  ,  
 \label{hamiloca}
\end{eqnarray}
on $ T^* \Ra^n_+ =
\Ra^+_r \times \Ra^{n-1}_{\theta} \times \Ra_{\rho} \times
\Ra^{n-1}_{\eta} $. Here $ q $ is an  homogeneous polynomial of degree $2$ w.r.t
$\eta$ and $w$ a positive function. Naturally, the motivation for the study of (\ref{hamiloca}) comes from the form of the principal
symbol of $ P $ given by  (\ref{symbolprincipaliota}).

 In subsection \ref{wpasquelconque}, we will assume that $w(r) = e^{-2r}$ but we start with more general cases in subsection \ref{wquelconque}. We emphasize that
$ p $ is defined on $ T^* \Ra^n_+ $ whereas $ p_{\iota} $ is only defined on a subset of the form $ T^* (R_{\mathcal K} , + \infty )
\times V_{\iota} $. The result of subsection \ref{wpasquelconque} will nevertheless hold for $ p_{\iota} $
 with no difficulty for we shall have good localization of the flow in the regions we consider (see Corollary \ref{dependencedudomaine}).

\subsection{A general result} \label{wquelconque} 
Let $ w = w (r) $ be a smooth function on $ \Ra^+ = (0,+ \infty) $
such that
\begin{eqnarray}
w > 0 , \qquad w^{\prime} < 0, \qquad \left( \frac{w^{\prime}}{w}
 \right)^{\prime} \geq 0,  \label{fourgene}
\end{eqnarray}
and, for some $ 0 < \gamma < 1 $,
\begin{eqnarray}
 \limsup_{r \rightarrow + \infty }
 \int_{r}^{(1+\gamma)r} \frac{w^{\prime}}{w}  \in [ - \infty , 0 ) . \label{fourgene2}
\end{eqnarray}
Note that $ \lim_{r \rightarrow + \infty} w (r) $ exists, by
(\ref{fourgene}), and that  (\ref{fourgene2}) imply that this
limit must be $ 0 $. Note also that, for all $ R > 0 $, we have
$$ w (r) \lesssim 1 \qquad \mbox{and} \qquad |w^{\prime}(r)| \lesssim w (r), \qquad \mbox{on} \ \ [ R , + \infty ) .  $$
 These assumptions are satisfied for instance by $ w (r) = r^{-2} $ or $ w (r) = e^{-2r} $.

We assume that $ q $ is an  homogeneous polynomial of degree $2$ w.r.t
$\eta$ of the form
\begin{eqnarray}
 q (r,\theta,\eta) = q_0 (\theta,\eta) + q_1 (r,\theta,\eta)
 \label{pourlasectionsuivante}
\end{eqnarray}
with $ q_0 , q_1 $ homogeneous polynomials of degree $2$ w.r.t
$\eta$ satisfying, for some $ 0 < \tau \leq 1 $,
\begin{eqnarray}
 |\partial_{\theta}^{\alpha} \partial_{\eta}^{\beta} q_0
(\theta,\eta) | & \lesssim & \scal{\eta}^{2-|\beta|} , \label{pourlasectionsuivante2} \\
|\partial_r^j
\partial_{\theta}^{\alpha}
\partial_{\eta}^{\beta} q_1 (r,\theta,\eta) | & \lesssim & \scal{r}^{-\tau-j} \scal{\eta}^{2-|\beta|}
, \label{longrange}
\end{eqnarray}
and, for some $ C > 0 $,
\begin{eqnarray}
C^{-1} |\eta|^2 \leq  q (r,\theta,\eta) \leq C |\eta|^2 ,
\label{secogene}
\end{eqnarray}
 for $
(r,\theta,\eta) \in \Ra^+ \times \Ra^{n-1} \times \Ra^{n-1} $. The
latter implies, by possibly increasing $ C $, that
\begin{eqnarray}
  C^{-1} |\eta|^2 \leq  q_0 (\theta,\eta) \leq C |\eta|^2, \qquad
  (\theta,\eta) \in \Ra^{n-1} \times \Ra^{n-1}.  \label{ellipiticiteuniformeinfini}
\end{eqnarray}
Setting $ q^{\prime} = \partial_r q  $ (ie $ \partial_r q_1 $), we
finally assume that,
\begin{eqnarray}
\frac{q^{\prime}}{q} \times \frac{w}{w^{\prime}} \longrightarrow 0
\qquad \mbox{as} \qquad r \rightarrow + \infty , \label{fiftgene}
\end{eqnarray}
 uniformly with respect to $\theta \in
\Ra^{n-1} $ and $ \eta \in \Ra^{n-1} \setminus 0 $.

The Hamiltonian flow $ \Phi^t = (r^t , \theta^t , \rho^t,\eta^t)
$, generated by $p$, is the solution to the system
\begin{eqnarray}
\begin{cases}
\ \dot{r} & = \  2 \rho   \\
\ \dot{\theta} & = \ w  \partial q /\partial \eta \\
\ \dot{\rho} & = \ - w^{\prime}  q - w q^{\prime}  \\
\ \dot{\eta} & = \  - w   \partial q /\partial \theta
\end{cases}
\label{Hamigene}
\end{eqnarray}
with initial condition
\begin{eqnarray}
(r^t,\theta^t,\rho^t,\eta^t)_{|t=0} = (r,\theta,\rho,\eta) .
\label{initgene}
\end{eqnarray}

Our main purpose is to show that, if $ \rho > - p^{1/2} $ (with $p
= p (r,\theta,\rho,\eta)$) and $ r $ is large enough, then
$\Phi^t$ is defined for all $ t \geq 0 $ and $ r_t \rightarrow +
\infty $ as $ t \rightarrow + \infty $ (we will obtain a similar
result for $ t \leq 0 $ provided $ \rho < p^{1/2}  $). This result
 relies mainly on the following remark: if $ \eta \ne 0 $, we can
write
$$
 - w^{\prime} q - w q^{\prime} = - \frac{w^{\prime}}{w} \left( p
- \rho^2 \right)\left(1 + \frac{w}{w^{\prime}} \times
\frac{q^{\prime}}{q} \right) .
$$
Using (\ref{fiftgene}) and the negativity of $ w^{\prime}/w $,
this shows that, for all $ \epsilon
> 0 $, we can find $ R > 0 $  such that
\begin{eqnarray}
 - w^{\prime} q - w q^{\prime} \geq - (1 - \epsilon) (p - \rho^2) \frac{w^{\prime}}{w} ,
 \qquad \mbox{on} \ \ [R, + \infty )_r
 \times \Ra^{n-1}_{\theta} \times  \Ra_{\rho} \times \Ra^{n-1}_{\eta}
\label{deriplus}
\end{eqnarray}
which we shall exploit to prove that $ \dot{\rho} \geq 0 $.

In the following lemma and in the sequel, we shall extensively use
the shorter notation
$$ p = p (r,\theta,\rho,\eta) . $$

\begin{lemm} \label{gene1} Denote by $ (-t_-,t_+) $ ($t_{\pm} \in (0,+\infty]$) the maximal interval on
which the solution of (\ref{Hamigene}), with initial condition
(\ref{initgene}), is defined. Then
$$ t_{\pm} \geq \frac{r}{2 p^{1/2}} . $$
Furthermore, either $ r_t \rightarrow 0 $ as $ t \rightarrow t_+ $
(resp. $ t \rightarrow - t_- $) or $ t_+ = + \infty $ (resp. $ t_-
= + \infty $).
\end{lemm}
Note that, if $ p (r,\theta,\rho,\eta) = 0 $, i.e. $ \rho = 0 $
and $ \eta = 0 $, then it is trivial that $ t_{\pm} = + \infty $.

\smallskip

\noindent {\it Proof.} We will only consider the case of $ t_+ $,
the one of $ t_- $ being similar.
 By the conservation of energy  we have $ |\rho^t| \leq p^{1/2}  $
 thus, for $ t \in [0, t_+ )
 $, $ \dot{r}^t $ is bounded,
\begin{eqnarray}
|r^t - r| \leq 2 t p^{1/2}, \label{bornetrivialeenr}
\end{eqnarray}
 and $  r^t \geq r - 2 t p^{1/2} $.
  We now argue by contradiction
and assume that $ t_+ < r / 2 p^{1/2} $ (in particular, that $
t_+$ is finite). Then $ r_+ : = r  - 2 t_+ p^{1/2} > 0 $ and $ r_t
\geq r_+ $ for all $ t \in [ 0, t_+ ) $. Furthermore, by
(\ref{secogene}), we have $ |w
\partial_{\eta} q | \leq C (w q + w ) \leq C
 (p
+w) $, with $ w $ bounded on $ [ r_+   , + \infty ) $, hence
  $ \dot{\theta}^t $ is bounded on $ [ 0,
t_+ ) $. One shows similarly that $ \dot{\rho}^t $ and $
\dot{\eta}^t $ are bounded on $ [0,t_+) $, using that $ |
w^{\prime} | \lesssim w $ on $[ r_+ , + \infty ) $ for $
\dot{\rho} $. This implies that $ \lim_{t \rightarrow t_+}
(r^t,\theta^t,\rho^t,\eta^t) $ exists and belongs to $ (0,+\infty)
\times \Ra^{n-1} \times \Ra \times \Ra^{n-1} $. The solution can
therefore be continued beyond $ t_+ $ which yields the
contradiction.

We now consider the second statement. Assume that $ t_+ < + \infty
$. We must show that $ r^t \rightarrow 0 $ as $ t \rightarrow t_+
$. Assume that this is wrong. Then there exists $ R > 0 $ small
enough and a sequence $ t_k \rightarrow t_+ $ such that $ r^{t_k}
\geq R $ for all $ k \geq 0 $. On the other hand, by energy
conservation, we have $ |r^t - r^s| \leq 2 p^{1/2}|t-s| $ for all
$ t,s \in [0,t_+) $, hence
$$ r^t \geq r^{t_k} - 2 p^{1/2}|t-t_k| \geq R/2 $$
 provided $ |t-t_k | \leq R / 4 p^{1/2}$. Since $t_k$ can be chosen as close to $ t_+ $ as we want,
 there exists  $\epsilon > 0$
 small enough such that $ r^t \geq R/2 $ for $ t \in [t_+ - \epsilon , t_+)
 $. Then, by the same argument as above, $ \lim_{t \rightarrow
t_+} (r^t,\theta^t,\rho^t,\eta^t) $ exists and belongs to $
(0,+\infty) \times \Ra^{n-1} \times \Ra \times \Ra^{n-1} $. The
solution can be continued beyond $ t_+ $ hence $ t_+ = + \infty $
which is a contradiction.
 \finpreuve

\smallskip

\begin{lemm} \label{gene2} Let $ 0 < \epsilon < 1 $. For any $ R > 0 $  such that
(\ref{deriplus})  holds, we have the following: if $ r^{t_0} \geq
R $ and  $ \rho^{t_0} > 0 $ (resp. $ \rho^{t_0} < 0 $) for some $
t_0 \in [0,t_+) $ (resp. $t_0 \in (-t_-,0]$), then $ t_+ = +
\infty $ (resp. $ - t_- = -  \infty $) and
$$ r^t \geq R , \qquad \rho^t \geq \rho^{t_0} \ \ \mbox{(resp.} \ \rho^t \leq \rho^{t_0} \mbox{)}
\qquad \forall \ t \geq t_0 \ \ \mbox{(resp.} \ t \leq
t_0\mbox{)}.
$$ Furthermore, $ r^t \geq r^{t_0} + 2 (t-t_0 )\rho^{t_0} $ for all $ t \geq t_0 $ (resp. $ t \leq t_0 $).
\end{lemm}
\noindent {\it Proof.} As in Lemma \ref{gene1}, we only consider
the case of $ t_+ $. It suffices to show that
\begin{eqnarray}
 r^t \geq R, \qquad \mbox{for all} \ t \in [t_0, t_+ )  . \label{cont2}
\end{eqnarray}
Indeed, if this is true, Lemma $ \ref{gene1} $ shows that $t_+ =
+ \infty $ and then, by (\ref{deriplus} ), we have $ \dot{\rho}^t
\geq 0 $ hence $ \rho^t \geq \rho^{t_0} $ and $ r^t - r^{t_0} \geq
2 \rho^{t_0}(t-t_0) $. Let us  prove (\ref{cont2}). We consider
the set
$$ I = \{ t \in [t_0, t_+ ) \ | \ r^s \geq R \
\mbox{and} \ \rho^s \geq \rho^{t_0} \ \ \forall \ s \in [t_0,t] \}
.$$ It is clearly an interval containing $ t_0 $ and we set $ T:=
\sup I $. By continuity, $ \rho^t \geq \rho^{t_0} / 2
> 0 $ for $ t  $ in a small neighborhood $ J $ of $ t_0 $. This
implies that $ \dot{r}^t
> 0 $ on  $ J $, hence that $ r^t \geq r^{t_0} \geq R $ on $ J \cap [t_0,t_+) $
and thus that $ \dot{\rho}^t \geq 0 $ on $ J \cap [t_0,t_+) $
which in turn shows that $ \rho^t \geq \rho^{t_0} $ on $ J \cap
[t_0,t_+) $. This proves that $ T > t_0 $. Then, on  $ [t_0,T) $,
we have
\begin{eqnarray}
 r^t \geq R , \qquad \rho^t \geq \rho^{t_0} . \label{mino1}
\end{eqnarray}
Now assume, by contradiction, that $ T < t_+ $. Then (\ref{mino1}) holds on $ [t_0,T] $ and in particular we have $
r^T \geq r^{t_0} + 2 (T-t_0) \rho^{t_0} > r^{t_0} $. Thus $ r^t
\geq R $ in a neighborhood of $ T$ and this implies that $
\dot{\rho}^t \geq 0 $ in this neighborhood. Hence there exists $
T^{\prime} > T $ such that  (\ref{mino1}) holds on $ [t_0,
T^{\prime}] $  yielding a contradiction. \finpreuve

\medskip

To state the next result, we define $ l \in (0,+\infty ] $ as
\begin{eqnarray}
 l = -  \limsup_{r \rightarrow + \infty }
 \int_{r}^{(1+\gamma)r} \frac{w^{\prime}}{w} \label{defilimi}
\end{eqnarray}
 and we choose an arbitrary $  \sigma \in \Ra $ such that
\begin{eqnarray}
 0 < \sigma < \begin{cases} - \frac{2}{l} + \left( \frac{4}{l^2} + 1 \right)^{1/2} & \ \mbox{if} \ l < +
 \infty \\
 1 & \ \mbox{if} \ l = +
 \infty
 \end{cases} . \label{choisigm}
\end{eqnarray}
Note that, if $ l $ is finite, $ 0 < - \frac{2}{l} +
(\frac{4}{l^2}+ 1)^{1/2} < 1 $ and that (\ref{choisigm}) is
equivalent to
$$ (1 - \sigma^2) l/2 > 2 \sigma
> 0 . $$
\begin{prop} \label{propgene} For any $ \sigma $ satisfying (\ref{choisigm}), there exists
$ R_{w,\gamma,\sigma} > 0 $ large enough such that:
 if $ r > R_{w,\gamma,\sigma} $  and  $ \rho > - \sigma p^{1/2} $ (resp $ \rho <  \sigma p^{1/2} $),
 then $ t_+ = + \infty $ (resp $ - t_- = - \infty $) and for all $  t \geq 0 $
(resp. $ t \leq 0 $) we have
\begin{eqnarray}
  r^t \geq \max \left((1-\gamma)r , (1-  \gamma - \sigma
\gamma) r + 2 \sigma p^{1/2} |t|  \right)  . \label{factgene}
\end{eqnarray}
\end{prop}

This proposition means that, by choosing an initial data with $ r
$ large enough and $ \rho > - \sigma p^{1/2}$ (resp. $ \rho <
\sigma p^{1/2} $), the forward (resp. backward) trajectory lies in
a neighborhood of infinity. In particular, the forward (resp.
backward) flow starting at $(r,\theta,\rho,\eta) $, with $  \rho
> - \sigma p^{1/2} $ (resp $ \rho < \sigma p^{1/2} $ ) depends only on the values of $ p $
on $ [(1-\gamma) r, + \infty ) \times \Ra^{n-1} \times \Ra \times
\Ra^{n-1}$.

\bigskip

\noindent {\it Proof.} We only consider the case where $ \rho > -
\sigma p^{1/2} $, the case where $ \rho <  \sigma p^{1/2} $ being
similar. If $ l < \infty $, (\ref{choisigm}) allows to choose  $ 0
< \epsilon < 1 $ such that
\begin{eqnarray}
 (1 - \epsilon)^2 (1-\sigma^2) l / 2 \geq 2 \sigma . \label{choiepso}
\end{eqnarray}
  If $ l =
\infty $, we choose an arbitrary $ \epsilon \in (0,1) $.
 We next choose $ R $
so that (\ref{deriplus}) holds with the above choice of $ \epsilon
$. If $ \rho \geq \sigma p^{1/2} $ (recall that $ p^{1/2} > 0 $
since $ \rho > - \sigma p^{1/2} $) and $ r \geq R $, then Lemma
\ref{gene2} shows that the result holds with $ R_{w,\gamma,\sigma}
= R $. We can therefore assume that $ \rho < \sigma p^{1/2} $. Let
us set
\begin{eqnarray}
 R_1 = (1 - \gamma)^{-1}R  \qquad \mbox{and} \qquad T = \gamma r
/ 2 p^{1/2} . \label{conditionsurT}
\end{eqnarray}
 By Lemma \ref{gene1}, we have $ t_+ >   T  $ and,
if $ r \geq R_1 $,
$$ r^t \geq r - 2 t p^{1/2} \geq (1-\gamma)r \geq R, \qquad \mbox{for} \ \ t \in [0,T] . $$
Using (\ref{deriplus}), this implies that $ \dot{\rho}^t \geq 0 $
on $ [0,T] $ and hence that $ \rho^t \geq - \sigma p^{1/2} $ for
all $ t \in [0,T] $. Let us now prove by contradiction that there
exists $ t \in [0,T] $ such that $ \rho^t \geq \sigma p^{1/2} $.
If this is wrong,  we have $ (\rho^t)^2 \leq \sigma^2 p $ on $
[0,T] $, thus (\ref{deriplus}) shows that, for all $ t \in [0,T]
$,
$$ \dot{\rho}^t \geq - (1 - \epsilon) (1 - \sigma^2) p \frac{w^{\prime}}{w}(r^t) \geq
- (1 - \epsilon) (1 - \sigma^2) p \frac{w^{\prime}}{w} \left( r +
2 t p^{1/2} \right) ,
$$ using the third estimate of (\ref{fourgene}) and the fact that $
r^t \leq r + 2 t p^{1/2} $ in the second inequality. By
integration over $ [0,T] $, we get
\begin{eqnarray}
\rho^T - \rho \geq - (1-\epsilon)(1-\sigma^2) p^{1/2} \frac{1}{2}
\int_{r}^{(1+\gamma)r } \frac{w^{\prime}}{w} , \label{clefgene}
\end{eqnarray}
using the second equality in (\ref{conditionsurT}). Let us now fix
$ R_2 $ such that, for all $ r
> R_2 $,
$$ - \int_{r}^{(1+\gamma) r}
\frac{w^{\prime}}{w}   \ \ > \begin{cases} (1- \epsilon) l  & \ \mbox{if} \ \ l < + \infty \\
 \frac{4 \sigma}{(1-\epsilon)(1-\sigma^2)} &  \ \mbox{if} \ \ l = + \infty \end{cases} . $$
With such a choice (and (\ref{choiepso}) if $ l $ is finite), we
see that, if $ r \geq \max (R_1,R_2) $, (\ref{clefgene}) implies
that $ \rho^T - \rho \geq 2 \sigma p^{1/2}  $ and hence that $
\rho^T \geq \sigma p^{1/2} $ which yields the expected
contradiction.

In summary, we have shown that for any $ r \geq \max( R_1, R_2 ) $
and any $ \rho
> - \sigma p^{1/2} $, there exists $ t_0 \in [0,T] $ such that $
\rho^{t_0} \geq \sigma p^{1/2} > 0 $ and $ r^{t_0} \geq R $, hence
$ t_+ = + \infty $ by Lemma  \ref{gene2}. Furthermore, $ r^t \geq
(1 - \gamma) r $ on $ [0,T] $ and $ r^t \geq r^T + 2(t-T)\sigma
p^{1/2} \geq (1 - (1+\sigma) \gamma ) r + 2 t \sigma p^{1/2} $ on
$ [T,+\infty) $. Since
$$ \max \left((1-\gamma)r , (1-  \gamma - \sigma \gamma) r + 2 \sigma p^{1/2} t
\right) =
\begin{cases} (1-\gamma)r & \mbox{if} \ \ t \in [0,T] \\ (1-  \gamma - \sigma \gamma) r + 2 \sigma p^{1/2} t
 & \mbox{if} \ \ t > T \end{cases} , $$
the result follows. \finpreuve

\subsection{The asymptotically hyperbolic case} \label{wpasquelconque}

In this part, we prove more precise estimates on the Hamiltonian
flow of $ p $  when
$$ w (r) = e^{-2r} . $$
In that case, the conditions (\ref{fourgene}),  (\ref{fourgene2})
and (\ref{fiftgene}) are fulfilled, with any $ 0 < \gamma < 1 $ in
(\ref{fourgene2}) and we have $l = + \infty $ in (\ref{defilimi}).

\medskip

In the sequel, we shall need the following improvement of
Proposition  \ref{propgene}.

\begin{prop} \label{localisationr} Let $ 0 < \sigma < 1 $. There exist $ R_{\sigma} > 0
$ and $ C_{\sigma} > 0 $ such that: if $ r \geq R_{\sigma} $ and
 $ \rho > - \sigma p^{1/2} $ (resp. $ \rho < \sigma p^{1/2}
$), then
$$ r^t \geq r + 2 \sigma p^{1/2}|t| - C_{\sigma} ,
\qquad \mbox{for all} \ \ t \geq 0 \ \ \mbox{(resp.}
\ t \leq 0 \mbox{)}.  $$
\end{prop}

The improvement   consists in replacing $ (1- \gamma -  \sigma
\gamma ) r $ in the estimate (\ref{factgene}) by $ r - C_{\sigma}
$.

\bigskip

\noindent {\it Proof.} Here again we only consider the case $ t
\geq 0 $.  By  Proposition \ref{propgene}, we may assume that $
r^t \geq R $ for all $t \geq 0$, with $ R $ large enough so that
(\ref{deriplus}) holds with $ \epsilon = 1/2 $. This implies that
\begin{eqnarray}
 \dot{\rho}^t =  2 e^{-2r^t} q (r^t,\theta^t,\eta^t) - e^{-2r^t} \partial_r q_1 (r^t,\theta^t,\eta^t) \geq
e^{-2r^t} q (r^t,\theta^t,\eta^t) = p - (\rho^t)^2 . \label{pourzoneintermediaire}
\end{eqnarray}
  If $ \rho
\geq \sigma p^{1/2} $, then the result follows from Lemma
\ref{gene2}  (with $ C_{\sigma} = 0 $). If $ \rho < \sigma p^{1/2}
$, we will show that, with $ T = 2 \sigma p^{-1/2} / (1 -
\sigma^2) $, there exists $ t \in [0,T] $ such that $ \rho^t \geq
\sigma p^{1/2} $. Assume that this is wrong. Then  $ (\rho^t)^2
\leq \sigma^2 p $ on $
 [0,T]$ and by integrating the above estimate on $ \dot{\rho}^t $, we get
$$ \rho^T - \rho \geq  T (1 - \sigma^2) p = 2 \sigma p^{1/2} . $$
This proves that $ \rho^T \geq \sigma p^{1/2} $ which is a
contradiction. Therefore, by Lemma \ref{gene2}, we see that $ r^t
- r^T \geq 2 \sigma p^{1/2} (t - T)  $ for $ t \geq T $. On the
other hand, we have $ r^t \geq r - 2 p^{1/2} t $ for $ t \in [0,T]
$. The latter implies that $ r^t \geq r + 2 \sigma p^{1/2}t - 2
p^{1/2}(1+\sigma) t \geq r + 2 \sigma p^{1/2} t - 4 \sigma / (1-
\sigma) $ for $ t \in [0,T] $. This holds in particular for $ t =
T $ and then for $ t \geq T $. Thus the results holds with $
C_{\sigma} = 4 \sigma / (1-\sigma)  $. \finpreuve

\bigskip

We have so far only studied some localization properties of $
\Phi^t $, the Hamiltonian flow of $p$. We shall now give estimates
on derivatives of $ \Phi^t $. We start with the following
proposition giving some rough estimates. They will serve as a
priori estimates for the proof of Proposition
\ref{estimeesflotprecises} below.

\begin{lemm} \label{lemmemotivezones} For all $ 0 < \sigma < 1  $, there exists $ R > 0 $
such that, for all $ (r,\theta,\rho,\eta) \in T^* \Ra^n_+ $
satisfying
\begin{eqnarray}
 r > R , \qquad \pm \rho > - \sigma p^{1/2}, \qquad p \in (1/4 , 4)
 ,
\label{motivezones}
\end{eqnarray}
and  all $ \pm t \geq 0 $, we have
$$ \left| e^{r|\beta|} \partial_{\eta}^{\beta} \partial_r^j \partial_{\theta}^{\alpha}
 \partial_{\eta}^k \left( \Phi^t - \Phi^0 \right)(r,\theta,\rho,\eta) \right| \lesssim \scal{t}  . $$
\end{lemm}

Note the $ e^{r |\beta|}  $ factor in front of the derivatives.

We will need two lemmas. The first one is a soft version of the
classical Faa Di Bruno formula.
\begin{lemm} \label{Faadibruno}  Let  $ \Omega_1 \subset \Ra^{n_1} $, $ \Omega_2 \subset \Ra^{n_2}
$ be open subsets, with $ n_1,n_2 \geq 1 $. Consider smooth maps $
y = (y_1, \ldots , y_{n_2}) : \Omega_1 \rightarrow \Omega_2 $ and
$ Z : \Omega_1 \times \Omega_2 \rightarrow \Ra^{n_3} $, with $ n_3
\geq 1 $. Then, for all $ | \gamma | \geq 1 $,
$$ \partial_x^{ \gamma } \left( Z (x,y(x)) \right) = (\partial_y Z)(x,y(x)) \partial_x^{ \gamma } y (x)
 + \left( \partial_x^{ \gamma } Z \right) (x,y(x)) + R_{ \gamma }(x)  $$
 where $ R_{\gamma } (x) \equiv 0 $ if $ | \gamma | = 1 $ and, otherwise,
 is a linear combination of
$$   \left( \partial^{\gamma - \gamma^{\prime}}_x \partial_y^{\nu} Z \right)
(x,y(x)) \left( \partial_x^{\gamma_1^1} y_1 (x) \ldots
\partial_x^{\gamma_{\nu_1}^1} y_1 (x) \right) \ldots \left( \partial_x^{\gamma_1^{n_2}}
y_{n_2}(x) \ldots
\partial_x^{\gamma_{\nu_{n_2}}^{n_2}} y_{n_2}(x) \right),
$$ with $ \gamma , \gamma^{\prime}, \gamma_j^k \in \Na^{n_1}  $, $
\nu = (\nu_1 , \ldots , \nu_{n_2}) \in \Na^{n_2} $ satisfying $
\gamma^{\prime} \ne 0 $, $ \nu \ne 0 $ and
$$  \gamma^{\prime} \leq \gamma , \qquad 2 \leq |\nu| + |\gamma -
\gamma^{\prime} | \leq |\gamma| , \qquad \gamma_1^1 + \cdots +
\gamma_{\beta_1}^1 + \cdots + \gamma_1^{n_2} + \cdots +
\gamma^{n_2}_{\beta_{n_2}} = \gamma^{\prime} ,
$$
 and using the convention that $ \partial^{\gamma^{k}_1}_x y_k (x)
\ldots
\partial_x^{\gamma^k_{\nu_k}} y_k (x) \equiv 1 $ if $ \nu_k = 0 $
(if $ \nu_k \ne 0 $ then  $\gamma^{k}_1, \ldots , \gamma^k_{\nu_k}
$ are all  non zero).
\end{lemm}

\noindent {\it Proof.} It follows by a direct induction.
\finpreuve

\bigskip

In the second lemma, we consider the linear differential equation
\begin{eqnarray}
 \dot{X} = A(t)X + Y(t), \label{edolineaire}
\end{eqnarray}
  with $ A (\cdot) \in C
([0,+\infty),{\mathcal M}_{N \times N}(\Ra)) $\footnote{${\mathcal
M}_{N \times N}(\Ra)$ the space of $ N \times N $ matrices with
real entries} and $ Y ( \cdot ) \in C ([0,+\infty), \Ca^N) $ for
some $ N \geq 1 $. We assume that $  A (\cdot)  $ belongs to a
subset $ {\mathcal B} \subset C ([0,+\infty),{\mathcal M}_{N
\times N}(\Ra)) $ for which there exist $ \delta_{\mathcal B} > 0
$ and $ C_{\mathcal B}
> 0 $ such that
$$ |||A(t)||| \leq  C_{\mathcal B} e^{- \delta_{\mathcal B} t}, \qquad t \geq
0 , \ \ A (\cdot) \in {\mathcal B}, $$ with $ ||| \cdot||| $ a
matrix norm associated  to the norm $ ||\cdot|| $ on $ \Ca^N $,
i.e. such that $ || M Z || \leq |||M||| || Z || $, for all $ M \in
{\mathcal M}_{N \times N}(\Ra) $ and $ Z \in \Ca^N $.
\begin{lemm} \label{Gronwallexponentiel} There exists $ C > 0 $
 such that, for all $ A (\cdot) \in {\mathcal B} $ and all $ Y (\cdot) $ satisfying
$$ \int_0^{\infty}||Y(t)||dt < \infty ,$$
the solutions $ X (\cdot) $ of  (\ref{edolineaire})  satisfy
\begin{eqnarray}
 || X (t) || \leq C \left( ||X(0)||+ \int_0^{\infty}||Y(s)||ds
 \right) , \qquad t \geq 0 .  \label{integreebor}
\end{eqnarray}
\end{lemm}

\noindent {\it Proof.} Fix first $ 0  < \delta  < \delta_{\mathcal
B}  $ and $ \varepsilon = \delta_{\mathcal B} -\delta $.  Choose $
T > 0 $ such that $ C_{\mathcal B}  e^{-\delta_{\mathcal B} t}
\leq \varepsilon $ for $ t \geq T $. By the Gronwall Lemma, we
have
$$ || X (t) || \leq \left( ||X(T) || +  \int_{T}^{\infty}||Y(s)||ds \right)
e^{\varepsilon(t-T)} , \qquad t \geq T , $$ and
$$  || X (t) || \leq \left( ||X(0) || +  \int_{0}^{T}||Y(s)||ds \right)
e^{C_{\mathcal B} T} , \qquad t \in [0, T ] . $$ These two
inequalities give, for some $ C $ depending only on $ C_{\mathcal
B} $, $ \delta_{\mathcal B} $, $ \delta $ and $ T $,
$$ || X (t) || \leq C \left( ||X(0) || +
\int_{0}^{\infty}||Y(s)||ds \right) e^{ \varepsilon t}, \qquad t
\geq 0 . $$ Used as an a priori estimate in (\ref{edolineaire}),
this yields
$$   ||\dot{X}(t) || \leq || Y (t) || + C C_{\mathcal B} e^{-\delta t} \left(
  ||X(0)|| + \int_0^{\infty}||Y(s)||ds  \right),
\qquad t \geq 0 ,  $$ which implies  (\ref{integreebor}).
\finpreuve

\bigskip

\noindent {\it Proof of Lemma \ref{lemmemotivezones}.} As
before, we only prove the result for $ t \geq 0 $. For $ |\beta| +
j + |\alpha| + k = 0 $, the result is a consequence of the motion
equations (\ref{Hamigene}) and the energy conservation. Indeed,
for $ r^t - r $, the estimate follows directly from
(\ref{bornetrivialeenr}).
 Next, the motion equation for $ \theta $,
(\ref{secogene}) and Proposition \ref{localisationr} show that
$$ |\dot{\theta}^t| \lesssim e^{-2 r^t} |\eta^t| \lesssim  e^{-2 r^t} \scal{\eta^t }^2 \lesssim 1 + p  $$
hence that $  |\theta^t - \theta| \lesssim \scal{t} $ by
integration. One similarly shows that $|\rho^t - \rho| + |\eta^t -
\rho| \lesssim \scal{t} $. We now consider the derivatives and
denote for simplicity $ \partial^{\gamma} =
\partial_{\eta}^{\beta}
\partial_r^j
\partial_{\theta}^{\alpha}
 \partial_{\eta}^k $. Denoting by $ H_p$ is the Hamiltonian vector field of $ p $ and applying
 $ \partial^{\gamma} $ to (\ref{Hamigene}), we obtain
 $$ e^{r|\beta|}\partial^{\gamma} \dot{\Phi}^t = (d H_p)(\Phi^t) e^{r|\beta|}\partial^{\gamma} \Phi^t + R (t) $$
 where, by Lemma \ref{Faadibruno},  $ R (t) = 0 $
 if $ |\gamma | = 1 $ or, if $ |\gamma| \geq 2 $, is  a
 linear combination of
\begin{eqnarray}
 (\partial^{\nu} H_p) (\Phi^t) e^{r|\beta|} \left( \partial^{
\gamma^{1}_{1}} r^t \cdots \partial^{\gamma^{1}_{\nu_1}} r^t
\right)
 \cdots \left( \partial^{ \gamma^{2n}_{1}} \eta^t_{n-1} \cdots \partial^{\gamma^{2n}_{\nu_{2n}}} \eta^t_{n-1}
 \right).
 \label{combinatoirederivee}
\end{eqnarray}
 Here $ \nu = (\nu_1 , \ldots , \nu_{2n}) $ is of length at least
$2$,  all the derivatives of $ \Phi^t $ involved in $ R (t) $ are
of strictly smaller order than $ \gamma $ (ie $ \gamma^i_{l_i}
\leq \gamma $ and $ \gamma^i_{l_i} \ne \gamma $) and
\begin{eqnarray}
2 \leq |\nu| \leq |\gamma|, \qquad \gamma^1_1 + \cdots +
\gamma^{2n}_{\nu_{2n}} = \gamma . \label{decompositiongamma}
\end{eqnarray}
Writing  $ d H_p  $ as a matrix, we have
$$ d H_p =  \left(
\begin{matrix}
0 & 0 & 2 & 0 \\
0 & 0 &  0 & 0 \\
0 & 0 & 0  & 0 \\
0 & 0 & 0  & 0
\end{matrix} \right) +
e^{-2r} \left(
\begin{matrix}
0 & 0 & 0  & 0 \\
  \partial^2_{r \eta} q_1 - 2 \partial_{\eta} q &  \partial^2_{\theta \eta} q & 0  & \partial^2_{\eta \eta} q \\
4 \partial_{r} q_1 - 4 q - \partial^2_{rr} q_1 & 2
\partial_{\theta} q -\partial^2_{\theta r} q_1 & 0 & 2
\partial_{\eta} q  - \partial^2_{ \eta r } q_1 \\
 2 \partial_{\theta} q -\partial^2_{r \theta } q_1
 & - \partial^2_{\theta \theta} q &
0  & - \partial^2_{\eta \theta} q
\end{matrix} \right). $$
Defining $ M $ as the first (constant) matrix of the right hand
side and using Proposition \ref{localisationr}, we have
$$  \left| d H_p (\Phi^t) - M \right| \lesssim e^{-2 r^t} \scal{\eta^t}^2 \lesssim e^{-2 r - 2 \sigma t}
 (\scal{\eta}^2 + \scal{t}^2) \lesssim e^{- \sigma t} , $$
using that $ 2 p^{1/2} \geq 1 $ and that $ e^{-2r} \scal{\eta}^2 $
is bounded, by (\ref{motivezones}). We then set
\begin{eqnarray}
 A (t) & =& e^{-tM} \left( d H_p (\Phi^t)  -M \right) e^{t M}, \nonumber \\  X (t) & =
& e^{-tM} e^{r|\beta|}
\partial^{\gamma}\Phi^t - e^{r|\beta|}
\partial^{\gamma}\Phi^0,  \nonumber \\ Y (t) & = & e^{-tM} R (t) + A (t) e^{r|\beta|}
\partial^{\gamma}\Phi^0 , \nonumber
\end{eqnarray}
so that
$$ \dot{X}(t) = A (t) X (t) + Y (t) , \qquad X (0) = 0 . $$
Noting that $ M^2 = 0 $, we have
\begin{eqnarray}
 \exp (\pm t M) = 1 \pm t M , \qquad | \exp (\pm tM)| \lesssim \scal{t} ,  \label{exponentielleM}
\end{eqnarray}
 thus
\begin{eqnarray}
  | A (t) | \lesssim e^{- \sigma t} \scal{t}^2 \lesssim
e^{-\sigma t / 2}  . \label{exponentielleentemps}
\end{eqnarray}
To estimate $ X (t) $ by Lemma \ref{Gronwallexponentiel}, we still
need to estimate $ Y (t) $. We first assume that $
\partial^{\gamma} = \partial^{\beta}_{\eta} $ with $ |\beta | = 1
$. We then have $ R (t) = 0 $ and
$$  A (t) e^{r|\beta|}
\partial^{\gamma}\Phi^0 = e^{-t M}  ( \partial_{\eta}^{\beta} H_p ) (\Phi^t) e^{r} $$
since $ M \partial_{\eta}^{\beta} \Phi^0 = 0 $. By Proposition
\ref{localisationr} and (\ref{motivezones}) again, we obtain
$$ | ( \partial_{\eta}^{\beta} H_p )
(\Phi^t) | \lesssim e^{-2 r - 2 \sigma t } \scal{\eta^t} \lesssim
e^{ -r - \sigma t} , $$ so that $ |Y(t)| \lesssim e^{-\sigma t /
2} $. Using (\ref{exponentielleentemps}) and Lemma
\ref{Gronwallexponentiel}, we get $ |X(t)| \lesssim 1 $. Since $ M
\partial_{\eta}^{\beta} \Phi^0 = 0 $, we can rewrite  $ X (t) = e^{-tM}
e^{r} \partial_{\eta}^{\beta}( \Phi^t - \Phi^0   ) $ and, using
(\ref{exponentielleM}), finally get
$$ \left| e^{r} \partial_{\eta}^{\beta}( \Phi^t - \Phi^0   ) \right| \lesssim \scal{t} . $$
The other first order derivatives of $ \Phi^t - \Phi^0 $ are
studied similarly (note that there is no $e^r$ factor then), by
showing that $ X (t) $ is bounded and using that  $ X (t) =
e^{-tM}
\partial^{\gamma}(\Phi^t - \Phi^0  ) + ( e^{-tM} - 1 )
\partial^{\gamma}\Phi^0 $ with (\ref{exponentielleM}) to get
$$ | \partial^{\gamma} ( \Phi^t - \Phi^0 ) | \lesssim \scal{t} .  $$
For higher order derivatives, $ \partial^{\gamma}\Phi^0 = 0 $ and
$
\partial^{\gamma}(\Phi^t - \Phi^0 ) = \partial^{\gamma}\Phi^t $.
Furthermore, since the derivatives of $ \Phi^t $ involved in $ R
(t) $ are of strictly smaller order than $ \gamma $, we can
proceed by induction. By writing $ x^t $ for  $
r^t,\rho^t,\theta^t$ and $ \partial^{\gamma^i_l} =
\partial^{\beta^i_l}_{\eta} \partial^{k^i_l}_r \partial^{\alpha^i_l}_{\theta} \partial^{j^i_l}_{\rho}
$ for the derivatives involved in (\ref{combinatoirederivee}),
with $ 1 \leq i \leq 2n  $ and $ 1 \leq l \leq \nu_i $ (recall
that, if $ \nu_i =0 $, the corresponding product in
(\ref{combinatoirederivee}) is $ 1 $), the induction hypothesis
yields
$$ | e^{|\beta^i_{l}|r} \partial^{\gamma^i_{l}} x^t | \lesssim \scal{t},  $$
since, if $ \beta^i_l \ne 0 $, $
\partial_{\eta}^{\beta^i_l} x^t =
\partial_{\eta}^{\beta^i_l} (x^t-x^0) $.  If $ n+2 \leq i \leq 2n $ (and $ \nu_i \ne 0
$), we also have
$$ | e^{|\beta^i_{l}|r} \partial^{\gamma^i_{l}} \eta^t_{i-n-1} | \lesssim \scal{t} , $$
unless $ \partial^{\gamma^i_{l}} =
\partial_{\eta}^{\beta^i_{l}} $ with $ |\beta^i_{l} | = 1 $,
in which case we only have $ |\partial^{\gamma^i_{l}}
\eta^t_{i-n-1} | \lesssim \scal{t} $. By setting
$$ {\mathcal E} =  \{ n+2 \leq i \leq 2n \ ; \ \exists 1 \leq l \leq \nu_i
\ \mbox{such that} \
\partial^{\gamma^i_{l}} =
\partial_{\eta}^{\beta^i_{l}} \ \mbox{with} \ |\beta^i_{l} | = 1
\} ,
$$
and $ N = \# {\mathcal E} $, we thus obtain
$$  |(\ref{combinatoirederivee})| \lesssim e^{N  r } | (\partial^{\nu} H_p )(
\Phi^t) | \scal{t}^{|\nu| - N} \prod_{\mathcal E} |
\partial^{\gamma^i_{l}} \eta^t_{i-n-1} | .
$$
Since the components of $ H_p $ are polynomial of degree $2$ with
respect to the last $n-1$ variables, we only need to consider the
case where $ N \leq 2 $, otherwise $ \nu_{n+2} + \cdots + \nu_{2n}
\geq 3 $ and $
\partial^{\nu} H_p \equiv 0$. Furthermore
$$ | (\partial^{\nu} H_p )(
\Phi^t) | \lesssim e^{-2 r^t} \scal{\eta^t}^{2 - \nu_{n+2} -
\cdots - \nu_{2n} } \lesssim e^{-2 r^t} \scal{\eta^t}^{2 - N } .
$$
For $ N \leq 2 $, we have $ \scal{\eta^t}^{2-N} \lesssim
\scal{\eta}^{2-N} + \scal{t}^{2-N} $ so, using that  $e^{Nr} e^{-2
r^t} \lesssim e^{- (2-N)r - 2 \sigma t}$, we see that $ e^{Nr}
e^{-2r^t} \scal{\eta^t}^{2-N}  \lesssim e^{- \sigma t} $ which
finally implies
$$ |(\ref{combinatoirederivee})| \lesssim  \scal{t}^{|\nu|} e^{- \sigma t}  \lesssim e^{-  \sigma t/2} . $$
Therefore $ | Y (t) | \lesssim \scal{t} e^{- \sigma t } $ and, by
Lemma \ref{Gronwallexponentiel}, $ |X(t)| $ is bounded . The
result follows then easily. \finpreuve

\bigskip

\begin{lemm} \label{lemmetempsdeltafini} For all $ 0  < \sigma  < 1 $, there exist $ R  > 0 $
and   $ C > 0 $   such that, for all $ (r,\theta,\rho,\eta) $ satisfying
(\ref{motivezones}),
\begin{eqnarray}
| \rho^t \mp p^{1/2} | \leq C e^{-|t|/C}, \qquad \pm t \geq 0 .
\label{ecartrt}
\end{eqnarray}
 In particular, $ \rho^t \rightarrow \pm
p^{1/2} $ as $ t \rightarrow \pm \infty $.
\end{lemm}

\noindent {\it Proof.} We consider the case where $ t \geq 0 $,
the case of negative times being similar. Using (\ref{deriplus}),
Proposition \ref{localisationr} and Lemma \ref{lemmemotivezones},
we can choose $ R $ large enough such that $ \dot{\rho}^t \geq 0 $
and
\begin{eqnarray}
  \dot{\rho}^t   \lesssim  e^{-2 r^t} |\eta^t|^2  \lesssim
e^{-2 r^t} ( |\eta| + \scal{t} )^2 \lesssim e^{-2 r - 2 \sigma t} ( |\eta| +
\scal{t} )^2  \lesssim e^{- \sigma t}
  , \label{estimatefrom}
\end{eqnarray}
using the fact that $ e^{-2r}|\eta|^2 \lesssim p  $ in
 the last estimate. Therefore, $ \rho^t $ has a limit as $ t \rightarrow + \infty $.
 By the energy conservation and the estimate on $ e^{-2 r^t}
|\eta^t|^2 $ given by (\ref{estimatefrom}), we have $ p =
(\rho^t)^2 + {\mathcal O} (e^{-\sigma t}) $, which  shows that $
(\rho^t)^2 \rightarrow p $.  Since $ \rho^t $ is non decreasing
and $  \rho^0 = \rho
> - p^{1/2} $, the limit must be $ p^{1/2} $. Then we get (\ref{ecartrt}) by integrating the motion equation for
 $ \rho^t $ between $t$ and $ + \infty
$, namely
\begin{eqnarray}
p^{1/2} - \rho^t = \int_t^{\infty} \dot{\rho}^s ds =
\int_t^{\infty} e^{-2 r^s} \left(2 q (r^s,\theta^s,\eta^s) -
(\partial_r q_1) (r^s,\theta^s,\eta^s) \right) ds
\label{utilisationrepetee}
\end{eqnarray}
where, by Proposition \ref{localisationr} and Lemma
\ref{lemmemotivezones}, the integrand is $ {\mathcal O}(e^{-2 r -
2 \sigma s } (\scal{s}+\scal{\eta})^2) $. \finpreuve
 
\bigskip

\begin{prop} \label{estimeesflotprecises} For all $ 0  < \sigma  < 1 $, there
  exists $ R  > 0 $
such that, for all $ j ,k\in \Na $, $ \alpha , \beta \in \Na^{n-1}
$, with the notation
 $$ D_{\rm hyp}^{j \alpha k \beta} =
e^{r|\beta|} \partial_{\eta}^{\beta} \partial_r^j \partial_{\theta}^{\alpha}
\partial_{\rho}^k  , $$ (introduced before Definition \ref{definhyp}) and $ (l)_+ = \max (0,l) $, we have
\begin{eqnarray*}
|  D_{\rm hyp}^{j \alpha k \beta}  ( r^t - r - 2 |t| p^{1/2} )  |
& \lesssim & \left( e^{-r} \scal{\eta / p^{1/2}} \right)^{(2 -
|\beta|)_+} p^{-(k+|\beta|)/2} ,
\\
|  D_{\rm hyp}^{j \alpha k \beta} ( \theta^t - \theta   )  | &
\lesssim & e^{-r} \left( e^{-r} \scal{\eta / p^{1/2}}
 \right)^{(1-|\beta|)_+} p^{-(k+|\beta|)/2} , \\
 | D_{\rm hyp}^{j \alpha k \beta} ( \rho^t - \rho  )  | + |  D^{j \alpha k \beta}_{\rm hyp} ( \eta^t - \eta  )  |
   & \lesssim &  \left( e^{-r} \scal{\eta /
p^{1/2}} \right)^{(2 - |\beta|)_+} p^{(1-k-|\beta|)/2} ,
\end{eqnarray*}
 and, for  all $ 0 < \varepsilon < 1 $,
\begin{eqnarray*}
 | D_{\rm
hyp}^{j \alpha k \beta} ( \rho^t \mp p^{1/2}  )  | & \lesssim &
\left( e^{-r} \scal{\eta / p^{1/2} } \right)^{(2 - |\beta|)_+}
e^{- 4 (1 - \varepsilon) |t| p^{1/2}} p^{(1-k-|\beta|)/2} ,
\end{eqnarray*}
uniformly with respect to $ (r,\theta,\rho,\eta) $ and $t$
satisfying
\begin{eqnarray}
 r   >  R , \qquad \pm \rho  >  - \sigma p^{1/2} , \qquad \pm t
\geq 0 . \label{invariantparechelonnage}
\end{eqnarray}
\end{prop}

We point out that, apart from the energy localization and the localization in $ \theta $, the conditions (\ref{invariantparechelonnage})
are the main ones that define outgoing/incoming areas according to Definition
\ref{definitsortantentrantiota}.

Note also that, 
if $ (r,\theta,\rho,\eta) $ are
restricted to a subset where $ p $ belongs to a compact subset of
$ (0,+\infty) $, the estimates of Proposition
\ref{estimeesflotprecises} read
\begin{eqnarray}
|  D_{\rm hyp}^{j \alpha k \beta} ( r^t - r - 2 | t | p^{1/2} )  |
+ | D_{\rm hyp}^{j \alpha k \beta} ( \rho^t - \rho  )  | + |
D_{\rm hyp}^{j \alpha k \beta} ( \eta^t - \eta ) |
  \lesssim  \left( e^{-r} \scal{\eta} \right)^{(2 - |\beta|)_+} ,  \label{flotrrhoeta} \\
|  D_{\rm hyp}^{j \alpha k \beta} ( \theta^t - \theta   )  |
\lesssim e^{-r} \left( e^{-r} \scal{\eta} \right)^{(1-|\beta|)_+}
, \label{flottheta} \qquad
\qquad \qquad \qquad \\
|  D_{\rm hyp}^{j \alpha k \beta} ( \rho^t \mp p^{1/2} )  |
\lesssim \left( e^{-r} \scal{\eta} \right)^{(2-|\beta|)_+ }  e^{-
4 (1 - \varepsilon) | t | p^{1/2}} . \label{flotrhop}
\end{eqnarray}
Actually the latter estimates are equivalent to Proposition
\ref{estimeesflotprecises} by the following elementary scaling
properties,
\begin{eqnarray}
( r^t , \theta^t ) (r,\theta,\rho,\eta) & = & (r^{\lambda t},
\theta^{\lambda t}) (r,\theta,\rho/\lambda,\eta/\lambda) ,
\label{scalingrtheta} \\
( \rho^t , \eta^t ) (r,\theta,\rho,\eta) & = & \lambda
(\rho^{\lambda t}, \eta^{\lambda t})
(r,\theta,\rho/\lambda,\eta/\lambda) , \label{scalingrhoeta}
\end{eqnarray}
for $ \lambda > 0 $. Note that the condition
(\ref{invariantparechelonnage}) is invariant under the scaling $
(t,\rho,\eta) \mapsto (\lambda t, \rho/ \lambda, \eta / \lambda)
$.

\bigskip

\noindent {\it Proof.} We only need to prove (\ref{flotrrhoeta}),
(\ref{flottheta}) and (\ref{flotrhop}) with  $ p \in (1/4,4) $
and, again, we only consider $ t \geq 0$ and $ \rho > - \sigma
p^{1/2} $. We first assume that $ j + |\alpha| + k + |\beta| = 0
$. By (\ref{Hamigene}), Proposition \ref{localisationr} and Lemma
\ref{lemmemotivezones}, we have
\begin{eqnarray*}
 | \dot{\theta}^t | & \lesssim &
e^{-2r - 2 \sigma t} ( |\eta| + \scal{t} ) \lesssim e^{-2 r -
\sigma t} \scal{\eta}, \\
| \dot{\eta}^t | & \lesssim & e^{-2r - 2 \sigma t} ( |\eta| +
\scal{t} )^2 \lesssim e^{-2 r - \sigma t} \scal{\eta}^2 ,
\end{eqnarray*}
hence $ |\eta^t - \eta| \lesssim e^{-2r} \scal{\eta}^2 $ and $
|\theta^t - \theta| \lesssim e^{-2r} \scal{\eta} $. In particular,
$ \eta^t - \eta $ and $ \theta^t - \theta $ are bounded. The
motion  equation for $ r^t $ yields
\begin{eqnarray}
 r^t - r - 2 t p^{1/2} & = & 2 \int_0^t (\rho^s - p^{1/2}) ds ,
 \label{flotintegresurr}
\end{eqnarray}
and, using (\ref{ecartrt}), we get $ | r^t - r - 2 t p^{1/2} |
\lesssim 1 $. The latter estimate, the boundedness  $ | \eta^t -
\eta | $ and (\ref{utilisationrepetee}) imply
\begin{eqnarray}
 | \rho^t - p^{1/2} | \lesssim e^{-2 r - 4 t p^{1/2}}
\scal{\eta}^2 . \label{estimationprecisee}
\end{eqnarray}
 Furthermore, since $ |p^{1/2} - \rho|  = |\rho^2
- p| / |\rho + p^{1/2} |\lesssim e^{-2r} |\eta|^2 $, we also have
$ |\rho^t - \rho |  \lesssim e^{-2r} \scal{\eta}^2 $. Putting
(\ref{estimationprecisee}) into (\ref{flotintegresurr}), we obtain
$ |r^t - r - 2 t p^{1/2}| \lesssim e^{-2r} \scal{\eta}^2 $ which
completes the proof of (\ref{flotrrhoeta}), (\ref{flottheta}) and
(\ref{flotrhop}) for $ j + |\alpha| + k + |\beta| = 0 $ (note that
we can choose $ \varepsilon = 0 $ in this case).

Let us now prove (\ref{flottheta}) when $ j + |\alpha| + k +
|\beta| \geq 1 $. We first note that, by Lemma
\ref{lemmemotivezones} and the boundedness of $ |r^t - r - 2 t
p^{1/2}| $, we have
\begin{eqnarray}
  | D_{\rm hyp}^{j^{\prime} \alpha^{\prime} k^{\prime} \beta^{\prime}} ( e^{-r^t}\eta^t ) | & \leq &
| D_{\rm hyp}^{j^{\prime} \alpha^{\prime} k^{\prime}
\beta^{\prime}} ( e^{-r^t}( \eta^t - \eta ) )| +
| D_{\rm hyp}^{j^{\prime} \alpha^{\prime} k^{\prime} \beta^{\prime}} ( e^{-r^t}\eta ) | , \nonumber \\
&\lesssim & e^{- 2 t p^{1/2}}
\scal{t}^{j^{\prime}+|\alpha^{\prime}|+k^{\prime}+|\beta^{\prime}|}
\left( e^{-r} + (e^{-r} |\eta|)^{(1-|\beta^{\prime}|)_+} \right) ,
\nonumber \\ & \lesssim & e^{- 2 t p^{1/2}}
\scal{t}^{j^{\prime}+|\alpha^{\prime}|+k^{\prime}+|\beta^{\prime}|}
\left( e^{-r} \scal{\eta} \right)^{(1-|\beta^{\prime}|)_+}
\label{estimeederiveesansdiff} ,
\end{eqnarray}
for all $ j^{\prime} + |\alpha^{\prime}| + k^{\prime} +
|\beta^{\prime}| \geq 0 $. By writing
$$ \theta^t - \theta = \int_0^t e^{-r^s}(\partial_{\eta} q) (r^s,\theta^s,e^{-r^s}\eta^s) ds, $$
and using  (\ref{estimeederiveesansdiff}), Lemma
\ref{lemmemotivezones} (ie $
 | D_{\rm hyp}^{j^{\prime \prime} \alpha^{\prime \prime} k^{\prime \prime} \beta^{\prime \prime}} r^t | +
 | D_{\rm hyp}^{j^{\prime \prime} \alpha^{\prime \prime} k^{\prime \prime} \beta^{\prime \prime}} \theta^t | \lesssim
 \scal{t}
 $ if $ j^{\prime \prime} + |\alpha^{\prime \prime}| + k^{\prime \prime} + |\beta^{\prime \prime}| \ne 0 $),
 the Leibniz formula
 and Lemma \ref{Faadibruno}, we  obtain (\ref{flottheta}).
 We obtain similarly (\ref{flotrhop}) and then (\ref{flotrrhoeta})
 (also  using that $ (e^{-r}\scal{\eta})^2 \lesssim
e^{-r}\scal{\eta} \lesssim 1 $). Note that, for $ r^t - r - 2 t
p^{1/2} $,  (\ref{flotrrhoeta})  follows directly from
(\ref{flotrhop}) and (\ref{flotintegresurr}). \finpreuve

\bigskip

\begin{coro} \label{dependencedudomaine} Let $ V \Subset V^{\prime} \Subset \Ra^{n-1} $ be two
relatively compact open subsets and let $ 0 < \sigma < 1 $. There
exists $ R > 0 $ and $ C > 0 $ such that the conditions
\begin{eqnarray}
 r > R, \qquad \theta \in V, \qquad \pm \rho > - \sigma p^{1/2},
 \label{zonequasisortante}
\end{eqnarray}
imply that, for all $ \pm t \geq 0 $,
$$ r^t > r - C, \qquad \theta^t \in V^{\prime} . $$
In particular, if (\ref{zonequasisortante}) holds, the flow $
\Phi^t(r,\theta,\rho,\eta) $ depend only on  $ p $ on $ T^* \left(
(r-C , + \infty) \times V^{\prime} \right)  $ for $ \pm t \geq 0
$.
\end{coro}

This corollary allows us to localize the estimates  of Proposition
\ref{estimeesflotprecises} in charts of asymptotically hyperbolic
manifolds.

\section{The Hamilton-Jacobi and transport equations} \label{sectionHJ}
\setcounter{equation}{0} In this section, we develop the analytical tools necessary for  the Isozaki-Kitada parametrix that will be
constructed in  Section \ref{IsozakiKitada}. 

As in subsection \ref{wpasquelconque}, we consider
\begin{eqnarray}
 p = p(r,\theta,\rho,\eta) := \rho^2 + e^{-2r} q (r,\theta,\eta)
, \label{symboleprincipalglobal}
\end{eqnarray}
 where $ q$ is defined by (\ref{pourlasectionsuivante}) with
$q_0, q_1 $ satisfying (\ref{pourlasectionsuivante2}),
(\ref{longrange}), (\ref{secogene}) and
(\ref{ellipiticiteuniformeinfini}). 
Although the functions $ q $ and $q_0 $ are well defined for $ \theta \in
\Ra^{n-1} $, we shall mainly work with $
\theta $ in a neighborhood  of a fixed bounded open set
\begin{eqnarray}
 V_0 \Subset \Ra^{n-1}_{\theta} . \nonumber
\end{eqnarray}
Here $ V_0 $ is
arbitrary. In the applications, $ p $ will be replaced by $ p_{\iota} = \rho^2 + e^{-2r} q_{\iota} (r,\theta,\eta) $, 
the principal symbol of $ P $ in the patch $ \Psi_{\iota}({\mathcal U}_{\iota}) $,
and $ V_0 $ will be replaced by $ V_{\iota} $ (see (\ref{margeouverts})) contained in  a coordinate
patch of the angular manifold. Although $ p_{\iota} $ is not defined on $ T^* \Ra^n_+ $, this does not cause any problem since, as proved by Corollary \ref{dependencedudomaine} and as will be clear from the analysis below, the constructions of phases, flows, etc... depend on  $ p_{\iota} $ only through
its values for  $ \theta $ in  an arbitrary neighborhood of $ V_{\iota} $ and $ r \gg 1 $. Therefore, there will be no loss of generality in assuming that $ p_{\iota} $ is the restriction to a neighborhood of $ [ R, + \infty ) \times \overline{V}_{\iota} \times \Ra^n $ of a symbol of the form (\ref{symboleprincipalglobal}). The construction of such a continuation of $ p_{\iota} $ is fairly standard using (\ref{longueportee}).




\subsection{Properties of outgoing, incoming and intermediate areas} \label{preuvesannexessortantes}





\begin{defi} For all $ R > 0 $, all relatively compact open subset $ V \Subset \Ra^{n-1} $, all $ \sigma \in (-1,1)
$ and all open interval $ I \Subset (0,+\infty) $, we define
\begin{eqnarray}
\Gamma^{\pm} (R,V,I,\sigma) = \{ (r,\theta,\rho,\eta) \in \Ra^{2n}
\ | \ r > R , \ \theta \in V , \ p \in I , \ \pm \rho > - \sigma
p^{1/2} \} .
\end{eqnarray}
 The open
set $ \Gamma^+ (R,V,I,\sigma) $ (resp. $ \Gamma^- (R,V,I,\sigma)
$) is called an outgoing (resp. incoming) area.
\end{defi}

 Definition \ref{definitsortantentrantiota} is of course a special case of the above one with $ p = p_{\iota} $,
 $ V \Subset V_{\iota}^{\prime} $ (see (\ref{margeouverts})) and $ R $ large enough.

Note that $ \Gamma^{\pm}(R,V,I,\sigma) \subset
\Gamma^{\pm}(R,V,I,|\sigma|) $. 
Therefore, by Proposition \ref{propgene}, outgoing (resp.
incoming) areas are regions where the forward (resp. backward)
flow of $p$ is well defined and non trapping in the future (resp.
in the past) in the sense that $ r^t \rightarrow + \infty $ as $ t
\rightarrow + \infty $ (resp $ t \rightarrow - \infty $).

We now prove two useful propositions on outgoing/incoming areas using the classes $ {\mathcal S}_{\rm hyp}(\Omega) $ introduced
in Definition \ref{definhyp}.

\begin{prop} \label{partitions} i) Assume that
\begin{eqnarray}
  R_1 > R_2, \qquad V_1 \Subset V_2, \qquad I_1 \Subset I_2, \qquad \sigma_1 < \sigma_2 . \label{strongerassumption}
\end{eqnarray}
Then we can find $ \chi_{1 \rightarrow 2}^{\pm} \in {\mathcal S}_{\rm hyp} \left( \Gamma^{\pm} (R_2,V_2,I_2,\sigma_2) \right) $ such that
$$ \chi_{1 \rightarrow 2}^{\pm} \equiv 1 \qquad \mbox{on} \ \ \Gamma^{\pm} (R_1,V_1,I_1,\sigma_1) . $$
ii) Any symbol $ a \in {\mathcal S}_{\rm hyp} \left( (R,+\infty) \times V \times \Ra^n \cap p^{-1}(I) \right) $ can be written
$$ a = a^+ + a^-, \qquad \mbox{with} \ \ \ a^{\pm} \in {\mathcal S}_{\rm hyp} \left(  \Gamma^{\pm} (R,V,I,1/2) \right) . $$
\end{prop}
One important point in this proposition is that $ \chi_{1 \rightarrow 2}^{\pm} $ and $ a^{\pm} $ can be chosen in $ {\mathcal S}_{\rm hyp} $.

\medskip

\noindent {\it Proof.} {\it i)} We may for instance choose
$$  \chi_{1 \rightarrow 2}^{\pm}(r,\theta,\rho,\eta) = \chi_{R_1 \rightarrow R_2}(r) \chi_{V_1 \rightarrow V_2}(\theta)
\chi_{I_1 \rightarrow I_2}(p) \chi_{\sigma_1 \rightarrow \sigma_2} ( \pm \rho / p^{1/2}),  $$
with $ \chi_{R_1 \rightarrow R_2} , \chi_{\sigma_1 \rightarrow \sigma_2} \in C^{\infty}(\Ra) $, $ \chi_{V_1 \rightarrow V_2} \in C_0^{\infty}(V_2) $ and $ \chi_{I_1 \rightarrow I_2} \in C_0^{\infty}(I_2) $ such that
$$ \mbox{supp}(\chi_{R_1 \rightarrow R_2}) \subset (R_2 , + \infty), \qquad \mbox{supp}(\chi_{\sigma_1 \rightarrow \sigma_2}) \subset (- \sigma_2 , + \infty) $$
and
$$ \chi_{R_1 \rightarrow R_2} \equiv 1 \ \ \mbox{on} \ (R_1 , + \infty), \ \ \ \chi_{V_1 \rightarrow V_2} \equiv 1 \ \ \mbox{on} \ V_1 , \ \ \
 \chi_{I_1 \rightarrow I_2} \equiv 1 \ \ \mbox{on} \ I_1 , \ \ \ \chi_{\sigma_1 \rightarrow \sigma_2} \equiv 1 \ \ \mbox{on} \ (-\sigma_1 , + \infty ). $$
 Notice that $  \rho / p^{1/2}$ is smooth on the support of $\chi_{I_1 \rightarrow I_2}(p) $. So defined $ \chi_{1 \rightarrow 2}^{\pm} $ is smooth on $ \Ra^{2n} $, supported in $ \Gamma^{\pm} (R_2,V_2,I_2,\sigma_2) $, $ \equiv 1 $ on $ \Gamma^{\pm} (R_1,V_1,I_1,\sigma_1) $ and one easily checks that it belongs to $ {\mathcal B}_{\rm hyp} \left(  \Gamma^{\pm} (R_2,V_2,I_2,\sigma_2) \right) $, using for instance Lemma \ref{lemmesymbole}.

\noindent {\it ii)} It is very similar to the first case. We may for instance choose
$$ a^{\pm} (r,\theta,\rho,\eta) = a (r,\theta,\rho,\eta) \chi^{\pm}_{ 1/2} (\rho / p^{1/2}) ,  $$
with $ \chi^{\pm}_{-1/2 \rightarrow 1/2} \in C^{\infty}(\Ra) $ such that
$$ \chi^{+}_{1/2} + \chi^{-}_{1/2} \equiv 1 , \qquad \mbox{supp}(\chi^{+}_{1/2}) \subset (-1/2 , + \infty ), \qquad \mbox{supp}(\chi^{+}_{1/2}) \subset (- \infty ,1/2 ) . $$
Here again $ \rho / p^{1/2} $ is smooth on the support of $a$ and $ a^{\pm} \in
{\mathcal B}_{\rm hyp} \left(  \Gamma^{\pm} (R,V,I,1/2) \right) $.
\finpreuve

\bigskip

Let us remark that, by the first part of Proposition \ref{partitions}, $ \Gamma^{\pm} (R_2,V_2,I_2,\sigma_2) $ is a neighborhood of the closure of $ \Gamma^{\pm} (R_1,V_1,I_1,\sigma_1)  $ under the assumption (\ref{strongerassumption}). In the following proposition, we make this remark more quantitative.
\begin{prop} \label{geometriebornee} Assume (\ref{strongerassumption}). There exists $ \varepsilon > 0 $ such that, for all $ (r^{\prime},\theta^{\prime},\rho^{\prime},\eta^{\prime}) \in \Ra^{2n} $  and all  $ (r,\theta,\rho,\eta) \in \Gamma^{\pm}
 (R_1,V_1,I_1,\sigma_1) $,
$$ | (r,\theta,\rho,\eta) - (r^{\prime},\theta^{\prime},\rho^{\prime},\eta^{\prime}) | \leq \varepsilon  \qquad \Rightarrow \qquad (r^{\prime},\theta^{\prime},\rho^{\prime},\eta^{\prime}) \in \Gamma^{\pm} (R_2,V_2,I_2,\sigma_2) . $$
\end{prop}

\noindent {\it Proof.} Choose first $ \varepsilon_0 > 0 $ such that, if $  |r-r^{\prime}| + |\theta - \theta^{\prime} | \leq \varepsilon_0 $,   $ r^{\prime} > R_2 $ and $ \theta^{\prime} \in V_2 $. Then, by writing
\begin{eqnarray}
 q (r^{\prime},\theta^{\prime},e^{-r^{\prime}}\eta^{\prime}) - q (r^{\prime},\theta^{\prime}, e^{-r} \eta) = e^{-2 r^{\prime}} q (r^{\prime},\theta^{\prime},\eta^{\prime}-\eta) + ( e^{2(r-r^{\prime})} - 1 ) q (r^{\prime},\theta^{\prime}, e^{-r} \eta)  ,
 \label{algebreqiota}
\end{eqnarray}
and using (\ref{pourlasectionsuivante2}), (\ref{longrange}) with the Taylor formula, we get
$$ | p (r^{\prime},\theta^{\prime},\rho^{\prime},\eta^{\prime}) -  p (r,\theta,\rho,\eta) | \leq |\rho^2 - \rho^{\prime 2}|
+ C|\eta^{\prime} - \eta|^2 + C (|r-r^{\prime}| + |\theta - \theta^{\prime}|) e^{-2r}|\eta|^2 , $$
where $ e^{-2r}|\eta|^2 $ is bounded, using (\ref{secogene}). Since $ \rho $ is bounded too, we obtain
$$ | p (r^{\prime},\theta^{\prime},\rho^{\prime},\eta^{\prime}) -  p (r,\theta,\rho,\eta) |  \leq C |(r,\theta,\rho,\eta)-(r^{\prime},\theta^{\prime},\rho^{\prime},\eta^{\prime})| , $$
provided that $ |(r,\theta,\rho,\eta)-(r^{\prime},\theta^{\prime},\rho^{\prime},\eta^{\prime})| \leq \varepsilon_0 $ and therefore,
\begin{eqnarray*}
 | p^{1/2} (r^{\prime},\theta^{\prime},\rho^{\prime},\eta^{\prime}) -  p^{1/2} (r,\theta,\rho,\eta) | \leq C |(r,\theta,\rho,\eta)-(r^{\prime},\theta^{\prime},\rho^{\prime},\eta^{\prime}) | , \\
 \left|\frac{\rho^{\prime}}{ p^{1/2} (r^{\prime},\theta^{\prime},\rho^{\prime},\eta^{\prime})} - \frac{ \rho }{ p^{1/2} (r,\theta,\rho,\eta)} \right| \leq C |(r,\theta,\rho,\eta)-(r^{\prime},\theta^{\prime},\rho^{\prime},\eta^{\prime}) | ,
\end{eqnarray*}
if  $ |(r,\theta,\rho,\eta)-(r^{\prime},\theta^{\prime},\rho^{\prime},\eta^{\prime}) | $ is small enough,  using that $ I_2 \Subset (0,+\infty) $. The conclusion is then easy. \finpreuve

\bigskip

Similarly to (\ref{RVepsiloniota}), we set
\begin{eqnarray}
 R (\epsilon)  =  1/ \epsilon,   \qquad
 V_{\epsilon}  =   \{ \theta \in \Ra^{n-1} \ | \ \mbox{dist}(\theta,V_0) < \epsilon^2 \} . \label{RVepsilon}
\end{eqnarray}

\begin{defi} \label{zonesstrongsansiota} For all $ 0 < \epsilon < 1/4 $, we set
$$ \Gamma^{\pm}_{\rm s} (\epsilon) := \Gamma^{\pm} (R(\epsilon),V_{\epsilon}, (1/4-\epsilon, 4+\epsilon),\epsilon^2-1)
. $$ The open set $ \Gamma^{+}_{\rm s}(\epsilon) $ (resp. $
\Gamma^{-}_{\rm s} (\epsilon) $) is called a strongly outgoing
(resp. incoming) area.
\end{defi}
Here again, Definition \ref{zonesstrong} is a special case of the latter. We also recall that,
for all $ \epsilon $ small enough and all $ (r,\theta,\rho,\eta) \in \Gamma^{\pm}_{\rm s} (\epsilon) $, (\ref{crucialdiffeo}) holds, as explained after Definition \ref{zonesstrong}.

\medskip

In the sequel, we shall need very often the following result.

\begin{prop} \label{libre}  For all $ M > 0 $, there exist $ \epsilon_M > 0 $
and $ C_M > 1 $ such that, for all $ 0 < \epsilon \leq \epsilon_M
$, the following holds: if
\begin{eqnarray}
 (r, \theta , \rho , \eta) \in \Gamma^{\pm}_{\rm s}
(\epsilon) , \label{conditionslibres}
\end{eqnarray}
and
\begin{eqnarray}
r^{\prime} - r \geq - M , \qquad |\theta^{\prime}- \theta| < M \epsilon^2 ,
\qquad |\rho^{\prime} - \rho| < M \epsilon^2 , \qquad |\eta^{\prime}-\eta| < M
\epsilon e^{1/\epsilon} , \label{conditionsperturbatives}
\end{eqnarray}
then, for all $ 0 \leq s \leq 1 $,
$$ (r^{\prime},\theta^{\prime},\rho^{\prime},s \eta^{\prime}) \in \Gamma^{\pm}_{\rm s} ( C_M \epsilon )  . $$
In particular, $  ( r^{\prime}, \theta^{\prime} , \rho^{\prime} , 0 ) \in \Gamma^{\pm}_{\rm
s} (C_M \epsilon) $.
\end{prop}

\noindent {\it Proof.} Using (\ref{crucialdiffeo}) and (\ref{algebreqiota}), we first note the existence of $ M^{\prime}
> 0 $ such that, for all $ 0 < \epsilon < 1/4 $, if
(\ref{conditionslibres}) and (\ref{conditionsperturbatives}) hold
then
$$ | p (r^{\prime},\theta^{\prime},\rho^{\prime},s \eta^{\prime}) - p (r,\theta,\rho,\eta) | \leq M^{\prime} \epsilon^2 , $$
 using in particular that $ s \eta^{\prime} - \eta = s (\eta^{\prime} - \eta) + (s-1)\eta
$. If $ C_M $ is large enough and $0 <
\epsilon C_M < 1/4$, we obtain
\begin{eqnarray}
0 < \frac{1}{4} -  C_M \epsilon  < \frac{1}{4} - \epsilon -
M^{\prime} \epsilon^2
  \leq p (r^{\prime},\theta^{\prime},\rho^{\prime},s \eta^{\prime}) \leq 4 + \epsilon + M^{\prime}
 \epsilon^2 < 4 +  C_M \epsilon  . \nonumber
\end{eqnarray}
If $ 0 < \epsilon \leq \epsilon_M $ with $ 
\epsilon_M  $ small enough, then $  p (r^{\prime},\theta^{\prime},\rho^{\prime},s \eta^{\prime}) / p
(r,\theta,\rho,\eta) = 1 + {\mathcal O}(\epsilon^2) $ so
that
 \begin{eqnarray}
 \frac{ \pm \rho^{\prime}}{p (r^{\prime},\theta^{\prime},\rho^{\prime},s \eta^{\prime})^{1/2}} & = &
\frac{\pm \rho}{p (r,\theta,\rho,\eta)^{1/2}}
\frac{p (r,\theta,\rho,\eta)^{1/2}}{p (r^{\prime},\theta^{\prime},\rho^{\prime},s
\eta^{\prime})^{1/2}} \pm \frac{\rho^{\prime} - \rho}{p (r^{\prime},\theta^{\prime},\rho^{\prime},s
\eta^{\prime})^{1/2}} , \nonumber \\
 & > &  1 - ( C_M \epsilon)^2 , \nonumber
\end{eqnarray}
by possibly increasing $ C_M $.
In addition, $ \mbox{dist}(\theta,V_0) \leq |\theta^{\prime}-\theta| +
\mbox{dist}(\theta,V_0) < (C_M \epsilon )^2 $, by possibly
increasing $ C_M $ again and decreasing $ \epsilon_M $. Finally, $ r^{\prime}
\geq r - M > e^{1/\epsilon} - M > e^{1/C_M \epsilon} $, for all
$ 0 < \epsilon \leq \epsilon_M $ by possibly  decreasing $
\epsilon_M $ again, so $ (r^{\prime},\theta^{\prime},\rho^{\prime}, s \eta^{\prime}) \in
\Gamma^{\pm}_{\rm s}(C_M \epsilon ) $. \finpreuve

\bigskip

\begin{prop} \label{sortantversfortementsortant} Fix $ \sigma \in (-1,1) $.
Then, there exists $ R_{\sigma}^{\prime} >  0  $ such that for all $ R \geq R_{\sigma}^{\prime} $ and   all $ \epsilon >  0  $ small enough, there exists $
t_{R,\epsilon} \geq 0 $ such that
$$ \Phi^{t} \left( \Gamma^{\pm}(R ,V_0,(1/4 - \epsilon, 4 + \epsilon),\sigma)  \right) \subset \Gamma^{\pm}_{\rm s}(\epsilon), \qquad
 \mbox{for all} \ \pm t \geq t_{R,\epsilon} , $$
 where $ \Phi^t $ is the Hamiltonian flow of $p$.
\end{prop}

In other words, one can reach a strongly outgoing (resp. incoming)
area from an outgoing (resp. incoming) one  in finite time, along
the geodesic flow.


\bigskip

\noindent {\it Proof.} We consider only the outgoing case. With no
loss of generality, we may assume that $ 0 < \sigma  <  1 $. By choosing $ R \geq R_{\sigma}^{\prime} $ large enough, we can use Proposition 
\ref{localisationr}
and Corollary \ref{dependencedudomaine}. 
By Proposition \ref{localisationr}, we
have $ r_t \geq r + c t - C $ for some $ C,c > 0 $, hence $ r_t >
R (\epsilon) $ for all $ t \geq t_{R,\epsilon} $, provided
\begin{eqnarray}
 c t_{R,\epsilon}  -  C + R  >  R(\epsilon) . \label{conditionR1}
\end{eqnarray}
 By  Proposition \ref{estimeesflotprecises}, we have $ |\theta^t - \theta | \lesssim e^{-r} $ hence
  $ \theta^t \in V_{\epsilon}
$, for $ \epsilon $ small enough and all $t  \geq 0$, since $
e^{-1/\epsilon } \ll \epsilon^2 $.
Using (\ref{flotrhop}) and the energy conservation, we shall have $  \rho^t / p^{1/2}(r^t,\theta^t,\rho^t ,\eta^t) > 1 - \epsilon^2 $  provided for instance that
\begin{eqnarray}
 e^{-  p^{1/2} t_{R,\epsilon} } \leq \epsilon^3 ,  \label{conditionR4}
\end{eqnarray}
with $ \epsilon $ small enough.  Choosing  $ t_{R,\epsilon} $ so that
(\ref{conditionR1}) and (\ref{conditionR4}) hold, we get the
result. \finpreuve

\bigskip

\bigskip

We conclude this part with the following explicit construction of
cutoffs.

In Section \ref{IsozakiKitada}, we will need a result similar to part {\it i)} of Proposition \ref{partitions}. This is the purpose of the following result.

\begin{prop} \label{construitcutoff} We can find $ 0 <\nu <1 $ and a family of cutoffs $ \chi^{\pm}_{\epsilon^2 \rightarrow
\epsilon}\in {\mathcal S}_{\rm hyp} (\Gamma^{\pm}_{\rm  s} ( \epsilon^{1+\nu} ))  $, defined for all $
\epsilon $ small enough, such that,
\begin{eqnarray}
  \chi^{\pm}_{\epsilon^2 \rightarrow \epsilon} = 1 \ \ \mbox{on} \
\ \Gamma^{\pm}_{\rm s}(\epsilon^2) ,
\label{supports}
\end{eqnarray}
and, uniformly on $ \Ra^{2n} $,
\begin{eqnarray}
 \left| e^{-2r}|\eta|^j  \partial_{r,\theta,\rho,\eta}
\chi^{\pm}_{\epsilon^2 \rightarrow \epsilon}  \right| + \left|
e^{-2r}|\eta|^2 \partial_{\rho,\eta}
\partial_{r,\theta} \chi^{\pm}_{\epsilon^2 \rightarrow \epsilon}
 \right| \lesssim \epsilon^{1/2}  , \qquad j = 1,2. \label{pourKuranishi}
\end{eqnarray}
\end{prop}

That we can find, for each $ \epsilon $, $ \chi^{\pm}_{\epsilon^2 \rightarrow \epsilon}
\in {\mathcal S}_{\rm hyp} (\Gamma^{\pm}_{\rm  s} (
\epsilon^{1+\nu} ))  $ satisfying (\ref{supports})  would follow directly from Proposition
\ref{partitions}. The important additional point here is the control with
respect to $ \epsilon $ given by (\ref{pourKuranishi}). Note also that
the power $ 1/2 $ is essentially irrelevant: we only mean that the
left hand side of (\ref{pourKuranishi}) is uniformly small as $
\epsilon \rightarrow 0 $. This rather technical point will only be
used in Section \ref{IsozakiKitada} to  globalize suitably certain
phase functions.

\bigskip

 \noindent {\it Proof.} For $ 0 <\delta < 1 $ to be chosen later, we consider the
characteristic functions $ \bar{\chi}^I_{\epsilon^{1 + \delta}} $
and $ \bar{\chi}^V_{\epsilon^{2+\delta}} $ of $ (1/4 -
\epsilon^{1+\delta},4 + \epsilon^{1+\delta}) $ and $ V + B
(0,\epsilon^{2+\delta}) $ respectively. Let us  choose $ \zeta^I
\in C_0^{\infty} (\Ra) $, $ \zeta^V \in C_0^{\infty} (\Ra^{n-1}) $
both equal to $ 1 $ near $ 0 $, such that $ \int \zeta^I = \int
\zeta^V = 1 $ and set
\begin{eqnarray}
\chi^I_{\epsilon^{1 + \delta}} (\lambda)  & = & \int
\bar{\chi}^I_{\epsilon^{1 + \delta}} (\mu) \zeta^I \left(
\frac{\lambda-\mu}{\epsilon^{1 + 2 \delta}} \right) \epsilon^{- 1
- 2 \delta } d \mu
, \nonumber \\
\chi^V_{\epsilon^{2 + \delta}} (\theta) & = &  \int
\bar{\chi}^V_{\epsilon^{2 + \delta}} (\vartheta) \zeta^V \left(
\frac{\theta - \vartheta}{\epsilon^{2 + 2 \delta}} \right)
\epsilon^{-(n-1)(2 + 2 \delta)} d \vartheta . \nonumber
\end{eqnarray}
One then easily checks that, if $ \epsilon $ is small enough,
\begin{eqnarray*}
\chi^I_{\epsilon^{1 + \delta}} \equiv 1 \ \ \mbox{on} \ \ (1/4 -
\epsilon^2 , 4 + \epsilon^2) , & \qquad & \chi^I_{\epsilon^{1 +
\delta}} \equiv 0 \ \ \mbox{outside} \ \ (1/4 - \epsilon^{1 +
\frac{\delta}{4}} , 4 + \epsilon^{1 + \frac{\delta}{4}}) , \\
\chi^V_{\epsilon^{2 + \delta}} (\theta) = 1 \ \ \mbox{if} \ \
\mbox{dist}(\theta,V) < \epsilon^4 , & \qquad &
\chi^V_{\epsilon^{2 + \delta}}(\theta) = 0 \ \ \mbox{if} \ \
\mbox{dist}(\theta,V) \geq \epsilon^{2 + \frac{\delta}{2}} .
\end{eqnarray*}
Choosing $ \omega \in C^{\infty} (\Ra) $ supported in $ ( 1/4 , +
\infty ) $ such that $ \omega  = 1 $ near $ [ 1/3 , \infty ) $, we
now define
$$ \chi^{\pm}_{\epsilon^2
\rightarrow \epsilon} (r,\theta,\rho,\eta) =  \omega ( r /R
(\epsilon^{3/2}) ) \chi^V_{\epsilon^{2+\delta}}
(\theta)\chi^I_{\epsilon^{1+\delta}} (p) \omega (\pm  \rho )
\zeta^I ( e^{-2r}|\eta|^2 / \epsilon^{4-\delta} ).  $$
 On the support of $
\chi^I_{\epsilon^{1+\delta}} (p) $, we have $ \rho^2 \geq
1/4 - {\mathcal O}(\epsilon )$  so the factor $ \omega (\pm  \rho) $ only determines
the sign of $ \rho $. By
(\ref{equivalencepratique}) and (\ref{crucialdiffeo}), one sees
that (\ref{supports}) holds with $ \nu = \delta/2 $, if $
\epsilon $ is small enough. Furthermore, $ \chi^{\pm}_{\epsilon^2
\rightarrow \epsilon} $ is supported in $
\Gamma^{\pm}_{\rm  s} ( \epsilon^{1 + \nu} ) $
and belongs to $  {\mathcal B}_{\rm hyp}
(\Gamma^{\pm}_{\rm  s} ( \epsilon^{1 + \nu} )) $.

Let us prove (\ref{pourKuranishi}). Since $ e^{-2r}|\eta|^2
\lesssim \epsilon^{4-\delta}  $ on the support of $
\chi^{\pm}_{\epsilon^2 \rightarrow \epsilon} $, the first order
derivatives satisfy
\begin{eqnarray*}
 | \partial_{r} \chi^{\pm}_{\epsilon^2
\rightarrow \epsilon} | & \lesssim & R (\epsilon^{3/2})^{-1} +
 e^{-2r}|\eta|^2 (\epsilon^{-1-2\delta} + \epsilon^{-4 + \delta}) \lesssim 1 , \\
 | \partial_{\rho}
\chi^{\pm}_{\epsilon^2 \rightarrow \epsilon} | & \lesssim &
\epsilon^{-1-2\delta}   , \\
 | \partial_{\theta}
\chi^{\pm}_{\epsilon^2 \rightarrow \epsilon} | & \lesssim &
\epsilon^{-2 - 2 \delta} +  \epsilon^{-1-2\delta} e^{-2r}|\eta|^2  \lesssim \epsilon^{-2 - 2 \delta} , \\
 | \partial_{\eta}
\chi^{\pm}_{\epsilon^2 \rightarrow \epsilon} | & \lesssim &
e^{-2r}|\eta| ( \epsilon^{-1-2\delta} +  \epsilon^{-4 + \delta } )
\ll e^{-\epsilon^{- 1/2}} ,
\end{eqnarray*}
using the fact that $ e^{-2r}|\eta| \lesssim e^{-r} \leq e^{-
\epsilon^{-1} } $ for the last estimate. Similarly
\begin{eqnarray*}
| \partial_{\rho} \partial_{r,\theta} \chi^{\pm}_{\epsilon^2
\rightarrow \epsilon} | & \lesssim &  \epsilon^{-2-2 \delta}
\times \epsilon^{-1-2 \delta} = \epsilon^{- 3 - 4 \delta} ,
\\
| \partial_{\eta} \partial_{r,\theta} \chi^{\pm}_{\epsilon^2
\rightarrow \epsilon} | & \lesssim & e^{- \epsilon^{-1/2} }  .
\end{eqnarray*}
Since $ e^{-2 r}|\eta|^2 e^{-3 - 4 \delta} \lesssim \epsilon^{1 -
5 \delta } $ and $ e^{-2 r}|\eta| \ll e^{- \epsilon^{-1/2}} $, the
result follows  with $ \delta = 1/10 $ (hence with $ \nu = 1/20
$). \finpreuve

\bigskip

We finally consider the case of intermediate areas.

\begin{defi} Given $ V \Subset \Ra^{n-1} $, $ \epsilon > 0 $, $ \delta > 0 $ and $ \sigma_0 , \ldots , \sigma_L $ 
satisfying (\ref{conditionsurlessigmas}), (\ref{recouvrementinterval}), (\ref{tailleintervalle}), we set
$$ \Gamma^{\pm}_{\rm inter} (\epsilon,\delta;l) := \left\{ (r,\theta,\rho,\eta) \in \Ra^{2n} \ | \ r > R(\epsilon), \ \theta \in V
, \ p \in I(\epsilon), \  \pm \frac{\rho}{p^{1/2}} \in (-\sigma_{l+1}, -\sigma_{l-1}) \right\}  , $$
for $ 1 \leq l \leq L-1 $. 
\end{defi}
As previously, this definition is nothing but Definition \ref{definitionintermediaireoriginale},
with $  p$ and $ V $ instead of $ p_{\iota} $ and $ V_{\iota} $.

\begin{prop} \label{deuxiemedecoupagesansiota} Fix $ \epsilon > 0 $ small enough, $ \delta > 0 $ and $ \sigma_0 , \ldots , \sigma_L $  satisfying (\ref{conditionsurlessigmas}), (\ref{recouvrementinterval}) and (\ref{tailleintervalle}).
Then, any symbol
$$ a^{\pm} \in {\mathcal S}_{\rm hyp} (  \Gamma^{\pm}_{\iota}(R(\epsilon),V_{\iota},I,1/2) )  $$
can be written
\begin{eqnarray}
 a^{\pm} = a^{\pm}_{\rm s} + a_{1,{\rm inter}}^{\pm} + \cdots + a_{L-1,{\rm inter}}^{\pm}, \nonumber
\end{eqnarray}
with
\begin{eqnarray}
 a^{\pm}_{\rm s} \in {\mathcal S}_{\rm hyp} (  \Gamma^{\pm}_{\iota,{\rm s}} ( \epsilon) )  ,
\qquad a_{l,{\rm inter}}^{\pm} \in {\mathcal S}_{\rm hyp} (  \Gamma^{\pm}_{\iota,{\rm inter}} (\epsilon,\delta;l) ) \nonumber .
\end{eqnarray}
\end{prop}

\noindent {\it Proof.} 
By (\ref{conditionsurlessigmas}) and (\ref{recouvrementinterval}), we can find
$ \chi_{- \infty }, \chi_{+\infty} \in C^{\infty}(\Ra)$ and $ \chi_{l} \in C_0^{\infty}(-\sigma_{l+1},-\sigma_{l-1})$, for $1 \leq l \leq L - 1 $,
such that
$$  \mbox{supp} ( \chi_{-\infty} )\subset (- \infty , - \sigma_{L-1} ), \qquad \mbox{supp} (\chi_{+\infty} ) \in ( 1 -  \epsilon^2 , + \infty ) , $$
and
$$ \chi_{+\infty} + \sum_{l = 1}^{L-1} \chi_{l} + \chi_{- \infty} \equiv 1 \qquad \mbox{on} \ \ \Ra. $$
This simply relies on the overlapping property of the intervals in (\ref{recouvrementinterval}). We then obtain the result by considering
\begin{eqnarray*}
 a^{\pm}_{\rm s}(r,\theta,\rho,\eta) & = & a^{\pm}(r,\theta,\rho,\eta) \times \chi_{+\infty} (\pm \rho / p^{1/2}) ,  \\
 a_{l,{\rm inter}}^{\pm} (r,\theta,\rho,\eta) & = & a^{\pm}(r,\theta,\rho,\eta) \times \chi_{l} (\pm \rho / p^{1/2}) , \qquad 1 \leq l \leq  L-2, \\
 a_{L-1,{\rm inter}}^{\pm} (r,\theta,\rho,\eta) & = & a^{\pm}(r,\theta,\rho,\eta) \times (\chi_{L-1} + \chi_{-\infty} ) (\pm \rho / p^{1/2}) .
\end{eqnarray*}
since, in the definition of $ a_{L-1,{\rm inter}}^{\pm} $, the cutoff guarantees that $ \pm \rho / p^{1/2} < - \sigma_{L-2} $ and $ a^{\pm} $ that $  \pm \rho / p^{1/2} > - 1/2 = - \sigma_L$.  \finpreuve

\bigskip

We conclude this subsection with the following proposition giving the main property of intermediate areas.

\begin{prop} \label{dynamiquedeltasansiota} Fix   $ \underline{t} > 0 $. Then for all $ \epsilon > 0 $ small enough,
 we can find $ \delta > 0 $ small enough such that, for any choice of $ \sigma_0, \ldots , \sigma_L $ satisfying (\ref{conditionsurlessigmas}), (\ref{recouvrementinterval}) and (\ref{tailleintervalle}), we have
\begin{eqnarray}
\Phi^t \left( \Gamma^{\pm}_{{\rm inter}} (\epsilon,\delta;l) \right) \cap \Gamma^{\pm}_{{\rm inter}} (\epsilon,\delta;l) = \emptyset, \qquad  \pm t \geq \underline{t}, \label{dynamiquedeltasansiotaeffectif}
\end{eqnarray}
for all $ 1 \leq l \leq L-1 $.
\end{prop}

\noindent {\it Proof.} We consider the outgoing case, the incoming one being similar. Using
 Corollary \ref{dependencedudomaine}, we may assume that, if $ \epsilon $ is small enough, (\ref{pourzoneintermediaire}) holds
for any initial condition such that  $ r > R(\epsilon) $, $ \theta \in V $ and $ \rho \geq - p^{1/2}/2 $. In particular $ t \mapsto \rho^t $ 
is non decreasing for $ t \geq 0 $.
  Assume that $1/2 \leq \rho/p^{1/2} \leq 1 - (\epsilon/2)^2  $ and set 
$$ t_{\epsilon} = t_{\epsilon} (r,\theta,\rho,\eta) := \sup \left\{ t \geq 0 \ | \ \frac{\rho^s}{p^{1/2}} < \frac{\rho}{p^{1/2}} + \epsilon^4 \ \ \mbox{for all} \ s \in [0,t] \right\} . $$
Notice that $ t_{\epsilon} $ is finite by Lemma \ref{lemmetempsdeltafini} and that $ \rho^{t_{\epsilon}} = \rho  + p^{1/2} \epsilon^4 $.
If $ 1 - (\epsilon/2)^2 + \epsilon^4 \geq 1/2 $, we have $ |\rho^t / p^{1/2}| \leq 1 - (\epsilon/2)^2 + \epsilon^4 $ on $ [0,t_{\epsilon}) $.
Thus, if $ \epsilon $ is small enough (independent of $ (r,\theta,\rho,\eta) $),
we have $ (\rho^t)^2/p \leq 1 - (\epsilon/2)^2 $ for all $ t \in [ 0 , t_{\epsilon} ) $ and then, by (\ref{pourzoneintermediaire}) again, we have
$ \dot{\rho}^t \geq (\epsilon/2)^2 p $ on $ [ 0 , t_{\epsilon} ] $ so that
$$ \rho^{t_{\epsilon}} - \rho \geq (\epsilon/2)^2 p t_{\epsilon} . $$
This shows that $ t_{\epsilon} \leq \epsilon^4 / (\epsilon/2)^2 p = 4 \epsilon^2 / p $. Then, for  $ \epsilon $ small enough such that $ 4 \epsilon^2 / p \leq \underline{t} $ for all $ (r,\theta,\rho,\eta) $ in
\begin{eqnarray}
 \left\{ (r,\theta,\rho,\eta) \in \Ra^{2n} \ | \ r > R(\epsilon), \ \theta \in V
, \ p \in I(\epsilon), \  -1/2 \leq \frac{\rho}{p^{1/2}} \leq 1 - (\epsilon/2)^2  \right\}  , \label{borneensemblisteintermediaire}
\end{eqnarray}
and with $ \delta = \epsilon^4 /2 $,
we have $ \rho^t - \rho \geq 2 \delta p^{1/2} $ for all $ t \geq \underline{t} $. This implies (\ref{dynamiquedeltasansiotaeffectif}) since, 
for any choice of $ \sigma_0 , \ldots , \sigma_L $ and any $l$, $ \Gamma^{\pm}_{{\rm inter}} (\epsilon,\delta;l) $ is contained in (\ref{borneensemblisteintermediaire}). \finpreuve

\subsection{Hyperbolic long/short range symbols}
In this short subsection, we introduce the definitions of short/long range hyperbolic symbols which will be useful
for the resolution of transport equations in Subsection \ref{Transport}. We prove in passing Proposition \ref{diffeoexact} below 
which will be used at several places, in particular in Subsection \ref{HamiltonJacobi}.

\begin{defi} \label{definitshortrange} A smooth function $ a_{\pm} $ on $  \Gamma^{\pm}_{\rm
      s} ( \epsilon )  $ is said to be of hyperbolic short range if
\begin{eqnarray}
| \partial_r^{j} \partial_{\theta}^{ \alpha } \partial_{\rho}^k
\partial_{\eta}^{\beta} a_{\pm} (r,\theta,\rho,\eta) | \lesssim \scal{r- \log
\scal{\eta}}^{-
  \tau - 1 - j  }   , \qquad  \ (r,\theta,\rho,\eta) \in \Gamma^{\pm}_{\rm
      s} ( \epsilon ) ,
 \label{shortrangehyperbolique}
\end{eqnarray}
and of hyperbolic long range if
\begin{eqnarray}
| \partial_r^{j} \partial_{\theta}^{ \alpha } \partial_{\rho}^k
\partial_{\eta}^{\beta} a_{\pm}  (r,\theta,\rho,\eta) | \lesssim \scal{r- \log
\scal{\eta}}^{-
  \tau  - j  }   , \qquad  \ (r,\theta,\rho,\eta) \in \Gamma^{\pm}_{\rm
      s} ( \epsilon ) .
  \label{longrangehyperbolique}
\end{eqnarray}
\end{defi}

\bigskip

Notice that in this definition, we do not assume that $ a \in {\mathcal B}_{\rm
hyp} (\Gamma^{\pm}_{\rm s} ( \epsilon )) $. However, this will be the case in the applications and we now give
a simple criterion to check that a symbol $ a \in {\mathcal B}_{\rm
hyp} (\Gamma^{\pm}_{\rm s} ( \epsilon )) $ is of hyperbolic
short/long range. 

For $ \epsilon $ small enough,
by restricting $a$ to a smaller area $ \Gamma^{\pm}_{\rm
s}(\epsilon / C) $, with $ C > 1 $ large enough (or to $
\Gamma^{\pm}_{\rm s}(\epsilon^2 ) $, $ \Gamma^{\pm}_{\rm
s}(\epsilon^3 ) $ as it will be the case in the applications),
using Lemma \ref{lemmesymbole} and Proposition \ref{libre}, we have
\begin{eqnarray}
 a (r,\theta,\rho,\eta) = a (r,\theta,\rho,0) + \int_0^1
(\partial_{\xi} \tilde{a}) (r,\theta,\rho, s \xi)_{|\xi =
e^{-r}\eta}  ds \cdot e^{-r}\eta , \label{Taylorutile}
\end{eqnarray}
 where $ \tilde{a} $ belongs to $ C^{\infty}_b ( F_{\rm hyp}
(\Gamma^{\pm}_{\rm s}(\epsilon)) ) $ and $ (r,\theta,\rho,s\eta)
\in \Gamma^{\pm}_{\rm s}(\epsilon) $ if $ (r,\theta,\rho,\eta) \in
\Gamma^{\pm}_{\rm s}(\epsilon / C) $. Since, for all $ N > 0 $,
\begin{eqnarray}
| \partial_r^j \partial_{\eta}^{\beta} e^{-r} \eta | \lesssim
\scal{r - \log \scal{\eta}}^{-N}, \qquad (r,\theta,\rho,\eta) \in
\Gamma^{\pm}_{\rm s}(\epsilon), \nonumber
\end{eqnarray}
we obtain that, for $ a \in {\mathcal B}_{\rm hyp}
(\Gamma^{\pm}_{\rm  s} ( \epsilon )) $,
\begin{eqnarray}
a   \ \mbox{is of hyperbolic  short/long range in } \
\Gamma^{\pm}_{\rm s}(\epsilon / C)  \ \Leftrightarrow  \ a_{|\eta
= 0} \ \mbox{is of usual short/long range}
\label{caracterisationslrange}
\end{eqnarray}
in the sense that
$$ | \left(  \partial_r^j \partial_{\theta}^{\alpha}
\partial_{\rho}^k a \right) (r,\theta,\rho,0) | \lesssim \scal{r}^{-
  \tau  - j  }    , \qquad (r,\theta,\rho,0) \in \Gamma^{\pm}_{\rm s}(\epsilon ) , $$
in the long range case (recall that $ 0 < \tau \leq 1 $) and
$$ | \left(  \partial_r^j \partial_{\theta}^{\alpha}
\partial_{\rho}^k a \right) (r,\theta,\rho,0) | \lesssim \scal{r}^{-
  \tau - 1 - j  }    , \qquad (r,\theta,\rho,0) \in  \Gamma^{\pm}_{\rm s}(\epsilon ) , $$
in the short range case.

\bigskip

 To calculate $ a_{| \eta = 0} $ in some
applications, we shall use the following elementary result.

\begin{prop} \label{diffeoexact} For all $ r > 0 $, all $ \theta \in \Ra^{n-1}  $ and all $ \pm \rho > 0
$, we have, for all $ \pm t \geq 0 $,
\begin{eqnarray}
(r^t,\theta^t,\rho^t,\eta^t )_{| \eta = 0} &  = & ( r + 2 t \rho,
\theta , \rho , 0 ), \label{sansderiver} \\
\partial_{\eta} (r^t,\theta^t,\rho^t,\eta^t )_{|\eta = 0} & = & \left(
0,\int_0^t   e^{-2r - 4 s \rho}  \emph{hess}_{\eta}[q] (r+s\rho,\theta) ds,0, \emph{Id} \right) .
\label{sansderiver2}
\end{eqnarray}
where $ \emph{hess}_{\eta}[ q ](r,\theta) $ is the Hessian matrix of
$q$ with respect to $ \eta $ (which is independent of $ \eta $).
\end{prop}

\noindent {\it Proof.}  One simply checks that the right hand side
of (\ref{sansderiver}) is a solution to (\ref{Hamigene}) (with $
 w(r) = e^{-2r}$) for $ \pm t \geq 0 $. Applying then
$
\partial_{\eta} $ to  (\ref{Hamigene}), one sees easily as well that
 the right hand side of (\ref{sansderiver2}) is a solution to the corresponding system. \finpreuve

\bigskip

\noindent {\bf Remark.} If  $ \epsilon  $ is small enough then, on
$ \Gamma^{\pm}_{\rm
      s} ( \epsilon )  $, we have
\begin{eqnarray}
  r- \log \scal{\eta} \geq 0   .  \label{controlelongueportee}
\end{eqnarray}
In particular, in this region, $ \scal{r- \log \scal{\eta}} $ is
equivalent to the weight
$$ \scal{r- \log \scal{\eta}}_+ := \max( 1 , r- \log
  \scal{\eta} )  $$
which was introduced by Froese-Hislop in \cite{FH}. For the study of global in time estimates, which hope to consider in a separate paper, 
the resolvent estimates proved in \cite{Boucletresolvente} suggest that these
hyperbolic short/long range conditions would play the same role as
the usual Euclidean short/long range conditions used in \cite{BoucletTzvetkov2}.

\subsection{The Hamilton-Jacobi equation} \label{HamiltonJacobi}

In this subsection, we use the results of subsection \ref{wpasquelconque} to solve the time independent Hamilton-Jacobi equations
giving the phases of the Isozaki-Kitada parametrix.

\begin{lemm} \label{preparationphase}  There exists  $ 0   <
\epsilon_0 <  1   $  such that, for all $ 0 < \epsilon \leq
\epsilon_0 $ and all $ \pm  t \geq 0 $,  the map
$$ \Psi_t^{\pm} : (r,\theta,\rho,\eta) \mapsto (r,\theta,\rho^t,\eta^t) $$
is a diffeomorphism from $ \Gamma^{\pm}_{\rm s}(\epsilon)  $ onto
its range and
\begin{eqnarray}
\Gamma^{\pm}_{\rm s} (\epsilon^3)   \subset \Psi_t^{\pm} \left(
\Gamma^{\pm}_{\rm s}(\epsilon)
 \right), \qquad \mbox{for all} \ \pm t \geq 0 . \label{inclusionmesuree}
\end{eqnarray}
\end{lemm}

\noindent {\it Proof.} See Appendix \ref{Diff}.

\bigskip

The power $ \epsilon^ 3 $ in (\ref{inclusionmesuree})
 is not very important. It is only a rough explicit
quantitative bound for the size of a strongly outgoing (resp.
incoming) area which is contained in $ \Psi_t^{+} (\Gamma^+_{\rm
s} (\epsilon))  $ (resp. $ \Psi_t^{-} (\Gamma^-_{\rm s}
(\epsilon)) $) for all $ t \geq 0 $ (resp. $ t \leq 0 $).

\bigskip

The components of the inverse map $ ( \Psi_t^{\pm} )^{-1} $ are of
the form $ ( r , \theta , \rho_t , \eta_t ) $ with
$$ \rho_t = \rho_t (r,\theta,\rho,\eta) , \qquad \eta_t = \eta_t (r,\theta,\rho,\eta) . $$
Here we omit the $ \pm $ dependence for notational simplicity.
We thus have
\begin{eqnarray}
 \rho^t (r,\theta,\rho_t,\eta_t) = \rho, \qquad \eta^t
(r,\theta,\rho_t,\eta_t) = \eta , \label{ffinverse}
\end{eqnarray}
 at least for all $
(r,\theta,\rho,\eta) \in \Gamma^{\pm}_{\rm s}( \epsilon_0^3 ) $
and $ \pm t \geq 0 $.

\bigskip

\noindent {\bf Remark.} It follows from the Proof of Lemma
\ref{preparationphase} and the scaling properties
(\ref{scalingrtheta}), (\ref{scalingrhoeta}) that $ \Psi_t^{\pm} $
is actually a diffeomorphism from the cone generated by $
\Gamma^{\pm}_{\rm s}(\epsilon_0) $ onto its range, the latter
range containing  the cone generated by $ \Gamma^{\pm}_{\rm
s}(\epsilon_0^3) $. Therefore
 $( \rho_t , \eta_t )$ is actually the restriction to $ \Gamma^{\pm}_{\rm s}(\epsilon_0^3) $ of a map defined on the cone
generated by $ \Gamma^{\pm}_{\rm s}(\epsilon_0^3) $ and, using
(\ref{scalingrhoeta}), we  have
\begin{eqnarray}
( \rho_t , \eta_t ) (r,\theta,\lambda \rho , \lambda \eta) =
\lambda (\rho_{\lambda t},\eta_{\lambda t}) (r,\theta,\rho,\eta) ,
\qquad   \ \pm t \geq 0 , \ (r,\theta,\rho,\eta) \in
\Gamma^{\pm}_{\rm s}(\epsilon_0^3), \label{scalingreciproque}
\end{eqnarray}
 for all $ \lambda > 0 $.

\begin{prop} \label{estimepreparephase}
There exists $ \epsilon_1 \leq \epsilon_0^3  $  such that,  for
all $ j,k \in \Na $, $ \alpha,\beta \in \Na^{n-1} $,
\begin{eqnarray}
  \big| D_{\rm hyp}^{j \alpha k \beta}(\rho_t - \rho) \big| + \big| D_{\rm hyp}^{j
\alpha k \beta} (\eta_t - \eta) \big|  \lesssim 1 , \qquad
(r,\theta,\rho,\eta) \in \Gamma^{\pm}_{\rm s}(\epsilon_1) , \ \
\pm t \geq 0 . \label{estimeinversepsit}
\end{eqnarray}
 In addition, if $ (r,\theta,\rho,0) \in \Gamma^{\pm}_{\rm s}(\epsilon_1) $, we have
\begin{eqnarray}
(\rho_t,\eta_t)_{| \eta = 0} = (\rho,0), \label{implicite1}  \\
\partial_{\eta} (\rho_t - \rho,\eta_t - \eta)_{| \eta = 0} = (0,0) . \label{implicite2}
\end{eqnarray}
\end{prop}

\noindent {\it Proof.}  By (\ref{inclusionmesuree}), any $
(r,\theta,\rho,\eta) \in \Gamma^{\pm}_{\rm s}(\epsilon_0^3)  $ can
be written $ \Psi_t^{\pm} (r,\theta,\tilde{\rho},\tilde{\eta})  $
with $ ( r,\theta,\tilde{\rho},\tilde{\eta} ) \in
\Gamma^{\pm}_{\rm s}(\epsilon_0)  $, hence
$$ \sup_{\Gamma^{\pm}_{\rm s}(\epsilon_0^3)}  |\rho_t - \rho| + |\eta_t - \eta| \leq
\sup_{\Gamma^{\pm}_{\rm s}(\epsilon_0)}  |\tilde{\rho}- \rho^t
(\tilde{r},\tilde{\theta},\tilde{\rho},\tilde{\eta})| +
|\tilde{\eta}- \eta^t
(\tilde{r},\tilde{\theta},\tilde{\rho},\tilde{\eta})| . $$ By
(\ref{flotrrhoeta}), the right hand side is bounded so   we obtain
(\ref{estimeinversepsit}) for $  j+ |\alpha|+ k+  |\beta| = 0  $.
Then, for $ \epsilon $ small enough, using Proposition
\ref{estimeesflotprecises}  and Lemma \ref{diffeoexact}, we
 remark that, for $ (r,\theta,\rho,\eta) \in \Gamma^{\pm}_{\rm
s}(\epsilon)  $,
$$ \left|  \partial_{\rho,\eta}  \left( \rho^t - \rho , \eta^t - \eta \right)
\right| \leq   \int_0^1   \left|  \partial_{\eta}
\partial_{\rho,\eta} (\rho^t,\eta^t)(r,\theta,\rho, s
  \eta)   \right|  d s |\eta| \lesssim | e^{-r} \eta | \lesssim \epsilon  ,
$$
since, by Proposition \ref{libre}, $ (r,\theta,\rho, s
  \eta) \in \Gamma^{\pm}_{\rm
s}(\epsilon_0) $ if $ (r,\theta,\rho,\eta) \in \Gamma^{\pm}_{\rm
s}(\epsilon) $ and $ \epsilon $ is small enough. Therefore, if $
\epsilon  $ is small enough,
\begin{eqnarray}
 |
\partial_{\rho,\eta}(\rho^t,\eta^t) - \mbox{Id}_n | \leq 1/2  , \qquad \mbox{
  on } \ \
 \Gamma^{\pm}_{\rm s} (\epsilon)  , \label{conditionShwartz}
\end{eqnarray}
 for all $ \pm t \geq 0  $.
Here $ |.|  $ is a matrix norm. We can  now  prove
(\ref{estimeinversepsit})  when $ j+ |\alpha|+ k+  |\beta| \geq 1
$. Assume first that $ D_{\rm hyp}^{j\alpha k \beta} = e^{r}
\partial^{\beta}_{\eta}  $, with $ |\beta| = 1 $, and denote for
simplicity
$$ \Xi_t (r,\theta,\rho,\eta) = (\rho_t,\eta_t)(r,\theta,\rho,\eta)  , \qquad
\Xi^t (r,\theta,\tilde{\rho},\tilde{\eta}) =
(\rho^t,\eta^t)(r,\theta,\tilde{\rho},\tilde{\eta}), \qquad \Xi  =
(\rho,\eta) , $$ when $ (r,\theta,\rho,\eta) \in \Gamma^{\pm}_{\rm
s}(\epsilon^3)  $,  $ (r,\theta,\tilde{\rho},\tilde{\eta}) \in
 \Gamma^{\pm}_{\rm s} (\epsilon)  $ and $ \pm t \geq 0  $. Applying $ e^{r} \partial^{\beta}_{\eta} $ to  (
 \ref{ffinverse}),  we get
$$ ( \partial_{\tilde{\rho},\tilde{\eta}} \Xi^t ) (r,\theta,\Xi_t)   e^{r}
\partial^{\beta}_{\eta}  \Xi_t = (0,e^{r}\partial^{\beta}_{\eta} \eta ) = e^{r} \partial_{\eta}^{\beta} \Xi ,   $$
and using that $ ( \partial_{\tilde{\rho},\tilde{\eta}} \Xi^t )
\partial^{\beta}_{\eta}   \Xi  =
\partial_{\tilde{\eta}}^{\beta} \Xi^t  $, we obtain
\begin{eqnarray*}
  ( \partial_{\tilde{\rho},\tilde{\eta}} \Xi^t ) (r,\theta,\Xi_t) e^{r}\partial^{\beta}_{\eta}  ( \Xi_t -  \Xi )
 =   e^{r}  \left( \partial_{\tilde{\eta}} (\Xi - \Xi^t) \right)_{|(r,\theta,\Xi_t)}
  ,
\end{eqnarray*}
where  the right hand side is bounded, by (\ref{flotrrhoeta}).
Using (\ref{conditionShwartz}), we see that $
e^{r}\partial_{\eta}^{\alpha} (\Xi_t - \Xi ) $ is bounded on $
\Gamma^{\pm}_{\rm s}(\epsilon_1)  $ for $ \pm t \geq 0  $, by
choosing  $ \epsilon_1 \leq \epsilon_0^{3} $ and such that
(\ref{conditionShwartz})  holds. The other first order derivatives
 are treated similarly and are simpler to handle since there is
no $ e^r $. When $ j + |\alpha| + k + | \beta | \geq 2 $, we
iterate this process using Lemma \ref{Faadibruno}. To complete
the proof of the proposition, we finally note  that
(\ref{implicite1})  and (\ref{implicite2}) are easy  consequences of
(\ref{ffinverse}) and  Proposition \ref{diffeoexact}. \finpreuve

\bigskip

\bigskip

By  Propositions \ref{preparationphase} and
\ref{estimepreparephase}, we can define
 $ r_t^s = r_t^s (r,\theta,\rho,\eta) $ and $ \theta_t^s =
\theta_t^s (r,\theta,\rho,\eta) $ on $ \Gamma^{\pm}_{\rm
s}(\epsilon_1) $ by
$$ r_t^s = r^s (r,\theta,\rho_t,\eta_t) , \qquad \theta_t^s = \theta^s
(r,\theta,\rho_t,\eta_t) ,  \qquad  \pm t \geq \pm s \geq  0  , $$
where $  \pm t \geq \pm s \geq  0 $ means more precisely that $  t
\geq s \geq 0  $ if $  (r,\theta,\rho,\eta)  \in \Gamma^+_{\rm
s}(\epsilon_1) $ and that $ t \leq s \leq 0 $ if $
(r,\theta,\rho,\eta) \in \Gamma^-_{\rm s}(\epsilon_1)  $. Here we
assume
 that $ \epsilon_1 $ is small enough so
that Proposition \ref{estimeesflotprecises} hold for $ r > R
(\epsilon_1)$ and $ \sigma = 1/2 $ (for instance), which justifies
that  $ r_t^s $ and $ \theta_t^s $ are well defined (and that
their derivatives can be estimated  using Proposition
\ref{estimeesflotprecises}).

 By  the classical  Hamilton-Jacobi
theory, the function $  \Sigma_{\pm} $ defined by
\begin{eqnarray}
 \Sigma_{\pm} (t,r,\theta,\rho,\eta)  = r_t^t \rho + \theta_t^t \cdot \eta - t \rho^2 - t
 e^{-2 r_t^t} q (r_t^t,\theta_t^t,\eta)  \label{timeexplicit}
\end{eqnarray}
solve the following time dependent  eikonal equation, for  $
(r,\theta,\rho,\eta) \in \Gamma^{\pm}_{\rm s} (\epsilon_1) $ and $
\pm t  \geq 0 $,
\begin{eqnarray}
 \partial_t \Sigma_{\pm} = p (r,\theta,\partial_r \Sigma_{\pm} , \partial_{\theta} \Sigma_{\pm}) ,
\qquad \Sigma_{\pm} |_{ t = 0} = r \rho + \theta \cdot \eta .
\label{HJtimedependent}
\end{eqnarray}
 To put it in a more standard way, note
that  (\ref{timeexplicit}) is obtained by defining $ \Sigma_{\pm}
$ via $ \Sigma_{\pm} (t,r,\theta,\rho^t,\eta^t) = r^t \rho^t +
\theta^t \cdot \eta^t - t p (r^t,\theta^t,\rho^t,\eta^t) $. Note
also that this simple expression uses the fact that $ p $ is
homogeneous of degree $ 2 $ in $ (\rho,\eta) $. Now assume for a
while that
\begin{eqnarray}
S_{\pm} (r,\theta,\rho,\eta) := r \rho + \theta \cdot \eta +
\int_0^{\pm
  \infty} \partial_t \left( \Sigma_{\pm} (t,r,\theta,\rho,\eta) - t \rho^2
\right) dt \label{defphaseind}
\end{eqnarray}
is well defined on $ \Gamma^{\pm}_{\rm s} (\epsilon_1)  $. Then,
at least formally,
\begin{eqnarray}
 \partial_{r,\theta} S_{\pm} (r,\theta,\rho,\eta) = \lim_{t \rightarrow \pm
  \infty}
\partial_{r,\theta} \Sigma_{\pm} (t,r,\theta,\rho,\eta) .  \label{limitederiveestemps}
\end{eqnarray}
The latter only uses the fact that the term $ t \rho^2 $ inside
the integral is independent of $ r,\theta $. If we know in
addition that
\begin{eqnarray}
 \lim_{t \rightarrow \pm
  \infty}
\partial_{\rho}  \Sigma_{\pm} (t,r,\theta,\rho,\eta) = + \infty   , \label{limitephasetemps}
\end{eqnarray}
 then, using that $ \Sigma^{\pm}  $ are  generating functions of
$ \Phi^t $, namely
\begin{eqnarray}
 \Phi^t (r,\theta,\partial_r \Sigma_{\pm},\partial_{\theta} \Sigma_{\pm} ) =
(\partial_{\rho} \Sigma_{\pm} , \partial_{\eta}
\Sigma_{\pm},\rho,\eta ) , \qquad \pm  t \geq 0 ,
\label{generatingfunction}
\end{eqnarray}
we obtain, on $ \Gamma^{\pm}_{\rm s} (\epsilon_1)  $,
$$  p (r,\theta,\partial_r S_{\pm},\partial_{\theta} S_{\pm}) = \lim_{t
  \rightarrow \pm \infty} p (\partial_{\rho} \Sigma_{\pm} , \partial_{\eta}
\Sigma_{\pm},\rho,\eta ) = \rho^2 .   $$ Let us state the
following proposition.
\begin{prop} \label{Phaseplusmoins} There exists $ 0 < \epsilon_2 \leq \epsilon_1 $  such that
we can find   $ S_{\pm} = S_{\pm}(r,\theta,\rho,\eta) $, defined
on $ \Gamma^{\pm}_{\rm s} (\epsilon_2) $, real valued, satisfying
\begin{eqnarray}
 p (r,\theta,\partial_r S_{\pm},\partial_{\theta} S_{\pm}) = \rho^2 , \qquad
\mbox{on} \ \ \Gamma^{\pm}_{\rm s} (\epsilon_2) ,
\label{HJtimeindependent}
\end{eqnarray}
 and such that
\begin{eqnarray}
  S_{\pm} (r,\theta,\rho,\eta) =  r \rho + \theta \cdot \eta +
\varphi_{\pm}(r,\theta,\rho,\eta), \label{formepratique}
\end{eqnarray}
 for some $ \varphi_{\pm} \in
{\mathcal B}_{\rm hyp} (\Gamma^{\pm}_{\rm s}(\epsilon_2)) $ satisfying,
when $ (r,\theta,\rho,0) \in \Gamma^{\pm}_{\rm s} (\epsilon_2) , $
\begin{eqnarray}
  \varphi_{\pm | \eta = 0} & = & 0 , \label{ordre0}  \\  e^{r} \partial_{\eta}
  \varphi_{\pm | \eta = 0} & = &  0 , \label{ordre1}  \\
 e^{2r} \emph{hess}_{\eta} [\varphi_{\pm } ]_{ | \eta = 0} & = &
\int_{0 }^{\pm \infty} e^{ - 4 t \rho} \emph{hess}_{\eta} [q] (r+2t \rho,\theta) dt . \label{ordre2}
\end{eqnarray}
\end{prop} 
It is convenient to note that, by possibly decreasing $ \epsilon_2
$ and by using Lemma \ref{lemmesymbole}, (\ref{Taylorutile}),
(\ref{ordre0}) and (\ref{ordre1}), we can write
\begin{eqnarray}
\varphi_{\pm} (r,\theta,\rho,\eta) = \sum_{|\beta|=2}
a^{\pm}_{\beta} (r,\theta,\rho,e^{-r}\eta) e^{-2 r}\eta^{\beta}
\label{crucialphasestationnaire},
\end{eqnarray}
with $ a^{\pm}_{\beta} \in C^{\infty}_b (
F_{\rm hyp} (\Gamma^{\pm}(\epsilon_2))) $.


\bigskip

\noindent {\it Proof.} We consider only the outgoing case. To
complete the proof of  (\ref{HJtimeindependent}), we have to prove
the missing details, namely the convergence of the integral in
(\ref{defphaseind}) (plus its derivability) and the limits
(\ref{limitederiveestemps}) and (\ref{limitephasetemps}).
  Defining $ (\rho_t^s,\eta_t^s) := (\rho^s , \eta^s )
(r,\theta,\rho_t,\eta_t)  $, the motion equations yield
\begin{eqnarray}
 r_t^t & = &  r + 2 \int_0^t \rho^s_t ds  , \nonumber \\
& = &  r + 2 t \rho - 2 \int_0^t  \int_s^t e^{-2 r^u_t} \left( 2 q
                          (r^u_t,\theta^u_t,\eta^u_t) - (\partial_r
                          q)(r^u_t,\theta^u_t,\eta^u_t)  \right) du ds
                        . \label{rtasymptotique}
\end{eqnarray}
By  Propositions \ref{estimeesflotprecises} and
\ref{estimepreparephase}, we have the following bounds on  $
\Gamma^+_{\rm s} (\epsilon_1) $, for $ s \geq 0 $ and $ t \geq 0
$,
\begin{eqnarray}
\left| D^{\rm hyp}_{j \alpha k \beta} \left( r^s_t - r  \right)
\right| \lesssim \scal{ s }, \qquad \left| D^{\rm hyp}_{j \alpha k
\beta} \left( \theta^s_t - \theta   \right) \right| \lesssim
e^{-r} , \qquad \left| D^{\rm hyp}_{j \alpha k \beta} \left(
\eta^s_t - \eta \right) \right| \lesssim   1 .
\label{borneutileredaction}
\end{eqnarray}
In addition, using  Proposition \ref{localisationr} and
(\ref{inclusionmesuree}), we have, for $ s \geq 0 $ and $ t \geq 0
$,
\begin{eqnarray}
 r^s_t \geq r + 2 (1-\epsilon^6) s p^{1/2}
 (r,\theta,\rho_t,\eta_t)- C \geq r +  s/4 -  C , \qquad \mbox{on} \ \Gamma^+_{\rm s} (\epsilon^3) ,
\label{borneutileredaction2}
\end{eqnarray}
with $ \epsilon  $ small enough such that, $ p^{1/2}
(r,\theta,\rho_t,\eta_t) \geq 1/4  $. Using
(\ref{rtasymptotique}), (\ref{borneutileredaction}),
(\ref{borneutileredaction2}),  with $ \epsilon_2 := \epsilon^3
\leq \epsilon_1 $ small enough, and Lemma \ref{Faadibruno}, we
obtain the existence of a bounded family $ (a_t)_{t
  \geq 0}  $ in $ {\mathcal B}_{\rm hyp} (\Gamma^+_{\rm s}(\epsilon_2)) $
such that
\begin{eqnarray}
 r_t^t = r + 2 t \rho + a_t (r,\theta,\rho,\eta) , \qquad t \geq 0 .  \label{famillebornee}
\end{eqnarray}
One shows similarly that $ (\theta_t^t - \theta ) \cdot \eta = e^r
(\theta_t^t - \theta) \cdot e^{-r} \eta  $ is bounded in  $
{\mathcal
  B}_{\rm hyp} (\Gamma^+_{\rm s}(\epsilon_2)) $ for $ t \geq 0 $ and hence
that
\begin{eqnarray}
  \Sigma_+ - ( r \rho + \theta \cdot \eta +  t \rho^2  )  \ \ \mbox{ is
    bounded in }  \ \    {\mathcal
  B}_{\rm hyp} (\Gamma^+_{\rm s}(\epsilon_2)) \ \ \mbox{for} \ t \geq 0, \label{phaseexplicitedemonstration}
\end{eqnarray}
 which  proves
(\ref{limitephasetemps}). Then, using (\ref{HJtimedependent}) and
(\ref{generatingfunction}), we note that
\begin{eqnarray}
 \partial_t \Sigma_+ - \rho^2 =
e^{-2 \partial_{\rho} \Sigma_+} q
(\partial_{\rho}\Sigma_+,\partial_{\eta}\Sigma_+,\eta) .
\label{clefintegrale}
\end{eqnarray}
Therefore, using (\ref{borneutileredaction2}),
(\ref{famillebornee}), (\ref{phaseexplicitedemonstration}) and
(\ref{clefintegrale}), we obtain the convergence of the integral
in (\ref{defphaseind}) and the limit (\ref{limitederiveestemps})
as well as the fact that $ S_{+} (r,\theta,\rho,\eta) - r \rho -
\theta \cdot \eta  $ belongs to $ {\mathcal B}_{\rm hyp}
(\Gamma^+_{\rm s} (\epsilon_2) )  $. Finally, the formulas
(\ref{ordre0}), (\ref{ordre1})  and (\ref{ordre2})  follow directly
from (\ref{clefintegrale}) combined with (\ref{implicite1}) and
(\ref{sansderiver}). \finpreuve

\bigskip

\noindent {\bf Remark 1.} We  point out that, by applying $
\partial_{\eta}  $ to (\ref{phaseexplicitedemonstration}), there
exists $ C $ such that, for all  $ (r,\theta,\rho,\eta) \in
\Gamma^{+}_{\rm s} (\epsilon_2) $ and $ t \geq 0 $,
\begin{eqnarray}
| \partial_{\eta} \Sigma_+ (t,r,\theta,\rho,\eta) - \theta  | \leq
C e^{-r}  \lesssim   e^{-R(\epsilon_2)} . \nonumber
\end{eqnarray}
This shows, in the spirit of Corollary \ref{dependencedudomaine},
that the proof above depends only on the definition of
$q(r,\theta,\eta)$ for $ \theta $ in an arbitrarily small
neighborhood of $ \overline{ V_0 } $, provided $ \epsilon_2 $ is
small enough.

\medskip

\noindent {\bf Remark 2.} Using (\ref{scalingrtheta}),
(\ref{scalingrhoeta}) and (\ref{scalingreciproque}), one sees that  $ S_{\pm}
$ is actually well defined  on the
conical area given by
$$ r > R (\epsilon_2), \qquad \theta \in V_{\epsilon_2}, \qquad \pm \rho > (1
- \epsilon^2_2) p^{1/2} , $$
and that
$$ \Sigma_{\pm} (t,r,\theta, \lambda \rho , \lambda \eta) = \lambda \Sigma_{\pm} (\lambda t , r , \theta , \rho , \eta),
 \qquad \lambda > 0 . $$
This implies that  $ S_{\pm}
$ is  the restriction to $
\Gamma^{\pm}_{\rm s}(\epsilon_2) $ of an homogeneous function of
degree $1$, with respect to $ (\rho,\eta) $.

\bigskip

We conclude this part with a useful result  to consider phases
globally defined on $ \Ra^{2n} $ when we shall construct Fourier
integral operators.

\begin{prop} \label{Phaseplusmoinsglobale} For some small enough $ \epsilon_3  > 0 $, there exists
 a family of functions $(S_{\pm,\epsilon})_{0 < \epsilon \leq \epsilon_3} $, globally
defined on $ \Ra^{2n} $,  such that
$$ {\varphi}_{\pm, \epsilon} (r, \theta,\rho,\eta) :=  S_{\pm,\epsilon}
(r,\theta,\rho,\eta) - r \rho - \theta \cdot \eta  $$ coincides
with $ \varphi_{\pm}  $ on  $ \Gamma^{\pm}_{\rm s} (\epsilon) $
and  satisfies
\begin{eqnarray}
\emph{supp} ( {\varphi}_{\pm,\epsilon} ) & \subset &  \Gamma^{\pm}_{\rm s} (\epsilon^{1/2}) , \label{supportstilde} \\
{\varphi}_{\pm,\epsilon}
 & \in & {\mathcal B}_{\rm hyp} (\Gamma^+_{\rm s} (\epsilon^{1/2}) ) \label{deriveestilde}
,
\\
 \left| \partial_{\rho,\eta} \otimes  \partial_{r,\theta}  {\varphi}_{\pm,\epsilon}
(r,\theta,\rho,\eta)
 \right| & \leq & 1/2 , \qquad  ( r,\theta,\rho,\eta ) \in \Ra^{2n}
 , \ \ 0 < \epsilon \leq \epsilon_3,
 \label{restepetitapplication}
\end{eqnarray}
with $ | \cdot | $ a matrix norm.
\end{prop}
In further applications, (\ref{restepetitapplication}) will also
be used under the equivalent form
\begin{eqnarray}
 \left|  \partial_{\rho,\eta} \otimes \partial_{r,\theta}
S_{\pm, \epsilon} (r,\theta,\rho,\eta) - \mbox{Id}_n
 \right| \leq 1/2 , \qquad  ( r,\theta,\rho,\eta ) \in \Ra^{2n} ,
 \ \ 0 < \epsilon \leq \epsilon_3 .  \label{diffeoKuranishi}
\end{eqnarray}

\bigskip \noindent {\bf Remark.} Although this proposition allows to assume that  they are
globally defined, the phases  $ S_{\pm} $  solve the
Hamilton-Jacobi equations on $ \Gamma^{\pm}_{\rm s} (\epsilon_2) $
only.

\medskip

\noindent {\it Proof.} We use Lemma \ref{construitcutoff} and
consider
\begin{eqnarray}
 S_{\pm,\epsilon} (r,\theta,\rho,\eta) : = r \rho + \theta \cdot \eta +
\chi_{\epsilon^{1/2} \rightarrow \epsilon}
(r,\theta,\rho,\eta)\varphi_{\pm}(r,\theta,\rho,\eta) ,
\label{vraiephasepourlescalculs}
\end{eqnarray}
 with $
\varphi_{\pm} $ defined in Proposition \ref{Phaseplusmoins}. We
have $ S_{\pm,\epsilon} = S_{\pm} $ on $ \Gamma^{\pm}_{\rm
s}(\epsilon) $ and, using(\ref{pourKuranishi})  and
(\ref{crucialphasestationnaire}),
$$ \left| \partial_{\rho,\eta} \otimes  \partial_{r,\theta}  S_{\pm,\epsilon}
(r,\theta,\rho,\eta) - \mbox{Id}_n
 \right| \lesssim \epsilon^{1/4}  , \qquad \mbox{on} \ \ \Ra^{2n} , $$
 since $ e^{-r}|\eta| \lesssim \epsilon^{1/2} $ on $ \Gamma^+_{\rm s} (\epsilon^{1/2})
 $. This yields the result if $ \epsilon $ is small. \finpreuve

\subsection{Fourier integral operators on $ \Ra^n $} \label{operateurintegraldeFourier}

In this subsection, we derive some basic properties of Fourier integral operators associated to the phases $ S_{\pm} $
obtained in Proposition \ref{Phaseplusmoins}.

For simplicity, we introduce the shorter notation
\begin{eqnarray}
  {\mathcal B}_{\rm s}^{\pm} (\epsilon) :=
 {\mathcal B}_{\rm hyp}(\Gamma^{\pm}_{\rm s}(\epsilon)) , \qquad {\mathcal S}_{\rm s}^{\pm} (\epsilon) := {\mathcal S}_{\rm hyp}^{\pm} (\Gamma^{\pm}_{\rm s}(\epsilon)), \label{simplificationclasses}
\end{eqnarray}
where the classes $ {\mathcal B}_{\rm hyp} $ and $ {\mathcal S}_{\rm hyp} $ were defined in Definition \ref{definhyp}.

\medskip

By Propositions  \ref{Phaseplusmoins} and
\ref{Phaseplusmoinsglobale},  for all $ h \in (0,1] $, all $
\epsilon $ small enough and all  $ a^{\pm} \in {\mathcal
S}^{\pm}_{\rm s} (\epsilon) $, we can define the operator
\begin{eqnarray}
 J^{\pm}_h(a^{\pm}) : {\mathcal S} (\Ra^n) \rightarrow {\mathcal
S} (\Ra^n)  , \label{OIFRn}
\end{eqnarray}
 as the operator with Schwartz kernel
$$ K_h^{\pm} (r,\theta,r^{\prime},\theta^{\prime}) =
 (2 \pi h)^{-n} \int e^{ \frac{i}{h} \left( S_{\pm}(r,\theta,\rho,\eta) - r^{\prime} \rho -
\theta^{\prime} \cdot \eta \right)  } a^{\pm} (r,\theta,\rho,\eta)
\ d\rho d \eta . $$
Since the symbol $ a^{\pm} $ is supported in $ \Gamma^{\pm}_{\rm s} (\epsilon)
 $, the phase $ S_{\pm} $ can be replaced by $ S_{\pm,\epsilon} $ which is globally defined (see Proposition \ref{Phaseplusmoinsglobale}).
Note also that
 $ J^{\pm}_h(a^{\pm}) $ maps clearly the
Schwartz space into itself since, for fixed $h$ say $h=1$, it can
be considered as the pseudo-differential operator with symbol $
e^{i \varphi_{\pm} } a^{\pm} = e^{i \varphi_{\pm,\epsilon} }
a^{\pm} $ which belongs to $ C^{\infty}_b (\Ra^{2n}) $.

To obtain the $ L^2 $ boundedness of such operators uniformly in $
h \in (0,1] $ as well as the factorization Lemma
\ref{factorisation} below, which are both consequences of the
usual Kuranishi trick,  we need a preliminary result.

Consider the maps   $
 (\underline{\rho}_{\pm,\epsilon} , \underline{\eta}_{\pm,\epsilon} ) :  \Ra^{3n}
 \rightarrow \Ra^n
$ defined by
\begin{eqnarray}
  (\underline{\rho}_{\pm,\epsilon},\underline{\eta}_{\pm,\epsilon})(r,\theta,r^{\prime},\theta^{\prime},\rho,\eta) : = \int_0^{1}
  \partial_{r,\theta} S_{\pm,\epsilon} (r^{\prime} + s
(r-r^{\prime}),\theta^{\prime} + s
(\theta-\theta^{\prime}),\rho,\eta) d s
\label{aevaluersurlinverse}
\end{eqnarray}
so that
\begin{eqnarray}
 (r-r^{\prime}) \underline{\rho}_{\pm,\epsilon} +
(\theta-\theta^{\prime})\cdot \underline{\eta}_{\pm,\epsilon}
 = S_{\pm,\epsilon}(r,\theta,\rho,\eta) - S_{\pm,\epsilon}(r^{\prime},\theta^{\prime},\rho,\eta)
 . \label{pourdiagonaliser}
\end{eqnarray}
\begin{lemm} \label{Kuranishi} For all $ (r,\theta,r^{\prime},\theta^{\prime}) \in \Ra^{2n}
$ and all $ 0 <  \epsilon  \leq \epsilon_3 $, the map $
(\rho,\eta) \mapsto
(\underline{\rho}_{\pm,\epsilon},\underline{\eta}_{\pm,\epsilon})
$ is a diffeomorphism  from $ \Ra^n $ onto itself. Denoting by $ (
\overline{\rho}_{\pm,\epsilon} , \overline{\eta}_{\pm,\epsilon} )
$ the corresponding inverse, we have, for all  $ 0 <  \epsilon
\leq \epsilon_3 $,
\begin{eqnarray}
 \left| \partial_{\eta}^{\beta}
\partial_{r}^j \partial_{r^{\prime}}^{j^{\prime}}
 \partial_{\theta}^{\alpha} \partial_{\theta^{\prime}}^{\alpha^{\prime}}
 \partial_{\rho}^k
  \left( (\overline{\rho}_{\pm,\epsilon},\overline{\eta}_{\pm,\epsilon}) - (\rho,\eta)  \right) \right| & \lesssim & 1
  ,\qquad \mbox{on} \   \Ra^{3n} .
  \label{controlesurlesderiveesKuranishi}
 \end{eqnarray}
Furthermore, there exists $ \epsilon_6  > 0 $ such that, for all $
0 < \epsilon  \leq \epsilon_6 $, we have
\begin{eqnarray}
(r,\theta, \rho , \eta ) \in \Gamma^{\pm}_{\rm s} (\epsilon ) &
\Rightarrow &  \left( r , \theta , \underline{\rho}_{\pm
,\epsilon}, \underline{\eta}_{\pm,\epsilon} \right)_{ |
r=r^{\prime},\theta=\theta^{\prime}} \in \Gamma^{\pm}_{\rm
s} (\epsilon^{1/3} ) , \qquad \qquad \label{pourlapropositionsuivante} \\
 \left( r , \theta , \underline{\rho}_{\pm ,\epsilon }, \underline{\eta}_{\pm,\epsilon} \right)_{ |
r=r^{\prime},\theta=\theta^{\prime}}  \in \Gamma^{\pm}_{\rm s}
(\epsilon^3 ) & \Rightarrow &   (r,\theta, \rho , \eta ) \in
\Gamma^{\pm}_{\rm s} (\epsilon ) ,
\label{pourlapropositionsuivante2}
\end{eqnarray}
and
\begin{eqnarray}
 \left| \partial_{\eta}^{\beta}
\partial_{r}^j \partial_{r^{\prime}}^{j^{\prime}}
 \partial_{\theta}^{\alpha} \partial_{\theta^{\prime}}^{\alpha^{\prime}}
 \partial_{\rho}^k
  \left( (\overline{\rho}_{\pm,\epsilon},\overline{\eta}_{\pm,\epsilon}) - (\rho,\eta)  \right)_{| r=r^{\prime},
  \theta = \theta^{\prime}} \right| & \lesssim & e^{-|\beta|  r },
  \qquad \mbox{on} \ \Gamma^+_{\rm s} (\epsilon^3 ) .
  \label{diagoanleKuranishi}
\end{eqnarray}
\end{lemm}

\noindent {\it Proof.}  The estimate (\ref{diffeoKuranishi})
implies directly that $ (\rho,\eta) \mapsto
(\underline{\rho}_{\pm,\epsilon},\underline{\eta}_{\pm,\epsilon})
$ is a diffeomorphism for all $
(r,\theta,r^{\prime},\theta^{\prime}) \in \Ra^{2n} $ and $ 0 <
\epsilon  \leq \epsilon_3  $. Evaluating
(\ref{aevaluersurlinverse}) at $
(r,\theta,r^{\prime},\theta^{\prime},
\overline{\rho}_{\pm,\epsilon},\overline{\eta}_{\pm,\epsilon} ) $,
namely
\begin{eqnarray}
 (\rho,\eta) = (
 \underline{\rho}_{\pm,\epsilon},\underline{\eta}_{\pm,\epsilon})(r,\theta,r^{\prime},
\theta^{\prime},\overline{\rho}_{\pm,\epsilon},\overline{\rho}_{\pm,\epsilon})
,  \label{pourlesderivees}
\end{eqnarray}
 yields
%
\begin{eqnarray}
(\rho,\eta) - (
\overline{\rho}_{\pm,\epsilon},\overline{\eta}_{\pm,\epsilon} ) =
\int_0^{1}   \partial_{r,\theta} \varphi_{\pm,\epsilon}
(r^{\prime} + s (r-r^{\prime}),\theta^{\prime} + s
(\theta-\theta^{\prime}),
\overline{\rho}_{\pm,\epsilon},\overline{\eta}_{\pm,\epsilon} ) d
s . \label{pourlafonction}
\end{eqnarray}
 By
(\ref{supportstilde}) and (\ref{deriveestilde}),  $
\varphi_{\pm,\epsilon} \in C^{\infty}_b (\Ra^{2n}) $, so
  $ (\overline{\rho}_{\pm,\epsilon},\overline{\eta}_{\pm,\epsilon}) -
(\rho,\eta) $ is bounded, for fixed $ \epsilon  $. For the
derivatives, we apply  $
\partial_{\eta}^{\beta}
\partial_{r}^j \partial_{r^{\prime}}^{j^{\prime}}
 \partial_{\theta}^{\alpha} \partial_{\theta^{\prime}}^{\alpha^{\prime}}
 \partial_{\rho}^k $ to the right hand
side of (\ref{pourlafonction}) and obtain
(\ref{controlesurlesderiveesKuranishi}) by induction, using Lemma
\ref{Faadibruno}.

To prove (\ref{pourlapropositionsuivante}), we simply notice that
$ \varphi_{\pm,\epsilon}  $ coincides with $ \varphi_{\pm}  $ on
$\Gamma^{\pm}_{\rm s} (\epsilon^3 ) $ so that
$$ \big| (\rho,\eta) - (\underline{\rho}_{\pm ,\epsilon},
\underline{\eta}_{\pm,\epsilon}  )_{ |
r=r^{\prime},\theta=\theta^{\prime}}  \big| = |\partial_{r,\theta}
\varphi_{\pm}(r,\theta,\rho,\eta) | \lesssim \epsilon^2 ,  $$
using (\ref{crucialdiffeo}) and (\ref{crucialphasestationnaire}).
The result follows from Proposition \ref{libre} and the fact that
$  \Gamma^{\pm}_{\rm s} ( C \epsilon ) \subset  \Gamma^{\pm}_{\rm
s} (\epsilon^{1/3} )  $ for $ \epsilon $ small enough. To get
(\ref{pourlapropositionsuivante2}), we use directly Proposition
\ref{propositiondouble}  proving that $ \Gamma^{\pm}_{\rm s}
(\epsilon^3) \subset \Psi^t ( \Gamma^{\pm}_{\rm
  s}(\epsilon) )  $ with
$$ \Psi^t (r,\theta,\rho,\eta) : = (r,\theta,\underline{\rho}_{\pm,\epsilon} ,
\underline{\eta}_{\pm,\epsilon}
)_{|r=r^{\prime},\theta=\theta^{\prime}} = (r,\theta,\partial_r
S_{\pm}(r,\theta,\rho,\eta),\partial_{\theta}
S_{\pm}(r,\theta,\rho,\eta)) ,  $$ which is actually independent
of $t$ and $ \epsilon  $.

By (\ref{controlesurlesderiveesKuranishi}),
(\ref{diagoanleKuranishi}) holds when $\beta =0$. Consider next
the first order derivatives when $ |\beta| = 1  $ and the other
multi-indices are $ 0 $. Applying $ \partial_{\eta}^{\beta}  $ to
(\ref{pourlesderivees}) and  evaluating at $ r= r^{\prime}  $,  $
\theta = \theta^{\prime}  $, we get
\begin{eqnarray*}
  \left( \partial_{\rho,\eta}
  (\underline{\rho}_{\pm,\epsilon},\underline{\eta}_{\pm,\epsilon}) \right)
  \partial_{\eta}^{\beta}
  \left(
    (\overline{\rho}_{\pm,\epsilon},\overline{\eta}_{\pm,\epsilon}) -  (\rho,\eta)   \right) = \partial_{\eta}^{\beta} \partial_{r,\theta}
\varphi_{\pm}
(r,\theta,\overline{\rho}_{\pm,\epsilon},\overline{\eta}_{\pm,\epsilon})
\end{eqnarray*}
where we have replaced $ \varphi_{\pm,\epsilon}  $ by $
\varphi_{\pm}  $  using (\ref{pourlapropositionsuivante2}). Since
$ \left( \partial_{\rho,\eta}
  (\underline{\rho}_{\pm,\epsilon},\underline{\eta}_{\pm,\epsilon}) \right)^{-1}  $
is uniformly bounded and $  e^{r \beta }\partial_{\eta}^{\beta}
\partial_{r,\theta} \varphi_{\pm}
(r,\theta,\overline{\rho}_{\pm,\epsilon},\overline{\eta}_{\pm,\epsilon})
$ is bounded,  using (\ref{pourlapropositionsuivante2}) again, we
get the result in this case. Higher order derivatives are obtained
similarly by induction, using Lemma \ref{Faadibruno}. \finpreuve

\bigskip

\begin{prop} \label{estimationsKuranishi} For all $ 0 < \epsilon \leq \epsilon_6 $ and all $ a^{\pm} ,
b^{\pm} \in {\mathcal S}^{\pm}_{\rm s} (\epsilon) $, we have
\begin{eqnarray}
 \big| \big|
 J^{\pm}_h (a^{\pm}) J^{\pm}_h (b^{\pm})^* - \sum_{ k \leq N} h^k c_k^{\pm} (r,\theta,hD_r,h D_{\theta})
  \big| \big|_{L^2 (\Ra^n) \rightarrow L^2(\Ra^n)} \leq C h^{N+1},
  \qquad h \in (0,1], \label{borneL2Kuranishi}
\end{eqnarray}
   where the constant $C$ can be chosen uniformly with respect to $a^{\pm}$
and $b^{\pm}$ when they vary in bounded subsets of $ {\mathcal
S}^{\pm}_{\rm s} (\epsilon) $ and where the symbols $ c_k^{\pm}
$ are given by
\begin{eqnarray}
 c^{\pm}_k = \sum_{j+|\alpha| = k} \frac{1}{j!\alpha!}\partial_{r^{\prime}}^j
 \partial_{\theta^{\prime}}^{\alpha} D_{\rho}^j D_{\eta}^{\alpha}
 \left(
 a (r,\theta,\overline{\rho}_{\pm,\epsilon},\overline{\eta}_{\pm,\epsilon})
 \overline{  b
 (r^{\prime},\theta^{\prime},\overline{\rho}_{\pm,\epsilon},\overline{\eta}_{\pm,\epsilon})}
  \emph{Jac}(\overline{\rho}_{\pm,\epsilon},\overline{\eta}_{\pm,\epsilon})
  \right)_{| r=r^{\prime}, \ \theta = \theta^{\prime}} , \label{partieprincipale}
\end{eqnarray}
with $
\emph{Jac}(\overline{\rho}_{\pm,\epsilon},\overline{\eta}_{\pm,\epsilon})
= |\emph{det} (\partial_{\rho,\eta}
(\overline{\rho}_{\pm,\epsilon},\overline{\eta}_{\pm,\epsilon})) |
$. In particular,
\begin{eqnarray}
  c_k^{\pm}  \in {\mathcal S}^{\pm}_{\rm s} ( \epsilon^{1/3}) .
  \label{localisationdusupport}
\end{eqnarray}
\end{prop}

\noindent {\it Proof.} The Schwartz kernel of $ J^{\pm}_h
(a^{\pm}) J^{\pm}_h (b^{\pm})^* $ takes the form
$$ (2 \pi h)^{-n} \int e^{\frac{i}{h} (S_{\pm,\epsilon}(r,\theta,\rho,\eta) - S_{\pm,\epsilon}(r^{\prime},\theta^{\prime},\rho,\eta) ) }
a(r,\theta,\rho,\eta) \overline{b
(r^{\prime},\theta^{\prime},\rho,\eta)} d\rho d\eta$$ and this can
be rewritten using the Kuranishi trick, ie
(\ref{pourdiagonaliser}) and Lemma \ref{Kuranishi}, as
\begin{eqnarray}
 (2 \pi h)^{-n} \int e^{\frac{i}{h} ((r - r^{\prime})\rho +
(\theta-\theta^{\prime}) \cdot \eta)  } a
(r,\theta,\overline{\rho}_{\pm,\epsilon},\overline{\eta}_{\pm,\epsilon})
 \overline{  b
 (r^{\prime},\theta^{\prime},\overline{\rho}_{\pm,\epsilon},\overline{\eta}_{\pm,\epsilon})}
  \mbox{Jac}(\overline{\rho}_{\pm,\epsilon},\overline{\eta}_{\pm,\epsilon}) d\rho d\eta
  .\label{pseudodouble}
\end{eqnarray}
By (\ref{controlesurlesderiveesKuranishi}), the symbol in
(\ref{pseudodouble}) belongs to $ C^{\infty}_b (\Ra^{3n}) $.
Therefore, the standard $h$-pseudo-differential calculus implies
that, with $c_k$ defined by (\ref{partieprincipale}), we obtain
the $ L^2 $ bound (\ref{borneL2Kuranishi}) by the
Calder\'on-Vaillancourt Theorem. In addition, by
(\ref{pourlapropositionsuivante}) (applied with $ (\rho,\eta) =
(\overline{\rho}_{\pm,\epsilon},\overline{\eta}_{\pm,\epsilon}
)_{| r=r^{\prime},\theta=\theta^{\prime}} $), we have $
\mbox{supp}(c_k^{\pm}) \subset \Gamma^+_{\rm s} (\epsilon^{1/3} )
$. One then checks that $ c_k^{\pm} \in {\mathcal B}^{\pm}_{\rm s}
( \epsilon^{1/3}) $, using (\ref{diagoanleKuranishi}).
 \finpreuve

\bigskip

We note in passing that this proposition shows that, for all $ 0 <
\epsilon \leq \epsilon_6  $ and all $ a^{\pm} \in {\mathcal
S}^{\pm}_{\rm s}(\epsilon) $,
\begin{eqnarray}
 || J^{\pm}_h (a^{\pm}) ||_{L^2 (\Ra^n) \rightarrow L^2 (\Ra^n)}
 \leq C , \qquad h \in (0,1] . \label{borneL2FIO}
\end{eqnarray}
More precisely, the constant $C$ can be chosen independently of $ a^{\pm} $ if, for $ \epsilon $ fixed, $a^{\pm} $ vary in a bounded
subset of $ {\mathcal
S}^{\pm}_{\rm s}(\epsilon) $.

\bigskip


\begin{prop} \label{factorisation} For all $ 0 < \epsilon \leq \epsilon_6  $, the following hold:
 if we are given
$$ a_0^{\pm}, \ldots , a_N^{\pm} \in
{\mathcal S}^{\pm}_{\rm s}(\epsilon)  , $$ with $ a^{\pm}_0 $
such that
\begin{eqnarray}
 a^{\pm}_0  \gtrsim  1 , \qquad \mbox{ on } \  \Gamma_{\rm s}^{\pm}
 (\epsilon^3) ,
\label{ellipticiteFourier}
\end{eqnarray}
then, for  all $ \chi^{\pm}_{\rm s} \in {\mathcal S}^{\pm}_{\rm
s} (\epsilon^9) $, we can find $ b_0^{\pm} , \ldots , b_N^{\pm}
\in {\mathcal S}^{\pm}_{\rm s} (\epsilon^3) $ such that, if we
set
$$  a^{\pm}(h) = a_0^{\pm} + \cdots + h^N a_N^{\pm}, \qquad
b^{\pm} (h) = b_0^{\pm} + \cdots + h^N b_N^{\pm} ,  $$ we have
$$ \left| \left| J_h^{\pm} (a^{\pm}(h)) J_h^{\pm} (b^{\pm}(h))^*
 - \chi^{\pm}_{\rm s}(r,\theta,hD_r,hD_{\theta}) \right| \right|_{L^2(\Ra^n) \rightarrow L^2 (\Ra^n)} \leq C h^{N+1},
\qquad h \in (0,1 ] . $$
\end{prop}

\noindent {\it Proof.} By Proposition \ref{estimationsKuranishi}
and the notation therein, we only need to find $ b_0^{\pm} ,
\ldots , b_N^{\pm} $ such that
\begin{eqnarray*}
c_0^{\pm} = \chi^{\pm}_{\rm s}, \qquad c_k^{\pm} = 0 , \ \ k = 1 ,
\ldots, N .
\end{eqnarray*}
Using Lemma \ref{Kuranishi} and (\ref{partieprincipale}), the
first equation, ie $ c_0^{\pm} = \chi^{\pm}_{\rm s}  $, is solved
explicitly by
$$  \overline{ b_0^{\pm}(r,\theta,\rho,\eta) } =
\left( \chi^{\pm}_{\rm s}
(r,\theta,\underline{\rho}_{\pm,\epsilon},\underline{\eta}_{\pm,\epsilon})
  \mbox{Jac}(\underline{\rho}_{\pm,\epsilon},\underline{\eta}_{\pm,\epsilon}) \right)_{
   | r^{\prime}=r, \ \theta^{\prime}=\theta}
  \times \frac{1}{
  a_0^{\pm} (r,\theta,\rho,\eta) } , $$
where $ 1/ a_0^{\pm} $ is well defined since $ \chi^{\pm}_{\rm s}
(r,\theta,\underline{\rho}_{\pm,\epsilon},\underline{\eta}_{\pm,\epsilon})_{
| r^{\prime}=r, \ \theta^{\prime}=\theta} $ is supported in $
\Gamma_{\rm s}^{\pm} (\epsilon^3) $ by
(\ref{pourlapropositionsuivante2}). Thus, $ b_0^{\pm} $ is well
defined, supported in $ \Gamma_{\rm s}^{\pm} (\epsilon^3) $ and
belongs to $ {\mathcal B}^{\pm}_{\rm s} (\epsilon^3) $ by
(\ref{aevaluersurlinverse}) and Proposition \ref{Phaseplusmoins}
(since $
(\underline{\rho}_{\pm,\epsilon},\underline{\eta}_{\pm,\epsilon})_{
| r^{\prime}=r, \ \theta^{\prime}=\theta} = \partial_{r,\theta}
S_{\pm} $ in $\Gamma_{\rm s}^{\pm} (\epsilon^3)$). Furthermore, $
b_0^{\pm} (r,\theta,\overline{\rho}_{\pm,\epsilon},
\overline{\eta}_{\pm,\epsilon})_{ | r^{\prime}=r, \
\theta^{\prime}=\theta}  $ is supported in $ \Gamma_{\rm s}^{\pm}
(\epsilon^9) $. We then find the other symbols by induction  for
we have a triangular system of equations. More precisely,  the $ k
$-th equation $ c_k \equiv 0  $ ($k \geq 1$), reads
$$
\left(
\overline{b_k^{\pm}(r,\theta,\overline{\rho}_{\pm,\epsilon},\overline{\eta}_{\pm,\epsilon})
}
a_0^{\pm}(r,\theta,\overline{\rho}_{\pm,\epsilon},\overline{\eta}_{\pm,\epsilon})
\mbox{Jac}(\overline{\rho}_{\pm,\epsilon},\overline{\eta}_{\pm,\epsilon})
\right)_{| r
  =r^{\prime},\theta = \theta^{\prime}} =  d_k^{\pm} (r,\theta,\rho,\eta)
  $$
where $ d_k^{\pm}  $
 is a linear combination of symbols of the form
$$  \overline{ (\partial^{\gamma} b_{k_2}^{\pm})(r,\theta,\overline{\rho}_{\pm,\epsilon},\overline{\eta}_{\pm,\epsilon})
}_{r=r^{\prime},\theta=\theta^{\prime}}
(\partial^{\gamma^{\prime}}a_{k_1}^{\pm})
(r,\theta,\overline{\rho}_{\pm,\epsilon},\overline{\eta}_{\pm,\epsilon})_{r=r^{\prime},\theta=\theta^{\prime}}
\delta_{k_1 k_2 \gamma \gamma^{\prime}}(r,\theta,\rho,\eta)
$$
with $ k_2 < k $ and $ \delta_{k_1 k_2 \gamma \gamma^{\prime}} $ a
product of derivatives of order $ \geq 1 $ of $
(\overline{\rho}_{\pm,\epsilon},\overline{\eta}_{\pm,\epsilon})(r,\theta,r^{\prime},\eta^{\prime},\rho,\eta)
$ evaluated at $ r = r^{\prime},\theta = \theta^{\prime} $. By the
induction assumption $ (\partial^{\gamma}
b_{k_2}^{\pm})(r,\theta,\overline{\rho}_{\pm,\epsilon},\overline{\eta}_{\pm,\epsilon}
)_{r=r^{\prime},\theta=\theta^{\prime}} $ is supported in $
\Gamma_{\rm s}^{\pm} (\epsilon^9)  $, so we have
$$
(r,\theta,\underline{\rho}_{\pm,\epsilon},\underline{\eta}_{\pm,\epsilon})_{r=r^{\prime},\theta=\theta^{\prime}}
\in \Gamma_{\rm s}^{\pm} (\epsilon^3) , $$ using
(\ref{pourlapropositionsuivante}). Therefore, $ \delta_{k_1 k_2
\gamma \gamma^{\prime}}
(r,\theta,\underline{\rho}_{\pm,\epsilon},\underline{\eta}_{\pm,\epsilon})_{r=r^{\prime},\theta=\theta^{\prime}}
$ belongs to  $ {\mathcal B}^{\pm}_{\rm s} (\epsilon^3) $ by
(\ref{diagoanleKuranishi}) and $ b_k^{\pm} $ satisfies the
expected properties. \finpreuve

\subsection{The transport equations} \label{Transport}

In this subsection, we solve the time independent transport
equations related to the phases constructed in Proposition
\ref{Phaseplusmoins}. If we  define $ (v^{\pm},w^{\pm}) =
(v^{\pm},w^{\pm}) (r,\theta,\rho,\eta)  $ by
\begin{eqnarray}
  \begin{pmatrix}
 v^{\pm} \\
 w^{\pm}
\end{pmatrix} :=  \begin{pmatrix}
  (\partial_{\rho} p)(r,\theta,\partial_r S_{\pm} , \partial_{\theta}S_{\pm}) \\
(\partial_{\eta} p)(r,\theta,\partial_r S_{\pm} ,
\partial_{\theta}S_{\pm})
\end{pmatrix} =  \begin{pmatrix}
  2 \partial_r S_{\pm}  \\
e^{-2r}(\partial_{\eta} q)(r,\theta , \partial_{\theta}S_{\pm})
\end{pmatrix}  ,  \label{champsdevecteurs}
\end{eqnarray}
these transport equations take the form
\begin{eqnarray}
 v^{\pm} \partial_r a^{\pm} + w^{\pm} \cdot \partial_{\theta} a^{\pm} +
y^{\pm} a^{\pm} = z^{\pm} ,  \label{formetransport}
\end{eqnarray}
where $ y^{\pm},z^{\pm}  $ are given and $ a^{\pm} $ is  the
unknown function of $ (r,\theta,\rho,\eta )$. Such equations arise
naturally in the construction of the Isozaki-Kitada parametrix
(see Section \ref{IsozakiKitada}). They can be solved standardly by
the method of characteristics and therefore, we start with the
study the integral curves of the vector field $ (v^{\pm},w^{\pm})
$.

Given $ (r,\theta,\rho,\eta) \in \Gamma^{\pm}_{\rm s} (\epsilon^2)
$, with $ \epsilon > 0 $ small enough (to be specified below), we
denote by  $$ r_t^{\pm} = r_t^{\pm} (r,\theta,\rho,\eta)  , \qquad
\theta_t^{\pm} = \theta_t^{\pm} (r,\theta,\rho,\eta) , $$ the
solution to
\begin{eqnarray}
\begin{cases}
 \dot{r}_t^{\pm} \ \  = & v^{\pm}
( r_t^{\pm},\theta_t^{\pm},\rho,\eta), \\
 \dot{\theta}_t^{\pm} \ \  = & w^{\pm}(r^{\pm}_t,\theta^{\pm}_t,\rho,\eta) ,
\end{cases} \label{systemetransport}
\end{eqnarray}
 with initial data
$$ r_0^{\pm} (r,\theta,\rho,\eta) = r \qquad
\theta_0^{\pm} (r,\theta,\rho,\eta) = \theta .  $$ In this
problem, $ \rho $ and $\eta $ are parameters. Equivalently,
$$ \phi_t^{\pm} =  \phi_t^{\pm} (r,\theta,\rho,\eta) := (r^{\pm}_t,
\theta^{\pm}_t , \rho , \eta  ) , $$
 is the flow of the autonomous vector field $ (v^{\pm},w^{\pm},0,0)  $.

\begin{prop} \label{transport} There exist  $ \epsilon_4 > 0 $
such that for all $ (r,\theta,\rho,\eta) \in \Gamma^{\pm}_{\rm s}
(\epsilon_4^2) $, the solution $ (r_t^+,\theta_t^+) $ (resp. $
(r_t^-,\theta_t^-) $) is globally defined on $ [0,+\infty) $
(resp. $(-\infty,0]$). There also exists $ C > 0 $  such that, for
all $ 0 < \epsilon \leq \epsilon_4 $ and all $
(r,\theta,\rho,\eta) \in \Gamma^{\pm}_{\rm s} (\epsilon^2) $, we
have
\begin{eqnarray}
 (r_t^{\pm},\theta^{\pm}_t
,\rho,\eta) \in \Gamma^{\pm}_{\rm s} (\epsilon) , \qquad  \pm t
\geq 0 , \label{localisationphasesepsilon}
\end{eqnarray}
and
\begin{eqnarray}
\big| r_t^{\pm} - r - 2 t \rho \big|  \leq  C \epsilon^2 \min( 1 , |t| ), \qquad \label{rplus}  \\
\big| \theta_t^{\pm} - \theta \big|  \leq  C e^{-r} . \qquad
\qquad \qquad \label{thetaplus}
\end{eqnarray}
Furthermore,
\begin{eqnarray}
 \big| D_{\rm hyp}^{j \alpha k \beta} (r_t^{\pm} - r - 2 t \rho) \big| +  \big| D_{\rm hyp}^{j \alpha k \beta}
(\theta_t^{\pm} - \theta) \big|  \leq  C_{j \alpha k \beta }.
\label{rthetauniftemps}
\end{eqnarray}
for  $ (r,\theta,\rho,\eta) \in \Gamma^{\pm}_{\rm s}(\epsilon^2_4)
$ and $ \pm t \geq 0 $.
\end{prop}


Since $  S_{\pm,\epsilon} =  S_{\pm} $ on $ \Gamma^{\pm}_{\rm
s}(\epsilon) $, the localization property
(\ref{localisationphasesepsilon}) shows that  $ \phi_t^{\pm}  $
still solves (\ref{systemetransport}) on $ \Gamma^{\pm}_{\rm
s}(\epsilon^2) $ if one replaces $ (v^{\pm},w^{\pm})  $ by $
(v^{\pm}_{\epsilon},w^{\pm}_{\epsilon}) $, the latter being
obtained by replacing $ S_{\pm} $  by $ S_{\pm,\epsilon}  $ in
(\ref{champsdevecteurs}).

\bigskip

\noindent {\it Proof.}   Here again we only consider the outgoing
case. By (\ref{crucialphasestationnaire}),  there exists $ C_0
\geq 1 $ such that, for all $ (r,\theta,\rho,\eta) \in
\Gamma^+_{\rm s} (\epsilon_2)  $,
\begin{eqnarray}
  | \partial_{r} S_{+} - \rho | \leq C_0 e^{-r}|\eta|,  \qquad
 |e^{-2r}(\partial_{\eta} q)(r,\theta , \partial_{\theta}S_{+}) | \leq C_0   e^{-2 r} |\eta |
 .  \label{estimeechamps}
\end{eqnarray}
By (\ref{crucialdiffeo}), there exists $ C_1 \geq 1 $ such
that, for all $ \epsilon > 0 $ small enough and all $
(r,\theta,\rho,\eta) \in \Gamma^+_{\rm s} (\epsilon) $,
\begin{eqnarray}
  e^{-r}|\eta| \leq C_1 \epsilon, \qquad e^{-2r}|\eta| \leq C_1 \epsilon^{2} , \label{estimeechamps2}
\end{eqnarray}
the last inequality following from $ e^{-R(\epsilon)} \leq \epsilon $. If
 $ \epsilon $  small enough, we may also assume that, for
all $ (r,\theta,\rho,\eta) \in \Gamma^+_{\rm s} (\epsilon) $,
$$  \rho > 1/8 . $$ Fix now $ M  = 5 C_0 C_1  $, and
 for   $ (r,\theta,\rho,\eta) \in
\Gamma^+_{\rm s } (\epsilon^2) $, consider $ {\mathcal T} := {\mathcal T}(r,\theta,\rho,\eta)  $
defined by
$$  {\mathcal T}  =
\{ t \geq 0 \ | \ (r_s^+ , \theta_s^+) \   \mbox{is defined and}
\ \  r^+_s \geq r + s/8, \ \ |\theta^+_s - \theta|
\leq M \epsilon^2 , \ \ \forall \ s \in [0,t]  \} . $$
The set $  {\mathcal T}  $ is clearly an interval containing $ 0 $ and,   if
$ \epsilon  $ is small  enough,
Proposition \ref{libre} shows that  $ (r_s^+ , \theta_s^+,\rho,\eta) \in \Gamma^+_{\rm
  s}(\epsilon) $ for all $ s \in {\mathcal T} $. Thus, by
(\ref{estimeechamps}) and (\ref{estimeechamps2}), we have
$$ |\dot{r}_s^+ - 2 \rho| \leq 2 C_0 C_1 \epsilon , \qquad |\dot{\theta}^+_s|
\leq C_0 C_1 \epsilon^2 , \qquad s \in {\mathcal T} , $$
and, by possibly assuming that  $ C_0 C_1
\epsilon < 1/8 $, we have $ \dot{r}_s^+ > 0  $ on $ {\mathcal T}  $.
Choosing $ C_M \geq 1 $ as in Proposition \ref{libre}, we now
claim that, if
$$ \epsilon < \epsilon_2 / C_M \qquad \mbox{and} \qquad  r > R (C_M \epsilon) , $$
 then $ T:= \sup {\mathcal T} = + \infty $. Assume that this is
wrong. Then $ T $ is finite, belongs to $ {\mathcal T} $ and,
on $ [0,T] $, we have
$$ r^+_s \geq r + s / 8 \geq r  , \qquad |\theta^+_s - \theta| \leq C_1 \epsilon^2 < M \epsilon^2 ,  $$
so, by  Proposition \ref{libre}, $ (r^+_s,\theta^+_s,\rho,\eta)
\in \Gamma^+_{\rm s} (C_M \epsilon) \subset \Gamma^+_{\rm s} (
\epsilon_2) $ and, by (\ref{estimeechamps}), (\ref{estimeechamps2}),
\begin{eqnarray}
 |r^+_T - r - 2 \rho T| \leq C_0 e^{-r}|\eta| \int_0^T e^{-s/8}
ds \leq C_0 e^{-r}|\eta| T  \leq C_0 C_1 \epsilon T , \label{rplusbis}  \\
|\theta_T^+ - \theta| \leq C_0 e^{-2r}|\eta| \int_0^T e^{-s/4} ds
\leq 4 C_0 e^{-2 r}|\eta| < 5 C_0 C_1 \epsilon^2 .
\label{thetaplusbis}
\end{eqnarray}
This implies that $ r_T^+ > r + T / 8 $ and that $
|\theta^+_T-\theta| < M \epsilon^2 $ so the flow can be continued beyond $ T
$, yielding a contradiction
with the definition of $ T $. The flow is thus well defined for $ t \geq 0 $.
Then, (\ref{rplus}) and
(\ref{thetaplus}) follow from the first inequalities of (\ref{rplusbis}) and
(\ref{thetaplusbis}) with an arbitrary $ t \geq 0 $ instead of $ T
$, since $ e^{-r}|\eta| \lesssim \epsilon^2  $ for $ (r,\theta,\rho,\eta) \in
 \Gamma^{+}_{\rm s}(\epsilon^2)  $.  If $ \epsilon  $ is small enough, Proposition \ref{libre},
 (\ref{localisationphasesepsilon}) shows that
is a direct
 consequence of (\ref{rplus}) and (\ref{thetaplus}), using that $ e^{-r} \ll \epsilon^4  $.

 It remains to prove (\ref{rthetauniftemps}) for $ j + |\alpha| +
k + |\beta| \geq 1 $.
 We consider  $ \overline{r}_t^+ := r^+_t
- 2 t \rho $ and $ \overline{\theta}^+_t := \theta^+_t $, which
satisfy
\begin{eqnarray}
 d \overline{r}_t^+ / d t   =  \overline{v} (t ,\overline{r}^+_t ,
 \overline{\theta}^+_t , \rho , \eta ) , \qquad
 d \overline{\theta}_t^+ / d t   =  \overline{w} (t ,\overline{r}^+_t ,
 \overline{\theta}^+_t , \rho , \eta ) , \label{appliquerAjkalphabeta}
\end{eqnarray}
with
\begin{eqnarray}
 \overline{v}(t,r,\theta,\rho,\eta) & = & (\partial_r
\varphi_+) (r + 2 t \rho , \theta , \rho ,  \eta)
 , \nonumber \\
 \overline{w}(t,r,\theta,\rho,\eta) & = & e^{-2 r - 4 t \rho} \left(
 \partial_{\eta} q \right) \left( r + 2 t \rho , \theta ,  \partial_{\theta} S_+ (r+2 t \rho,\theta,
 \rho,\eta) \right) . \nonumber
\end{eqnarray}
Using (\ref{crucialphasestationnaire}), we have, for all $
j^{\prime},\alpha^{\prime},k^{\prime},\beta^{\prime} $,
\begin{eqnarray}
| D_{\rm hyp}^{j^{\prime} \alpha^{\prime} k^{\prime}
\beta^{\prime}} (\overline{v},\overline{w}) | \lesssim
\scal{t}^{k^{\prime}} e^{-4 t \rho} \lesssim e^{-2 t \rho}, \qquad
t \geq 0, \ \ \mbox{on} \ \Gamma^{+}_{\rm s} (\epsilon_2 / C),
\label{decroissanceexponentielleflot2}
\end{eqnarray}
with $ C $ such that if $ (r,\theta,\rho,\eta) \in
\Gamma^{+}_{\rm s}(\epsilon_2 /C) $ then $ (r + 2 t \rho ,
\theta , \rho,\eta ) \in \Gamma^{+}_{\rm s}(\epsilon_2 /C)
$. Note also that if $ \epsilon $ is small enough and $
(r,\theta,\rho,\eta) \in
 \Gamma^{+}_{\rm s}(\epsilon^2) $, we have $ ( \overline{r}_t^+ , \overline{\theta}_t^+,\rho,\eta )
 \in \Gamma^{+}_{\rm s}(\epsilon_2 /C) $, using (\ref{rplus}), (\ref{thetaplus}) and
 Proposition \ref{libre}. We then obtain (\ref{rthetauniftemps}) by induction by applying $ D_{\rm
hyp}^{j \alpha k \beta} $ to (\ref{appliquerAjkalphabeta}).
Indeed, using Lemma \ref{Faadibruno} and
(\ref{decroissanceexponentielleflot2}), we have
$$ \frac{d}{d t} D_{\rm
hyp}^{j \alpha k \beta} \left( \overline{r}_t^+ ,
\overline{\theta}_t^+ \right) = \left(
\partial_{r,\theta} (\overline{v} , \overline{w} ) \right)   D_{\rm
hyp}^{j \alpha k \beta} ( \overline{r}_t^+ ,
\overline{\theta}_t^+ ) + {\mathcal O}(e^{-2 \rho t}),
$$
where $ {\mathcal O}(e^{-2 \rho t}) = 0 $ for first order
derivatives and, otherwise, follows from the induction assumption.
Since $ |\partial_{r,\theta} (\overline{v} , \overline{w} )| \lesssim
e^{- 2 \rho t } $,
 Lemma \ref{Gronwallexponentiel} yields the result. \finpreuve

\bigskip

We now come to the resolution of (\ref{formetransport}) in a way
suitable to further purposes.

\begin{prop} \label{suffisanttransport} There exists $ \epsilon_5 > 0 $ such that, for all $
0 < \epsilon \leq \epsilon_5 $ and all $ y^{\pm} \in {\mathcal
B}_{\rm hyp} (\Gamma^{\pm}_{\rm s}(\epsilon))$  of hyperbolic
short range in $ \Gamma^{\pm}_{\rm s}(\epsilon) $,  the function
\begin{eqnarray}
 a^{\pm}_{\rm hom} = \exp \left(  \int_0^{\pm \infty} y^{\pm} \circ
    \phi^{\pm}_s d s   \right) , \nonumber
\end{eqnarray}
solves (\ref{formetransport}) on $ \Gamma^{\pm}_{\rm
s}(\epsilon^2) $ with $ z^{\pm} \equiv 0 $ , belongs to $
{\mathcal B}_{\rm hyp} (\Gamma^{\pm}_{\rm s}(\epsilon^2)) $ and
$$  a^{\pm}_{\rm hom} - 1 \ \mbox{is of
hyperbolic long range in} \ \ \Gamma^{\pm}_{\rm s}(\epsilon^2) .
$$
In addition, for all $z^{\pm} \in
    {\mathcal B}_{\rm hyp} \left( \Gamma^{\pm}_{\rm
      s} ( \epsilon )   \right) $, of hyperbolic
short range in $ \Gamma^{\pm}_{\rm s}(\epsilon) $, the function
\begin{eqnarray}
a^{\pm}_{\rm inhom} = - \int_0^{\pm \infty} z^{\pm} \circ
\phi^{\pm}_s
   \exp \left( \int_0^s y^{\pm} \circ
    \phi^{\pm}_u d u   \right) ds , \nonumber
\end{eqnarray}
solves (\ref{formetransport}) on $ \Gamma^{\pm}_{\rm
s}(\epsilon^2) $, belongs to $ {\mathcal B}_{\rm hyp}
(\Gamma^{\pm}_{\rm s}(\epsilon^2)) $ and
$$  a^{\pm}_{\rm inhom}  \ \mbox{is of
hyperbolic long range in} \ \ \Gamma^{\pm}_{\rm s}(\epsilon^2) .
$$
\end{prop}


 We need the
following lemma.

\begin{lemm} \label{lemmelongueportee}  There exists $ \epsilon_5 > 0 $ such
that, for all $ j ,\alpha, k,\beta $ and all $ N \geq 0 $,
\begin{eqnarray}
 \big| \partial_r^j \partial_{\theta}^{\alpha} \partial_{\rho}^k \partial_{\eta}^{\beta}
  \big( r_t^{\pm} - r - 2 t \rho \big)  \big|
+  \big| \partial_r^j \partial_{\theta}^{\alpha} \partial_{\rho}^k
\partial_{\eta}^{\beta} (\theta_t^{\pm} - \theta )  \big| &  \lesssim &
\scal{r-\log
  \scal{\eta}}^{-N} ,
\nonumber
\end{eqnarray}
on  $  \Gamma^{\pm}_{\rm s}(\epsilon_5)  $, uniformly with respect
to $ \pm t \geq 0 $.
\end{lemm}

\noindent {\it Proof.} By Proposition \ref{libre}, there exists $
C > 0 $ such that, for all $ \epsilon $ small enough and all $ s
\in [0,1] $,
\begin{eqnarray}
(r,\theta,\rho,\eta) \in \Gamma^{\pm}_{\rm s}(\epsilon^2)  \
\Rightarrow \ (r,\theta,\rho,s \eta) \in \Gamma^{\pm}_{\rm s}(C
\epsilon^2) .
\end{eqnarray}
Therefore, if $ C \epsilon^2 \leq \epsilon_4^2 $ and if  we set $
X_t^{\pm} (r,\theta,\rho,\eta) =( r^{\pm}_t - r - 2 t \rho ,
\theta^{\pm}_t - \theta) $, we can write
$$ X_t^{\pm} (r,\theta,\rho,\eta) = X_t^{\pm} (r,\theta,\rho,0) + \int_0^{1} (e^{r}\partial_{\eta}
X_t^{\pm}) (r,\theta,\rho,s \eta) ds \cdot e^{-r}\eta ,  $$ on $
\Gamma^{\pm}_{\rm s}(\epsilon^2)  $. The crucial remark is that $
X_t^{\pm} (r,\theta,\rho,0) = 0 $. Indeed, by
(\ref{formepratique}) and (\ref{ordre0}), we have $ \partial_r
S_{\pm} \equiv \rho $ and $
\partial_{\theta} S_{\pm} \equiv 0 $  at $ \eta = 0  $  (notice that $ (r,\theta,\rho,0) \in
\Gamma^{\pm}_{\rm s}(\epsilon_2) $ if  $ C \epsilon^2 \leq
\epsilon_2 $),  so the solution to (\ref{systemetransport}) is
simply $ (r + 2 t \rho , \theta)  $ in this case. In addition, by
(\ref{rthetauniftemps}), $ (X_{t}^{\pm})_{t \geq 0} $ is bounded
in $ {\mathcal B}_{\rm hyp} ( \Gamma^{\pm}_{\rm s}(\epsilon^2) )$.
Thus, for all $ N \geq 0 $,
$$ \big| \partial_r^j \partial_{\theta}^{\alpha} \partial_{\rho}^k
\partial_{\eta}^{\beta} X_t^{\pm} (r,\theta,\rho,\eta)  \big|   \lesssim
e^{-r} \scal{\eta} \lesssim \scal{r-\log \scal{\eta}}^{-N} ,
\qquad \pm t \geq 0, \ (r,\theta,\rho,\eta) \in \Gamma^{\pm}_{\rm
s}(\epsilon^2),
$$ which yields the result.  \finpreuve

\bigskip

\noindent {\it Proof of Proposition \ref{suffisanttransport}.}
 For simplicity we set $ \partial^{\gamma} = \partial_r^j
\partial_{\theta}^{\alpha} \partial_{\rho}^k
\partial_{\eta}^{\beta} $. Then, using Lemma \ref{Faadibruno} with $ |\gamma| \geq 1 $, 
$  \partial^{\gamma}  \left( y^{\pm} \circ \phi_s^{\pm} \right) $ is the sum of
\begin{eqnarray}
(\partial_{r} y^{\pm})\circ \phi^{\pm}_s \partial^{\gamma}
 r^{\pm}_s + (\partial_{\theta} y^{\pm})\circ \phi^{\pm}_s \cdot
 \partial^{\gamma} \theta^{\pm}_s + \delta_{j 0} \delta_{\alpha 0}
 (\partial^k_{\rho} \partial_{\eta}^{\beta} y^{\pm}) \circ
 \phi^{\pm}_s \label{avecKronecker}
\end{eqnarray}
and of a linear combination of
\begin{eqnarray}
 (\partial^{k-k^{\prime}}_{\rho}  \partial_{\eta}^{\beta-\beta^{\prime}} \partial_{r,\theta}^{\nu}
  y^{\pm}) \circ \phi^{\pm}_s  \left( \partial^{\gamma_1^1} r^{\pm}_s  \ldots
\partial^{\gamma_{\nu_1}^1} r^{\pm}_s  \right) \ldots \left( \partial^{\gamma_1^{n}}
(\theta^{\pm}_s)_{n-1} \ldots
\partial_x^{\gamma_{\nu_{n}}^{n}} (\theta^{\pm}_s)_{n-1} \right),
\label{sansKronecker}
\end{eqnarray}
where $ (\theta^{\pm}_s)_{1}, \ldots, (\theta^{\pm}_s)_{n-1} $ are
the components of $ \theta^{\pm}_s $, $(0,0,k^{\prime},
\beta^{\prime} ) + \sum \gamma_i^j = \gamma $, using the
convention and the notation of Lemma \ref{Faadibruno}. By (\ref{rthetauniftemps}),
we have
$$  | (\partial_{r} y^{\pm})\circ \phi^{\pm}_s \partial^{\gamma}
 r^{\pm}_s | \lesssim \scal{r^{\pm}_s - \log \scal{\eta}}^{-\tau-2} e^{-r|\beta|} \scal{s}^{\kappa} ,  $$
 where $ \kappa = 1 $ if $ k = 1 $ and $ j +|\alpha|+|\beta| = 0
 $, and $ \kappa = 0  $ otherwise. On the other hand, by Lemma
 \ref{lemmelongueportee}, we have
$$  | (\partial_{r} y^{\pm})\circ \phi^{\pm}_s \partial^{\gamma}
 r^{\pm}_s | \lesssim \scal{r^{\pm}_s - \log \scal{\eta}}^{-\tau-2}
 \scal{r-\log\scal{\eta}}^{- \tilde{j}}  \scal{s}^{\kappa} ,  $$
with the same $ \kappa $ as above and $ \tilde{j} = j $ if $ j
\geq 2 $, or $\tilde{j} = 0 $ for $ j \leq 1 $. Similarly, we also
have
$$ | (\partial_{\theta} y^{\pm})\circ \phi^{\pm}_s \cdot
 \partial^{\gamma} \theta^{\pm}_s | \lesssim \scal{r^{\pm}_s - \log
 \scal{\eta}}^{-\tau-1}\times
\min \left( e^{-|\beta|r} , \scal{r - \log \scal{\eta}}^{-j}
\right) ,
$$ while, for the last term of (\ref{avecKronecker}), we have
$$ | \delta_{j 0} \delta_{\alpha 0}
(\partial^k_{\rho} \partial_{\eta}^{\beta} y^{\pm}) \circ
\phi^{\pm}_s  | \lesssim \min \left( e^{-|\beta| r}
e^{-2|\beta||\rho s|} , \scal{r^{\pm}_s - \log
 \scal{\eta}}^{-\tau-1-j}   \right) , $$
since $ e^{-|\beta|r^{\pm}_s} \lesssim e^{-|\beta| r}
e^{-2|\beta||s|} $ for $ r_s^{\pm} - r - 2 \rho s $ is  bounded
from below and $ \rho s \geq 0 $. Now, we remark that
$$ \big| \left( \partial^{\gamma_1^1} r^{\pm}_s  \ldots
\partial^{\gamma_{\nu_1}^1} r^{\pm}_s  \right) \big| \lesssim \scal{s}^{\tilde{\nu}_1}
\scal{r-\log \scal{\eta}}^{-N_0} , $$ where $ \tilde{\nu}_1 $ is
the number of $ \partial^{\gamma^1_{l}} =  \partial_r^{j^1_l}
\partial_{\theta}^{\alpha_l^1} \partial_{\rho}^{k^1_l}
\partial_{\eta}^{\beta_l^1}  $ for which $ j^1_l = 0 $, $ N_0 = 0 $
if  $ j^1_l \leq 1 $ for all $ l  $ and $ N_0 $ is any positive
number if $ j^1_l \geq 2 $ for at least one $ l $. We therefore
obtain, if $ \beta = \beta^{\prime} $,
$$ |(\ref{sansKronecker})| \lesssim \scal{r^{\pm}_s - \log
 \scal{\eta}}^{-\tau-1-\nu_1}  \min \left( e^{-r|\beta|} \scal{s}^{\nu_1} , \scal{r - \log
 \scal{\eta}}^{\nu_1-\tilde{\nu}_1 - j} \scal{s}^{\tilde{\nu_1}}  \right) , $$
 since $ \nu_1-\tilde{\nu}_1 - j \geq 0  $ in the case where no
 $r$ derivative fall on the components of $ \theta^{\pm}_s $ and
 only $ r$  derivatives of order at most $1$ fall on $ r^{\pm}_s
 $. If $ \beta \ne \beta^{\prime}
 $, we have
$$ |(\ref{sansKronecker})| \lesssim \min \left( e^{-2 |\beta - \beta^{\prime}| |\rho s|} e^{- |\beta|r} \scal{s}^{\nu_1} ,
 \scal{r^{\pm}_s - \log
 \scal{\eta}}^{-\tau-1-\nu_1} \scal{r - \log
 \scal{\eta}}^{\nu_1 - \tilde{\nu}_1-j} \scal{s}^{\tilde{\nu_1}}
    \right)  .$$
Since $ r^{\pm}_s - r - 2 \rho s $ is bounded from below, $ \rho s
\geq 0 $ (with $ |\rho | \gtrsim 1 $) and using
(\ref{controlelongueportee}), we have
$$  \scal{r^{\pm}_s - \log
 \scal{\eta}}^{-\tau-1-\nu_1} \lesssim \scal{r - \log
 \scal{\eta} +  |s|}^{-\tau-1-\nu_1} . $$
 All this implies that
 \begin{eqnarray*}
 | D_{\rm hyp}^{j \alpha k \beta} \left( y^{\pm} \circ
\phi_s^{\pm} \right) | & \lesssim & \scal{s}^{-\tau-1}, \\
 |
\partial_r^j
\partial_{\theta}^{\alpha} \partial_{\rho}^k
\partial_{\eta}^{\beta} \left( y^{\pm} \circ
\phi_s^{\pm} \right)|  & \lesssim & \scal{r - \log
 \scal{\eta} +  |s|}^{-\tau-1} \scal{r - \log
 \scal{\eta} }^{- j} ,
\end{eqnarray*}
and since
$$ \int_0^{+\infty}  \scal{r-\log \scal{\eta}+|s|}^{- \tau - 1} d s \lesssim
\scal{r-\log \scal{\eta}}^{- \tau } , $$
(using  (\ref{controlelongueportee}) on strongly outgoing/incoming areas),
 we see that the function
$ \int_{0}^{\pm \infty} y^{\pm} \circ \phi^{\pm}_s ds $ belongs to
$ {\mathcal B}_{\rm hyp}(\Gamma^{\pm}_{\rm s}(\epsilon^2)) $ and
is of hyperbolic long range. This implies easily that the same
holds for $ a^{\pm}_{\rm hom} -1 $. One then checks that $
a^{\pm}_{\rm hom} $ solves the homogeneous transport equation by
computing $ d (a^{\pm}_{\rm hom} \circ \phi^{\pm}_t ) / d t $ at $
t = 0^{\pm} $. One studies similarly the case of $ a^{\pm}_{\rm
inhom} $. \finpreuve

\section{An Isozaki-Kitada type parametrix} \label{IsozakiKitada}
\setcounter{equation}{0}

In this section, we prove an  approximation of $ e^{-ithP} \widehat{O \! p}_{\iota}(\chi_{\rm s}^{\pm}) $ when $ \chi_{\rm s}^{\pm} $
is supported in the strongly outgoing ($+$)/ incoming ($-$) region $ \Gamma^{\pm}_{\iota,{\rm s}}(\epsilon) $ (see Definition \ref{zonesstrong} 
for these areas and Definition \ref{definitOpchapeau} for
$ \widehat{O \! p}_{\iota} (\cdot) $). We recall that $ \iota $ is an arbitrary index corresponding to the chart at infinity we consider and where  the symbols are supported (see (\ref{Diffeoaveciota}) and (\ref{presquepartition})).

Here we will prove an $ L^2 $ approximation, valid for times such that $ 0 \leq \pm t  \lesssim h^{-1} $.
Basically, we will show that, for any $ N $, $ e^{-ithP} \widehat{O \! p}_{\iota}(\chi_{\rm s}^{\pm}) $ is the sum of a Fourier integral operator
and of a term of order $ h^N $ in the operator norm of $ L^2 ({\mathcal M}, \widehat{dG}) $, uniformly for $ 0 \leq t \lesssim h^{-1} $.

We will therefore essentially prove half of Proposition \ref{sousIK}, namely the estimate (\ref{resteL2IsozakiKitada}).
The dispersion estimate (\ref{dispersionIsozakiKitada}), following from a stationary phase
argument on the Fourier integral operator, will be proved in Section \ref{sectiondispersion}.

\bigskip

In Theorem \ref{IsozakiKitadaansatz} below, we use the classes of symbols
$ {\mathcal S}_{\rm hyp}(\cdot) $ introduced in
Definition \ref{definhyp} and the Fourier integral operators  (\ref{OIFRn}) defined in Subsection \ref{operateurintegraldeFourier}.
For these operators, the phases are associated to the Hamiltonian $ p = p_{\iota} $, the principal symbol of $ P $ in the $ \iota$-th chart, 
as explained in the beginning of Section \ref{sectionHJ}.

In the sequel, we fix $ \iota \in {\mathcal I} $ (see (\ref{cover})).

\begin{theo} \label{IsozakiKitadaansatz}  For all $ N \geq 0  $, there exists  $
  \epsilon (N)  > 0 $ such that, for all $ 0  < \epsilon \leq
  \epsilon (N) $,
the following holds:   there exists  $ a^{\pm} (h) = a_0^{\pm}
+ \cdots + h^N a_N^{\pm} $ with
$$ a^{\pm}_0 , \ldots , a^{\pm}_N  \in {\mathcal S}_{\rm hyp} \left( \Gamma^{\pm}_{\iota,{\rm s}}( \epsilon)  \right) , $$
such that for all
\begin{eqnarray}
 \chi^{\pm}_{\rm s} \in {\mathcal S}_{\rm hyp} \left( \Gamma^{\pm}_{\iota,{\rm s}}( \epsilon^9)  \right)  ,
\end{eqnarray}
we can find  $ b^{\pm}(h) =  b_0^{\pm} + \cdots + h^N b_N^{\pm} $,
with
\begin{eqnarray}
b^{\pm}_0 , \ldots , b^{\pm}_N \in    {\mathcal S}_{\rm hyp} \left( \Gamma^{\pm}_{\iota,{\rm s}}( \epsilon^3)  \right)   ,
\end{eqnarray}
such that, for all $
T > 0 $,   there exists $ C > 0 $ such that
\begin{eqnarray}
 \left| \left|  \ e^{-it h P
} \widehat{{O \! p}}_{\iota} \left( \chi^{\pm}_{\rm s}
 \right)
 - \Psi_{{\iota}}^* \left( J^{\pm}_h \left( a^{\pm} (h) \right) e^{- i t h
D_r^2} J^{\pm}_h \left( b^{\pm } (h) \right)^* \right) \big(
\Psi_{\iota}^{-1} \big)^* \right| \right|_{L^2 ( \widehat{dG}) \rightarrow L^2 ( \widehat{dG})}
\leq C h^{N-1} , \nonumber
\end{eqnarray}
provided
$$ 0 \leq  \pm t \leq T h^{-1} , \qquad h \in (0,1 ] . $$
\end{theo}

We emphasize that, by (\ref{pourLp}), (\ref{margeouverts}) and (\ref{RVepsiloniota}), the symbols $ a^{\pm}(h) $ and $ b^{\pm} (h) $
are supported in $ (\epsilon^{-1},+\infty) \times V_{\iota, \epsilon} \times \Ra^n \subset 
 (R_{\mathcal K}+1,+\infty) \times V_{\iota}^{\prime} \times \Ra^n $, for $ \epsilon $ small. Therefore 
 the Schwartz kernel of  the operator $ J^{\pm}_h \left( a^{\pm} (h)
\right) e^{- i t h D_r^2} J^{\pm}_h \left( b^{\pm } (h) \right)^*
$ is supported in $ \left( (R_{\mathcal K}+1,+\infty) \times V_{\iota}^{\prime}
\right)^2 $ and hence
$$ \Psi_{{\iota}}^* \left( J^{\pm}_h \left( a^{\pm} (h) \right) e^{- i t h
D_r^2} J^{\pm}_h \left( b^{\pm } (h) \right)^* \right) \big(
\Psi_{\iota}^{-1})^* $$ is well defined on the whole manifold (by the implicit requirement that its kernel vanishes outside the coordinate
 patch $ {\mathcal U}_{\iota} \times {\mathcal U}_{\iota} $ of $ {\mathcal M} \times {\mathcal M} $).

We also remark that $ \epsilon(N)  $  could certainly  be chosen
independently of $ N  $. However this is useless for the applications we have
in mind and
we will not consider this refinement.

\bigskip

Before starting the proof, we fix or recall some notation. We set 
$$ P_{\iota} = \big( \Psi_{\iota}^{-1}
\big)^* P \big( \Psi_{{\iota}} \big)^*  = p (r,\theta,D_r,D_{\theta}) + p_1 (r,\theta,D_r,D_{\theta}) + p_2 (r,\theta) $$
with $ p $ the principal symbol  and $ p_{k} $ of degree $2-k$ in
$ (\rho,\eta) $ for $ k = 1,2 $.

Accordingly, we drop the indices $ \iota $ from the notation of strongly outgoing/incoming areas . In other words, we consider the areas
given by Definition \ref{zonesstrongsansiota} with  $ V_0 = V_{\iota} $. The explanation at the beginning of Section \ref{sectionHJ} shows this should no cause any ambiguity.
 
For simplicity, we also use the notation (\ref{simplificationclasses}).

Recall finally that, for some fixed $ \epsilon_{\iota} > 0 $ small enough, Proposition \ref{Phaseplusmoinsglobale} proves the existence of $
S_{\pm} $  solving
\begin{eqnarray}
 p (r,\theta,\partial_r S_{\pm} , \partial_{\theta} S_{\pm} ) =
\rho^2 , \qquad (r,\theta,\rho,\eta) \in \Gamma^+_{\rm
s}(\epsilon_{\iota}) . \label{HamiltonJacobieffectif}
\end{eqnarray}

\bigskip

\noindent {\it Proof of Theorem \ref{IsozakiKitadaansatz}.} We
denote for simplicity
$$  A_{\pm} = J_h^{\pm} (a^{\pm}(h)) ,
\qquad B_{\pm} = J_h^{\pm} (b^{\pm}(h)) . $$ By the Duhamel
formula,
we have
\begin{eqnarray}
e^{-i t h P} \Psi_{{\iota}}^* A_{\pm} = \Psi_{{\iota} }^* A_{\pm}
e^{-it h D_r^2} - \frac{i}{h} \int_0^t e^{-i(t-s)h P}
\Psi_{{\iota}}^* (h^2 P_{\iota} A_{\pm} - A_{\pm} h^2 D_r^2)
e^{-ishD_r^2} ds . \label{Duhamelclef}
\end{eqnarray}
Multiplying (\ref{Duhamelclef}) by $ B^*_{\pm} (\Psi_{\iota}^{-1})^* $
and  denoting
\begin{eqnarray}
 C_{\pm} & := &  \chi^{\pm}_{\rm s} (r,\theta,h D_r,h D_{\theta})( \tilde{\kappa} \otimes \tilde{\kappa}_{\iota} )  -
  A_{\pm} B^*_{\pm} , \\
 D_{\pm}(s) & : = & (h^2 P_{\iota} A_{\pm} - A_{\pm} h^2 D_r^2) e^{-ishD_r^2}
 B_{\pm}^*,
\end{eqnarray}
(where $ \tilde{\kappa} $ and $ \tilde{\kappa}_{\iota} $ are the cutoffs used in Definition \ref{definitOpchapeau}),
we obtain 
\begin{eqnarray*}
  e^{-it h P}
\widehat{O \! p}_{\iota} (\chi^{\pm}_{\rm s}) & = &
\Psi_{{\iota}}^* A_{\pm} e^{-ithD_r} B^*_{\pm} (\Psi_{\iota}^{-1})^* \\
& & \ \ \ \  + e^{-it
h P} \Psi_{{\iota}}^* C_{\pm}  (\Psi_{\iota}^{-1})^* - \frac{i}{h}
\int_0^t e^{-i(t-s)h P} \Psi_{{\iota}}^* D_{\pm} (s)
(\Psi_{\iota}^{-1})^* ds .
\end{eqnarray*}
Using (\ref{equivalenceLp}) with $ q= 2 $, the theorem will then be proved if we find $ a^{\pm}(h) $ and
$ b^{\pm}(h) $ such that
\begin{eqnarray}
 || C_{\pm} ||_{L^2(\Ra^n) \rightarrow L^2(\Ra^n) } \lesssim h^N,
\qquad || D_{\pm} (s) ||_{L^2(\Ra^n) \rightarrow L^2(\Ra^n) }
\lesssim h^{N+1} , \qquad h \in (0,1] , \label{conditionIK}
\end{eqnarray}
 uniformly with respect to
$ 0 \leq \pm s \leq T h^{-1}$ for $  D_{\pm} (s) $.

 For simplicity we only consider the outgoing case but the
incoming one is of course completely similar.

\medskip

\noindent {\bf Construction of $  a^+ (h) $.} 
We first define $ (v^{+},w^{+})  $ by (\ref{champsdevecteurs}) and also set
\begin{eqnarray}
y^{+} & := & p (r,\theta,\partial_r , \partial_{\theta}) S_+
  +  p_1 (r,\theta,\partial_r,\partial_{\theta}) S_{+}. \label{vraiy}
\end{eqnarray}
\begin{lemm} \label{vraitransport} There exists $  \tilde{\epsilon}_{\iota} \leq \epsilon_{\iota} $ such that $
  y^+  $ belongs to  $  {\mathcal B}_{\rm
    hyp}(\Gamma^+_{\rm s}(\tilde{\epsilon}_{\iota}))  $ and is of hyperbolic short range on
  $ \Gamma^+_{\rm s}(\tilde{\epsilon}_{\iota})  $.
\end{lemm}

\noindent {\it Proof.} It follows from (\ref{coefficientsconcrets}) and
(\ref{caracterisationslrange}) since Proposition
\ref{Phaseplusmoins} shows that $ y^+_{| \eta = 0} \equiv 0 $.
\finpreuve

\bigskip

Elementary computations show that, for all $ a_0^+ , \ldots ,
a_N^+ \in {\mathcal S}_{\rm hyp}^{+} ( \epsilon) $ and $ a^+ (h) =
a_0^+ + \cdots + h^N a_N^+ $,
$$ h^2 P_{\iota}
J_h^+ (a^+ (h)) - J_h^+ (a^+ (h)) h^2 D_r^2 = \sum_{l=0}^{N+2} h^l
J_h^+ (d_l^+),
$$ where the symbols are given by
\begin{eqnarray}
d^+_l  & = &  \left( p (r, \theta , \partial_r S_+ ,
\partial_{\theta} S_+ ) - \rho^2 \right) a_{l}^{+} - i \left(
v^+ \partial_r a_{l-1}^+ + w^+ \cdot \partial_{\theta} a_{l-1}^+ +
 y^+ a_{l-1}^+ \right) + P_{\iota} a_{l-2}^+ \nonumber \\
& = &  - i \left( v^+ \partial_r a_{l-1}^+ + w^+ \cdot
\partial_{\theta} a_{l-1}^+ +
 y^+ a_{l-1}^+ \right) + P_{\iota} a_{l-2}^+ ,
\end{eqnarray}
using (\ref{HamiltonJacobieffectif}) and assuming $ \epsilon \leq
\epsilon_{\iota} $. Here, we have
 $ 0 \leq l \leq N+2 $ and the convention that $ a^{+}_{-2} =
a^{+}_{-1} = a^+_{N+1} = a^{+}_{N+2} \equiv 0 $. In particular,
the first three terms  are given by
\begin{eqnarray}
d_0^+   & = &   0 , \\
i d_1^+ &  = &  v^+ \partial_r a_{0}^+ + w^+ \cdot
\partial_{\theta} a_{0}^+ +
 y^+ a_{0}^+, \\
i d_2^+   & = & v^+ \partial_r a_{1}^+ + w^+ \cdot
\partial_{\theta} a_{1}^+ +
 y^+ a_{1}^+      + i P_{\iota} a_0^+ .
\end{eqnarray}
Using Proposition
\ref{suffisanttransport}, Lemma \ref{vraitransport} and assuming $
\hat{\epsilon}_{\iota} \leq \min (\tilde{\epsilon}_{\iota}^2,\epsilon_5)  $  we can define
$$ \hat{a}^+_0 (r,\theta,\rho,\eta) = \exp \left( \int_0^{+\infty} y^{+} \circ
  \phi^+_s(r,\theta,\rho,\eta) ds  \right) , \qquad (r,\theta,\rho,\eta) \in
\Gamma^{+}_{\rm s}(\hat{\epsilon}_{\iota})
 $$ so that $ \hat{a}^+_0 \in {\mathcal B}_{\rm hyp} ( \Gamma^{+}_{\rm
   s}(\hat{\epsilon}_{\iota}) )  $, $ \hat{a}^+_0 - 1  $ is of hyperbolic
 long range in $ \Gamma^{+}_{\rm s}(\hat{\epsilon}_{\iota})  $ and
\begin{eqnarray}
  v^+ \partial_r \hat{a}^+_0 + w^+ \cdot \partial_{\theta} \hat{a}^+_0 +
 y^+ \hat{a}^+_0 \equiv 0 \qquad \mbox{on} \ \  \Gamma^{+}_{\rm
   s}(\hat{\epsilon}_{\iota}) . \nonumber
\end{eqnarray}
Since the function $ \int_0^{\infty} y^+ \circ \phi^+_s ds  $ is
bounded on $  \Gamma^{+}_{\rm
   s}(\hat{\epsilon}_{\iota}) $ (see the proof of Proposition
 \ref{suffisanttransport}), we also have
\begin{eqnarray}
\hat{a}_0^+ (r,\theta,\rho,\eta) \gtrsim 1   , \qquad
(r,\theta,\rho,\eta) \in \Gamma^+_{\rm s} (\hat{\epsilon}_{\iota} ) .
\label{verifieellipticiteFourier}
\end{eqnarray}
Using (\ref{coefficientsconcrets}) and the fact that $ \hat{a}_0^+
- 1  $ is of hyperbolic long range, it is easy to check that $
P_{\iota} \hat{a}_0^+  $ is of hyperbolic short range
in $ \Gamma^+_{\rm s} (\hat{\epsilon}_{\iota}^2 )   $. By
Proposition \ref{suffisanttransport}, we can then define
$$  \hat{a}^+_1  = i \int_0^{+\infty}
(P_{\iota}  \hat{a}_0^+) \circ \phi_s^+ \exp \left( \int_0^{s} y^{+} \circ
  \phi^+_u du  \right) ds  , \qquad \mbox{on} \ \ \Gamma^{+}_{\rm
  s}(\hat{\epsilon}_{\iota}^2) , $$
which belongs to $   {\mathcal B}_{\rm hyp} ( \Gamma^{+}_{\rm
   s}(\hat{\epsilon}_{\iota}^2) )  $,  is of hyperbolic
 long range in $ \Gamma^{+}_{\rm s}(\hat{\epsilon}_{\iota})  $ and satisfies
\begin{eqnarray}
  v^+ \partial_r \hat{a}^+_1 + w^+ \cdot \partial_{\theta} \hat{a}^+_1 +
 y^+ \hat{a}^+_1 \equiv - i P_{\iota} \hat{a}_0^+
\qquad \mbox{on} \ \  \Gamma^{+}_{\rm s}(\hat{\epsilon}_{\iota}^2) . \nonumber
\end{eqnarray}
More generally, for $ 1 \leq l \leq N $, we can define iteratively
$$  \hat{a}^+_l  = i \int_0^{+\infty}
(P_{\iota}  \hat{a}_{l-1}^+) \circ \phi_s^+ \exp \left( \int_0^{s} y^{+} \circ
  \phi^+_u du  \right) ds  , \qquad \mbox{on} \ \ \Gamma^{+}_{\rm
  s}(\hat{\epsilon}_{\iota}^{2^l}) , $$
which belongs to $   {\mathcal B}_{\rm hyp} ( \Gamma^{+}_{\rm
   s}(\hat{\epsilon}_{\iota}^{2^l}) )  $,  is of hyperbolic
 long range in $ \Gamma^{+}_{\rm s}(\hat{\epsilon}_{\iota}^{2^l})  $ and satisfies
\begin{eqnarray}
  v^+ \partial_r \hat{a}^+_l + w^+ \cdot \partial_{\theta} \hat{a}^+_l +
 y^+ \hat{a}^+_l \equiv - i P_{\iota} \hat{a}_{l-1}^+
\qquad \mbox{on} \ \  \Gamma^{+}_{\rm s}(\hat{\epsilon}_{\iota}^{2^l}) , \nonumber
\end{eqnarray}
using Proposition  \ref{suffisanttransport} and the fact that $ P_{\iota} \hat{a}_{l-1}^+  $
is of hyperbolic short range if $  \hat{a}_{l}^+  $ is of hyperbolic long range.
 Hence, using Lemma \ref{construitcutoff} with $ \epsilon \leq \hat{\epsilon}_{\iota}^{2^N}  $ and setting
$$  a_l^+ = \chi^{+}_{\epsilon^2 \rightarrow \epsilon} \hat{a}^+_l , \qquad 0 \leq l \leq N $$
with the $ \hat{a}_l^+ $ defined above, we have constructed $
a_0^+ , \ldots , a_N^+ \in {\mathcal S}^+_{\rm hyp}(\epsilon) $
with $ a_0^+  $ satisfying (\ref{ellipticiteFourier}), by
(\ref{verifieellipticiteFourier}). Furthermore,
\begin{eqnarray}
 d_l^+  \in  {\mathcal S}^{+}_{\rm hyp}(\epsilon)  \ \ \mbox{ for } \  0 \leq
l \leq N+2, \label{nonstationnaire0} \\
 d_l^+ \equiv 0 \ \
\mbox{on} \ \Gamma^{+}_{\rm s} (\epsilon^2) \ \ \mbox{for} \ 0
\leq l \leq N . \label{nonstationnaire1}
\end{eqnarray}

\medskip

\noindent {\bf Construction of $ b^+ (h)  $.} Given $ \chi^+_{\rm
s} \in {\mathcal
  S}_{\rm hyp}^+ (\epsilon^9)  $, we then simply
choose the symbols $ b_0^+, \ldots , b_N^+  $
according to the Proposition \ref{factorisation}, with $ \epsilon \leq \min
(\hat{\epsilon}_{\iota}^{2^N},\epsilon_6  )  $.

\medskip

\noindent {\bf Justification of the parametrix. } Since $
\tilde{\kappa} \otimes \tilde{\kappa}_{\iota} \equiv 1 $ near the
support of $ \chi^{+}_{\rm s} $, we have
$$  \big| \big| \chi^{+}_{\rm s}
(r,\theta,h D_r,h D_{\theta}) - \chi^{+}_{\rm s} (r,\theta,h D_r,h
D_{\theta})(\tilde{\kappa} \otimes \tilde{\kappa}_{\iota}) \big|
\big|_{L^2(\Ra^n)\rightarrow L^2(\Ra^n)}  \lesssim  h^M, \qquad h
\in (0,1] , $$ for all $ M $, using the standard symbolic  calculus,
the Calder\'on-Vaillancourt Theorem and the fact that $ {\mathcal
S}^+_{\rm hyp}(\epsilon) \subset C^{\infty}_b (\Ra^{2n}) $. Using
Proposition \ref{factorisation}, we therefore obtain
$$  \big| \big| C^+
\big| \big|_{L^2(\Ra^n)\rightarrow L^2(\Ra^n)}  \lesssim  h^{N+1},
\qquad h \in (0,1]  . $$ It remains to consider $ D_+ (s) $ which
reads
$$  D_+ (s) = \sum_{ l = 0}^{N+2} \sum_{m=0}^N h^{l+m} J_h^+ (d_l^+) e^{-i s h D_r^2} J_h^+ (b_m^+)^* . $$
By (\ref{borneL2FIO}) and (\ref{nonstationnaire0}), the part of
the sum where $ l \geq N + 1 $, has an $ L^2 $ operator norm of
order $ h^{N+1} $. Once divided by $h$ and integrated over an
interval of size at most $ h^{-1} $, the corresponding operator
norm is $ {\mathcal O}(h^{N-1}) $. The control of the other terms
of the sum will follow from the next result.
\begin{prop} \label{entrantsortant} If $ \epsilon  $ is small enough, then,  for
all $ 0 \leq l , m \leq N  $ and all $ M \geq 0  $, we have
$$
\left| \left| J_h^+ (d_l^+) e^{-i s h D_r^2} J_h^+ (b_m^+)^*
\right| \right|_{L^2 (\Ra^n) \rightarrow L^2 (\Ra^n)} \leq
C_{\epsilon} h^M , \qquad h \in (0,1], \ \  0 \leq s \leq T h^{-1}
. $$
\end{prop}

The proof is based on a fairly elementary non stationary phase
argument. To control the operator norms of the kernels obtained
after integrations by parts, we need the following rough lemma.

\begin{lemm} \label{lemmrough} For $ a \in C^{\infty}_b (\Ra^{3n}) $ compactly supported with respect to $ \rho $, let us set
$$ [a]^{+}_h (r,\theta,r^{\prime},\theta^{\prime}) = (2 \pi h)^{-n} \int \! \! \int e^{ \frac{i}{h}(
S_{+,\epsilon}(r,\theta,\rho,\eta)- s \rho^2 -
S_{+,\epsilon}(r^{\prime},\theta^{\prime},\rho,
\eta))}a(r,\theta,r^{\prime},\theta^{\prime},\rho,\eta) d \rho
d\eta , $$ using $ S_{+,\epsilon} $ defined in Proposition
\ref{Phaseplusmoinsglobale}.
 Denote by $ {\mathcal A}_h^{+} : L^2(\Ra^n)
\rightarrow L^2(\Ra^n) $ the operator with Schwartz kernel $
[a]^{+}_h $. Then, there exists $ n_0 (n) \geq 0 $ such that, for
all $ \epsilon $ small enough,
$$  ||  {\mathcal A}_h^{+} ||_{L^2(\Ra^n)
\rightarrow L^2(\Ra^n)} \leq C_{\epsilon} h^{-n_0} \scal{s}^{n_0}
\max_{|\gamma| \leq n_0} \sup_{\Ra^{3n}} || \partial^{\gamma} a
||_{\infty} , $$ for all $ h \in (0,1] $, all $ s \in \Ra $ and
all $ a \in C^{\infty}_b (\Ra^{3n}) $ satisfying
$$ \emph{supp}(a) \subset \{
|\rho| \leq 10 \} . $$
\end{lemm}

\noindent {\it Proof.} It is a simple consequence of the
Calder\'on-Vaillancourt Theorem by interpreting $ {\mathcal A}_h^{+}
$ as the pseudo-differential operator with symbol
$$ e^{\frac{i}{h}(  \varphi_{+,\epsilon} (r,\theta,\rho,\eta) - s \rho^2 -
\varphi_{+,\epsilon} (r^{\prime},\theta^{\prime},\rho,\eta) ) }
a(r,\theta,r^{\prime},\theta^{\prime},\rho,\eta) , $$ where $
\varphi_{+,\epsilon} $ is defined in Proposition
\ref{Phaseplusmoinsglobale}. \finpreuve

\bigskip

\noindent {\it Proof of Proposition \ref{entrantsortant}.} We
notice first that, by Proposition \ref{Phaseplusmoins} and
(\ref{crucialphasestationnaire})
\begin{eqnarray}
\partial_{\rho} \left( S_{+}(r,\theta,\rho,\eta) - s \rho^2 - S_{+}(r^{\prime},\theta^{\prime},\rho,\eta)
\right) = r - r^{\prime} - 2 s \rho + {\mathcal O}(\epsilon^2),
\label{deriveephaserho} \\
\partial_{\eta} \left( S_{+}(r,\theta,\rho,\eta) - s \rho^2 - S_{+}(r^{\prime},\theta^{\prime},\rho,\eta)
\right) = \theta - \theta^{\prime}  + {\mathcal O}(e^{-1/\epsilon
} ) , \label{deriveephaseeta}
\end{eqnarray}
on the support of $ d_l^+ (r,\theta,\rho,\eta) b_m^+
(r^{\prime},\theta^{\prime},\rho,\eta) $. On the other hand, by
construction, we have
$$ d_l^{+} = i^{-1}\left( v^+ \partial_r \chi_{\epsilon^2
\rightarrow \epsilon} + w^+ \cdot \partial_{\theta}
\chi_{\epsilon^2 \rightarrow \epsilon} \right) \hat{a}^+_{l-1} +
P_{\iota} ( \chi_{\epsilon^2 \rightarrow \epsilon}
\hat{a}^+_{l-2} ) - \chi_{\epsilon^2 \rightarrow \epsilon}
P_{\iota} \hat{a}^+_{l-2} , $$ (with the convention
that $ \hat{a}^+_{-2} = \hat{a}^+_{-1} \equiv 0 $). Using in particular that
$$ w^+ = e^{-r} (\partial_{\eta} q)(r,\theta,e^{-r} \partial_{\theta}S_+) ,
$$
 $ d_l^+ $ is a sum of
terms of the form $ c(r,\theta,\rho,\eta) \partial_r^j
(e^{-r}\partial_{\theta})^{\alpha} \chi^+_{\epsilon^2 \rightarrow
\epsilon} $, with $ j + |\alpha| \geq 1 $ and $ c \in {\mathcal
B}_{\rm s}^+ (\epsilon) $. 
Using the form of $ \chi^+_{\epsilon^2 \rightarrow \epsilon} $
given by Lemma \ref{construitcutoff}, we see that, on the support
of such terms, at least one of the following properties hold
\begin{eqnarray}
 \epsilon^{-1} \leq  r  \leq \epsilon^{-2} , \qquad \qquad \qquad \qquad \label{phasenonstat1} \\
 p (r,\theta,\rho,\eta) \leq 1/4 - \epsilon^2 \ \  \mbox{or} \ \  p (r,\theta,\rho,\eta) \geq 4 +
 \epsilon^2,  \label{phasenonstat2} \\
 \epsilon^{4 - 2 \kappa}   \lesssim  e^{-2r}|\eta|^2  \lesssim
 \epsilon^2,  \label{phasenonstat3} \qquad \qquad \qquad \\
 \mbox{dist}(\theta,V_{\iota})  \geq \epsilon^4 , \qquad \qquad \qquad \qquad  \label{phasenonstat4}
\end{eqnarray}
for some fixed  $ 0 < \kappa  < 1 $ in (\ref{phasenonstat3}).
   For terms such that  (\ref{phasenonstat1})
holds on  their supports, we have
\begin{eqnarray}
 (\ref{deriveephaserho}) \leq \epsilon^{-2} - \epsilon^{-3}
- 2 s \rho + C \leq - 1 - 2  s \rho . \label{minorationphaserho}
\end{eqnarray}
for $ \epsilon $ small enough and integrate by parts w.r.t $ \rho
$. For those satisfying (\ref{phasenonstat2}) on their supports,
then we must have
$$ \rho^2 - 1/4 \leq  - \epsilon^2 \qquad \mbox{or} \qquad \rho^2 - 4 \gtrsim \epsilon^2 , $$
since $ e^{-2r}|\eta|^2 \lesssim \epsilon^2 $ in any case, whereas
on the support of $ b_l^+ $, where $ p
(r^{\prime},\theta^{\prime},\rho,\eta) \in (1/4 - \epsilon^3 , 4 +
\epsilon^3) $ and $ e^{-2 r^{\prime}}|\eta|^2 \lesssim \epsilon^6
$,
$$ \rho^2 - 1/4 \gtrsim - \epsilon^3 \qquad \mbox{and} \qquad \rho^2 - 4 \leq \epsilon^3 , $$
so that the amplitude vanishes identically, again if $ \epsilon $
is small enough. For those satisfying (\ref{phasenonstat3}) on
their supports, we have $ e^{r}|\eta|^{-1} \lesssim
\epsilon^{\kappa - 2} $. Since $  e^{- r^{\prime}}|\eta| \lesssim
\epsilon^3 $, we get
$$ e^{r - r^{\prime}} \leq C + (1+\kappa) \ln \epsilon \ll 0 , $$ which
implies again that $ (\ref{deriveephaserho})  \leq - 1 - 2  s \rho
$, if $ \epsilon $ is small enough. Thus on the supports of  terms
satisfying either (\ref{phasenonstat1}) or (\ref{phasenonstat2})
or (\ref{phasenonstat3}), we have $ |(\ref{deriveephaserho})|
\gtrsim \scal{s} $. By standard integrations by parts, the kernel
of corresponding operator can be written, for all $ M $, as in
Lemma \ref{lemmrough} with amplitudes of order $ (h / \scal{s})^M
$ in $ C^{\infty}_b (\Ra^{3n}) $. Hence, their $ L^2 $ operator
norms are of order $ (h / \scal{s})^{M-n_0 } $ with an arbitrary $
M $.

For the remaining terms satisfying (\ref{phasenonstat4}) on their
supports, we remark that $ |\theta^{\prime} - \theta| \geq
\epsilon^5 $ (otherwise $ \mbox{dist}(\theta,V_{\iota}) \leq |\theta -
\theta^{\prime} | + \mbox{dist}(\theta^{\prime},V_{\iota}) < \epsilon^5
+ \epsilon^6 \ll \epsilon^4 $) hence
$$ |(\ref{deriveephaseeta}) | \gtrsim \epsilon^5 . $$
Thus, for all $ M \geq 0 $, the kernel of the corresponding
operators can be written as in Lemma \ref{lemmrough} with
amplitudes of order $ h ^M $ in $ C^{\infty}_b (\Ra^{3n}) $. Since
$ M $ is arbitrary,  their $ L^2 $ operator norms are of order $
h^M $ if $ |s | \lesssim h^{-1} $.  \finpreuve

\bigskip

This completes the proof of Theorem \ref{IsozakiKitadaansatz}.
\finpreuve

\section{Geometric optics and Egorov Theorem on AH manifolds} \label{WKBsection}
\setcounter{equation}{0}

\subsection{Finite time WKB approximation} 

In this subsection, we give a short time
parametrix of $ e^{-ith P} \widehat{O \!
p}_{\iota}(\chi^{\pm}) $ when $ \chi^{\pm} $ is supported in an outgoing ($ + $) or an
incoming ($ - $) area. This parametrix is the standard geometric
optic (or WKB) approximation which is basically well known. Nevertheless,
in the literature, one mostly find local versions (ie with $ \chi
\in C_0^{\infty} $) or versions in $ \Ra^n $ for elliptic
operators. Here we are neither in a relatively compact set nor in
the uniformly elliptic setting  so we recall the construction with some details. 

Analogously to Section \ref{IsozakiKitada}, we prove here an $ L^2 $ approximation. The related dispersion estimates leading to
(\ref{weightWKB})
 will be derived in Section \ref{sectiondispersion}.

We also emphasize that, although we shall prove this approximation with a specified
time orientation (ie $ t \geq 0 $ for $ \chi^+ $ and $ t \leq 0 $
for $ \chi^- $), this result has nothing to do with
outgoing/incoming areas; in principle we should be able to state a
similar result for any $ \chi $ supported in $ p^{-1} (I) $ and
for times $ | t |\ll 1 $. We restrict the sense of time for only two reasons: firstly,
because it is sufficient for our purpose and, secondly,  because we
can use directly Proposition \ref{estimeesflotprecises} (we should
otherwise give a similar result for the geodesic flow for $ t $ in
an open neighborhood of $ 0 $).

In this subsection $ \iota$ is fixed (ie we work in a fixed coordinate patch).  We drop the subscript $ \iota $ from the notation of symbols, phases, etc... For instance, we replace $ p_{\iota} $ by $ p $ everywhere (see the beginning of Section  \ref{sectionHJ}). 

 Fix then  
$$ I_1 \Subset I_2 \Subset I_3 \Subset (0,+\infty), $$
three relatively compact open subsets of $ V^{\prime}_{\iota} $ (see (\ref{margeouverts})),
$$ V_1 \Subset V_2 \Subset V_3 \Subset V_{\iota}^{\prime} , $$
and three real numbers
$$ -1 < \sigma_1 < \sigma_2 < \sigma_3 < 1 .$$
For some $ R_3 $ large enough to be fixed below, we also choose arbitrary $  R_1 , R_2$ real numbers such that
$$ R_1 > R_2 > R_3 . $$
\begin{theo} \label{theoremWKB} For all $ R_3 $ large enough, there exists $ t_{\rm WKB} > 0 $ and a function
$$ \Sigma \in C^{\infty}\left( [0,\pm t_{\rm WKB}] \times \Ra^{2n}, \Ra \right) , $$
such that, for any
\begin{eqnarray}
  \chi^{\pm}
\in {\mathcal S}_{\rm hyp} \left(
\Gamma^{\pm}_{\iota}(R_1,V_1,I_1,\sigma_1)  \right) ,
\label{categorieini}
\end{eqnarray}
we can find
$$ a_0^{\pm}(t), \ldots , a_N^{\pm}(t)  \in {\mathcal S}_{\rm hyp} \left(
\Gamma^{\pm}_{\iota}(R_2,V_2,I_2,\sigma_2)  \right) , $$
depending smoothly on $t$ for $ 0 \leq \pm t \leq t_{\rm WKB} $, and such that, if we set 
$$ a^{\pm}_N (t,h) = a_0^{\pm}(t) + \cdots + h^N a_N^{\pm}(t) , $$
the operator defined on $ 
C_0^{\infty}(\Ra^n) $ by the following kernel
$$ \left[ {\mathcal J}^{\pm}_h \left( t, a^{\pm}_N (t,h) \right) \right] (t,r,\theta,r^{\prime},\theta^{\prime}) =
(2 \pi h)^{-n}  \int \! \! \int e^{\frac{i}{h} \left(
\Sigma^{\pm} (t,r,\theta,\rho,\eta) - r^{\prime}\rho - \theta^{\prime} \cdot \eta \right)} a_N^{\pm} (t,h,r,\theta,\rho,\eta)
 d \rho d\eta  ,  $$
satisfies, with $ {\bf 1}_{\iota} $  the
characteristic function of $ (R_3 , +\infty) \times V_3 $,
\begin{eqnarray}
 \left| \left|  \ e^{-it h P } \widehat{O \! p}_{\iota} \left( \chi^{\pm}
 \right)
 - \Psi_{{\iota}}^*  {\mathcal J}^{\pm}_h \left( t, a^{\pm}_N(t,h)  \right) {\bf 1}_{\iota} \left( \Psi_{{\iota}}^{-1} \right)^*
  \right| \right|_{L^2 ({\mathcal
M},\widehat{dG}) \rightarrow L^2 ({\mathcal M},\widehat{dG})}
\leq C h^{N+1} , \label{vraieparametrixeWKB}
\end{eqnarray}
for
$$ 0 \leq  \pm t \leq t_{\rm WKB} , \qquad h \in (0,1 ] . $$
In addition, the functions $ \Sigma^{\pm} $ are of the form 
$$ \Sigma^{\pm}(t,r,\theta,\rho,\eta) = r \rho + \theta \cdot \eta  + 
\left( \Sigma^{\pm}_0 (t,r,\theta,\rho,\eta) - r\rho - \eta \cdot \eta \right) \chi_{2 \rightarrow 3}^{\pm}(r,\theta,\rho,\eta) , $$
with $ \chi_{2 \rightarrow 3}^{\pm} \in {\mathcal S}_{\rm hyp} (\Gamma^{\pm}_{\iota}(R_3,V_3,I_3,\sigma_3)) $ such that $ \chi^{\pm}_{2 \rightarrow 3} \equiv 1 $ on $ \Gamma^{\pm}_{\iota}(R_2,V_2,I_2,\sigma_2) $, and some bounded family
$  (\Sigma_0^{\pm}(t))_{0 \leq \pm t \leq t_{\rm WKB}} $ in $ {\mathcal B}_{\rm hyp} (\Gamma^{\pm}_{\iota}(R_3,V_3,I_3,\sigma_3)) $ satisfying
\begin{eqnarray} \begin{cases}
\partial_t \Sigma^{\pm}_0 + p (r,\theta,\partial_r \Sigma^{\pm}_0,\partial_{\theta}
\Sigma^{\pm}_0)= 0 , \\ \Sigma^{\pm}_0(0,r,\theta,\rho,\eta) = r
\rho + \theta \cdot \eta \end{cases} , \qquad \qquad \label{eikonalequation} 
\end{eqnarray}
and
\begin{eqnarray}
 | D_{\rm hyp}^{j \alpha k \beta}
\left( \Sigma^{\pm}_0(t,r,\theta,\rho,\eta) - r\rho - \theta \cdot
\eta - t p (r,\theta,\rho,\eta) \right) |
  \leq  C_{j \alpha k \beta} t^2 , \label{developpementphaseWKB}
\end{eqnarray}
both for 
$$ 0 \leq \pm t \leq t_{\rm WKB} \ \ \ \mbox{and} \ \ \ (r,\theta,\rho,\eta) \in \Gamma^{\pm}_{\iota} (R_3,V_3,I_3,\sigma_3) . $$ We also have
\begin{eqnarray}
 \left( \Sigma^{\pm}(t,r,\theta,\rho,\eta) - r\rho - \theta \cdot
\eta  \right)_{0 \leq \pm t \leq t_{\rm WKB}} \ \ \mbox{bounded in } \ \ {\mathcal S}_{\rm hyp} \left(
\Gamma^{\pm}_{\iota}(R_3,V_3,I_3,\sigma_3)  \right) . \label{globalisationphaseWKB}
\end{eqnarray}
Finally, for all $ 0 \leq j \leq N $,
\begin{eqnarray}
\left( a_j^{\pm}(t) \right)_{0 \leq \pm t \leq t_{\rm WKB}} \ \ \mbox{is bounded in} \ \ {\mathcal S}_{\rm hyp} \left(
\Gamma^{\pm}_{\iota}(R_2,V_2,I_2,\sigma_2)  \right). \label{supporttransport}
\end{eqnarray}
\end{theo}

Notice that $ V_1 \Subset V_{\iota}^{\prime} $ so it makes sense to consider $ \widehat{ O \! p }_{\iota} (\chi^{\pm}) $ (see (\ref{symbolpseudo})).

In principle it is not necessary to have  $ R_3 $ large to get such a lemma, but this will be sufficient for our applications.
The interest of choosing $ R_3 $ large is simply to allow to use directly Proposition \ref{estimeesflotprecises}.

Note also that, by (\ref{supporttransport}), 
the kernel of $  {\mathcal J}^{\pm}_h \left(t, a^{\pm}_N(t,h)  \right)
{\bf 1}_{\iota}  $ is supported in $ \left(
(R_3,+\infty) \times V_3 \right)^2 $.

\medskip

 The rest of the subsection is devoted to the proof of 
Theorem \ref{theoremWKB}. We  need to find $ \Sigma_{\pm} $ and $ a^{\pm}_N
(t,h) $ such that
\begin{eqnarray}
 {\mathcal J}^{\pm}_h \left(0, a^{\pm}_N (0,h) \right) & = & \chi^{\pm} (r,\theta,h D_r , h
 D_{\theta}), \label{initialWKB} \\
 \left( h D_t + h^2 P_{\iota} \right)  {\mathcal J}^{\pm}_h \left(t, a^{\pm}_N (t,h)
 \right) & = & h^{N+2} R_N^{\pm}(t,h), \label{equationWKB}
\end{eqnarray}
where $ P_{\iota} = (\Psi_{\iota}^{-1})^* P
\Psi_{\iota}^* $ and
\begin{eqnarray}
|| R_N^{\pm}(t,h) ||_{L^2 (\Ra^n) \rightarrow L^2 (\Ra^n)} \leq C,
\qquad h \in (0,1], \ 0 \leq \pm t \leq t_{\rm WKB} .
\label{resteWKB}
\end{eqnarray}
Indeed, if (\ref{initialWKB}), (\ref{equationWKB}) and
(\ref{resteWKB}) hold  then
\begin{eqnarray*}
 \Psi_{{\iota}}^*  {\mathcal J}^{\pm}_h \left( t, a^{\pm}_N (t,h)
\right) {\bf 1}_{\iota} \left( \Psi_{{\iota}}^{-1} \right)^* -
e^{-ith P} \Psi_{{\iota}}^* \chi^{\pm} (r,\theta,h
D_r, h D_{\theta}) {\bf 1}_{\iota} \left( \Psi_{{\iota}}^{-1}
\right)^* & = & \\  i h^{N+1} \int_0^t e^{-i(t-s)h \widetilde{P}}
\Psi_{{\iota}}^* R_N (s,h) {\bf 1}_{\iota} \left(
\Psi_{{\iota}}^{-1} \right)^* ds
\end{eqnarray*}
will yield (\ref{vraieparametrixeWKB}) since, for all $ M >0 $,
$$ || \Psi_{{\iota}}^* \chi^{\pm} (r,\theta,h D_r, h
D_{\theta}) {\bf 1}_{\iota} \left( \Psi_{{\iota}}^{-1} \right)^* -
\widehat{O \! p}_{\iota}(\chi^{\pm}) ||_{ L^2 ({\mathcal M},\widehat{dG}) \rightarrow L^2 ({\mathcal M},\widehat{dG}) } \leq C_M h^M ,
$$
by standard off diagonal decay (see Definition \ref{definitOpchapeau} for $ \widehat{O \! p}_{\iota} $), since 
$ {\bf 1}_{\iota} \equiv 1 $ near $ \Pi_{r,\theta} \left( \mbox{supp}(\chi^{\pm}) \right) $.

To get the conditions to be satisfied by $ \Sigma^{\pm} $
and $ a_0^{\pm}, \ldots , a^{\pm}_N $ we observe that
\begin{eqnarray}
   \left( h D_t + h^2 P_{\iota} \right)  {\mathcal J}^{\pm}_h \left(t, a^{\pm}_N (t,h)
 \right) = \sum_{j = 0}^{N+2} h^j  {\mathcal J}^{\pm}_h \left( t,b_j^{\pm}(t) 
 \right) , \label{compositionWKB}
\end{eqnarray}
 where, if we additionally set $ a^{\pm}_{-2} = a^{\pm}_{-1} = a^{\pm}_{N+1} = a^{\pm}_{N+2} \equiv 0
 $,
\begin{eqnarray}
b_j = (\partial_t \Sigma^{\pm} + p (r,\theta,\partial_r
\Sigma^{\pm}, \partial_{\theta} \Sigma^{\pm})) a_j^{\pm} + i^{-1}
(\partial_t + {\mathcal T}^{\pm}) a_{j-1}^{\pm} + P a_{j-2}^{\pm}, \label{calculcompositionWKB}
\end{eqnarray}
with
\begin{eqnarray}
 {\mathcal T}^{\pm} = 2 \partial_r \Sigma^{\pm} \partial_r +
(\partial_{\eta} q ) (r,\theta, e^{-r} \partial_{\theta}
\Sigma^{\pm} )\cdot e^{-r}
\partial_{\theta} + (p+ p_1)
(r,\theta,\partial_r,\partial_{\theta})\Sigma^{\pm} ,
\label{motivvfield}
\end{eqnarray}
where $q$ is defined in (\ref{symboleprincipalglobal}) and $ p_1 $ is the homogeneous part of degree $1$ of the full
symbol of $ P_{\iota} $. To obtain (\ref{initialWKB}),
(\ref{equationWKB}) and (\ref{resteWKB}) it will therefore be
sufficient to solve the eikonal equation (\ref{eikonalequation}),
then the transport equations
\begin{eqnarray}
(\partial_t + {\mathcal T}^{\pm} ) a^{\pm}_0 = 0 , \qquad
a^{\pm}_0 (0,.) = \chi^{\pm}(.), \label{transportWKB0} \\
(\partial_t + {\mathcal T}^{\pm} ) a^{\pm}_k = -i
P_{\iota} a^{\pm}_{k-1} , \qquad a^{\pm}_k (0,.) = 0 ,
\label{transportWKBk}
\end{eqnarray}
for $ 1 \leq k \leq N $, and finally to get an $ L^2 $ bound for
Fourier integral operators of the form $ {\mathcal J}^{\pm}_h \left( t,a
 \right) $ (using the Kuranishi trick).

\medskip

To solve (\ref{eikonalequation}), we need the following lemma for which we recall that $ (r^t,\theta^t,\rho^t,\eta^t) $ is the Hamiltonian flow of $p$.

\begin{lemm} \label{injectivityradius} For all $ - 1 < \sigma_{\rm eik} < \sigma_{\rm eik}^{\prime} <  1 $, 
all open intervals $ I_{\rm eik} \Subset I_{\rm eik}^{\prime} \Subset (0,+\infty) $, all open subsets $ V_{\rm eik} \Subset V_{\rm eik}^{\prime} \Subset V_{\iota}^{\prime} $ and all $
R_{\rm eik}
> R_{\rm eik}^{\prime}  $  large enough, there exists  $ t_1 > 0 $ small enough such that
$$ {\Psi}^t_{\pm} : (r,\theta,\rho,\eta) \mapsto (r^t,\theta^t,\rho,\eta) $$
is a diffeomorphism from $ \Gamma^{\pm} (R_{\rm eik}^{\prime},V_{\rm eik}^{\prime},I_{\rm eik}^{\prime},\sigma_{\rm eik}^{\prime}) $
onto its range for all $ 0 \leq \pm t < t_1 $ and
$$ \Gamma^{\pm} (R_{\rm eik},V_{\rm eik},I_{\rm eik},\sigma_{\rm eik}) \subset 
{\Psi}^t_{\pm} \left( \Gamma^{\pm} (R_{\rm eik}^{\prime},V_{\rm eik}^{\prime},I_{\rm eik}^{\prime},\sigma_{\rm eik}^{\prime}) \right)
\qquad \mbox{for all} \ 0 \leq \pm t < t_1 . $$
\end{lemm}

\noindent {\it Proof.} Let us choose first $ \sigma_{\rm eik}^{\prime \prime}  \in \Ra$, $
I_{\rm eik}^{\prime \prime} $ open interval and $ V_{\rm eik}^{\prime \prime} $ open set  such that
$$ \sigma_{\rm eik}^{\prime} <
\sigma_{\rm eik}^{\prime \prime} < 1 , \qquad I_{\rm eik}^{\prime} \Subset I_{\rm eik}^{\prime \prime} \Subset (0,+\infty) \ \ \
\mbox{and} \ \ \ V_{\rm eik}^{\prime} \Subset V_{\rm eik}^{\prime \prime} \Subset V_{\iota}^{\prime} . $$ We also choose $
R_{\rm eik}^{\prime \prime}  > 0 $ large enough such that Proposition \ref{estimeesflotprecises} holds with $
\sigma = |\sigma_{\rm eik}^{\prime \prime}|  $ and $ R = R_{\rm eik}^{\prime \prime}  $. We then choose arbitrary $ R_{\rm eik} $
and $ R_{\rm eik}^{\prime }  $ such that
$$ R_{\rm eik} > R_{\rm eik}^{\prime }  > R_{\rm eik}^{\prime \prime}  , $$
and then $ \chi_{ \prime \rightarrow  \prime \prime }^{\pm} \in
{\mathcal S}_{\rm hyp} \left(  \Gamma^{\pm} (R_{\rm eik}^{\prime \prime} ,V_{\rm eik}^{\prime \prime} ,
I_{\rm eik}^{\prime \prime} ,\sigma_{\rm eik}^{\prime \prime}  ) \right) $ such $
\chi_{\prime \rightarrow  \prime \prime}^{\pm} \equiv 1 $ on $  \Gamma^{\pm}
(R_{\rm eik}^{\prime } ,V_{\rm eik}^{\prime } ,I_{\rm eik}^{\prime } ,\sigma_{\rm eik}^{\prime} ) $. The existence of such a function
follows from Proposition \ref{partitions} {\it i)}. In particular, $ \chi_{\prime \rightarrow \prime \prime}^{\pm} $ and $
\partial_{r,\theta,\rho,\eta} \chi_{\prime \rightarrow \prime \prime}^{\pm} $ are
bounded on $ \Ra^{2n} $. For $ \pm t \geq 0 $, consider  the map
\begin{eqnarray}
 {\varepsilon}^t_{\pm} : \Ra^{2n} \ni (r,\theta,\rho,\eta)
\mapsto \left( \int_0^t 2 \rho^s ds, \int_0^t e^{-r^s}
(\partial_{\eta}q)(r^s,\theta^s,e^{-r^s}\eta^s) ds \right) \chi_{\prime
\rightarrow \prime \prime}^{\pm}(r,\theta,\rho,\eta) \in \Ra^n,
\end{eqnarray}
so that, by the motion equations,  
$$ {\Psi}^t_{\pm} =
\mbox{Id}_{\Ra^{2n}} + \left({\varepsilon}^t_{\pm} , 0 \right)
\ \ \ \mbox{  on } \ \   \Gamma^{\pm} (R_{\rm eik}^{\prime },V_{\rm eik}^{\prime },I_{\rm eik}^{\prime },\sigma_{\rm eik}^{\prime }) . $$ 
By Proposition \ref{estimeesflotprecises} we have
 $ |
\partial_{r,\theta,\rho,\eta} {\varepsilon}^t_{\pm} | \lesssim |t|
$, hence $ \mbox{Id}_{\Ra^{2n}} +
\left({\varepsilon}^t_{\pm} , 0 \right) $ is a diffeomorphism
from $ \Ra^{2n} $ onto itself,  for all $ \pm t \geq 0 $ small
enough. Therefore, it remains to show that, if $ t $ is small enough and $
(r,\theta,\rho,\eta) \in
\Gamma^{\pm} (R_{\rm eik},V_{\rm eik},I_{\rm eik},\sigma_{\rm eik}) $ is of the form
$$ (r,\theta,\rho,\eta) = (r^{\prime},\theta^{\prime},\rho^{\prime},\eta^{\prime}) +
 ({\varepsilon}^t_{\pm} (r^{\prime},\theta^{\prime},\rho^{\prime},\eta^{\prime}),0) , $$
then $ (r^{\prime},\theta^{\prime},\rho^{\prime},\eta^{\prime}) \in \Gamma^{\pm}
(R_{\rm eik}^{\prime},V_{\rm eik}^{\prime},I_{\rm eik}^{\prime},\sigma_{\rm eik}^{\prime}) $. We have trivially $ \rho = \rho^{\prime} $
and $ \eta = \eta^{\prime} $. By Proposition \ref{estimeesflotprecises}, 
$| {\varepsilon}^t_{\pm} | \lesssim |t|$ on $ \Ra^{2n} $, so $
|r-r^{\prime}|+|\theta-\theta^{\prime}| \lesssim |t| $, hence $ r^{\prime}
> R_{\rm eik}^{\prime} $ and $ \theta^{\prime} \in V_{\rm eik}^{\prime} $ if $ t $ is small enough.
Furthermore, by writing
$$ q (r,\theta,e^{-r}\eta) - q (r^{\prime},\theta^{\prime},e^{-r^{\prime}}\eta) =
 q (r,\theta,e^{-r}\eta) - q (r^{\prime},\theta^{\prime},e^{-r}\eta) + (1-e^{-2(r^{\prime}-r)})
  q (r^{\prime},\theta^{\prime},e^{-r}\eta) , $$
 we see that $$ |p(r,\theta,\rho,\eta)-p (r^{\prime},\theta^{\prime},\rho,\eta)| \lesssim |t|
 $$ using the boundedness of $ |e^{-r^{\prime}}\eta | $ and the Taylor formula. Hence
 $$ p(r^{\prime},\theta^{\prime},\rho,\eta) \in I_{\rm eik}^{\prime} \qquad \mbox{and} \qquad \pm \rho > - \sigma_{\rm eik}^{\prime} p(r^{\prime},\theta^{\prime},\rho,\eta)^{1/2}  $$
 if $ t$ is small enough, since  $ p(r,\theta,\rho,\eta)
  \in I_{\rm eik} $ and $ \pm \rho > - \sigma_{\rm eik} p(r,\theta,
 \rho,\eta)^{1/2} $. This completes the proof. \finpreuve

\bigskip


Let us now fix $ I_{\rm eik} , I_{\rm eik}^{\prime} $, $ V_{\rm eik} , V_{\rm eik}^{\prime} $, and $ \sigma_{\rm eik} ,
\sigma_{\rm eik}^{\prime} $ as in Lemma \ref{injectivityradius} with the additional conditions 
$$ V_{\rm eik} = V_3 , \qquad I_{\rm eik} = I_3 , \qquad \sigma_{\rm eik} = \sigma_3 . $$
We  denote by $ \Psi^{\pm}_t $ the inverse of $ \Psi_{\pm}^t
$ and define $ ( r_t , \theta_t ) = (r_t ,
\theta_t)(r,\theta,\rho,\eta) $ by 
$$ \Psi^{\pm}_t (r,\theta,\rho,\eta) = (r_t,\theta_t,\rho,\eta) \in \Gamma^{\pm}
(R_{\rm eik}^{\prime},V_{\rm eik}^{\prime},I_{\rm eik}^{\prime},\sigma_{\rm eik}^{\prime}) , $$
if $ (r,\theta,\rho,\eta) \in \Gamma^{\pm} (R_{\rm eik},V_{\rm eik},I_{\rm eik},\sigma_{\rm eik})
$ and $ 0 \leq \pm t < t_1 $. Here $ t_1 $, $ R_{\rm eik} $ and $ R_{\rm eik}^{\prime} $ are those given by Lemma \ref{injectivityradius}.

\begin{prop} For all $ R_{3} > R_{\rm eik} $, there exists $ t_{\rm eik} > 0 $
such that
$$ \Sigma^{\pm}_0 ( t,r,\theta,\rho,\eta) := r_t \rho + \theta_t \cdot \eta  + t p (r_t,\theta_t,\rho,\eta) ,  $$
solves (\ref{eikonalequation}) on $ \Gamma^{\pm}(R_3 ,V_3,I_3,\sigma_3) $ for $ 0 \leq \pm t \leq
t_{\rm eik} $, 
and such that
\begin{eqnarray}
 (\Sigma^{\pm}_0 (t,r,\theta,\rho,\eta) - r \rho -
\theta \cdot \eta )_{0 \leq \pm t \leq t_{\rm eik} }  \ \ \mbox{
is bounded in }  \ \ {\mathcal B}_{\rm hyp} \left(
\Gamma^{\pm}(R_3,V_3,I_3,\sigma_3)
\right) . \label{borneeikonale}
\end{eqnarray}
\end{prop}

\noindent {\it Proof.} That $ \Sigma^{\pm}_0 $ solves the eikonal
equation is standard so we only have to show
(\ref{borneeikonale}). Since
$$ \Sigma^{\pm}_0 (t,r,\theta,\rho,\eta) = r \rho +
\theta \cdot \eta + (r_t -r) \rho +  e^{r} (\theta_t -\theta) \cdot e^{-r}\eta + t e^{-2(r_t-r)} q (r_t,\theta_t,e^{-r}\eta), $$
(\ref{borneeikonale}) would follow from the following
estimates
\begin{eqnarray}
| D_{\rm hyp}^{j \alpha k \beta} ( r_t - r )| + | D_{\rm hyp}^{j
\alpha k \beta} \left( e^{r} (\theta_t - \theta) \right) | \leq
C_{j \alpha k \beta} , \label{pourborneeikonale}
\end{eqnarray}
for $ 0 \leq \pm t \leq \pm t_{\rm eik} $ and $
(r,\theta,\rho,\eta) \in \Gamma^{\pm}(R_3,V_3,I_3,\sigma_3) $. 
The motion equations yield
\begin{eqnarray}
 r^t = r + \int_0^{t} 2 \rho^s ds , \qquad \theta^t = \theta + \int_0^t
e^{-r^s} (\partial_{\eta}q)(r^s,\theta^s,e^{-r^s}\eta^s) ds ,
\label{equamouv}
\end{eqnarray}
so, by Proposition \ref{estimeesflotprecises} with $ R^{\prime}_{\rm eik} $ of Lemma \ref{injectivityradius} and by choosing $
t_{\rm eik } $ small enough, we see that, for $ 0 \leq \pm t \leq
t_{\rm eik} $, we have
$$ | \partial_{r,\theta} (r^t , \theta^t) - \mbox{Id}_n | \leq 1/2, \qquad  \mbox{on} \ \ \Gamma^{\pm}(R_{\rm eik}^{\prime},V_{\rm eik}^{\prime},I_{\rm eik}^{\prime},\sigma_{\rm eik}^{\prime}) ,  $$
where $ |.| $ is a matrix norm. Therefore, by differentiating the
identity $ ( r^t , \theta^t )(r_t,\theta_t,\rho,\eta) = (r,\theta)
$ one obtains, similarly  to Proposition \ref{estimepreparephase},
\begin{eqnarray}
  | D_{\rm hyp}^{j \alpha k \beta} ( r_t - r )| + | D_{\rm
hyp}^{j \alpha k \beta}  (\theta_t - \theta)  | \leq C_{j \alpha k
\beta} , \label{intermediairebis}
\end{eqnarray}
for $ 0 \leq \pm t \leq t_{\rm eik} $ and $ (r,\theta,\rho,\eta)
\in \Gamma^{\pm}(R_3,V_3,I_3,\sigma_3) $. This proves the expected estimates for $ r_t - r $. The second equation
of (\ref{equamouv}) evaluated at $ (r_t,\theta_t,\rho,\eta) $
yields
\begin{eqnarray}
 e^r \left( \theta - \theta_t \right) = \int_0^t e^{ r - r^s_t }
(\partial_{\eta} q) (r^s_t,\theta^s_t, e^{-r^s_t}\eta^s_t) ds ,
\label{mouvexpo}
\end{eqnarray}
where $ x_t^s = x^s (r_t,\theta_t,\rho,\eta) $ for $ x=
r,\theta,\eta $. Combining (\ref{intermediairebis}) and
Proposition \ref{estimeesflotprecises}, we have, on $
\Gamma^{\pm}(R_3,V_3,I_3,\sigma_3) $,
$$  | D_{\rm hyp}^{j \alpha k \beta} ( r_t^s - r )| + | D_{\rm
hyp}^{j \alpha k \beta}  (\theta_t^s - \theta)  | + | D_{\rm
hyp}^{j \alpha k \beta}  (\eta_t^s - \eta)  |  \leq C_{j \alpha k
\beta} , \qquad 0 \leq \pm t , \pm s \leq t_{\rm eik} ,  $$ from
which the estimate of the second term of (\ref{pourborneeikonale})
follows  using (\ref{mouvexpo}). \finpreuve

\bigskip

We now solve the transport equations. By (\ref{motivvfield}), we
have to
 consider the time dependent vector field $ (v_t^{\pm} ,
w_t^{\pm}) $ defined on $ \Gamma^{\pm} (R_3,V_3,I_3,\sigma_3) $,
for $ 0 \leq \pm  t \leq t_{\rm eik} $, by
\begin{eqnarray}
  \begin{pmatrix}
 v^{\pm}_t \\
 w^{\pm}_t
\end{pmatrix} :=  \begin{pmatrix}
  (\partial_{\rho} p)(r,\theta,\partial_r \Sigma^{\pm}_0 , \partial_{\theta} \Sigma^{\pm}_0) \\
(\partial_{\eta} p)(r,\theta,\partial_r \Sigma^{\pm}_0 ,
\partial_{\theta} \Sigma^{\pm}_0)
\end{pmatrix} =  \begin{pmatrix}
  2 \partial_r \Sigma^{\pm}_0  \\
e^{-2r}(\partial_{\eta} q)(r,\theta , \partial_{\theta}
\Sigma^{\pm}_0)
\end{pmatrix}  .  \label{champsdevecteurstemps}
\end{eqnarray}
We then denote by $ \phi^{\pm}_{s \rightarrow t} $ the flow, from
time $s$ to time $t$, of $ (v^{\pm}_t,w^{\pm}_t,0_{\Ra^n}) $
namely the solution to
\begin{eqnarray}
\partial_t \phi^{\pm}_{s \rightarrow t} = ( v^{\pm}_t (\phi^{\pm}_{s \rightarrow t}) ,
w^{\pm}_t(\phi^{\pm}_{s \rightarrow t}) , 0 ), \qquad
\phi^{\pm}_{s \rightarrow s} (r,\theta,\rho,\eta) =
(r,\theta,\rho,\eta). \label{flotaderiverenprincipe}
\end{eqnarray}

\begin{lemm} \label{lemmeflottransportWKB} For all $ I_{\rm tr} $ open interval, $ \sigma_{\rm tr} \in \Ra
$ and $ V_{\rm tr} \subset \Ra^{n-1} $ open subset such that
$$ R_{\rm tr} > R_3 , \qquad V_{\rm tr} \Subset V_3 , \qquad I_{\rm tr} \Subset I_3 ,  \qquad - 1 < \sigma_{\rm tr} < \sigma_3 ,  $$
 there exists $ 0 < t_2 \leq t_{\rm
eik} $ small enough such that:
\begin{eqnarray}
\phi^{\pm}_{s \rightarrow t} \ \ \mbox{is well defined on} \ \
\Gamma^{\pm} (R_{\rm tr},V_{\rm tr},I_{\rm tr},\sigma_{\rm tr}), \ \ \mbox{for all} \ \ 0 \leq
\pm s \leq t_2 , \ 0 \leq \pm t \leq t_2 , \label{suppositionimplicite} 
\end{eqnarray}
and
\begin{eqnarray}
\left| D^{\rm hyp}_{j \alpha k \beta} \left( \phi^{\pm}_{s
\rightarrow t} - \emph{Id} \right) \right| \lesssim 1 , \qquad
\mbox{on} \ \ \Gamma^{\pm} (R_{\rm tr},V_{\rm tr},I_{\rm tr},\sigma_{\rm tr}), \ \mbox{for} \
\  0 \leq \pm s, \pm t \leq t_2 . \label{restatable}
\end{eqnarray}
\end{lemm}
By (\ref{suppositionimplicite}), we  mean in particular that
\begin{eqnarray}
 \phi_{s
\rightarrow t} \left( \Gamma^{\pm} (R_{\rm tr},V_{\rm tr},I_{\rm tr},\sigma_{\rm tr}) \right)
\subset \Gamma^{\pm} (R_3,V_3,I_3,\sigma_3), \qquad 0 \leq \pm s,
\pm t \leq t_2 .
\end{eqnarray}
 The estimate (\ref{restatable}) can be restated by saying  that
the components of $ \phi^{\pm}_{s \rightarrow t} - \mbox{Id} $ are
bounded families of $ {\mathcal B}_{\rm hyp} \left(
\Gamma^{\pm}(R_{\rm tr},V_{\rm tr},I_{\rm tr},\sigma_{\rm tr}) \right) $ for $ 0 \leq \pm s,
\pm t \leq t_2 $.

\bigskip

\noindent {\it Proof.}  For all $
\delta > 0 $ small enough, we have
\begin{eqnarray}
 |r-r^{\prime}|+|\theta-\theta^{\prime}| \leq \delta  \ \mbox{and}
 \ (r,\theta,\rho,\eta) \in \Gamma^{\pm} (R_{\rm tr},V_{\rm tr},I_{\rm tr},\sigma_{\rm tr})
\Rightarrow (r^{\prime},\theta^{\prime},\rho,\eta) \in
\Gamma^{\pm} (R_3,V_3,I_3,\sigma_3)  \label{stabilitedomaine}
\end{eqnarray}
by Proposition \ref{geometriebornee}. Denoting by $ (r^{\pm}_{s
\rightarrow t }, \theta^{\pm}_{s \rightarrow t},\rho,\eta) $ the
components of $ \phi^{\pm}_{s \rightarrow t} $, they must be solutions
 of the problem
$$ (r^{\pm}_{s \rightarrow
t }, \theta^{\pm}_{s \rightarrow t}) = (r,\theta) + \int_s^t
(v_{\tau}^{\pm},w_{\tau}^{\pm})(r^{\pm}_{s \rightarrow \tau },
\theta^{\pm}_{s \rightarrow \tau},\rho,\eta) d \tau . $$ By
(\ref{borneeikonale}), we have
\begin{eqnarray}
 |(v_{\tau}^{\pm},w_{\tau}^{\pm})| + | \partial_{r,\theta} ( v_{\tau}^{\pm} ,
w_{\tau}^{\pm} ) | \leq C  , \label{clefPicard}
\end{eqnarray}
 on $ \Gamma^{\pm} (R_3,V_3,I_3,\sigma_3)
$, for $ 0 \leq \pm \tau \leq t_{\rm eik} $. Therefore, the sequence
$ u_n^{\pm}(t)= u_n^{\pm}(t,s,r,\theta,\rho,\eta) $ defined by
\begin{eqnarray}
u_0^{\pm}(s) =  (r,\theta) , \qquad u_{k+1}^{\pm}(t) = (r,\theta)
+ \int_s^t
(v_{\tau}^{\pm},w^{\pm}_{\tau})(u_k^{\pm}(\tau),\rho,\eta) d\tau
\nonumber ,
\end{eqnarray}
is a Cauchy sequence in $ C^0 ([0,\pm t_2], \Ra^n) $  for all $
(r,\theta,\rho,\eta) \in \Gamma^{\pm}(R_{\rm tr},V_{\rm tr},I_{\rm tr},\sigma_{\rm tr}) $ and
$ 0 \leq \pm s \leq t_2$, for some $t_2$ small enough independent
of $ (r,\theta,\rho,\eta) $. Indeed, using
(\ref{stabilitedomaine}) and choosing $t_2$ small enough so that $
\sum_{k \geq 0 } (C t_2)^{k+1} \leq \delta $,  a standard
induction using (\ref{clefPicard}) shows that
$$  |u_{k+1}^{\pm}(t)- u_k^{\pm} (t)| \leq  (C t_2)^{k+1} ,  $$
which makes the sequence well defined and convergent. This proves
(\ref{suppositionimplicite}). We then obtain (\ref{restatable}) by
induction by differentiating the equations in
(\ref{flotaderiverenprincipe}). This proof is completely similar
to the one of the estimate (\ref{rthetauniftemps}) in Proposition
\ref{transport} (and much simpler for it is local in time) so we
omit the details. \finpreuve

\bigskip
 
In the next proposition, we denote by $ q_t^{\pm} = q_t^{\pm}(r,\theta,\rho,\eta) $ the
function defined on $ [0,\pm t_{\rm eik}] \times
\Gamma^{\pm}(R,V,I,\sigma) $ by
$$ q_t^{\pm} := (p+ p_1)
(r,\theta,\partial_r,\partial_{\theta})\Sigma^{\pm}_0 , $$ which is
involved in (\ref{motivvfield}).

\begin{prop} Choose $ R_{\rm tr},V_{\rm tr}, I_{\rm tr} $ and $ \sigma_{\rm tr} $ such that
$$ R_2 > R_{\rm tr} > R_3 , \qquad V_2 \Subset V_{\rm tr} \Subset V_3 , \qquad I_2 \Subset I_{\rm tr} \Subset I_3 ,  \qquad \sigma_2 < \sigma_{\rm tr} < \sigma_3 $$
Then, there exists $ t_{\rm tr} > 0 $ small enough such that, for
all $ \chi^{\pm} $ satisfying (\ref{categorieini}), the  functions 
$$ a_0^{\pm}, \ldots , a_{N}^{\pm} : [0,\pm t_{\rm tr}] \times
\Ra^{2n} \rightarrow \Ca , $$ 
defined iteratively by
$$ a_0^{\pm}(t) := \chi^{\pm} \circ \phi^{\pm}_{t \rightarrow 0}
\exp \left( \int_0^t q^{\pm}_s \circ \phi^{\pm}_{t \rightarrow s}
\right) $$ and
$$ a_k^{\pm}(t) := - \int_0^t  i (P_{\iota} a^{\pm}_{k-1}) (s_1, \phi^{\pm}_{t \rightarrow
s_1}) \exp \left( \int_{s_1}^t q^{\pm}_{s_2} \circ \phi^{\pm}_{t \rightarrow s_2} d
s_2 \right) d s_1 , \qquad 1 \leq k \leq N , $$ on $ \Gamma^{\pm}
(R_2,V_2,I_2,\sigma_2) $ and by $ 0
$ outside are smooth and solve (\ref{transportWKB0}) and
(\ref{transportWKBk}). Furthermore, for all $ 0 \leq k \leq N $,
\begin{eqnarray}
 (a_k^{\pm}(t))_{0 \leq \pm t \leq t_{\rm tr}} \ \ \mbox{is bounded in}  \ \ {\mathcal S}_{\rm hyp} \left(
\Gamma^{\pm}(R_2,V_2,I_2,\sigma_2) \right) . \label{bonneclasseWKB}
\end{eqnarray}
\end{prop}

\noindent {\it Proof.} Let us fix $R^{\prime}_{\rm tr} , V^{\prime}_{\rm tr} , I^{\prime}_{\rm tr} $ and
 $ \sigma^{\prime}_{\rm tr} $
such that
$$R_2 > R^{\prime}_{\rm tr} > R_{\rm tr}, \qquad V_2 \Subset V^{\prime}_{\rm tr} \Subset V_{\rm tr} , \qquad I_2 \Subset I^{\prime}_{\rm tr}
 \Subset I_{\rm tr} , \qquad \sigma_2 < \sigma^{\prime}_{\rm tr} < \sigma_{\rm tr} . $$
By choosing $ 0 < t_{\rm tr} \leq t_2 $ small enough, we then have, for
all $ 0 \leq \pm s , \pm t \leq t_{\rm tr} $,
\begin{eqnarray}
\phi^{\pm}_{s \rightarrow t} \left( \Gamma^{\pm}(R_1,V_1,I_1,\sigma_1)
\right) & \subset & \Gamma^{\pm}(R_2,V_2,I_2,\sigma_2), \label{conditiontransport0}
\\
 \phi^{\pm}_{s \rightarrow t} \left(  \Gamma^{\pm}(R_2,V_2,I_2,\sigma_2) \right) & \subset &
\Gamma^{\pm}(R^{\prime}_{\rm tr},V^{\prime}_{\rm tr},I^{\prime}_{\rm tr},\sigma^{\prime}_{\rm tr}),
\label{conditiontransport1}
\\
\phi^{\pm}_{s \rightarrow t} \left(
\Gamma^{\pm}(R^{\prime}_{\rm tr},V^{\prime}_{\rm tr},I^{\prime}_{\rm tr},\sigma^{\prime}_{\rm tr}) \right) & \subset &
\Gamma^{\pm}(R_{\rm tr},V_{\rm tr},I_{\rm tr},\sigma_{\rm tr}) . \label{conditiontransport2}
\end{eqnarray}
This follows from Proposition \ref{geometriebornee} and the fact that $ | \phi^{\pm}_{t \rightarrow s} - \mbox{Id} | \lesssim |t-s| $, which comes  from the integration of (\ref{flotaderiverenprincipe}) between $s$ and $t$, using (\ref{restatable}).
By Lemma \ref{lemmeflottransportWKB}, the flow is well defined on $
\Gamma^{\pm}(R_{\rm tr},V_{\rm tr},I_{\rm tr},\sigma_{\rm tr}) $, therefore the condition
(\ref{conditiontransport2}) ensures that we have the pseudo-group
property
\begin{eqnarray}
 \phi^{\pm}_{t \rightarrow u} \circ \phi^{\pm}_{s \rightarrow t}
 = \phi^{\pm}_{s \rightarrow u} ,
\qquad 0 \leq \pm s , \pm t , \pm u \leq t_{\rm tr} ,
\label{pseudogroup2}
\end{eqnarray}
 on $ \Gamma^{\pm}(R^{\prime}_{\rm tr},V^{\prime}_{\rm tr},I^{\prime}_{\rm tr},\sigma^{\prime}_{\rm tr})
 $. In particular, $
\phi^{\pm}_{t \rightarrow s} \circ \phi^{\pm}_{s \rightarrow t} =
\mbox{Id} $ on this set. Therefore,
 by
(\ref{conditiontransport1}), we have
$$   \Gamma^{\pm}(R_2,V_2,I_2,\sigma_2)   \subset \phi^{\pm}_{t
\rightarrow s} \left( \Gamma^{\pm}(R^{\prime}_{\rm tr},V^{\prime}_{\rm tr},I^{\prime}_{\rm tr},\sigma^{\prime}_{\rm tr}) 
\right) . $$ This implies that the map
$$ (t,r,\theta,\rho,\eta) \mapsto (t,\phi_{s \rightarrow t
}^{\pm}(r,\theta,\rho,\eta)) $$ is a diffeomorphism from $ (0,\pm
t_{\rm tr}) \times \Gamma^{\pm}(R^{\prime}_{\rm tr},V^{\prime}_{\rm tr},I^{\prime}_{\rm tr},\sigma^{\prime}_{\rm tr})  $
onto its range and that this range contains $ (0,\pm t_{\rm tr}) \times
\Gamma^{\pm}(R_2,V_2,I_2,\sigma_2)
$. Restricted to the latter set, the inverse is given by $ (t,
\phi^{\pm}_{t \rightarrow s}) $ which shows that $ \phi^{\pm}_{t
\rightarrow s} $ is smooth with respect to $ t$. Furthermore, by
differentiating in $t$ the relation $ \phi^{\pm}_{t \rightarrow s}
\circ \phi^{\pm}_{s \rightarrow t} = \mbox{Id}  $, one
obtains
\begin{eqnarray}
\partial_t \phi^{\pm}_{t \rightarrow s} + (\partial_{r,\theta} \phi^{\pm}_{t \rightarrow
s}) \cdot (v^{\pm}_t,w^{\pm}_t) = 0 , \qquad \mbox{on} \ \  \Gamma^{\pm}(R_2,V_2,I_2,\sigma_2), \nonumber
\end{eqnarray}
for $ 0 < \pm t < t_{\rm tr} $. Using this relation, one easily checks that $ a_0^{\pm} $ solves (\ref{transportWKB0}) on $ \Gamma^{\pm}(R_2,V_2,I_2,\sigma_2) $. In
addition, if $$ (r,\theta,\rho,\eta) \in
\Gamma^{\pm}(R^{\prime}_{\rm tr},V^{\prime}_{\rm tr},I^{\prime}_{\rm tr},\sigma^{\prime}_{\rm tr}) 
\setminus \Gamma^{\pm}(R_2,
V_2, I_2, \sigma_2) , $$ we have $
\phi^{\pm}_{t \rightarrow 0} (r,\theta,\rho,\eta) \notin \mbox{supp}(\chi^{\pm}) $
otherwise $ (r,\theta,\rho,\eta) \in  \Gamma^{\pm}(R_2, V_2, I_2, \sigma_2) $ by (\ref{categorieini}),
(\ref{conditiontransport0}) and (\ref{pseudogroup2}). This shows
that, extended by $ 0 $ outside $ \Gamma^{\pm}(R_2, V_2,
I_2, \sigma_2) $, $ a_0^{\pm} $ is smooth. The
property (\ref{bonneclasseWKB}) for $k=0$ is then a direct consequence of (\ref{restatable}). We
note in passing that we have
$$ \mbox{supp}(a_0^{\pm}(t)) \subset \phi^{\pm}_{0 \rightarrow t} (\mbox{supp}(\chi^{\pm})) . $$
The proof for the higher order terms $ a_k^{\pm} $, $ k \geq 1 $,
is then obtained similarly by induction using  that $ \mbox{supp}
(P_{\iota} a^{\pm}_{k-1} (s_1)) \subset
\phi^{\pm}_{0 \rightarrow s_1} (\mbox{supp}(\chi^{\pm})) $ for all
$ s_1 $. \finpreuve

\bigskip

\noindent {\bf Proof of Theorem \ref{theoremWKB}.} It remains to prove (\ref{developpementphaseWKB}), to globalize $ \Sigma^{\pm}_0 $, to prove (\ref{globalisationphaseWKB}) and the bound (\ref{resteWKB}).
By Proposition \ref{partitions}, we can choose 
$$ \chi^{\pm}_{2 \rightarrow 3} \in {\mathcal S}_{\rm hyp}(\Gamma^{\pm}(R_3,V_3,I_3,\sigma_3)) \ \ \mbox{such that} \ \
 \chi^{\pm}_{2 \rightarrow 3} \equiv 1 \ \ \mbox{on} \ \ \Gamma^{\pm}(R_2,V_2,I_2,\sigma_2) . $$
 We set
$$ \Sigma^{\pm}(t,r,\theta,\rho,\eta) = r \rho + \theta \cdot \eta + \chi^{\pm}_{2 \rightarrow 3}(r,\theta,\rho,\eta)
\times \left( \Sigma^{\pm}_0(t,r,\theta,\rho,\eta) - r \rho - \theta
\cdot \eta \right).
$$
It coincides with $ \Sigma^{\pm}_0 $ on $ [0,\pm t_{\rm eik}] \times \Gamma^{\pm} (R_2,V_2,I_2,\sigma_2) $ so it is a solution to the eikonal equation on $ [0,\pm t_{\rm WKB}] \times \Gamma^{\pm} (R_2,V_2,I_2,\sigma_2) $, for any $ 0 < t_{\rm WKB} \leq t_{\rm eik} $. Furthermore, (\ref{borneeikonale}) implies (\ref{globalisationphaseWKB}) and, by using
\begin{eqnarray}
\Sigma^{\pm}_0
(t,r,\theta,\rho,\eta) = r \rho + \theta \cdot \eta + \int_0^t p
(r,\theta,\partial_r \Sigma^{\pm}_0(s) , \partial_{\theta} \Sigma^{\pm}_0(s)) ds , \label{pourTaylorWKB}
\end{eqnarray}
 we get  (\ref{developpementphaseWKB}) since (\ref{borneeikonale}) and (\ref{pourTaylorWKB}) itself show that
 the components of $ (\partial_r \Sigma^{\pm}(s)-\rho,\partial_{\theta} \Sigma^{\pm} (s) -\eta ) $ are $ {\mathcal O}(s) $ in 
 $ {\mathcal B}_{\rm hyp} (\Gamma^{\pm}(R_3,V_3,I_3,\sigma_3)) $.

To prove (\ref{resteWKB}), we use the Kuranishi trick which is as follows. By the Taylor formula, we can write
$$  \Sigma^{\pm}(t,r,\theta,\rho,\eta) - \Sigma^{\pm}(t,r^{\prime},\theta^{\prime},\rho,\eta)
 = (r-r^{\prime}) \tilde{\rho}^{\pm}(t,r,\theta,r^{\prime},\theta^{\prime},\rho,\eta) + (\theta-\theta^{\prime})
 \cdot \tilde{\eta}^{\pm}(t,r,\theta,r^{\prime},\theta^{\prime},\rho,\eta) .
$$
Using again (\ref{pourTaylorWKB}) 
and (\ref{borneeikonale}), we obtain
\begin{eqnarray}
| \partial_r^j \partial_{\theta}^{\alpha}
\partial_{r^{\prime}}^{j^{\prime}}\partial_{\theta^{\prime}}^{\alpha^{\prime}}
\partial_{\rho}^k \partial_{\eta}^{\beta} \left( (\tilde{\rho}^{\pm}, \tilde{\eta}^{\pm})(t,r,\theta,r^{\prime},\theta^{\prime},\rho,\eta)
- (\rho,\eta) \right) | \leq C_{j\alpha j^{\prime} \alpha^{\prime}
k \beta } |t| ,
\end{eqnarray}
for $ ( r,\theta,r^{\prime},\theta^{\prime},\rho,\eta \in \Ra^{3n}
) $ and $ 0 \leq \pm t \leq t_{\rm eik} $. The latter implies
that, for  all $0 \leq  \pm t \leq t_{\rm WKB} $ small enough and all $ (r,\theta,r^{\prime},\theta^{\prime}) \in \Ra^{2n}
$, the map
$$ (\rho,\eta) \mapsto
(\tilde{\rho}^{\pm} , \tilde{\eta}^{\pm} ) , $$ is a
diffeomorphism from $ \Ra^n $ onto itself. Furthermore, proceeding
similarly to the proof of (\ref{controlesurlesderiveesKuranishi})
in Lemma \ref{Kuranishi}, we see that its inverse $
(\tilde{\rho},\tilde{\eta}) \mapsto (\rho^{\pm} , \eta^{\pm}) $
satisfies
\begin{eqnarray}
| \partial_r^j \partial_{\theta}^{\alpha}
\partial_{r^{\prime}}^{j^{\prime}}\partial_{\theta^{\prime}}^{\alpha^{\prime}}
\partial_{\tilde{\rho}}^k \partial_{\tilde{\eta}}^{\beta} \left( (\rho^{\pm}, \eta^{\pm})
(t,r,\theta,r^{\prime},\theta^{\prime},\tilde{\rho},\tilde{\eta})
- (\tilde{\rho},\tilde{\eta}) \right) | \leq C_{j\alpha j^{\prime}
\alpha^{\prime} k \beta }  ,
\end{eqnarray}
on $ \Ra^{3n} $, uniformly with respect to $ 0 \leq \pm t \leq t_{\rm WKB} $. 
Thus, for any bounded family $ (a^{\pm}(t))_{0 \leq \pm t \leq t_{\rm WKB}} $ in $ {\mathcal S}_{\rm hyp}(\Gamma^{\pm}(R_2,V_2,I_2,\sigma_2)) $, the kernel of $  {\mathcal J}^{\pm}_h \left(t, a^{\pm}(t) \right)   {\mathcal J}^{\pm}_h \left( t, a^{\pm} (t)
\right)^* $, which reads
\begin{eqnarray}
(2 \pi h)^{-n} \int e^{\frac{i}{h}
(\Sigma^{\pm}(t,r,\theta,\rho,\eta) -
\Sigma^{\pm}(t,r^{\prime},\theta^{\prime},\rho,\eta))}
a^{\pm} (t,r,\theta,\rho,\eta) \overline{ a^{\pm}
(t,r^{\prime},\theta^{\prime},\rho,\eta) } d \rho d \eta , 
\end{eqnarray}
can be written
 \begin{eqnarray} (2 \pi
h)^{-n} \int e^{\frac{i}{h} ( (r - r^{\prime} ) \tilde{ \rho } +
(\theta - \theta^{\prime}) \cdot \tilde{\eta}))} B
(t,r,\theta,r^{\prime},\theta^{\prime}, \tilde{\rho},\tilde{\eta})
d \tilde{\rho} d \tilde{\eta}, \label{pseudoWKB}
\end{eqnarray}
with $ B (t,.) $ bounded in $ C^{\infty}_b (\Ra^{3n})  $ as $ 0
\leq \pm t \leq t_{\rm WKB} $. By the Calder\'on-Vaillancourt
Theorem we obtain the uniform boundedness of the operator given by
(\ref{pseudoWKB}) hence
$$ || {\mathcal J}^{\pm}_h \left( t,a^{\pm}
(t) \right)   ||_{L^2 (\Ra^n) \rightarrow L^2 (\Ra^n)} \leq C,
\qquad 0 \leq \pm t \leq t_{\rm WKB}, \ h \in (0,1 ] , $$ where $
C $ depends only a finite number of semi-norms of $a^{\pm}(t)$ in $
C^{\infty}_b (\Ra^{2n}) $.  Using (\ref{compositionWKB}), (\ref{calculcompositionWKB}) (with $ a_{k}^{\pm}(t) $ solutions to the transport equations)  and (\ref{bonneclasseWKB}), the bound above yields (\ref{resteWKB}) which
completes the proof of Theorem \ref{theoremWKB}. \finpreuve

\subsection{Proof of Proposition \ref{propositionfinale} }


To prove Proposition \ref{propositionfinale}, we first need a version of the semi-classical Egorov Theorem in the asymptotically hyperbolic setting.

 As in the previous subsection, $ \iota $ is fixed and we drop the subscript $ \iota $ from most notations. In particular $ p = p_{\iota} $ and $ \Phi^t =(r^t,\theta^t,\rho^t,\eta^t) $ denotes the Hamiltonian flow of $p$
studied in subsection \ref{wpasquelconque}.

Fix $ V \Subset V_{\iota}^{\prime} $ an open subset, $ I \Subset (0,+\infty) $ an open interval and $ -1 < \sigma < 1 $.
\begin{theo} \label{EgorovAH}  If $ R > 0 $ is large enough  the following holds: for all $ T > 0 $, all $ N \geq 0 $
and all
\begin{eqnarray}
  a \in {\mathcal S}_{\rm hyp} \left( \Gamma^{\pm}_{\iota} (R,V,I,\sigma) \right) , \label{conditionaEgorov}
\end{eqnarray}  
we can find 
\begin{eqnarray}
 a_{0}(t), \ldots , a_N (t) \in {\mathcal S}_{\rm hyp} \left( \Phi^t \left( \emph{supp}(a) \right) \right) 
, \qquad  0 \leq \pm t \leq T ,
 \label{supportEgorov}
\end{eqnarray}
such that, 
\begin{eqnarray}
 \left| \left|  e^{-ith P} \widehat{O \! p}_{\iota}(a) e^{ithP} - \sum_{k = 0}^N h^k  \widehat{O \! p}_{\iota}(a_k(t))
   \right| \right|_{L^2 ({\mathcal M},\widehat{dG}) \rightarrow L^2 ({\mathcal M},\widehat{dG}) } 
   \leq C_{N,T,a} h^{N+1} , \label{sensEgorov}
\end{eqnarray}
for all $ 0 \leq \pm t \leq T $ and all $ 0 < h \leq 1 $.   
\end{theo}

This theorem is basically well known. Here the main point is to check 
(\ref{supportEgorov}), namely that $ a_0 (t) , \ldots , a_N (t) \in {\mathcal B}_{\rm hyp} \left( \Phi^t \left( \mbox{supp}(a) \right) \right) $.
Notice that, by Corollary \ref{dependencedudomaine}, $ \Phi^t \left( \mbox{supp}(a) \right) $ is contained in the same chart
as $ a $ in which it is therefore sufficient to work.


 Using the group property, it is sufficient to prove the result when $ T $ is small enough 
(depending only on $ V, I , \sigma $). To check this point, we choose open sets $ V_1,V_2 $ such that $ V \Subset V_1
\Subset V_2 \Subset V_{\iota}^{\prime} $. Then for some $ C > 0 $ and all $ R $ large enough 
$$ \Phi^t \left( \Gamma^{\pm}_{\iota} (R,V,I,\sigma) \right) \subset \Gamma^{\pm}_{\iota} (R-C,V_1,I,\sigma) , \qquad \pm t \geq 0 , $$
and
$$ \Phi^t \left( \Gamma^{\pm}_{\iota} (R-C,V_1,I,\sigma) \right) \subset \Gamma^{\pm}_{\iota} (R-2C,V_2,I,\sigma) ,
\qquad \pm t \geq 0 . $$
This follows from Corollary \ref{dependencedudomaine} and the fact that $ \rho^t $ can be assumed to be non decreasing, using (\ref{pourzoneintermediaire}). Thus, it is sufficient to prove (\ref{sensEgorov})  for $ 0 \leq \pm t \leq \varepsilon $ with  $\varepsilon > 0 $ small enough independent of
 $ a \in  {\mathcal S}_{\rm hyp} \left( \Gamma^{\pm}_{\iota} (R-C,V_1,I,\sigma) \right)  $. Indeed if this holds, it holds for
 $ a $ satisfying (\ref{conditionaEgorov}) and 
$$ e^{i \varepsilon h P} \widehat{O \! p}_{\iota}(a) e^{-i\varepsilon hP} - \sum_{k = 0}^N h^k  \widehat{O \! p}_{\iota}(a_k(\varepsilon)) 
+ h^{N+1} R_N(h,\varepsilon)$$
with $ R_N (h,\varepsilon) $ uniformly bounded on $ L^2 ({\mathcal M},\widehat{dG}) $ and  $ a_k (\varepsilon) \in {\mathcal S}_{\rm hyp} (
\Gamma^{\pm}_{\iota} (R-C,V_1,I,\sigma) ) $, with $ a_k(\varepsilon) $ supported in $  \Phi^{\varepsilon} \left( \mbox{supp}(a) \right) $ more precisely. Conjugating the expression above by $ e^{- i \varepsilon h P} $ and then applying the same result with $ a_k(\varepsilon) $ instead of $a$ we can write
$$ e^{i 2 \varepsilon h P} \widehat{O \! p}_{\iota}(a) e^{- 2 i\varepsilon hP} - \sum_{k = 0}^N h^k  \widehat{O \! p}_{\iota}(a_k( 2 \varepsilon)) 
+ h^{N+1} R_N(h, 2 \varepsilon) , $$
where $  a_k(2\varepsilon) $ is supported in $ \Phi^{2 \varepsilon} \left( \mbox{supp}(a) \right) $ which is still contained in $  
\Gamma^{\pm}_{\iota} (R-C,V_1,I,\sigma)  $ and thus allows to iterate the procedure. 

The interest of considering small times is justified by the following lemma.

\begin{lemm} \label{controleflotinverseEgorov} Fix $ V_1, I , \sigma $ as above, then for some $ R_1 > 0 $ large enough and $  \varepsilon > 0 $ small enough,
$$ \left| D_{\rm hyp}^{j \alpha k \beta} \left( (\Phi^{t})^{-1} - \emph{Id}_{2n} \right) \right| \leq C_{j\alpha k \beta} \ \ \ \mbox{on} \
\ \ \Phi^{t} \left( \Gamma^{\pm}_{\iota}(R_1,V_1,I,\sigma) \right) , $$
for all $ 0 \leq \pm t \leq \varepsilon $.
\end{lemm}

\noindent {\it Proof.} Using the identity 
$$ d \left( \Phi^t - \mbox{Id}_{2n} \right) = \int_0^t d H_p (\Phi^s) d \Phi^s ds $$ and
Proposition \ref{estimeesflotprecises}, we have $ | d \left( \Phi^t - \mbox{Id}_{2n} \right) | \lesssim |t| $ hence 
$ | (d \Phi^t)^{-1}| \lesssim 1 $ on $\Gamma^{\pm}_{\iota} (R_1,V_1,I,\sigma)$ if $ R_1 $ is large enough and $t$ is small enough. We then obtain the result
by applying $ D_{\rm hyp}^{j \alpha k \beta } $ to $ \Phi^t \circ (\Phi^t)^{-1} $ and using the Faa Di Bruno formula. 
For instance, if $ j = k = | \alpha | =0 $ and
$ | \beta | = 1 $, we have
$$ d \Phi^t_{|(\Phi^t)^{-1}} e^{r} \partial_{\eta}^{\beta} \left( (\Phi^t)^{-1} - \mbox{Id}_{2n} \right) = ( \mbox{Id}_{2n}
 - d \Phi^t_{|(\Phi^t)^{-1}} )
e^r \partial_{\eta}^{\beta}  \mbox{Id}_{2n}  
 $$
where, using Proposition \ref{estimeesflotprecises}, the right hand side is bounded for this is simply
$ e^{r} \partial_{\eta}^{\beta} \left( \mbox{Id}_{2n} - \Phi^t \right) $ evaluated at $ (\Phi^t)^{-1}  $.
Higher order derivatives are studied similarly by iteration using Lemma \ref{Faadibruno}.
\finpreuve

\bigskip

Naturally, $ (\Phi^{t})^{-1} $ is the reverse Hamiltonian flow, namely flowing
$ \Phi^t \left( \Gamma^{\pm}_{\iota} (R_1,V_1,I,\sigma) \right) $ back to $  \Gamma^{\pm}_{\iota} (R_1,V_1,I,\sigma)  $. 
More precisely, for $ 0 \leq \pm t \leq \varepsilon $,
\begin{eqnarray}
 \frac{d}{dt}(\Phi^{t})^{-1}(r,\theta,\rho,\eta) = - H_p \left( (\Phi^{t})^{-1}(r,\theta,\rho,\eta) \right), \qquad 
(r,\theta,\rho,\eta) \in \Phi^{\pm \varepsilon} \left( \Gamma^{\pm}_{\iota} (R_1,V_1,I,\sigma) \right) . \label{equationinverseEgorov}
\end{eqnarray}
We prefer to keep the notation $ (\Phi^t)^{-1} $ on $ \Phi^t (\Gamma^{\pm}_{\iota}(R_1,V_1,I,\sigma)) $ 
rather than using $ \Phi^{-t} $ since we have only studied $ \Phi^t $  
for $ t \geq 0 $ on outgoing areas and $ t \leq 0 $ on incoming areas.

We have essentially all the necessary tools to solve the transport equations considered in the following lemma.
\begin{lemm} \label{LemmetransportEgorov} There exists $ C > 0 $ such that, for all $ R $ large enough, the following holds: for all $ a_{\rm ini} \in {\mathcal S}_{\rm hyp} \left( \Gamma^{\pm}_{\iota} (R,V,I,\sigma) \right) $ and
 $$ (f (t))_{0 \leq \pm t \leq \varepsilon} \ \ \mbox{a bounded family of} \ \ 
 {\mathcal S}_{\rm hyp}  \left( \Gamma^{\pm}_{\iota}(R -C,V_1,I,\sigma ) \right) , $$
 smooth with respect to $t$ and such that
 $$ \emph{supp}(f(t)) \subset \Phi^t (\emph{supp}(a_{\rm ini})) , $$
 the function defined   for $ 0 \leq \pm t \leq \varepsilon $ by
$$ a (t) := \begin{cases} a_{\rm ini} \circ (\Phi^{t})^{-1} + \int_0^{t} f (s)  \circ \Phi^s \circ(\Phi^{t})^{-1} ds & 
\mbox{on} \ \Phi^t \left( \emph{supp}(a) \right), \\
0 & \mbox{outside} , \end{cases} $$
is smooth and satisfies
\begin{eqnarray}
 \partial_t a(t)  + \{ p , a(t) \} = f(t), \qquad a (0) = a_{\rm ini} .  \label{transportEgorov}
\end{eqnarray}
Furthermore
\begin{eqnarray}
 ( a (t))_{ 0 \leq \pm t \leq \varepsilon} \ \ \mbox{is bounded in } 
{\mathcal S}_{\rm hyp}  \left( \Gamma^{\pm}_{\iota}(R -C,V_1,I,\sigma ) \right) . \label{bornetransportEgorov}
\end{eqnarray}
\end{lemm}

In (\ref{bornetransportEgorov}), we consider $ \Gamma^{\pm}_{\iota} (R -C,V_1,I,\sigma ) $ for it is independent of $t$ but, by construction, 
$ a (t) $ is supported in the smaller region $ \Phi^t (\mbox{supp}(a)) $.

\bigskip

\noindent {\it Proof.} To check the smoothness of $ a_0 (t) $ it is sufficient to see that $ a_{\rm ini} \circ (\Phi^{t})^{-1} $ and
$ f (s)  \circ (\Phi^{t-s})^{-1} $
are defined and smooth in a neighborhood of $ \Phi^t \left( \mbox{supp}(a) \right)  $ where they vanish on the complement of $ \Phi^t \left( \mbox{supp}(a) \right) $. Indeed $ (\Phi^t)^{-1} $ is defined on $ \Phi^t \left( \Gamma^{\pm}_{\iota} (R-C,V_1,I,\sigma) \right) $ and  if 
 $ (r,\theta,\rho,\eta) $ belongs to $ \Phi^t \left( \Gamma^{\pm}_{\iota} (R-C,V_1,I,\sigma) \right) $ but doesn't belong to $ \Phi^t \left( \mbox{supp}(a) \right) $, then
  $  a_{\rm circ} \circ (\Phi^t)^{-1}(r,\theta,\rho,\eta) = 0$ otherwise, $ (\Phi^t)^{-1}(r,\theta,\rho,\eta) $ should belong to $ \mbox{supp}(a) $
  and thus $(r,\theta,\rho,\eta) $ should belong to $\Phi^t ( \mbox{supp}(a) ) $. Similarly $ \int_0^{t} f (s)  \circ \Phi^s \circ(\Phi^{t})^{-1}
  (r,\theta,\rho,\eta) ds $ must vanish otherwise there would be an $s$ between $0$ and $t$ such that $ \Phi^s \circ (\Phi^t)^{-1}(r,\theta,\rho,\eta) \in \Phi^s (\mbox{supp}(a)) $
    implying that $ (r,\theta,\rho,\eta) \in \Phi^t (\mbox{supp}(a)) $.
  Then (\ref{transportEgorov})
  follows directly from (\ref{equationinverseEgorov}) and (\ref{bornetransportEgorov}) follows from Lemma \ref{controleflotinverseEgorov}.
\finpreuve

\bigskip

\noindent {\bf Proof of Theorem \ref{EgorovAH}.} By Lemma \ref{LemmetransportEgorov}, 
the solutions of the transport equations (\ref{transportEgorov}) belongs to ${\mathcal S}_{\rm hyp}  \left( \Gamma^{\pm}_{\iota}(R -C,V_1,I,\sigma ) \right) $. 
The proof is then standard  (see for instance \cite{Robebook}). \finpreuve

\bigskip

\noindent {\bf Proof of Proposition \ref{propositionfinale}.} We start by choosing $ \epsilon > 0 $ and $ \delta > 0  $ according to Proposition 
\ref{dynamiquedeltasansiota}
with 
$ \underline{t} = t_{\rm WKB} $. Then, using (\ref{borneL2gratuite}), (\ref{compositionpseudodifferentiel}), (\ref{adjointpseudodifferentiel}) and 
Theorem \ref{EgorovAH}, it is straightforward that, for all $ T  \geq t_{\rm WKB} $ and all $ N \geq 0 $, 
$$ || \widehat{O \! p}_{\iota} (b_{l,{\rm inter}}^{\pm}) e^{-ithP}
\widehat{O \! p}_{\iota} (b_{l,{\rm inter}}^{\pm})^* ||_{L^2
(\widehat{dG}) \rightarrow L^2 (\widehat{dG})} \leq C_{T,l,N} h^N , \qquad h \in (0,1 ] , \ t_{\rm WKB} \leq \pm t \leq T  . $$ 
It is therefore sufficient to show the existence of $ T  $ large enough such that, 
\begin{eqnarray}
 || \widehat{O \! p}_{\iota} (b_{l,{\rm inter}}^{\pm}) e^{-ithP}
\widehat{O \! p}_{\iota} (b_{l,{\rm inter}}^{\pm})^* ||_{L^2
(\widehat{dG}) \rightarrow L^2 (\widehat{dG})} \leq C_{l,N} h^N , \qquad h \in (0,1], \  T \leq \pm t \leq 2 h^{-1} . \label{tempspuissance}
\end{eqnarray}
For simplicity we consider positive times and set $ B = \widehat{O \! p}_{\iota} (b_{l,{\rm inter}}^{+}) $. For $ T $ to be chosen, we write
$$ e^{-it h P} B^* = e^{-i(t-T)hP} B(T)^* e^{-iThP}, \qquad B (T) = e^{-iThP} B e^{i T h P} . $$
As above, we may write
$$  B (T)^* = \sum_{k \leq N} h^k \widehat{O \! p}_{\iota} (b^*_{k}(T)) + h^{N+1} B_N (h) $$
with $ B_N (h) $ uniformly bounded on $ L^2 ({\mathcal M},\widehat{dG}) $ and
$$ b^*_{k}(T) \in {\mathcal S}_{\rm hyp} \left( \Phi^T \left( \mbox{supp}(b_{l,{\rm inter}}^+) \right) \right) \subset {\mathcal S}_{\rm hyp} \left( \Phi^T \left(  \Gamma^{+}_{\iota,{\rm inter}} (\epsilon,\delta;l) \right) \right) .  $$
By Proposition \ref{sortantversfortementsortant}, for all $ \tilde{\epsilon} > 0 $, we can choose $ T_{\tilde{\epsilon}}  $ large enough such that
$$ \Phi^T \left(  \Gamma^{+}_{\iota,{\rm inter}} (\epsilon,\delta;l) \right) \subset \Gamma^{+}_{\iota,s}(\tilde{\epsilon}^9) . $$
Thus, if $ \tilde{\epsilon} $ is small enough, Theorem \ref{IsozakiKitadaansatz} allows to write, for $ t \geq T_{\tilde{\epsilon}} $,
$$ e^{-i(t-T_{\tilde{\epsilon}})hP} \widehat{O \! p}_{\iota} (b^*_{k}(T_{\tilde{\epsilon}})) = \Psi_{\iota}^* \left( J_h^+ (\tilde{a}_{\tilde{\epsilon}}(h)) e^{-i(t-T_{\tilde{\epsilon}})hD_r^2} J_h^+ (\tilde{b}_{\tilde{\epsilon}}(h))^* \right) (\Psi_{\iota}^{-1})^* + h^N R_N (t,h) , $$
with $ R_N (t,h) $ uniformly bounded on $ L^2 ({\mathcal M},\widehat{dG}) $ for $ h \in (0,1] $ and $ 0 \leq t - T_{\tilde{\epsilon}} \leq 2 h^{-1} $, and
$$ \tilde{a}_{\tilde{\epsilon}}(h) \in {\mathcal S}_{\rm hyp}  \left(  \Gamma^{+}_{\iota,s}(\tilde{\epsilon})  \right) . $$
We will therefore get (\ref{tempspuissance}) with $ T = T_{\tilde{\epsilon}} $ if we choose $ \tilde{\epsilon} $ small enough such that, for all $N$,
$$ || b_{l,{\rm inter}}^{+} (r,\theta,hD_r,hD_{\theta}) J_h^+ (\tilde{a}_{\tilde{\epsilon}}(h)) ||_{L^2(\Ra^n) \rightarrow L^2 (\Ra^n)} \leq C_N h^N . $$
By the standard composition rule between  pseudo-differential and Fourier integral operators (see \cite{Robebook}), $ b_{l,{\rm inter}}^{+} (r,\theta,hD_r,hD_{\theta}) J_h^+ (\tilde{a}_{\tilde{\epsilon}}(h)) $ is the sum of an operator with norm  of order $ h^N $ and of Fourier integral operators with amplitudes vanishing outside the support of
$$ b_{l,{\rm inter}}^{+} (r,\theta,\partial_r S_+,\partial_{\theta}S_+) \tilde{a}_{\tilde{\epsilon}}(r,\theta,\rho,\eta,h) , $$
where $ S_+ = S_+ (r,\theta,\rho,\eta) $ is the phase defined in Proposition \ref{Phaseplusmoins}. 
It is therefore sufficient to show that, for $ \tilde{\epsilon} $ small enough, the support of the amplitude above is empty.
Indeed, on this support we have
\begin{eqnarray}
 \frac{\partial_r S_+ }{p(r,\theta,\partial_r S_+,\partial_{\theta}S_+)^{1/2}} \leq 1 - (\epsilon/2)^2 \qquad \frac{\rho}{p(r,\theta,\rho,\eta)^{1/2}} > 1 - \tilde{\epsilon}^2 . \label{deuxiemecondition} 
 \end{eqnarray}
Furthermore, by Proposition  \ref{Phaseplusmoins}, we also have
$$ |\partial_r S_+ - r | + |\partial_{\theta} S_+ - \eta| \lesssim \tilde{\epsilon}^2 ,$$
on $ \Gamma^{+}_{\iota,s}(\tilde{\epsilon}) $ where $ \tilde{a}_{\tilde{\epsilon}} (h) $ is supported. Since $p$ is bounded from above and from below on $ \Gamma^{+}_{\iota,s}(\tilde{\epsilon}) $, we obtain
that, for all $ \tilde{\epsilon} $ small enough,
$$ \frac{\rho}{p(r,\theta,\rho,\eta)^{1/2}} \leq 1 - (\epsilon/2)^2 + C \tilde{\epsilon}^2 \leq 1 - (\epsilon/4)^2 ,$$
which is clearly incompatible with the second condition of (\ref{deuxiemecondition}).  \finpreuve

\section{Dispersion estimates} \label{sectiondispersion}
\setcounter{equation}{0}
In this section, we prove Propositions \ref{sousIK} and \ref{sousWKB}, using respectively the parametrices given in Theorems \ref{IsozakiKitadaansatz} and \ref{theoremWKB}. For simplicity, we drop the index $ \iota $ from the notation.

\subsection{Stationary and non stationary phase estimates}
The dispersion estimates will basically follow from the Stationary Phase Theorem in the WKB or the Isozaki-Kitada parametrices. In both cases,
we have to consider oscillatory integrals of the form
\begin{eqnarray}
(2 \pi h )^{-n} \int \! \! \int e^{\frac{i}{h} \Phi^{\pm} (t,r,\theta,r^{\prime},\theta^{\prime},\rho,\eta)} A^{\pm}(t,r,\theta,r^{\prime},\theta^{\prime},\rho,\eta) 
d\rho d \eta . \label{sameform} 
\end{eqnarray}
For the Isozaki-Kitada parametrix, the amplitude is independent of $t$ and of the form
$$  A_{\rm IK}^{\pm}(t,r,\theta,r^{\prime},\theta^{\prime},\rho,\eta) =  a^{\pm} (r,\theta,\rho,\eta) 
\overline{ b^{\pm}(r^{\prime},\theta^{\prime},\rho,\eta)} , $$
with
\begin{eqnarray}
 a^{\pm} \in {\mathcal S}_{\rm hyp} \left( \Gamma^{\pm}_{\rm s}(\epsilon) \right), \qquad b^{\pm} \in {\mathcal S}_{\rm hyp} \left( 
\Gamma^{\pm}_{\rm s}(\epsilon^3) \right) , \label{amplitudeIKmolle}
\end{eqnarray}
with $ \epsilon > 0 $  small to be fixed.  The phase reads
$$ \Phi_{\rm IK}^{\pm} (t,r,\theta,r^{\prime},\theta^{\prime},\rho,\eta) =
S_{\pm,\epsilon} (r,\theta,\rho,\eta) -
 t \rho^2 - S_{\pm,\epsilon} (r^{\prime},\theta^{\prime},\rho,\eta) , $$
where $ S_{\pm,\epsilon} $ is defined in Proposition \ref{Phaseplusmoinsglobale}. We recall that it coincides with $ S_{\pm} $ on $ \Gamma^{\pm}_{\rm s}(\epsilon) $ (hence on 
$ \Gamma^{\pm}_{\rm s}(\epsilon^3) $ too), where $ S_{\pm} $ is given by  Proposition \ref{Phaseplusmoins}. We can therefore freely replace $ S_{\pm,\epsilon} $
by $ S_{\pm} $, or more generally by any other continuation of $ S_{\pm} $ outside $ \Gamma^{\pm}_{\rm s}(\epsilon) $. Here we have $ 0 \leq \pm t \leq 2 h^{-1} $. 
The integral (\ref{sameform}) is well defined for $ (r,\theta,r^{\prime},\theta^{\prime}) \in \Ra^{2n} $ but, using (\ref{amplitudeIKmolle}), we can assume that
\begin{eqnarray}
 r \geq \epsilon^{-1}, \qquad \theta \in V_{\epsilon}, \qquad r^{\prime} \geq \epsilon^{-3}, \qquad \theta^{\prime} \in V_{\epsilon^{3}} . \label{parametresIK1} 
\end{eqnarray}
The first purpose of this section is to prove that, if $ \epsilon $ is small enough, we can use stationary phase estimates.

The second purpose is to show a similar result for the WKB parametrix, using $ t_{\rm WKB} $ as small parameter (see Theorem \ref{theoremWKB}). In this case, we have to consider
$$   A_{\rm WKB}^{\pm}(t,r,\theta,r^{\prime},\theta^{\prime},\rho,\eta) =  a^{\pm} (t,r,\theta,\rho,\eta)  ,  $$ 
where, for $ V_2 \Subset \psi_{\iota}(U_{\iota}) $, $ I_2 \Subset (0,+\infty) $, $ \sigma_2 \in (-1,1) $, some $ R_2 > 0 $ large enough and some $ t_{\rm WKB} > 0 $, we have
\begin{eqnarray}
(a^{\pm}(t))_{0 \leq \pm t \leq t_{\rm WKB}} \qquad \mbox{bounded in} \ {\mathcal S}_{\rm hyp} \left( \Gamma^{\pm} (R_2,V_2,I_2,\sigma_2) \right) . \label{amplitudeWKBmolle}
\end{eqnarray}
In particular, we can assume that
\begin{eqnarray}
 r \geq R_2 , \qquad \theta \in V_2 . \label{parametresWKB1}
\end{eqnarray}
 The phase is of the form 
\begin{eqnarray}
\Phi_{\rm WKB}^{\pm} (t,r,\theta,r^{\prime},\theta^{\prime},\rho,\eta) =
  \Sigma^{\pm} (t,r,\theta,\rho,\eta ) -r^{\prime} \rho - \theta^{\prime} \cdot \eta , \label{noyauWKB}
\end{eqnarray}
and we refer to Theorem \ref{theoremWKB} for more details.
We only recall here that the phase $ \Sigma^{\pm} $ is  defined on $ [0,\pm t_{\rm WKB}] \times \Ra^{2n} $ and solves the eikonal equation (\ref{eikonalequation}) on  $ [0,\pm t_{\rm WKB}] \times \Gamma^{\pm} (R_3,V_3,I_3,\sigma_3) $, with $ \Gamma^{\pm} (R_2,V_2,I_2,\sigma_2) \subset \Gamma^{\pm} (R_3,V_3,I_3,\sigma_3) $. Here again, the condition (\ref{amplitudeWKBmolle}) implies that we can freely modify $ \Sigma^{\pm} $ outside $ \Gamma^{\pm}(R_2,V_2,I_2,\sigma_2) $.

Below, we will use the notation $ \Phi^{\pm} $ (resp. $ A^{\pm} $) either for $ \Phi_{\rm IK}^{\pm} $ or $ \Phi_{\rm WKB}^{\pm} $ 
(resp. $ A_{\rm IK}^{\pm} $ or $ A_{\rm WKB}^{\pm} $), as long as a single analysis for both cases will be possible. For convenience we also define
$$ 0 \leq \pm t \leq T (h) := \begin{cases} 2 h^{-1} & \mbox{for Isozaki-Kitada}, \\ t_{\rm WKB} & \mbox{for WKB}
\end{cases} . $$



In the next lemma, we summarize the basic properties of $ A^{\pm} $ and $ \Phi^{\pm} $ needed to get a first non stationary phase result. For simplicity, we set $ \partial^{\gamma} = \partial_{r}^{j} \partial_{\theta}^{\alpha} \partial_{r^{\prime}}^{j^{\prime}} 
\partial_{\theta^{\prime}}^{\alpha^{\prime}} \partial_{\rho}^{k} \partial_{\eta}^{\beta} $.
\begin{lemm} In each case,  for all $ |\gamma| \geq 0 $, the amplitude satisfies
\begin{eqnarray}
 | \partial^{\gamma}  A^{\pm} (t,r,\theta,r^{\prime},\theta^{\prime},\rho,\eta) | \leq
C_{\gamma},    \label{estimationamplitude}
\end{eqnarray}
 for all
\begin{eqnarray}
  (r,\theta,r^{\prime},\theta^{\prime},\rho,\eta) \in \Ra^{3n}, \qquad h \in (0,1 ], \qquad 0 \leq \pm  t \leq T (h) , \label{conditionvariables}
\end{eqnarray}
and we may assume that the phase satisfies, 
\begin{eqnarray}
| \partial^{\gamma} \left( \Phi^{\pm} (t,r,\theta,r^{\prime},\theta^{\prime},\rho,\eta) -  (r-r^{\prime}) \rho - (\theta - \theta^{\prime}) \cdot \eta \right) |
\leq C_{\gamma} \scal{t}  , \label{pourCalderonVaillancourtstat}
\end{eqnarray}
under the condition (\ref{conditionvariables}) too. In particular,  for all $ |\gamma| \geq 1 $, 
\begin{eqnarray}
| \partial^{\gamma}  \partial_{\rho} \Phi^{\pm}  (t,r,\theta,r^{\prime},\theta^{\prime},\rho,\eta)  |
 \leq C_{\gamma} \scal{t} , \label{pourCalderonVaillancourtstat2}
\end{eqnarray}
under the condition (\ref{conditionvariables}).
\end{lemm}

\noindent {\it Proof.} If $ A^{\pm} = A^{\pm}_{\rm IK} $, (\ref{estimationamplitude}) follows easily from  Definition \ref{definhyp}, (\ref{amplitudeIKmolle}), (\ref{amplitudeWKBmolle}) and the time independence of $ A^{\pm}_{\rm IK} $.
 If $ A^{\pm} = A^{\pm}_{\rm WKB} $, (\ref{estimationamplitude}) is a direct consequence of (\ref{amplitudeWKBmolle}). 
 For the phase, Proposition \ref{Phaseplusmoinsglobale} (resp. Lemma \ref{asympphaseWKB}) show that $ \Phi_{\rm IK}^{\pm} - (r-r^{\prime})\rho - (\theta-\theta^{\prime})\cdot \eta  $ (resp. $ \Phi_{\rm WKB}^{\pm} - (r-r^{\prime})\rho - (\theta-\theta^{\prime})\cdot \eta  $)
 is the sum of a function  $ f \in C^{\infty}_b (\Ra^{3n}) $
and  $- t \rho^2 $ (resp. $- t p (r,\theta,\rho,\eta) $).  Since the amplitude is compactly supported with respect to $ \rho $ and $ p (r,\theta,\rho,\eta) $, we may  replace $ \Phi_{\rm IK}^{\pm} $ (resp $ \Phi_{\rm WKB}^{\pm} $) by 
$ (r-r^{\prime})\rho - (\theta-\theta^{\prime})\cdot \eta + f - t \rho^2 \chi_1 (\rho) $ (resp. by $(r-r^{\prime})\rho - (\theta-\theta^{\prime})\cdot \eta + f - t p (r,\theta,\rho,\eta) \chi_1 (p(r,\theta,\rho,\eta))$) for some $ \chi_1 \in C_0^{\infty}(\Ra) $. This implies (\ref{pourCalderonVaillancourtstat}) and completes the proof. \finpreuve


\bigskip

Let us choose now $ \chi_1 \in C_0^{\infty}(-1,1), \chi_2 \in C_0^{\infty}(\Ra^{n-1}) $, both equal to $1$ near $ 0 $ and define, for any $ c_1,c_2 > 0 $,
$$ A_{c_1,c_2}^{\pm}  = \chi_1 \left( \frac{\partial_{\rho} \Phi^{\pm} }{c_1 \scal{t} } \right) \chi_2 \left( \frac{\partial_{\eta} \Phi^{\pm}}{c_2} \right) A^{\pm}  . $$
We denote by $ E^{\pm} (t,h) $ the operator with Schwartz kernel (\ref{sameform}) and by $ E_{c_1,c_2}^{\pm}(t,h) $
the operator with Schwartz kernel
\begin{eqnarray}
 (2 \pi h)^{-n} \int \! \! \int e^{\frac{i}{h} \Phi^{\pm} (t,r,\theta,r^{\prime},\theta^{\prime},\rho,\eta) } A_{c_1,c_2}^{\pm}(t,r,\theta,r^{\prime},\theta^{\prime},\rho,\eta)  d\rho d \eta, \label{noyauFPS}
\end{eqnarray}
for $ h \in (0,1] $ and $ 0 \leq \pm t \leq T(h) $.
\begin{prop}[Semi-classical finite speed of propagation] \label{speed} For all $ c_1 ,c_2 > 0 $ and all $ N \geq 0 $, we have
\begin{eqnarray}
|| E^{\pm} (t,h) - E_{c_1,c_2}^{\pm}(t,h) ||_{L^2 (\Ra^n ) \rightarrow L^2 (\Ra^n)} \leq C_{N,A,\Phi,c_1,c_2} h^N , \qquad h \in (0,1], \ \ 0 \leq \pm t \leq T (h) . \label{premierevraiereduction}
\end{eqnarray}
In addition, if $c_1 $ is small enough, there exists $ C \geq 0 $, independent of $ \pm t \in [0,T(h)] $ and of $ c_2 > 0 $, such that
\begin{eqnarray}
 r^{\prime} - r \leq C  , \label{rcompares}
\end{eqnarray}
on the support of $ A_{c_1,c_2}^{\pm} $.
\end{prop}

\noindent {\it Proof.} The kernel of $ E^{\pm} (t,h) - E_{c_1,c_2}^{\pm}(t,h) $ is an oscillatory integral similar to (\ref{noyauFPS}) with amplitude 
$$ A^{\pm} - A_{c_1,c_2}^{\pm} = \left(1 - \chi_1 \left( \frac{\partial_{\rho} \Phi^{\pm} }{c_1 \scal{t} } \right) \right) \chi_2 \left( \frac{\partial_{\eta} \Phi^{\pm}}{c_2} \right) A^{\pm} + \left(1 - \chi_2 \left( \frac{\partial_{\eta} \Phi^{\pm}}{c_2} \right) \right) A^{\pm} . $$ 
On the support of the second term of the right hand side, we integrate by part $ M $ times with 
$$   \frac{h}{i |\partial_{\eta} \Phi^{\pm}|^2} 
\partial_{\eta} \Phi^{\pm} \cdot \partial_{\eta} . $$
Here we note that all derivatives of $ \partial_{\eta} \Phi^{\pm} / | \partial_{\eta} \Phi^{\pm} |^2 $ are bounded since $t$ is bounded in the WKB case and since $ \partial^{\gamma} \partial_{\eta} \Phi^{\pm}_{\rm IK} $ is independent of $t$ and bounded for $ |\gamma| \geq 1 $.
On the support of the first term, we integrate by part $ M $ times with 
$$ \frac{h}{i\partial_{\rho} \Phi^{\pm}} \partial_{\rho}  . $$ 
Using (\ref{pourCalderonVaillancourtstat2}), we have, on the support of the first term,  $ | \partial^{\gamma} (1/\partial_{\rho} \Phi^{\pm}) | \lesssim 1 $, for all $ \gamma $.
Thus, using also (\ref{estimationamplitude}), we end up in both cases  with an integral of the form 
$$ h^{M-n} \int \! \! \int e^{\frac{i}{h} \Phi^{\pm}(t,r,\theta,r^{\prime},\theta^{\prime},\rho,\xi)} B^{\pm}(t,r,\theta,r^{\prime},\theta^{\prime},\rho,\xi) d \rho d \xi $$ with $ B^{\pm}(t,.) $ bounded in $ C^{\infty}_b (\Ra^{3n}) $, for
$ 0 \leq \pm t \leq T (h) $.  We then  interpret this integral as the kernel of a pseudo-differential operator with symbol
$ h^M \exp \{ i(\Phi^{\pm} - (r-r^{\prime})\rho - (\theta-\theta^{\prime})\cdot \eta)/h \} B^{\pm} $ (in the spirit of Lemma \ref{lemmrough}). By the Calder\'on-Vaillancourt Theorem and (\ref{pourCalderonVaillancourtstat}), its operator norm is of order $ h^M (\scal{t}/h)^{n_0} $, for some universal $ n_0 $ (depending only on $n$). We therefore get (\ref{premierevraiereduction}) by choosing $ M = N + 2 n_0 $.

To prove the second statement, we consider separately the two cases. For the WKB parametrix,  $ t$ is bounded. Thus, by (\ref{pourCalderonVaillancourtstat}),  $ \partial_{\rho} \Phi^{\pm}_{\rm WKB} - (r-r^{\prime}) $ is bounded and since 
$ |\partial_{\rho} \Phi^{\pm}_{\rm WKB} | \lesssim c_1 \scal{t}  $,  on the support of $ A_{{\rm WKB},c_1,c_2}^{\pm} $, $ r-r^{\prime} $ must be bounded too. 
For the Isozaki-Kitada parametrix, as long as $t$ belongs to a bounded set the same argument holds. We may therefore assume that $  \pm t \geq T  $ with $ T > 0 $ a fixed large constant. We  then exploit two facts: first, for some $c > 0$, we have $ c < \pm \rho < c^{-1} $ and $ t \rho \geq 0 $ on the support of $ A_{\rm IK}^{\pm} $. Second,  $ f^{\pm}:= \Phi^{\pm}_{\rm IK} - (r-r^{\prime}) \rho - (\theta - \theta^{\prime}) \cdot \eta + t \rho^2 $ is independent of $t$ and bounded, together with all its derivatives on the support of
$ A^{\pm}_{\rm IK} $. Then
$$ \partial_{\rho}\Phi^{\pm}_{\rm IK} = r-r^{\prime} - 2 t \rho + \partial_{\rho} f^{\pm} ,  $$
hence, on the support of $ \chi_1 (\partial_{\rho} \Phi^{\pm}_{\rm IK} / c_1 \scal{t}) $, we have
$$ r-r^{\prime} \geq - c_1 \scal{t} + 2 t \rho - \partial_{\rho} f^{\pm} . $$
If $c_1$ is small enough and $ T $ large enough, we have $ 2t \rho - c_1 \scal{t} \geq  0 $ for $ t \geq T $ and this completes the proof. 
\finpreuve

\bigskip

\noindent {\bf Remark.} It is clear from the proof that the constant $C $ in (\ref{rcompares}) is uniform  to $ \epsilon > 0 $ small in the Isozaki-Kitada case (recall that the amplitudes depend respectively on $t$ and $ \epsilon $ for the WKB and the IK parametrices). 

\bigskip

From now on, we fix $c_1 > 0$ small enough such that (\ref{rcompares}) holds.

\bigskip

\begin{prop}[Dispersion estimate for times $ \leq h $] \label{tempsdordreh} For all $ c_2 > 0   $, we have
$$ || e^{- \gamma_n r} E_{c_1,c_2}^{\pm}(t,h) e^{- \gamma_n r} ||_{ L^1 (\Ra^n ) \rightarrow L^{\infty} (\Ra^n) } \leq 
C_{A,\Phi,c_2} |ht|^{-n/2}, \qquad 0 < \pm t \leq  \min (T(h),h) , $$
where we recall that $ \gamma_n = (n-1)/2 $.
\end{prop}

 Note that
the condition  $ \pm t \leq \min (T(h),h)$ is essentially the condition $ \pm t \leq h $. We have put it under this form to only because of those  $h$ such that $ h \geq t_{\rm WKB} $. The latter will not modify the rest of the analysis. Furthermore, such $h$ correspond to bounded frequencies and their contribution to the Strichartz estimates can be treated by Sobolev embeddings.

\bigskip

\noindent {\it Proof.} In the Isozaki-Kitada case, both $ e^{-r^{\prime}} \eta = \xi $ and $ e^{-r} \eta $ are supported in a compact set. In the WKB one, $ e^{-r} \eta $ is compactly supported but, using (\ref{rcompares}), this also implies that $ e^{-r^{\prime}}\eta $
is compactly supported. Therefore, in both cases, the change of variable $ e^{-r^{\prime}} \eta = \xi $ shows that the kernel of $ E^{\pm}_{c_1,c_2}(t,h) $ is an integral of the form
$$ h^{-n} e^{(n-1)r^{\prime}} \int e^{\frac{i}{h} \Phi^{\pm} (t,r,\theta,r^{\prime},\theta^{\prime},\rho,e^{r^{\prime}}\xi)} B^{\pm}(t,\theta,r^{\prime},\theta^{\prime},\rho,\xi) d \rho d \xi , $$
with $ B^{\pm} $ bounded on $ [0,\pm T(h)] \times \Ra^{3n} $ and supported in a region where $ |\rho|+|\xi| \lesssim 1 $. The kernel of $ 
e^{- \gamma_n r} E^{\pm}_{c_1,c_2}(t) e^{- \gamma_n r} $  is then simply obtained by multiplying the above integral
by $ e^{-\gamma_n (r+r^{\prime})} $ so its modulus is controlled by
$$ h^{-n} e^{\gamma_n (r^{\prime}-r)} \lesssim |ht|^{-n/2} ,$$
using (\ref{rcompares})  and the fact that $ 0 < \pm t \leq h $. This completes the proof. \finpreuve

\bigskip

To prove the dispersion estimates for $ h \leq \pm t \leq T (h) $ we need to analyze more precisely the phases.

In the following lemma and its proof, we shall use the notation (\ref{pourlasectionsuivante}).

\begin{lemm} \label{respIKlemme} For all $ \tilde{\epsilon} > 0 $ small enough (fixed), we 
can find a family of real valued functions $  (\varphi_{\pm, \epsilon}^{\rm st})_{0 < \epsilon \ll 1} $ such that,
\begin{eqnarray}
\varphi_{\pm, \epsilon}^{\rm st} = \varphi_{\pm} = \varphi_{\pm, \epsilon} \qquad \mbox{on} \ \Gamma^{\pm}_{\rm s} (\epsilon),  \label{coincideIKst} \\
\varphi_{\pm,\epsilon}^{\rm st} \in {\mathcal S}_{\rm hyp} \left(  \Gamma_{\rm s}(\tilde{\epsilon}) \right), \label{typefonctionIK}
\end{eqnarray}
and, if we set
$$ R_{\pm,\epsilon} (r,\theta,\rho,\eta) = \varphi_{\pm, \epsilon}^{\rm st} (r,\theta,\rho,\eta) - \frac{q_0(\theta,e^{-r}\eta)}{4 \rho}, $$
the following holds  for $ j + |\alpha| \leq 1 $: if $ k + |\beta|
 \leq 2 $
\begin{eqnarray}
 \sup_{(r,\theta,\eta) \in \Ra^{2n-1},\atop \pm \rho \in [1/4,4]}  \left| (e^{r} \partial_{\eta})^{\beta} \partial_r^{j} 
 \partial_{\theta}^{\alpha} \partial_{\rho}^k  \
 R_{\pm, \epsilon} (r,\theta,\rho,\eta) \right| \leq C \epsilon^{\tau/2} ,
 \label{regimeperturbatifIK}
\end{eqnarray} 
and if $ k + |\beta| \geq 3 $,
$$ \sup_{(r,\theta,\eta) \in \Ra^{2n-1},\atop \pm \rho \in [1/4,4]}  \left| (e^{r} \partial_{\eta})^{\beta} \partial_r^{j} 
 \partial_{\theta}^{\alpha} \partial_{\rho}^k  \
 R_{\pm, \epsilon} (r,\theta,\rho,\eta) \right| \leq C_{\epsilon j \alpha k \beta} , $$
where $ \tau $, the decay rate in (\ref{longueportee}), satisfies (\ref{tauxlongueportee}).
\end{lemm}

\noindent {\it Proof.} Using (\ref{ordre0}), (\ref{ordre1}), (\ref{ordre2}) and Taylor's formula, we can write
$$ \varphi_{\pm}(r,\theta,\rho,\eta) = \int_{0}^{\pm \infty} e^{-4t\rho} q (r+ 2 t \rho, \theta,e^{-r}\eta) dt + \sum_{|\gamma|=3} a_{\gamma}(r,\theta,\rho,\eta) ( e^{-r} \eta )^{\gamma} , $$
with $ a_{\gamma} \in {\mathcal B}_{\rm hyp} (\Gamma^{\pm}_{\rm s}(\epsilon_0)) $ for some fixed $ \epsilon_0 > 0 $.
Therefore, 
\begin{eqnarray}
 \varphi_{\pm}(r,\theta,\rho,\eta) - \frac{q_0 (\theta,e^{-r}\eta)}{4\rho}  =  
 \int_{0}^{\pm \infty} e^{-4t\rho} q_1 (r+ 2 t \rho, \theta,e^{-r}\eta) dt + \sum_{|\gamma|=3} a_{\gamma}(r,\theta,\rho,\eta) ( e^{-r} \eta )^{\gamma} ,
 \label{righthandsidecontinuation}
\end{eqnarray}
with $ q_1 $ satisfying (\ref{longrange}). Denote by $ R (r,\theta,\rho,\eta) $ the right hand side of (\ref{righthandsidecontinuation}) and choose $ \chi_1 \in C_0^{\infty}(\Ra), \chi_2 \in C_0^{\infty}(\Ra^{n-1}) $ both equal to $1$ near $0$. For some $ \tilde{\epsilon} > 0 $ to be fixed below, we also choose $ \chi^{\pm}_{\tilde{\epsilon}}  $ such that
$$ \chi^{\pm}_{\tilde{\epsilon}} \in {\mathcal S}_{\rm hyp}(\Gamma^{\pm}_{\rm s}(\tilde{\epsilon})), \qquad 
\chi^{\pm}_{\tilde{\epsilon}} \equiv 1 \ \ \mbox{on} \ \ \Gamma^{\pm}_{\rm s}(\tilde{\epsilon}^2) , $$ 
using  Proposition \ref{partitions} (we don't need Proposition \ref{construitcutoff} here, since $ \epsilon^0 $ will be fixed). We then claim that,
if $ \tilde{\epsilon} $ is small enough (and fixed) and $ \epsilon $ with $ \tilde{\epsilon}^{\prime} $ is small enough too, the function
$$ \varphi_{\pm, \epsilon}^{\rm st} (r,\theta,\rho,\eta) := \frac{q_0 (\theta,e^{-r}\eta)}{4\rho} + R (r,\theta,\rho,\eta) 
\chi_{\tilde{\epsilon}}^{\pm}(r,\theta,\rho,\eta)  \chi_2 (e^{-r}\eta/\epsilon^{1/2}) (1 - \chi_1)(\epsilon^{1/2} r) , $$
satisfies  (\ref{coincideIKst}), (\ref{typefonctionIK}) and (\ref{regimeperturbatifIK}). Indeed, by choosing $ \tilde{\epsilon}  $ small enough, we have $ \pm \rho \approx 1 $ on the support of $ \chi_{\tilde{\epsilon}}^{\pm} $ so the integral in (\ref{righthandsidecontinuation}) is exponentially convergent. Furthermore, since 
$$ | ( e^{r} \partial_{\eta} )^{\beta} \partial_r^j \left(  (e^{-r}\eta)^{\gamma} \chi_1 (e^{-r}\eta/\epsilon^{1/2}) \right) | \leq C  ( 
\epsilon^{1/2} )^{|\gamma|-|\beta|} , $$
 for all $ \gamma  $, and using the fact that, if $ t \rho \geq 0 $ and $ r \geq 0 $,
$$| (e^{r} \partial_{\eta} )^{\beta}\partial_r^{j} \partial_{\theta}^{\alpha} \partial_{\rho}^k q_1 (r+t \rho,\theta,e^{-r}\eta) | \leq C |t|^k
\scal{r}^{-\tau}|e^{-r} \eta|^{2-|\beta|}, $$
we get the estimate (\ref{regimeperturbatifIK}).
Finally, since $ e^{r}|\eta| \lesssim \epsilon $ and $ r \geq \epsilon $ on $ \Gamma^{\pm}_{\rm s}(\epsilon) $, we have (\ref{coincideIKst})
for all $ \epsilon $ small enough. The property (\ref{typefonctionIK}) is clear thanks to $ \chi_{\tilde{\epsilon}}^{\pm} $. \finpreuve

\bigskip

In the following lemma, we use  the notation of Theorem \ref{theoremWKB}.

\begin{lemm} \label{asympphaseWKB}   We can find a family of real valued functions 
$  (\Sigma^{\pm}_{\rm st}(t) )_{0 \leq \pm t \leq t_{\rm WKB}} $ such that
\begin{eqnarray}
\Sigma^{\pm}_{\rm st}(t) = \Sigma^{\pm}(t) \qquad \mbox{on} \ \Gamma^{\pm} (R_2 , V_2,I_2,\sigma_2),
\end{eqnarray}
and, for all $ k , \beta $,
\begin{eqnarray}
 \sup_{\Ra^{2n}}  \left| (e^{r} \partial_{\eta})^{\beta}  \partial_{\rho}^k  \left(
 \Sigma^{\pm}_{\rm st} (t,r,\theta,\rho,\eta) - r \rho - \theta \cdot \eta - t p (r,\theta,\rho,\eta)  \right) \right| \leq
 C_{k \beta}t^2 . \label{regimeperturbatifWKB}
\end{eqnarray}
\end{lemm}

\noindent {\it Proof.} Using the function $ \chi^{\pm}_{2 \rightarrow 3} $ of Theorem \ref{theoremWKB},  the result is straightforward by considering
\begin{eqnarray*}
  \Sigma^{\pm}_{\rm st} (t,r,\theta,\rho,\eta) & = &  \chi^{\pm}_{2 \rightarrow 3}(r,\theta,\rho,\eta) 
  \left( \Sigma^{\pm}_{\rm 0} (t,r,\theta,\rho,\eta) - r \rho - \theta \cdot \eta - t p (r,\theta,\rho,\eta) \right) + \\ 
 & & r \rho + \theta \cdot \eta + t p (r,\theta,\rho,\eta) ,
 \end{eqnarray*}
 and using (\ref{developpementphaseWKB}). \finpreuve

\bigskip

Let us remark that $ \Sigma^{\pm}  $ satisfies (\ref{globalisationphaseWKB}) whereas $ \Sigma^{\pm}_{\rm st} $ does not. This was the reason for considering
$ \Sigma^{\pm} $ first, since the property (\ref{globalisationphaseWKB}) is convenient to prove $ L^2 $ bounds for Fourier integral operators.

\bigskip

The estimates (\ref{regimeperturbatifIK}) and (\ref{regimeperturbatifWKB}) show  that we have good asymptotics for the phases in certain regimes, namely
$ \epsilon \rightarrow 0 $ for the Isozaki-Kitada parametrix and $ t \rightarrow 0 $ for the WKB parametrix.
Using Lemma \ref{respIKlemme} (resp. Lemma \ref{asympphaseWKB}), we replace $ \varphi_{\pm,\epsilon} $ (resp. $\Sigma^{\pm}$) by $ \varphi_{\pm,\epsilon}^{\rm st} $ 
(resp. by $ \Sigma_{\rm st}^{\pm} $) in the
expression of $ \Phi^{\pm}_{\rm IK} $ (resp. $ \Phi^{\pm}_{\rm WKB} $).

 To use a single formalism for both cases, we introduce the parameter
$$ \lambda_{\rm st} := \begin{cases} \epsilon & \mbox{for the Isozaki-Kitada parametrix} \\ t_{\rm WKB}^{\rm st} & 
\mbox{for the WKB parametrix} \end{cases} , $$
where $  t_{\rm WKB}^{\rm st} > 0  $ will denote the size of the time interval where $ t$ will be allowed to live.
Using the change of variable $ \xi = e^{-r^{\prime}} \eta $ and factorizing by $ t$ in the phase, the integral (\ref{noyauFPS}) can be written
$$ (2 \pi h)^{-n} e^{2 \gamma_n r^{\prime}} \int e^{ i \frac{t}{h} \widetilde{\Phi}_{\lambda_{\rm st}}^{\pm} (z,\rho,\xi)} 
\widetilde{A}_{c_1,c_2\lambda_{\rm st}}^{\pm} (z,\rho,\xi) d \rho d \xi , $$
where $ h \in (0,1] $,
\begin{eqnarray}
\widetilde{\Phi}^{\pm}_{\lambda_{\rm st}} (y,\rho,\xi) = \frac{1}{t} \Phi^{\pm} (t,r,\theta,r^{\prime},\theta^{\prime},\rho,e^{r^{\prime}}\xi), \\
\widetilde{A}^{\pm}_{c_1,c_2,\lambda_{\rm st}} (y,\rho,\xi) = A_{c_1,c_2} (t,r,\theta,r^{\prime},\theta^{\prime},\rho,e^{r^{\prime}}\xi) ,
\label{referencebis}
\end{eqnarray}
and 
\begin{eqnarray}
 y = (h, t,r,\theta,r^{\prime}, \theta^{\prime} ), \label{referenceone}
\end{eqnarray}
with $ r,r^{\prime} $ satisfying (\ref{rcompares})  and 
$$   0 < \pm t \leq T (h,\lambda_{\rm st}) :=  \begin{cases} 2 h^{-1} & \mbox{for the Isozaki-Kitada parametrix} \\ t_{\rm WKB}^{\rm st} & 
\mbox{for the WKB parametrix} \end{cases} , $$
 The kernel of $ e^{-\gamma_n r} E_{c_1,c_2}^{\pm}(t,h) e^{- \gamma_n r} $ then reads
$$  (2 \pi h)^{-n} e^{ \gamma_n ( r^{\prime}-r)} \int e^{ i \frac{t}{h} \widetilde{\Phi}^{\pm}_{\lambda_{\rm st}} (y,\rho,\xi)} 
\widetilde{A}^{\pm}_{c_1,c_2,\lambda_{\rm st}} (y,\rho,\xi) d \rho d \xi .  $$

\begin{prop}[Non stationary phase] \label{nonstationnaireIKWKB} There exists $ C^{\prime} > 0 $ such that the condition
\begin{eqnarray}
\left|\frac{r-r^{\prime}}{t} \right| + e^{r^{\prime}} \left|\frac{\theta-\theta^{\prime}}{t} \right| \geq C^{\prime},  \label{grandparametrenonstationnaire}
\end{eqnarray}
imply that for all $ c_2 > 0 $, all $ N \geq 0 $ and all $ 0 < \lambda_{\rm st} \ll 1 $, we can find $ C_{c_2,N,\lambda_{\rm st}} $ such that, for all
$$   h \in (0,1 ], \qquad \pm t \in [h , T (h,\lambda_{\rm st}) ] , \qquad \omega \geq 1 , \qquad (r,\theta,r^{\prime},\theta^{\prime}) \in \Ra^{2n} ,$$
 with $ r,r^{\prime} $ satisfying (\ref{rcompares}), we have
\begin{eqnarray}
\left| (2 \pi h)^{-n} e^{ \gamma_n (r^{\prime}-r)}  \int e^{i \omega \widetilde{\Phi}_{\lambda_{\rm st}}^{\pm}(y,\rho,\xi)} 
\widetilde{A}^{\pm}_{c_1,c_2,\lambda_{\rm st}} (y,\rho,\xi) d\rho d \xi    \right| \leq C_{c_2,N,\lambda_{\rm st}} h^{-n} \omega^{-N} . \nonumber
\end{eqnarray}
\end{prop}

\noindent {\it Proof.} For $ t \ne 0 $, we define $
 \widetilde{\Phi}_t^{\rm free} : = \frac{r-r^{\prime}}{t} \rho + e^{r^{\prime}} \frac{\theta-\theta^{\prime}}{t} \cdot \xi  $ for which
 $$ \nabla_{\rho,\xi} \widetilde{\Phi}_t^{\rm free} = \left(\frac{r-r^{\prime}}{t}, e^{r^{\prime}} \frac{\theta-\theta^{\prime}}{t} \right) . $$ We then start with the case
 of $ \Phi^{\pm}_{\rm WKB} $. By Lemma \ref{asympphaseWKB} and (\ref{rcompares}), $ \nabla_{\rho,\xi} ( \widetilde{\Phi}_{\lambda_{\rm st}} - \widetilde{\Phi}_t^{\rm free} )$ is a function of $ (t,r,\theta,r^{\prime},\rho,\xi) $ which is bounded on the support of the amplitude, as well as all its derivatives in $ \rho,\xi $, uniformly with respect to $ (t,r,\theta,r^{\prime}) $. Therefore, if $ C^{\prime} $ is large enough, we have
\begin{eqnarray}
\left| \nabla_{\rho,\xi}  \widetilde{\Phi}_{\lambda_{\rm st}} \right| \gtrsim \left|\frac{r-r^{\prime}}{t} \right| + e^{r^{\prime}} \left|\frac{\theta-\theta^{\prime}}{t} \right|, \label{areutiliserpourIK}
\end{eqnarray}
and the result follows from standard integrations by parts. Note that, here, we have not used the smallness of $ \lambda_{\rm st} $ (ie of $t$).
We shall use it  for the case of $ \Phi^{\pm}_{\rm IK} $ which we now consider. Since $ \pm \rho \in [1/4,4] $ on the support of the amplitude if $ \epsilon = \lambda_{\rm st} $ is small enough, Lemma \ref{respIKlemme} and Taylor's formula imply that
$$ \nabla_{\rho,\xi} ( \widetilde{\Phi}_{\lambda_{\rm st}} - \widetilde{\Phi}_t^{\rm free} ) = (-2 \rho,0) + \nabla_{\rho,\xi} \frac{q_0 (\theta,e^{r^{\prime}-r}\xi)-q_0 (\theta^{\prime},\xi)}{t \rho}+ \varepsilon_{\epsilon}(y,\rho,\xi) \left( \frac{r-r^{\prime}}{t} , \frac{\theta-\theta^{\prime}}{t} \right) , $$ 
where  $ \varepsilon_{\epsilon}(y,\rho,\xi)  $ and all its derivatives in $ \rho,\xi $  go to $ 0 $
as $ \epsilon \rightarrow 0 $, 
uniformly with respect to $ y  $ (see (\ref{referenceone})) with $ r,r^{\prime} $ satisfying (\ref{rcompares}) and $ (\pm \rho,\xi) \in [1/4,4] \times \Ra^{n-1} $. Furthermore,
using (\ref{rcompares}) and the fact that $ |\xi| \lesssim \epsilon^{3} $ on the support of the amplitude, we have
$$  \left|\nabla_{\rho,\xi} \frac{q_0 (\theta,e^{r^{\prime}-r}\xi)-q_0 (\theta^{\prime},\xi)}{t \rho} \right| \lesssim \epsilon^{3}\left| \left(\frac{r-r^{\prime}}{t} , \frac{\theta-\theta^{\prime}}{t} \right) \right|  $$
 thus, using that $ r^{\prime} \geq 0 $ on the support of the amplitude, we have (\ref{areutiliserpourIK}) if $ \epsilon $ is small enough.
In addition, for all $ k + |\beta| \geq 2 $, we also have
$$ \left| \partial_{\rho}^k \partial_{\xi}^{\beta} \widetilde{\Phi}_{\lambda_{\rm st}} \right| \lesssim \left| \left(\frac{r-r^{\prime}}{t} , \frac{\theta-\theta^{\prime}}{t} \right) \right|  $$
on the support of the amplitude, using  (\ref{rcompares}). The result then follows again from integrations by parts. \finpreuve

\bigskip

We next state a convenient form of the Stationary Phase Theorem
with parameters.
\begin{prop}[Stationary Phase Theorem] \label{Hormanderrevu}
 Let $ \Omega $ be a set and $f$ be a  function
$$ f : \Ra^n \times \Omega  \ni (x,y) \mapsto f (x,y) \in \Ra $$
smooth with respect to $x$, such that :
\begin{eqnarray}
\emph{Hess}_x [f](x,y) = S (y) + R (x,y), \qquad (x,y) \in \Ra^n
\times \Omega \label{independentdex},
\end{eqnarray}
with $ S (y) $ a symmetric non singular matrix such that
\begin{eqnarray}
|S(y)^{-1}| \lesssim 1 , \qquad y \in \Omega ,
 \label{borneinferieurematrice}
\end{eqnarray}
and $ R (x,y) $ a symmetric matrix such that
\begin{eqnarray}
 ||S(y)^{-1}R(x,y)|| \leq 1/2 \qquad (x,y) \in \Ra^n
\times \Omega , \label{petiteperturbation}
\end{eqnarray}
where  $|| \cdot || $ is the Euclidean matrix norm. Then there
exists  $ N \geq 0 $ such that, for all $ K \Subset \Ra^n $, there
exists $ C_{K}
> 0 $ satisfying
\begin{eqnarray}
 \left| \int e^{i \omega f(x,y)} u (x) dx \right|
\leq C_K \omega^{-n/2} \sup_{ |\alpha| \leq N} ||\partial^{\alpha}
u||_{L^\infty(K)}  \sup_{2 \leq |\alpha| \leq N} \left( \sup_{x
\in K} |\partial^{\alpha} f (x,y)| + 1 \right)^N , \nonumber
\end{eqnarray}
for all  $ y \in \Omega $, all $ u \in C_0^{\infty}(K) $ and all $
\omega \geq 1 $.
\end{prop}

\noindent {\it Proof.} It is a simple adaptation of the proof of
Theorem 7.7.5 in \cite{Horm1}. We give a proof in Appendix
\ref{Phasestationnaire} for completeness. \finpreuve

\bigskip

For the WKB parametrix, we shall use this proposition fairly directly by considering
$$ \Omega_{\rm WKB}^{\pm} \left( t_{\rm WKB}^{\rm st} \right) = \left\{ (h,t,r,\theta,r^{\prime},\theta^{\prime}) \ | \ h \in (0,1 ], \ \ \left|\frac{r-r^{\prime}}{t} \right|  \leq C^{\prime}, \ \  h \leq 
\pm t \leq t_{\rm WKB}^{\rm st}   \right\} . $$ 
Notice that, since $ t_{\rm WKB}^{\rm st} $ is bounded (it will be chosen small enough), $ r - r^{\prime} $ is bounded on $ \Omega_{\rm WKB} \left( t_{\rm WKB}^{\rm st} \right) $. 

\begin{prop}[Dispersion estimate for the WKB parametrix] \label{dispersionpourWKB} Fix $c_2 > 0 $. There exists $ t_{\rm WKB}^{\rm st} > 0 $ small enough such that, for all
$ y =  (h,t,r,\theta,r^{\prime},\theta^{\prime}) \in \Omega_{\rm WKB}^{\pm} \left( t_{\rm WKB}^{\rm st} \right)$ and all $ \omega \geq 1 $, 
we have
\begin{eqnarray}
\left| (2 \pi h)^{-n} e^{ \gamma_n (r^{\prime}-r)}  \int e^{i \omega \widetilde{\Phi}_{t_{\rm WKB}^{\rm st}}^{\pm}(y,\rho,\xi)} \widetilde{A}^{\pm}_{c_1,c_2,t_{\rm WKB}^{\rm st}} (y,\rho,\xi) 
d\rho d \xi  \right| \lesssim \omega^{-n/2}  . \nonumber
\end{eqnarray}
\end{prop}

\noindent {\it Proof.} It is a straightforward application of Proposition \ref{Hormanderrevu} since, using (\ref{regimeperturbatifWKB}), we have 
$$ \mbox{Hess}_{\rho,\xi} [ \widetilde{\Phi}_{t_{\rm WKB}^{\rm st}} ] = \left(
\begin{matrix}
2 & 0  \\
0 & \mbox{Hess}_{\eta}(q)
\end{matrix} \right) + {\mathcal O}(t_{\rm WKB}^{\rm st}) , $$
where the first matrix of the right hand side satisfies (\ref{borneinferieurematrice}) by the uniform ellipticity of $q$. The conclusion is then clear since all derivatives, in $ \rho, \xi $, of $ \widetilde{A}_{t_{\rm WKB}^{\rm st}}^{\pm} $ are bounded, as well as those   of $ \widetilde{\Phi}_{t_{\rm WKB}^{\rm st}} $ of order at least $2$, on the support of the amplitude. \finpreuve

\bigskip

To be in position to use Proposition  \ref{Hormanderrevu} for the Isozaki-Kitada parametrix, we still need two lemmas.

\begin{lemm}[Sharper localization for IK] \label{sharper} Let $ \chi_0 \in C_0^{\infty}(\Ra) $ be equal to $1$ near $ 0 $ and set 
\begin{eqnarray}
 \chi_{\epsilon} (y,\rho) = \chi_0 \left( \epsilon^{-\tau/4} \left( 2 \rho - \frac{r-r^{\prime}}{ t} \right) \right) . \label{referencereference}
\end{eqnarray}
Then, for all $ \epsilon > 0 $ small enough, all $ N \geq 0 $ and all $ c_2 > 0 $, there exists $ C_{c_2,N,\epsilon} $ such that, for all
$$ h \in (0,1 ], \qquad \pm h \leq t \leq 2 h^{-1},  \qquad \omega \geq 1 , $$
and all $ (r,\theta,r^{\prime},\theta^{\prime}) \in \Ra^{2n} $ satisfying (\ref{parametresIK1}) and such that
\begin{eqnarray}
 \left|\frac{r-r^{\prime}}{t} \right| + e^{r^{\prime}} \left|\frac{\theta-\theta^{\prime}}{t} \right| \leq C^{\prime} , \label{zonestationnaire}
\end{eqnarray} 
we have
\begin{eqnarray}
\left| (2 \pi h)^{-n} e^{ \gamma_n( r^{\prime}-r)}  \int e^{i \omega \widetilde{\Phi}_{\epsilon}^{\pm}(y,\rho,\xi)} (1 - \chi_{\epsilon}(y,\rho)) \widetilde{A}_{c_1,c_2,\epsilon}^{\pm} (y,\rho,\xi) 
d\rho d \xi    \right| \leq C_{c_2,N,\epsilon} h^{-n} \omega^{-N} . \nonumber
\end{eqnarray}
\end{lemm}

\noindent {\it Proof.} By the same analysis as in the proof of Proposition \ref{nonstationnaireIKWKB}, using Lemma \ref{respIKlemme} and (\ref{zonestationnaire}), we may write
$$ \widetilde{\Phi}_{\epsilon}^{\pm}(y,\rho,\xi) = \frac{(r-r^{\prime})}{t} \rho - \rho^{2} + R_{\epsilon}^{\pm} (y,\rho,\xi) $$
where, on the support of the amplitude, 
$$ | \partial_{\rho} R_{\epsilon}^{\pm} | \lesssim \epsilon^{\tau/2} , \qquad | \partial_{\rho}^k \partial_{\xi}^{\beta} R_{\epsilon}^{\pm} | 
\lesssim 1 , $$
for $ k +| \beta| \geq 1 $. On the other hand, on the support of $ (1 - \chi_{\epsilon}(y,\rho)) $ we also have, for some $c > 0$,
$$ \frac{r-r^{\prime}}{t} - 2 \rho \geq c \epsilon^{\tau/4} \qquad \mbox{or} \qquad   \frac{r-r^{\prime}}{t} - 2 \rho \leq - c \epsilon^{\tau/4} . $$
Therefore, if $ \epsilon $ is small enough, 
$$ | \partial_{\rho} \widetilde{\Phi}_{\epsilon}^{\pm}(y,\rho,\xi) | \gtrsim \epsilon^{\tau/4} , $$ 
on the support of the amplitude and the result follows from integrations by parts in $ \rho $. \finpreuve

\bigskip

Basically, the interest of the localization (\ref{newamplitude}) is to replace $1/4\rho$ in (\ref{regimeperturbatifIK}) by $ 2 t / (r-r^{\prime})  $ up to a small error. We implement this idea as follows.
By  Lemma  \ref{sharper}, we can replace $ \widetilde{A}^{\pm}_{c_1,c_2,\epsilon} (y,\rho,\xi)  $ in (\ref{referencebis}) by 
\begin{eqnarray}
  \chi_{\epsilon}(y,\rho) \widetilde{A}^{\pm}_{c_1,c_2,\epsilon} (y,\rho,\xi)  . \label{newamplitude}
 \end{eqnarray}
If $ \epsilon $ is small enough, $ \pm \rho \in [1/4,4] $  on the support of $ \widetilde{A}^{\pm}_{c_1,c_2,\epsilon} $ hence,  for some $ c > 0 $,
\begin{eqnarray}
 c |t| \leq  r - r^{\prime} \leq  c^{-1} |t| , \label{controlerr}
\end{eqnarray}
on the support of (\ref{newamplitude}), which is stronger than (\ref{rcompares}). 
Furthermore, the condition (\ref{zonestationnaire}) together with (\ref{parametresIK1}) implies that we may assume that $|\theta-\theta^{\prime}| \leq C^{\prime} e^{- \epsilon^{-3}}|t| $. We fix from now on 
$$ c_2 = \epsilon . $$ 
Thus, by writing
$$ \partial_{\eta} \Phi^{\pm}_{\rm IK} = \theta - \theta^{\prime
} + \partial_{\eta} \varphi_{\pm}(r,\theta,\rho,\eta) - \partial_{\eta} \varphi_{\pm}(r^{\prime},\theta^{\prime},\rho,\eta) , $$ with $ \varphi_{\pm} 
\in {\mathcal B}_{\rm hyp} (\Gamma^{\pm}(\epsilon_2))$, we have$ |\partial_{\eta} \varphi_{\pm}(r,\theta,\rho,\eta) | \lesssim e^{-r} $ and
$ |\partial_{\eta} \varphi_{\pm}(r^{\prime},\theta^{\prime},\rho,\eta) | \lesssim e^{-r^{\prime}} $ on the support of the amplitude. By (\ref{parametresIK1}), we have for instance $ | \partial_{\eta} \Phi^{\pm}_{\rm IK} - (\theta - \theta^{\prime}) | \leq \epsilon^2 $ 
if $ \epsilon $ is small enough. We may therefore assume  that
\begin{eqnarray}
  |\theta-\theta^{\prime}| \leq C^{\prime \prime} \epsilon \frac{|t|}{\scal{t}} .  \label{controlethetatheta}
\end{eqnarray}
To be set of parameters for the stationary phase theorem, we will thus choose
$$ \Omega_{\rm IK}^{\pm}(\epsilon) = \left\{ (h,t,r,\theta,r^{\prime},\theta^{\prime}) \ | \ h \in (0,1], \ \pm t \in [h,2 h^{-1}] \ \mbox{and  (\ref{parametresIK1}), (\ref{zonestationnaire}), (\ref{controlerr}), (\ref{controlethetatheta}) hold} \right\} . $$
Before applying Proposition \ref{Hormanderrevu}, we still need to modify the phase $ \widetilde{\Phi}_{\epsilon}^{\pm} $ outside the support of the new amplitude (\ref{newamplitude}).

\begin{lemm} \label{phasestationaireexplicite} We can find $ \Psi_{\epsilon}^{\pm} $ smooth and real valued such that, on the support of (\ref{newamplitude}),
$$ \Psi_{\epsilon}^{\pm}(y,\rho,\xi) = \widetilde{\Phi}_{\epsilon}^{\pm} (y,\rho,\xi) ,$$
and, 
\begin{eqnarray}
  \Psi_{\epsilon}^{\pm}(y,\rho,\xi) =  \frac{r-r^{\prime}}{t} \rho + \frac{\theta-\theta^{\prime}}{t}e^{r^{\prime}}\xi - \rho^2 - \frac{1-e^{2(r^{\prime}-r)}}{2(r-r^{\prime})}  q_0 (\theta^{\prime},\xi)  + \psi_{\epsilon}^{\pm}(y,\rho,\xi) , \label{expansionphasepourst}
\end{eqnarray}
where, for all $ k + |\beta| \leq 2 $,
\begin{eqnarray}
 \sup_{(\rho,\xi) \in \Ra^n , \atop
y \in \Omega_{\rm IK}^{\pm}(\epsilon)} | \scal{t}^{|\beta|/2} \partial_{\rho}^k \partial_{\xi}^{\beta} \psi_{\epsilon}^{\pm}(y,\rho,\xi) | \rightarrow 0 , \qquad \epsilon \rightarrow 0  , \label{ordre2phase}
\end{eqnarray}
and, $ |k| + |\beta | \geq 3 $,
\begin{eqnarray}
 \sup_{(\rho,\xi) \in \Ra^n , \atop
y \in \Omega_{\rm IK}^{\pm}(\epsilon)} | \ \partial_{\rho}^k \partial_{\xi}^{\beta} \psi_{\epsilon}^{\pm}(y,\rho,\xi) | \leq C_{\epsilon,k,\beta}  .
\label{deriveesdordreaumoins3}
 \end{eqnarray}
\end{lemm}

\noindent {\it Proof.} We shall basically combine (\ref{regimeperturbatifIK}) with the fact that 
\begin{eqnarray}
 | 2 \rho - (r-r^{\prime})/t | \lesssim \epsilon^{\tau/4} , \label{pourindependenceenrho}
 \end{eqnarray}
  on the support of (\ref{newamplitude}). By Lemma \ref{respIKlemme}, the phase reads
 $$  \frac{r-r^{\prime}}{t} \rho + \frac{\theta-\theta^{\prime}}{t}e^{r^{\prime}}\xi - \rho^2 - \frac{  q_0 (\theta^{\prime},\xi) - e^{2(r^{\prime}-r)} q_0 (\theta,\xi)}{4 \rho t}  + \frac{R_{\pm,\epsilon}(r,\theta,\rho,e^{r^{\prime}}\xi)-R_{\pm,\epsilon}(r^{\prime},\theta^{\prime},\rho,e^{r^{\prime}}\xi)}{t} . $$
The last term of this sum satisfies the estimates (\ref{ordre2phase}) and (\ref{deriveesdordreaumoins3}): for $ 0 < \pm t \leq 1 $,
it follows from Taylor's formula using (\ref{zonestationnaire}) and Lemma \ref{respIKlemme} with $ j + |\alpha| = 1 $, and for $ \pm t \geq 1 $ it follows from Lemma \ref{respIKlemme}
 with $ j + |\alpha| = 0 $. For the term involving $ q_0 $ we write
 $$ \frac{1}{4 \rho t} = \frac{1}{2(r-r^{\prime})} + \left( \frac{1}{4 \rho t} - \frac{1}{2(r-r^{\prime})} \right) \chi_1 \left( \frac{2 \rho - (r-r^{\prime})/t}{\epsilon^{\tau/8}} \right) ,$$
 using (\ref{pourindependenceenrho}) with $ \epsilon $ small enough  and $ \chi_1 \in C_0^{\infty}(\Ra^{n-1}) $ equal to $ 1 $ near $ 0 $, and
$$ q_0 (\theta,e^{r^{\prime}-r}\xi)  = e^{2(r^{\prime}-r)}q_0 (\theta^{\prime},\xi) + e^{2(r^{\prime}-r)} ( q_0 (\theta,\xi) - q_0 (\theta^{\prime},\xi) ) \chi_2 (\xi) , $$
with $ \chi_2 \in C_0^{\infty}(\Ra^{n-1}) $ equal to $ 1 $ near $ 0 $. We obtain the estimates (\ref{ordre2phase}) and (\ref{deriveesdordreaumoins3}) for
$$ \frac{1}{4 \rho t} e^{2(r^{\prime}-r)} ( q_0 (\theta,\xi) - q_0 (\theta^{\prime},\xi) ) \chi_2 (\xi) , $$
using (\ref{controlethetatheta}), and for 
$$ ( 1 - e^{2(r^{\prime}-r)}) q_0 (\theta^{\prime},\xi) \left( \frac{1}{4 \rho t} - \frac{1}{2(r-r^{\prime})} \right) \chi_1 \left( \frac{2 \rho - (r-r^{\prime})/t}{\epsilon^{\tau/8}} \right) $$
using (\ref{controlerr}). In both cases, we can freely multiply the functions by a compactly support cutoff in $ \rho $ using that $ \pm \approx 1 $
on the support of the amplitude. This completes the proof. \finpreuve

\bigskip

\begin{prop}[Bounded times] \label{tempsbornesIK} There exists $ \epsilon_{\rm st} > 0 $ such that,
 for all $ T > 0 $, all $ 0 < \epsilon \leq \epsilon_{\rm st} $, there exists $ C_{\epsilon,T} $ such that, for  all
 \begin{eqnarray}
  h \in (0,1], \qquad h \leq \pm  t \leq T, \qquad (r,\theta,r^{\prime},\theta^{\prime}) \ \ \mbox{satisfying (\ref{parametresIK1}), (\ref{controlerr}) and (\ref{controlethetatheta}) }  , \label{parametrestempsbornes}
\end{eqnarray}
 we have
 \begin{eqnarray}
\left| (2 \pi h)^{-n} e^{ \gamma_n (r^{\prime}-r)}  \int e^{i \frac{t}{h} \widetilde{\Phi}_{\epsilon}^{\pm}(y,\rho,\xi)} \chi_{\epsilon}(y,\rho) \widetilde{A}_{c_1,\epsilon,\epsilon} (y,\rho,\xi) 
d\rho d \xi  \right| \leq C_{\epsilon,T} |ht|^{-n/2}  . \label{pasdetempsexponentiel}
\end{eqnarray}
\end{prop}

\noindent {\it Proof.} By Lemma \ref{phasestationaireexplicite}, we can replace $ \widetilde{\Phi}_{\epsilon}^{\pm} $ by 
$ \Psi_{\epsilon}^{\pm} $. We then have
$$ \mbox{Hess}_{\rho,\xi} [ \Psi_{\epsilon}^{\pm} ] = \left(
\begin{matrix}
2 & 0  \\
0 & \frac{1-e^{2(r^{\prime}-r)}}{2 (r-r^{\prime})} \mbox{Hess}_{\eta}(q_0)
\end{matrix} \right) + o(1) , $$
where $ o(1) \rightarrow 0 $ as $ \epsilon \rightarrow 0 $ , uniformly with respect to $ (\rho,\xi) \in \Ra^{n} $ and to the parameters satisfying
(\ref{parametrestempsbornes}). Using the upper bound in (\ref{controlerr}) and the boundedness of $t$, the positive number
$$ \frac{1-e^{2(r^{\prime}-r)}}{2 (r-r^{\prime})} $$
 belongs to a compact subset of $ (0,\infty) $, yielding the condition (\ref{borneinferieurematrice}). We conclude by applying Proposition \ref{Hormanderrevu}. \finpreuve

\bigskip

Notice that, to obtain (\ref{pasdetempsexponentiel}), we have  used the boundedness of $ e^{\gamma_n (r^{\prime}-r)} $, since $ |r - r^{\prime} |$ was bounded. In principle, the condition (\ref{controlerr}) implies that $ e^{\gamma_n (r^{\prime}-r)} $ decays exponentially in time. We shall exploit the latter below.

\bigskip

\begin{prop}[Large times] \label{grandstempsIK} There exists $ T > 0 $ and $ \epsilon_{\rm st}^{\prime} $ such that, for all $ 0 < \epsilon \leq 
\epsilon_{\rm st}^{\prime} $, there exists $ C_{\epsilon} $ such that, for all
\begin{eqnarray}
 h \in (0,1], \qquad T \leq \pm  t \leq 2 h^{-1}, \qquad (r,\theta,r^{\prime},\theta^{\prime}) \ \ \mbox{satisfying (\ref{parametresIK1}), (\ref{controlerr}) and (\ref{controlethetatheta}) }  ,  \label{parametrestempsnonbornes}
\end{eqnarray}
 we have
 \begin{eqnarray}
\left| (2 \pi h)^{-n} e^{ \gamma_n (r^{\prime}-r)}  \int e^{i \frac{t}{h} \widetilde{\Phi}_{\epsilon}^{\pm}(y,\rho,\xi)} \chi_{\epsilon}(y,\rho) \widetilde{A}_{c_1,\epsilon,\epsilon} (y,\rho,\xi) 
d\rho d \xi  \right| \leq C_{\epsilon} |ht|^{-n/2}  . \nonumber
\end{eqnarray}
\end{prop}

\noindent {\it Proof.} Choose $T$ large enough such that, for $t\geq T$ and $ r,r^{\prime} $ satisfying (\ref{controlerr}), we have $ e^{2(r^{\prime}-r)} \leq 1/2   $. To compensate the factor $ 1/(r-r^{\prime}) $ in (\ref{expansionphasepourst}) (of order $1/|t|$ by (\ref{controlerr})), we consider the new variable $ |t|^{1/2} \zeta = \xi $. By (\ref{ordre2phase}), if $ \epsilon  $ is small enough, this new phase satisfies the assumptions of Proposition \ref{Hormanderrevu}. In the corresponding estimate given by Proposition \ref{Hormanderrevu},  derivatives of the new amplitude as well as derivatives of the new phase of order at least $3$ will grow at most polynomially with respect to $t$. This gives a polynomial growth in $t$ of the coefficient in the stationary phase estimate of Proposition \ref{Hormanderrevu}  but such a growth is controlled by the exponential decay of $ e^{\gamma_n (r^{\prime}-r)} \lesssim e^{-c|t|} $. This completes the proof. \finpreuve

\subsection{Proof of Proposition \ref{sousIK}}

By (\ref{adjointpseudodifferentiel}), up to a remainder of operator norm of size $h^n$ (uniformly in time), we may replace $ \widehat{O \! p}_{\iota}(a_{\rm s}^{\pm})^* $ by a linear combination of operators of the form $ \widehat{O \! p}_{\iota}(\tilde{a}_{\rm s}^{\pm}) $
with $ \mbox{supp} (\tilde{a}_{\rm s}^{\pm} ) \subset \mbox{supp}(a_{\rm s}^{\pm}) $. We next apply Theorem \ref{IsozakiKitadaansatz} to order $  n + 1 $ and are left with the study of the Fourier integral operator part. By Proposition \ref{speed}, the amplitude can be modified so that, up to a remainder of operator norm of order
 $ h^n $ uniformly in time, we are left with an operator whose kernel $ K^{\pm} (r,\theta,r^{\prime},\theta^{\prime},t,h) $ satisfies
 $$ | e^{- \gamma_n r} K^{\pm} (r,\theta,r^{\prime},\theta^{\prime},t,h) e^{- \gamma_n r^{\prime}} | 
 \lesssim |ht|^{-n/2}, \qquad h \in (0,1 ], 0 < \pm t \leq 2h .  $$
 Indeed, for $ t \leq h $, this follows from Proposition \ref{tempsdordreh} and for $ t \geq h $, from Propositions  \ref{tempsbornesIK} and \ref{grandstempsIK} with $ \omega = \pm t/h $ and also from Proposition \ref{nonstationnaireIKWKB} and Lemma \ref{sharper}
with $ N \geq n / 2 $. \finpreuve
 
\subsection{Proof of Proposition \ref{sousWKB}}
It is completely similar to the one of Proposition \ref{sousIK} by considering times $ 0 \leq \pm t \leq  t_{\rm WKB}^{\rm st} $ with  $ t_{\rm WKB}^{\rm st}$ small enough to be in position to use both Theorem \ref{theoremWKB} and Proposition \ref{dispersionpourWKB}. \finpreuve

\appendix

\section{Control on the range of some diffeomorphisms}
\label{Diff}
 \setcounter{equation}{0}
In this section, we prove a proposition implying
 Lemma \ref{preparationphase} and (\ref{pourlapropositionsuivante2})  in Lemma \ref{Kuranishi}.
 For simplicity, we consider  the outgoing case only but the symmetric result holds
 in the incoming one.

Let us define the following conical subset of $ T^* \Ra^n_+
\setminus 0 $,
\begin{eqnarray}
\Gamma^{+}_{\rm s-con} (\epsilon) & = & \left\{
(r,\theta,\rho,\eta)  \ | \ r
> R(\epsilon) , \ \theta \in V_{\epsilon} , \  \rho
> (1-\epsilon^2) (\rho^2 + q (r,\theta,e^{-r}\eta))^{1/2}
 \right\} , \label{coniqueoblige}
\end{eqnarray}
which  is the cone generated by $ \Gamma^{+}_{\rm s}(\epsilon)  $.

\begin{prop} \label{propositiondouble}  Assume that, for some $ 0
< \bar{\epsilon} < 1/4  $, we are given a family of maps $
(\Psi^t)_{t \geq 0} $ defined on $ \Gamma^{+}_{\rm
s-con}(\bar{\epsilon}) $, of the form
$$ \Psi^t (r,\theta,\rho,\eta) = (r,\theta,
\underline{\rho}^t(r,\theta,\rho,\eta),\underline{\eta}^t
(r,\theta,\rho,\eta)) \in \Ra^{2n} ,  $$ satisfying, for all   $ r
> R (\bar{\epsilon})$, $\theta \in V_{\bar{\epsilon}} $, $ \rho
> (1-\bar{\epsilon}^2) p^{1/2} $, $ t \geq 0 $ and $ \lambda > 0
$,
\begin{eqnarray}
( \underline{\rho}^t ,\underline{\eta}^t ) (r,\theta,\lambda \rho
, \lambda \eta) & = & \lambda ( \underline{\rho}^{\lambda t},
\underline{\eta}^{\lambda t} ) (r,\theta,\rho,\eta)  , \label{pseudohomogene} \\
 ( \underline{\rho}^t , \underline{\eta}^t ) (r,\theta,\rho,0) & = &(\rho,0) , \label{restrictionidentite}
\end{eqnarray}
 and such that,
\begin{eqnarray}
(\underline{\rho}^t - \rho)_{t \geq 0} \ \mbox{and (the components
of)}\ (\underline{\eta}^t - \eta)_{t \geq 0} \ \mbox{ are bounded
in } \ {\mathcal B}_{\rm hyp} (\Gamma^{+}_{\rm s}(\bar{\epsilon}))
. \label{borneuniforme}
\end{eqnarray}
Then, there exists  $0 < \tilde{\epsilon} \leq \bar{\epsilon} $
such that,  for all $ t \geq 0 $ and all $ 0 < \epsilon \leq
\tilde{\epsilon} $, $ \Psi^t $ is a diffeomorphism from $
\Gamma^{+}_{\rm s}(\epsilon) $ onto its range and
$$ \Gamma^{+}_{\rm s}(\epsilon^3) \subset \Psi^t \left( \Gamma^{+}_{\rm s}(\epsilon) \right) , \qquad t \geq 0 , \ \
0 < \epsilon \leq \tilde{\epsilon} . $$
\end{prop}

 Lemma \ref{preparationphase} is indeed a consequence of Proposition  \ref{propositiondouble}
since Proposition \ref{estimeesflotprecises},
(\ref{scalingrhoeta}) and (\ref{sansderiver}) show that
(\ref{pseudohomogene}), (\ref{restrictionidentite}) and
(\ref{borneuniforme}) hold with $
(\underline{\rho}^t,\underline{\eta}^t ) = (\rho^t , \eta^t) $.
Similarly, for Lemma \ref{Kuranishi}, we consider
$$ (\underline{\rho}^t,\underline{\eta}^t)
(r,\theta,\rho,\eta) := (\underline{\rho}_+,\underline{\eta}_+)
(r,\theta,r,\theta,\rho,\eta) $$ which is independent of $t$ and
satisfies the assumptions (\ref{pseudohomogene}),
(\ref{restrictionidentite}), (\ref{borneuniforme}) by
(\ref{aevaluersurlinverse}), Proposition \ref{Phaseplusmoins} and
Remark 2 after Proposition \ref{Phaseplusmoins}.

\bigskip

To prove the proposition, we need another conical subset of $ T^*
\Ra^n_+ \setminus 0 $:
\begin{eqnarray}
 \widetilde{\Gamma}^{+}_{\rm s-con} (\epsilon ) & = & \left\{
(r,\theta,\rho,\eta)  \ | \ r > R(\epsilon) , \ \theta \in
V_{\epsilon} , \  \rho
> (1-\epsilon^2) (\rho^2 + |\eta|^2)^{1/2}
 \right\} . \nonumber
\end{eqnarray}
Using the diffeomorphism $F_{\rm hyp}$ defined by (\ref{Fhyp}), we have
\begin{eqnarray}
 F_{\rm hyp}^{-1} \left( \widetilde{\Gamma}^{+}_{\rm s-con} (\epsilon) \right)  =  \left\{
(r,\theta,\rho,\eta) \ | \ r > R(\epsilon) , \ \theta \in
V_{\epsilon} , \  \rho
> (1-\epsilon^2) (\rho^2 + |e^{-r}\eta|^2)^{1/2}
 \right\} .
\end{eqnarray}
The latter is of interest in view of the following lemma.
\begin{lemm} \label{lemmeconique} There exists  $
C > 1 $  such that, for all $  \epsilon > 0  $ small enough,
$$ \Gamma^{+}_{\rm s-con} ( \epsilon/C  )
\subset  F_{\rm hyp}^{-1} \left( \widetilde{\Gamma}^{+}_{\rm s-con} (\epsilon
) \right) \subset \Gamma^{+}_{\rm s-con} (C \epsilon ) . $$
\end{lemm}

\bigskip

\noindent {\it Proof.} By  (\ref{secogene}),
we have, for some $ 0 < c < 1  $,
$$ c e^{-2r}|\eta|^2 \leq q (r,\theta,e^{-r}\eta) \leq c^{-1}
|e^{-r}\eta|^2 , \qquad  r > R(\epsilon) , \ \theta \in
V_{\epsilon} , \ \eta \in \Ra^{n-1} . $$ Using
(\ref{equivalencepratique}), it suffices to show the existence of
$C> 1$ satisfying, for all $ \epsilon $ small enough,
\begin{eqnarray}
 c^{-1} (1- (\epsilon /C)^2)^{-2} \left( 1 - (1 -
(\epsilon/C)^2)^2 \right) \leq  (1 - \epsilon^2)^{-2} \left( 1 -
(1-\epsilon^2)^2 \right) , \label{developpement1}
\end{eqnarray}
and
\begin{eqnarray}
  (1 - \epsilon^2)^{-2} \left( 1 - (1-\epsilon^2)^2 \right)
\leq c (1- (C \epsilon )^2)^{-2} \left( 1 - (1 - (C \epsilon)^2)^2
\right) . \label{developpement2}
\end{eqnarray}
For $ \epsilon \rightarrow 0 $, the left hand side of
(\ref{developpement1}) is equivalent to $ 2 c^{-1} (\epsilon /
C)^2 $ and the right hand side to $ 2 \epsilon^2 $. Therefore,
(\ref{developpement1}) holds if $ c^{-1} / C^2 < 1 $ and $
\epsilon $ is small enough. We get (\ref{developpement2})
similarly.
 \finpreuve


\bigskip

Let us now consider  $ ( 1,0) = ( 1 , 0 , \ldots , 0) \in \Ra^n
\setminus 0  $. For all  $ 0 < \epsilon < 1  $, let us denote by $
{\mathcal C}^{+} (\epsilon) $ the cone generated by $ B ((
1,0),\epsilon) $, namely
$$   {\mathcal C}^{+} (\epsilon) =  \{ (\lambda \rho , \lambda \eta) \ | \ \lambda > 0 , \
(\rho - 1)^2 + |\eta|^2 < \epsilon^2  \}  . $$  Since $ \rho
> 1 - \epsilon > 0 $ and  $ \rho^2 / (\rho^2 + |\eta|^2) > 1 - \epsilon^2 / (1 - \epsilon)^2 $
on $ B ((1,0),\epsilon) $,
  it is then not hard to check that, for all $ \epsilon $ small
  enough,
$$ {\mathcal C}^{+} (\epsilon^2/4) \subset
 \{  \rho > (1-\epsilon^2) (\rho^2 + |\eta|^2)^{1/2}  \}
   , $$
and
$$ \{  \rho > (1-\epsilon^2) (\rho^2 + |\eta|^2)^{1/2}  \}
 \subset {\mathcal C}^{+} (2 \epsilon ), $$
since, if $ \rho > (1-\epsilon^2) (\rho^2 + |\eta|^2)^{1/2} $ then
$(1,\eta/\rho) \in B ((1,0),2 \epsilon)$, using that $ 1 - (1 -
\epsilon^2)^2 < 4 \epsilon^2 (1 - \epsilon^2)^2 $ for $ \epsilon $
small enough. In particular,  we obtain
\begin{eqnarray}
 (R(\epsilon),+\infty) \times V_{\epsilon} \times {\mathcal C}^{+}(\epsilon^2/4)
\subset   \widetilde{\Gamma}^{+}_{\rm s-con}(\epsilon) \subset
(R(\epsilon),+\infty) \times V_{\epsilon} \times {\mathcal C}^{+}
(2 \epsilon ) .  \label{lemmeconique0}
\end{eqnarray}

\bigskip

We next recall a standard lemma the simple proof of which we omit.
\begin{lemm} \label{lemmenondemontre} Let $ x_0 \in \Ra^n $, $ \varepsilon > 0 $ and
$ f : B (x_0,\varepsilon) \rightarrow \Ra^n $ such that $ f (x_0)
= x_0 $ and $ f- \emph{id} $ is $ 1/2 $ Lipschitz, ie $ | f(x)-x +
y - f (y) | \leq |x-y|/2 $, on $ B (x_0,\varepsilon) $. Then $ f $
is injective on $ B (x_0,\varepsilon) $ and
$$ B (x_0,\varepsilon/2) \subset f (B(x_0,\varepsilon)) . $$
\end{lemm}

\noindent  {\it Proof of Proposition \ref{propositiondouble}.} Let
us  set
$$ f_{r,\theta,t} (\rho,\xi) = \left( \underline{\rho}^t (r,\theta,\rho,e^{r}\xi) , e^{-r}
  \underline{\eta}^t (r,\theta,\rho,e^{r}\xi) \right) . $$
By  Lemma \ref{lemmesymbole} and (\ref{borneuniforme}), we have, for
$ k + |\beta| = 2  $,
\begin{eqnarray}
 | \partial_{\rho}^k \partial^{\beta}_{\eta} f_{r,\theta,t} (\rho,\eta) | \lesssim 1  , \qquad
 \ t \geq 0, \ (r , \theta , \rho , \xi ) \in F_{\rm hyp} ( \Gamma^+_{\rm
s}(\bar{\epsilon})) ,  \label{descaling}
\end{eqnarray}
and, by choosing $ \bar{\epsilon}_1  $ small enough, we also have
$$ r > R ( \bar{\epsilon} ), \ \ \theta \in V_{\bar{\epsilon} } , \ \ (\rho , \xi) \in B ((1,0),\bar{\epsilon}_1)
 \ \ \Rightarrow
(r , \theta , \rho , \xi ) \in F_{\rm hyp} ( \Gamma^+_{\rm
s}(\bar{\epsilon})) . $$
 By
(\ref{restrictionidentite})  $ \partial_{\rho,\xi} f_{r,\theta,t}
(\rho,0) = \mbox{Id}_n $ , so (\ref{descaling}) implies that $
f_{r,\theta,t} - \mbox{Id}_n  $ is $ 1/2 $-Lipschitz on $ B ((1,0),2
\epsilon)$ for all $ \epsilon $ small enough, all $ t \geq 0  $, $
r  > R (\bar{\epsilon})$, and $\theta \in V_{\bar{\epsilon}} $.
Therefore, by Lemma \ref{lemmenondemontre},
$$ B ( (1,0),\epsilon ) \subset f_{t,r,\theta} \left( B ( (1,0), 2 \epsilon)
\right) , \qquad t \geq 0 , \ r > R (\tilde{\epsilon}), \ \theta
\in V_{\tilde{\epsilon}} .  $$ Using (\ref{pseudohomogene}), we
can replace the balls in the above inclusion by the cones they
generate  and, using  Lemma \ref{lemmeconique} with
(\ref{lemmeconique0}), we get
\begin{eqnarray}
  \Gamma^+_{\rm s-con} (\epsilon /2 C )   \subset
 \Psi^t \left( \Gamma^+_{\rm s-con}  (2 \sqrt{2} C \epsilon^{1/2}) \right) , \qquad t \geq 0, \label{versionconique}
\end{eqnarray}
 for all $ \epsilon  $ small
enough, with the $ C > 1 $ of Lemma \ref{lemmeconique}. Since $
f_{r,\theta,t} - \mbox{Id}_n  $ is $ 1/2 $-Lipschitz on $ B ((1,0),2
\epsilon)$ for all $ t \geq 0 $, (\ref{pseudohomogene}) implies
that it is also $ 1/2 $-Lipschitz on the cone generated by $ B
((1,0),2 \epsilon) $ so
 $  f_{r,\theta,t}  $ is injective on this cone. Thus, for all $ \epsilon
 $ small enough and $ t \geq 0 $, $ \Psi^t  $ is injective on $ \Gamma^+_{\rm s-con}  (\epsilon)  $
and is a diffeomorphism onto its range. By (\ref{versionconique}),
we have
$$\Gamma^+_{\rm s}(\epsilon^3) \subset \Gamma^+_{\rm
s-con}(\epsilon^3) \subset \Psi^t \left( \Gamma^+_{\rm
s-con}(\epsilon) \right), $$
 for all $t \geq 0 $ and all $ \epsilon $ small enough, so
 the proof will be complete by showing that, for all $ \epsilon $ small enough and all $ t
\geq 0 $, the following implication
holds
\begin{eqnarray}
 (r,\theta,\rho,\eta) = \Psi^t (r,\theta,\rho_1,\eta_1) \in
\Gamma^+_{\rm s}(\epsilon^3) \ \  \mbox{with} \ \
(r,\theta,\rho_1,\eta_1) \in \Gamma^+_{\rm s-con}(\epsilon)
\nonumber \qquad
\qquad \\
\Rightarrow \qquad p (r,\theta,\rho_1,\eta_1) \in (1/4 -
\epsilon, 4 + \epsilon) . \label{implicationholds}
\end{eqnarray}
Assume the first line of (\ref{implicationholds}). Using
(\ref{restrictionidentite}) at $ (\rho_1,0) $ and the fact that $
f_{t,r,\theta} - \mbox{Id}_n $ is $ 1/2 $-Lipschitz, we have
\begin{eqnarray}
 |(\rho,e^{-r}\eta)-(\rho_1,e^{-r}\eta_1)| =  | f_{t,r,\theta} (\rho_1,e^{-r}\eta_1) - (\rho_1,e^{-r}\eta_1) | \leq
|e^{-r} \eta_1| / 2 . \label{autiliser2fois}
\end{eqnarray}
Therefore $  |e^{-r}\eta-e^{-r}\eta_1| \leq |e^{-r} \eta_1| / 2
$ and we get $ |\eta_1| \leq 2 |\eta|  $. Since $ e^{-r}|\eta|
\lesssim \epsilon^3  $, (\ref{autiliser2fois}) shows that $ |\rho
- \rho_1|+|e^{-r}(\eta-\eta_1)|  \lesssim \epsilon^3  $ hence that
$$ |  p (r,\theta,\rho_1,\eta_1) - p (r,\theta,\rho,\eta) | \lesssim
\epsilon^3 .  $$ Since $  p (r,\theta,\rho,\eta) \in  (1/4 -
\epsilon^3, 4 + \epsilon^3) $, the latter yields
(\ref{implicationholds}) for $ \epsilon  $ small enough.
\finpreuve

\section{Proof of Lemma \ref{Hormanderrevu}} \label{Phasestationnaire}
\setcounter{equation}{0}

 Note first that, for all $ y \in \Omega $, the map
$$ \Ra^n \ni x  \mapsto \nabla_{x} f (x,y) \in \Ra^n $$
is a diffeomorphism since, by considering $ F (x,y) := S (y)^{-1}
\nabla_{x} f (x,y) $ and using (\ref{independentdex}),
(\ref{petiteperturbation}) and (\ref{borneinferieurematrice}), $ x
\mapsto F (x,y) -x $ is $ 1/2 $ Lipschitz. For all $ y \in \Omega
$, we then denote by $ x_0 = x_0
(y) $ the unique solution to
$$ \nabla_{x} f (x_0,y) = 0  . $$
Let us now consider
$$ g (x,y) = f(x,y) - f (x_0,y) -
\left\langle \mbox{Hess}_x [f](x_0,y)(x-x_0),x-x_0 \right\rangle /
2,
$$
and, for all $ s \in [0,1] $,
$$ f_s (x,y) = f (x_0,y) + \left\langle \mbox{Hess}_x [f](x_0,y)(x-x_0),x-x_0
\right\rangle / 2 + s g (x,y) . $$ Notice that $ f_1 = f $, that $
f_0 - f (x_0,y) $ is quadratic with respect to $ x - x_0 $ and
that
$$ \nabla_{x} f_s (x,y) = \left\{ S (y) + s \int_0^1 R (x_0 + t (x-x_0),y) dt
+ (1-s)R(x_0,y) \right\}(x-x_0), $$
 by the Taylor formula and (\ref{independentdex}).
 By (\ref{borneinferieurematrice}), there exists $
c
> 0 $ such that $ |S (y)X| \geq 2 c |X| $, for all $ X \in \Ra^n $
and all $ y \in \Omega $ hence (\ref{petiteperturbation}) implies
that
\begin{eqnarray}
 | \nabla_{x} f_s (x,y) | \geq c |x-x_0 (y)| ,
 \qquad s \in [0,1], \ (x,y) \in \Ra^n
\times \Omega . \label{controlepointcritique}
\end{eqnarray}

\begin{lemm} \label{nonstationnaireHormander} For all $ K \Subset \Ra^n $ and all integer $ k \geq 1 $,
 there exists $ C > 0$ and $ N > 0 $ such that,
for all $s \in [0,1] $, all $ y \in \Omega $ and all $u$ such that
\begin{eqnarray}
 u \in C_0^{2k-1}(K) \cap
C^{2k}(\Ra^n \setminus \{ x_0 (y) \}),
\label{regularitenoncinifinie1} \\ \partial^{\alpha}_x u (x_0(y))
= 0 , \qquad |\alpha| < 2 k,
 \label{regularitenoncinifinie2} \\
\partial_x^{\alpha }u \in L^{\infty}(\Ra^n), \qquad |\alpha| = 2 k, \label{regularitenoncinifinie3}
\end{eqnarray}
we have
$$  \left|  \int e^{i \omega f_s (x,y)} u (x) d x \right| \leq
C \omega^{-k} \max_{|\alpha| \leq 2k}
||\partial^{\alpha}u||_{L^{\infty}(K)}  \max_{2 \leq |\alpha| \leq
2k} \left( 1 + \sup_{x \in K} |\partial^{\alpha} f_s|
\right)^N , \qquad \omega \geq 1 . $$
\end{lemm}

Notice that the assumption (\ref{regularitenoncinifinie3}) is only
a condition  near $ x_0 (y) $. It guarantees the boundedness of $
\partial^{\alpha}u(x)/|x-x_0|^{2k - |\alpha|} $.

\bigskip

\noindent {\it Proof.} We proceed by induction and consider first
$k=1$. We would like to integrate by part using the operator $
|\nabla_x f_s|^{-2} \nabla_x f_s \cdot
\nabla_x $ but, since $ \nabla_x f_s $ may vanish on the support
of $u$, we consider $ L_{\delta } := (|\nabla_x f_s|^2 +
\delta)^{-1} \nabla_x f_s \cdot \nabla_x $ which satisfies
$$ i \omega \int e^{i \omega f_s (x,y)} u (x) d x = \lim_{\delta \downarrow 0}
\int  ( L_{\delta} e^{i \omega f_s (x,y)} ) u (x) d x .
$$
We then integrate by part at fixed $ \delta > 0 $, using that
$$ ^t \! L_{\delta} =
- \frac{1}{|\nabla_x f_s|^2 + \delta} \left\{ \nabla_x f_s \cdot \nabla_x + \Delta_x f_s -
\frac{2}{|\nabla_x f_s|^2 + \delta} \left\langle \mbox{Hess}_x
[f_s] \nabla_x f_s , \nabla_x f_s \right\rangle \right\} . $$
Since $ | \Delta_x f_s (x,y) u(x)| \lesssim \max_{|\alpha| = 2} || \Delta_x f_s(.,y) ||_{L^{\infty}(K)}||\partial^{\alpha}u||_{L^{\infty}} |x-x_0(y)|^2 $ and using (\ref{controlepointcritique}), by
letting $ \delta \downarrow 0 $ we get
$$ \left| i \omega \int e^{i \omega f_s (x,y)} u (x) d x \right| \leq
C \max_{|\alpha| \leq 2} ||\partial^{\alpha}u||_{L^{\infty}(K)}
\max_{|\alpha|=2} \left( 1 + \sup_{x \in \Ra^n} |\partial^{\alpha}
f_s| \right).
$$ Here the constant $C$
is independent of $y$, $u$, $s$ and $ \omega $; it depends only on
$ K  $ and the constant $c$ in (\ref{controlepointcritique}). The
result then follows by induction using that
$$ |\nabla_x f_s|^{-2} \langle \nabla_x f_s, \partial_x u \rangle , \qquad |\nabla_x f_s|^{-2} (\Delta_x f_s) u ,
\qquad |\nabla_x f_s|^{-4} \left\langle \mbox{Hess}_x [f_s]
\nabla_x f_s , \nabla_x f_s \right\rangle  u  $$ satisfy the
assumptions (\ref{regularitenoncinifinie1}),
(\ref{regularitenoncinifinie2}) and
(\ref{regularitenoncinifinie3}) if $ u $ does for $ k+1 $.
\finpreuve

\bigskip

\noindent {\it End of the proof of Lemma \ref{Hormanderrevu}.} We
next consider $ I (s) = \int e^{i \omega f_s (x,y)} u (x) dx $ so
that, for all $j \in \Na$, we have
$$ I^{(2j)} (s) = (i \omega)^{2j} \int e^{i \omega f_s (x,y)} g (x,y)^{2j} u(x) dx . $$
Since $ \partial^{\alpha}_x \left( g (x,y)^{2j} \right)_{| x = x_0
(y)} = 0 $ for all $ |\alpha| < 6 j $, Lemma
\ref{nonstationnaireHormander} yields, with $ k = 3 j \geq n / 2
$,
$$ |I^{(2j)}(s)| \leq C \omega^{-n/2} \max_{|\alpha| \leq 6 j} || \partial^{\alpha} u ||_{L^{\infty}(K)}
  \max_{2 \leq |\alpha| \leq
6 j} \left( 1 + \sup_{x \in \Ra^n} |\partial^{\alpha} f_s|
\right)^N , \qquad s \in [0,1] . $$ Since $ I (1) = \int e^{i
\omega f(x,y)} u(x) dx $, the estimate
$$ | I (1) - \sum_{l < 2 j} I^{(l)}(0) / l! | \leq \sup_{s \in [0,1]} |I^{(2j)}(s)|/ (2j)! , $$
reduces the proof to estimating the integrals $ I^{(l)}(0) $ whose
common phase $f_0$ is quadratic, up to a constant term and whose
amplitude is $ u (x) g (x,y)^{l} $. By Taylor's formula $ g (x,y) $ is of order $ |x -x_0(y)|^2 $
so the derivatives of $ u (x) g (x,y)^{l} $ may be of order $
\scal{x_0(y)}^{2l} $ on which we have no control. By choosing  $ \widetilde{K} $  a bounded neighborhood of $ K $ and
 applying Lemma \ref{nonstationnaireHormander} to the subset of $ \Omega$ on which $ x_0 (y) \notin \widetilde{K} $, we can assume that we consider those $y$ for which $ x_0 (y) \in \widetilde{K}  $.  
We then use the  Lemma 7.7.3 of \cite{Horm1} on  oscillatory integrals with quadratic phases, observing  that
 $ || \partial^{\alpha}_x  g
(.,y)^{l} ||_{L^{\infty}(K_x)} $ is controlled by (products
of)  of norms $ ||
\partial^{\beta}_x f (.,y) ||_{L^{\infty}(K_x)} $ with $ |\beta| \geq
2$, since $x$ is bounded on the support of $u$ and $ x_0 (y) $
remains bounded. \finpreuve


\begin{thebibliography}{99}

\bibitem{Banica} {\sc V. Banica}, {\it The non linear Schr\"odinger equation
    on hyperbolic space}, Comm. PDE, vol. 32, 10, 1643-1677 (2007).

\bibitem{BanicaCarles}{\sc V. Banica, R. Carles, G. Staffilani}, {\it
    Scattering theory for radial nonlinear Schr\"odinger equations on
    hyperbolic space}, GAFA, to appear.

\bibitem{BanicaDuyckaerts} {\sc V. Banica, T. Duyckaerts}, {\it Weighted Strichartz estimates for radial Schr\"odinger
equation on non-compact manifolds}, preprint.

\bibitem{Boucletresolvente} {\sc J.-M. Bouclet}, {\it Resolvent estimates for the Laplacian on asymptotically hyperbolic manifolds},
Ann. Henri Poincar\'e, no. 7, 527-561 (2006).

\bibitem{BoucletLpCF} {\sc \name}, {\it Semi-classical functional calculus on manifolds with ends and weighted $ L^p $
estimates}, preprint.

\bibitem{BoucletLP} {\sc \name}, {\it Littlewood-Paley decompositions on manifolds with ends}, preprint.

\bibitem{BoucletTzvetkov} {\sc J.-M. Bouclet, N. Tzvetkov}, {\it Strichartz estimates for long range
perturbations}, Amer. J. Math. vol 129, no. 6 (2007).

\bibitem{BoucletTzvetkov2} {\sc \name}, {\it On global Strichartz estimates
    for non trapping metrics}, J. Funct. Analysis, to appear.

\bibitem{Bourgain} {\sc J. Bourgain}, {\it Fourier transform restriction   phenomena for certain lattice subsets and applications
to nonlinear evolution equations I. Schr\"odinger equations.} GAFA, 3, 107-156 (1993).

\bibitem{BGT} {\sc N. Burq, P. G\'erard, N. Tzvetkov}, {\it Strichartz inequalities and the non linear 
Schr\"odinger equation on compact manifolds}, Amer. J. Math., 126, 569-605 (2004).

\bibitem{CardosoVodev} {\sc F. Cardoso, G. Vodev}, {\it Uniform estimates of the resolvent of the Laplace-Beltrami operator on infinite volume manifolds, II}, Ann. Henri Poincar\'e, 3, 673-691 (2002).

\bibitem{Doi} {\sc S. Doi}, {\it Smoothing effects of Schr\"odinger evolution groups on Riemannian manifolds}, Duke Math. J. 82, 679-706 (1996).

\bibitem{FH} {\sc R. G. Froese, P. D. Hislop}, {\it Spectral analysis of second-order elliptic operators on non-compact manifolds},
Duke Math. J. 58, 103-129 (1989).

\bibitem{GinibreVelo} {\sc J. Ginibre, G. Velo}, {\it The global Cauchy problem for the non linear Schr\"odinger equation}, Ann. IHP-Analyse non lin\'eaire 2, 309-327 (1985).

\bibitem{HTW} {\sc A. Hassell, T. Tao, J. Wunsch}, {\it Sharp Strichartz estimates on non-trapping asymptotically conic manifolds}, 
Amer. J.  Math., 128, 963-1024 (2006).

\bibitem{Horm1} {\sc L. H\"ormander},
{\it The Analysis of Linear Partial Differential Operators},
Springer (1983).

\bibitem{KeelTao}{\sc M. Keel, T. Tao}, {\it Endpoint Strichartz estimates},
  Amer. J. Math. vol. 120, 955-980 (1998).



\bibitem{MaMe1} {\sc R. Mazzeo, R.B. Melrose},
{\it Meromorphic extension of the resolvent on complete spaces
with asymptotically constant negative curvature}, J. Funct. Anal.
75, No. 2, 260-310 (1987).

\bibitem{Melr0} {\sc R.B. Melrose}, {\it Geometric scattering
theory}, Stanford lectures, Cambridge Univ. Press (1995).

\bibitem{Pier1}{\sc V. Pierfelice}, {\it Weighted Strichartz estimates for the radial perturbed Schr\"odinger equation on the hyperbolic
space}, Manuscripta Mathematica, vol 120, 4, 377-389 (2006).

\bibitem{Pier2}{\sc \name}, {\it  Weighted Strichartz estimates for the Schr\"odinger and Wave equations on Damek-Ricci
spaces}, preprint.



\bibitem{RZ} {\sc L. Robbiano, C. Zuily}, {\it  Strichartz estimates for Schr\"odinger equations with variable coefficients}, 
M\'em. SMF, No. 101-102 (2005).

\bibitem{Robebook} {\sc D. Robert}, {\it Autour de l'approximation semi-classique}, Birkh\"auser  (1987).


\bibitem{Strichartz} {\sc R. S. Strichartz}, {\it Restriction of Fourier transform to quadratic surfaces and decay of solutions of Wave equations}, Duke Math. J. 44, 705-774 (1977).

\bibitem{StTa} {\sc G. Staffilani, D. Tataru}, {\it   Strichartz estimates for a Schr\"odinger operator with non smooth coefficients}, Comm. PDE 27, 1337-1372 (2002).

\bibitem{Tatahype} {\sc D. Tataru}, {\it Strichartz estimates for the wave equation in the hyperbolic space and global existence for the semilinear Schr\"odinger equation}, Trans. Amer. Math. Soc. 353, no. 2, 795-807 (2001).


\end{thebibliography}
\end{document}